\documentclass[final,1p,nonatbib]{elsarticle}
\pdfoutput=1
\usepackage{subfiles}
\usepackage[shortlabels]{enumitem}

\makeatletter
\let\c@author\relax
\makeatother
\let\bibfont\relax

\makeatletter
\def\ps@pprintTitle{%
 \let\@oddhead\@empty
 \let\@evenhead\@empty
 \def\@oddfoot{\centerline{\thepage}}%
 \let\@evenfoot\@oddfoot}
\makeatother

\usepackage{\subfix{double-theory-macros}}
\numberwithin{equation}{section}

\usepackage[breaklinks,bookmarksnumbered,colorlinks]{hyperref}
\usepackage[capitalize,noabbrev,nameinlink]{cleveref}
\usepackage{doi} %

\usepackage{bookmark}
\usepackage{thmtools}
\makeatletter
\newcommand*{\theorembookmark}[1]{%
  \ifx\@currentlabelname\@empty
  \else
  \bookmark[dest=\@currentHref,rellevel=1,keeplevel]{%
    #1\space\csname the\thmt@envname\endcsname :\space
    \@currentlabelname
  }
  \fi
}
\makeatother
\declaretheorem[within=section,postheadhook={\theorembookmark{Thm}}]{theorem}
\declaretheorem[sibling=theorem,postheadhook={\theorembookmark{Prop}}]{proposition}
\declaretheorem[sibling=theorem,postheadhook={\theorembookmark{Cor}}]{corollary}
\declaretheorem[sibling=theorem]{lemma}
\declaretheorem[sibling=theorem,style=definition,postheadhook={\theorembookmark{Def}}]{definition}
\declaretheorem[sibling=theorem,style=remark,qed=\qedsymbol]{remark}
\declaretheorem[sibling=theorem,style=remark,qed=\qedsymbol]{example}
\declaretheorem[sibling=theorem,style=remark,postheadhook={\theorembookmark{Th}},qed=\qedsymbol]{theory}

\usepackage[style=trad-abbrv,maxbibnames=10,sorting=nyt]{biblatex}
\renewbibmacro{in:}{}

\usepackage{quiver}
\usepackage{nicematrix}
\NiceMatrixOptions{cell-space-limits=3pt}
\newenvironment{dblArray}[1]{\begin{NiceArray}{#1}[hvlines]}{\end{NiceArray}}

\usepackage[titles]{tocloft}
\begin{document}

\begin{frontmatter}

\title{Cartesian double theories: \\%
  A double-categorical framework for categorical doctrines}

\author[1]{Michael Lambert}
\ead{michael.james.lambert@gmail.com}

\author[2]{Evan Patterson}
\ead{evan@epatters.org}

\affiliation[1]{organization={University of Massachusetts-Boston},
  city={Boston, MA}, country={United States}}
\affiliation[2]{organization={Topos Institute},
  city={Berkeley, CA}, country={United States}}

\begin{abstract}
  The categorified theories known as ``doctrines'' specify a category equipped
  with extra structure, analogous to how ordinary theories specify a set with
  extra structure. We introduce a new framework for doctrines based on double
  category theory. A cartesian double theory is defined to be a small double
  category with finite products and a model of a cartesian double theory to be a
  finite product-preserving lax functor out of it. Many familiar categorical
  structures are models of cartesian double theories, including categories,
  presheaves, monoidal categories, braided and symmetric monoidal categories,
  2-groups, multicategories, and cartesian and cocartesian categories. We show
  that every cartesian double theory has a unital virtual double category of
  models, with lax maps between models given by cartesian lax natural
  transformations, bimodules between models given by cartesian modules, and
  multicells given by multimodulations. In many cases, the virtual double
  category of models is representable, hence is a genuine double category.
  Moreover, when restricted to pseudo maps, every cartesian double theory has a
  virtual equipment of models, hence an equipment of models in the representable
  case. Compared with 2-monads, double theories have the advantage of being
  straightforwardly presentable by generators and relations, as we illustrate
  through a large number of examples.
\end{abstract}

\end{frontmatter}

\tableofcontents

\section{Introduction}
\label{sec:introduction}

The passage from algebraic theories as syntactical objects to Lawvere theories
with their functorial semantics \cite{lawvere1963} transformed universal algebra
and spurred the development of a wide-reaching assimilation of logic and type
theory into category theory. In categorical logic, a logical system is defined
by a 2-category of categories with extra structure. Theories within the logic
are objects of the 2-category; models of a theory are structure-preserving
functors out of the theory; and model homomorphisms are natural transformations.
In some cases, the logic defined by a 2-category can be identified with a
preexisting system given in the traditional syntactic style. Lawvere theories
correspond to single-sorted algebraic theories, categories with finite products
to multi-sorted algebraic theories, categories with finite limits to essentially
algebraic theories, and cartesian closed categories to simply typed lambda
calculus with product types. Much work in categorical logic has been devoted to
building this dictionary \cite{crole1993,johnstone2002}. In other cases, the
correspondence is less clear and category theory provides a toolbox for
generating new logics. Symmetric monoidal categories are a relatively weak,
resource-sensitive logic whose type theory is still a subject of active
investigation \cite{shulman2021}. Markov categories are a logic of
nondeterminism whose models include probabilistic and statistical models
\cite{fritz2020,patterson2020}.

The purpose of any logic is to precisely specify and interpret a uniform class
of structures. For a category theorist confronted by the landscape of
categorical logics, the inevitable next step is to categorify the study of logic
itself, giving a means to precisely specify and interpret a uniform class of
\emph{logics}. The higher theories which have logics as models are known as
\emph{doctrines}.\footnote{The word ``doctrine'' apparently first appears in the
  category theory literature through Lawvere \cite{lawvere1969}, but Lawvere
  attributes the usage to Jon Beck.} Whereas an ordinary, one-dimensional theory
is interpreted to give a category of models, a doctrine is interpreted to give
at least a 2-category of models, with a model being a category (or another
category-like object, such as a multicategory or a polycategory) with extra
structure.

The word ``doctrine,'' like the word ``theory,'' connotes a general idea that
can be made precise in different ways. By far the best known formalism,
extensively developed by Kelly and collaborators
\cite{kellystreet1974,kelly1974,blackwell1989} \cite[\S{4}]{lack2010}, takes a
doctrine to be a 2-monad on $\Cat$ or a similar 2-category. For example, there
are 2-monads on $\Cat$ whose algebras are strict or weak monoidal categories
and, in either case, whose strict, pseudo, and lax morphisms of algebras are
strict, strong, and lax monoidal functors. The view of doctrines as 2-monads on
$\Cat$ categorifies the view of theories as monads on $\Set$. Since finitary
monads and Lawvere theories are equivalent \cite[\mbox{Theorem
  A.37}]{adamek2010}, one might wonder whether there is a formalism for
doctrines that categorifies the concept of a Lawvere theory or, more generally,
a finite products theory. It is the aim of this work to develop such a framework
for categorical doctrines.

An advantage of our approach is the ability to present doctrines by generators
and relations. From its earliest beginnings, formal logic has been understood as
a finitary and mechanistic calculus. Thus, while it is important that a Lawvere
theory is an invariant description of a theory, it is equally important that a
Lawvere theory can be presented by generators and relations, much like a group,
another kind of invariant object, can be presented by generators and relations.
In this paper we will present many different doctrines by generators and
relations, from categories and presheaves on categories to monoidal categories
and multicategories to cartesian and cocartesian categories. In contrast,
presentations of 2-monads are less elementary and direct.\footnote{To present a
  2-monad, one starts with free 2-monads and then glues them together using
  2-categorical colimits, such as coproducts, co-inserters, and co-equifiers, in
  the 2-category of 2-monads \cite{kelly1993} \cite[\S{5}]{lack2010}.}

To realize a notion of doctrine based on two-dimensional functorial semantics,
we begin with Bénabou's famous observation that a lax functor
\begin{equation*}
  F: \cat{1} \to \bicat{Span}
\end{equation*}
from the terminal 2-category to the bicategory of sets, spans, and
(feet-preserving) maps of spans is equivalent to a small category
\cite{benabou1967}. The unique object of $\cat{1}$ is sent to the set of objects
of the category, the identity on the unique object is sent to the set of
morphisms equipped with their source and target maps, and the laxators and
unitors of the lax functor define the composition and identities of the
category. This observation can be generalized in various ways. A lax functor
$F: \cat{1} \to \bicat{B}$ into an arbitrary bicategory $\bicat{B}$ is a monad
in $\bicat{B}$. A lax functor $F: \coDisc(X) \to \bicat{B}$ from the codiscrete
2-category on a set $X$ into $\bicat{B}$ is a category enriched in $\bicat{B}$
with set of objects $X$. For Bénabou, who invented bicategories and lax
functors, such facts were a key motivation to study lax functors in the first
place \cite[\S{3.1}]{lack2010}.

Despite the appeal and long history of these correspondences, lax functors
between bicategories cannot play the role of models of a doctrine, for a simple
reason: they give the wrong maps between models. The well-behaved
transformations
\begin{equation*}
  \begin{tikzcd}
    {\cat{1}} & {\bicat{Span}}
    \arrow[""{name=0, anchor=center, inner sep=0}, "F", curve={height=-18pt}, from=1-1, to=1-2]
    \arrow[""{name=1, anchor=center, inner sep=0}, "G"', curve={height=18pt}, from=1-1, to=1-2]
    \arrow[shorten <=5pt, shorten >=5pt, Rightarrow, from=0, to=1]
  \end{tikzcd}
\end{equation*}
between lax functors are called \emph{icons} \cite{lack2010icons}. In the
correspondence between span-valued lax functors and categories, icons are
\emph{identity-on-objects} functors, which are far too restrictive.

As is now generally appreciated, the solution to such problems is to move from
bicategories to double categories. A lax functor $F: \dbl{1} \to \Span$ from the
terminal double category to the double category of sets, functions, spans, and
maps of spans is again equivalent to a small category. But now a natural
transformation
\begin{equation*}
  \begin{tikzcd}
    {\dbl{1}} & \Span
    \arrow[""{name=0, anchor=center, inner sep=0}, "F", curve={height=-18pt}, from=1-1, to=1-2]
    \arrow[""{name=1, anchor=center, inner sep=0}, "G"', curve={height=18pt}, from=1-1, to=1-2]
    \arrow[shorten <=5pt, shorten >=5pt, Rightarrow, from=0, to=1]
  \end{tikzcd}
\end{equation*}
between lax functors is precisely a functor between the corresponding
categories. Moreover, with the right choice of 2-morphisms, which turn out to be
\emph{modulations} (between identity modules), the resulting 2-category
$\LaxFun(\dbl{1},\Span)$ of lax functors is equivalent to $\Cat$, the 2-category
of small categories. Modulations were introduced by Paré \cite{pare2011},
building on the definition by Cockett and others in the bicategorical setting
\cite{cockett2003}.

Having recovered the 2-category of categories, the most basic of all doctrines,
we are encouraged to make the leap that any small, strict double category
$\dbl{T}$ defines a doctrine through the 2-category $\LaxFun(\dbl{T},\Span)$ of
span-valued lax functors on it. In this capacity, the double category $\dbl{T}$
will be called a ``double theory,'' specifically a \define{simple double theory}
to emphasize that no additional double-categorical structure is assumed. So, a
simple double theory $\dbl{T}$ has the 2-category $\LaxFun(\dbl{T},\Span)$ as
its 2-category of models.

This proposal is justified by the fact, made no less remarkable by the ease of
its proof, that a lax double functor contains within it all of the most
fundamental definitions of category theory. A lax functor $F: \dbl{D} \to \Span$
sends each object in $\dbl{D}$ to a category, each arrow in $\dbl{D}$ to a
functor, and each cell in $\dbl{D}$ bounded by identity proarrows to a natural
transformation, in a 2-functorial way. Moreover, it sends each proarrow in
$\dbl{D}$ to a profunctor and each cell in $\dbl{D}$ to map of profunctors,
which is again a kind of natural transformation. More precise and general
versions of these statements are proved in \cref{sec:lax-functors}.

To obtain from these promising observations a viable framework for doctrines,
two obstacles must be overcome. First, although their 2-categories of models
include several interesting examples, such as adjunctions and monads
(\cref{sec:simple-double-theories}), simple double theories are too inexpressive
to present the most familiar categorical doctrines, beginning with monoidal
categories. Second, although natural transformations between lax double functors
correctly capture the \emph{strict} maps between models, it is the pseudo, lax,
or oplax maps that are often more important. For example, the default notion of
map between monoidal categories is most commonly taken to be \emph{strong}
monoidal functors. Let us consider these problems in turn.

Just as categorical logic has its origins in cartesian categories with their
close connection to algebraic theories, a natural point of departure for
doctrines is \emph{cartesian} double categories \cite{aleiferi2018}. These are
double categories with finite products, in the sense of Grandis and Paré's
theory of limits in double categories \cite{grandis1999}. In a cartesian double
category $\dbl{D}$, both underlying categories $\dbl{D}_0$ and $\dbl{D}_1$ have
finite products, which are preserved by the source and target functors
$s,t: \dbl{D}_1 \rightrightarrows \dbl{D}_0$, as well as by the external
composition and identities. The prototypical example of a cartesian double
category is $\Span$ itself.

We say that a \define{cartesian double theory} is a small, strict cartesian
double category $\dbl{T}$ and that a \define{model} of the theory is a cartesian
lax functor $F: \dbl{T} \to \Span$, meaning that both underlying functors
$F_0: \dbl{T}_0 \to \Set$ and
$F_1: \dbl{T}_1 \to \Set^{\{\bullet \leftarrow \bullet \rightarrow \bullet\}}$
preserve finite products. As a motivating example, the \define{theory of a
  pseudomonoid} is generated by a single object $x$, arrows $\otimes: x^2 \to x$
and $I: 1 \to x$, and cells
\begin{equation*}
  \begin{tikzcd}
    {x^3} & {x^3} \\
    {x^2} & {x^2} \\
    x & x
    \arrow[""{name=0, anchor=center, inner sep=0}, "{\mathrm{id}_x^3}", "\shortmid"{marking}, from=1-1, to=1-2]
    \arrow["{\otimes \times 1_x}"', from=1-1, to=2-1]
    \arrow["\otimes"', from=2-1, to=3-1]
    \arrow["{1_x \times \otimes}", from=1-2, to=2-2]
    \arrow[""{name=1, anchor=center, inner sep=0}, "{\mathrm{id}_x}"', "\shortmid"{marking}, from=3-1, to=3-2]
    \arrow["\otimes", from=2-2, to=3-2]
    \arrow["\alpha"{description}, draw=none, from=0, to=1]
  \end{tikzcd}
  \qquad\qquad
  \begin{tikzcd}
    x & x \\
    {x^2} \\
    x & x
    \arrow["{I \times 1_x}"', from=1-1, to=2-1]
    \arrow["\otimes"', from=2-1, to=3-1]
    \arrow[""{name=0, anchor=center, inner sep=0}, "{\mathrm{id}_x}"', "\shortmid"{marking}, from=3-1, to=3-2]
    \arrow[""{name=1, anchor=center, inner sep=0}, "{\mathrm{id}_x}", "\shortmid"{marking}, from=1-1, to=1-2]
    \arrow[Rightarrow, no head, from=1-2, to=3-2]
    \arrow["\lambda"{description}, draw=none, from=1, to=0]
  \end{tikzcd}
  \qquad\qquad
  \begin{tikzcd}
    x & x \\
    {x^2} \\
    x & x
    \arrow["{1_x \times I}"', from=1-1, to=2-1]
    \arrow["\otimes"', from=2-1, to=3-1]
    \arrow[""{name=0, anchor=center, inner sep=0}, "{\mathrm{id}_x}"', "\shortmid"{marking}, from=3-1, to=3-2]
    \arrow[""{name=1, anchor=center, inner sep=0}, "{\mathrm{id}_x}", "\shortmid"{marking}, from=1-1, to=1-2]
    \arrow[Rightarrow, no head, from=1-2, to=3-2]
    \arrow["\rho"{description}, draw=none, from=1, to=0]
  \end{tikzcd}
\end{equation*}
representing the associators and left and right unitors, as well as their
inverses, subject to equations expressing the usual coherence axioms (see
\cref{th:pseudomonoid} for details). A model of the theory of pseudomonoids is
precisely a (weak) monoidal category.

Natural transformations between models of the theory of pseudomonoids are
\emph{strict} monoidal functors. To loosen the maps, we introduce a notion of
\emph{lax} natural transformation $\alpha: F \To G$ between lax double functors
$F, G: \dbl{D} \to \dbl{E}$ whose data includes, for every arrow $f: x \to y$ in
$\dbl{D}$, a naturality comparison cell
\begin{equation*}
  \begin{tikzcd}[row sep=scriptsize]
    Fx & Fx \\
    Gx & Fy \\
    Gy & Gy
    \arrow[""{name=0, anchor=center, inner sep=0}, "{\mathrm{id}_{Fx}}", "\shortmid"{marking}, from=1-1, to=1-2]
    \arrow["{\alpha_x}"', from=1-1, to=2-1]
    \arrow["Gf"', from=2-1, to=3-1]
    \arrow["Ff", from=1-2, to=2-2]
    \arrow["{\alpha_y}", from=2-2, to=3-2]
    \arrow[""{name=1, anchor=center, inner sep=0}, "{G(\mathrm{id}_y)}"', "\shortmid"{marking}, from=3-1, to=3-2]
    \arrow["{\alpha_f}"{description}, draw=none, from=0, to=1]
  \end{tikzcd}
\end{equation*}
in $\dbl{E}$, subject to several axioms (\cref{def:lax-transformation}). When
$\dbl{D}$ and $\dbl{E}$ are cartesian, we say that a lax transformation is
\define{cartesian} when it is strictly natural with respect to projection maps
in $\dbl{D}_0$. We show that there is a 2-category
$\CartLaxFun_\ell(\dbl{D},\dbl{E})$ of cartesian lax functors
$\dbl{D} \to \dbl{E}$, cartesian lax transformations, and modulations. When
$\dbl{T}$ is the theory of pseudomonoids, the 2-category of models
$\CartLaxFun_\ell(\dbl{T},\Span)$ is equivalent to the 2-category of monoidal
categories, lax monoidal functors, and monoidal natural transformations. Oplax
and pseudo maps between models of double theories are constructed similarly.

A double-categorical framework for doctrines might be expected to produce not
just a 2-category but a double category of models. This is true for double
theories with the caveat that the double category of models is in general only
\emph{virtual} \cite{burroni1971,leinster2004,cruttwell2010}. Extending Paré's
definition of a module between lax double functors \cite{cockett2003,pare2011},
we define a \emph{cartesian} module between cartesian lax functors to be a
module that suitably preserves finite products. We then show that for any
cartesian double categories $\dbl{D}$ and $\dbl{E}$, there is a virtual double
category $\vCartLaxFun_\ell(\dbl{D},\dbl{E})$ of cartesian lax functors,
cartesian lax natural transformations, cartesian modules, and multimodulations.
In particular, every cartesian double theory $\dbl{T}$ has a virtual double
category of models, $\vCartLaxFun_\ell(\dbl{T},\Span)$. This virtual double
category always has units but, due to obstructions to composing modules between
lax functors \cite{pare2013}, it is not always representable as a double
category. Some sufficient conditions for representability are known
\cite{pare2013} but the general situation is not well understood.

At least when our theories are purely 2-categorical---that is, when the double
theories have only trivial proarrows---we can give a more satisfactory answer:
the virtual double categories of models of such theories are always
representable as double categories. In more detail, any 2-category defines a
strict double category whose proarrows are all identities. Several of the
theories that we present are of this type, including the simple double theories
of adjunctions (\cref{th:adjunctions}) and monads (\cref{th:monad}) and the
cartesian double theories of monoids (\cref{th:monoid}), pseudomonoids
(\cref{th:pseudomonoid}), and cartesian monoidal categories
(\cref{th:cart-mon-cat-I}). In each of these cases, as shown in
\cref{prop:lax-func-form-double-category}, we obtain a genuine double category
of models whose proarrows are familiar profunctor-like structures between the
models. For example, when $\dbl{T}$ is the theory of pseudomonoids, the double
category of models $\vCartLaxFun_\ell(\dbl{T},\Span)$ is equivalent to the
double category of monoidal categories, lax monoidal functors, monoidal
profunctors, and monoidal natural transformations.

Two-dimensional theories based on 2-categories have been studied in other
contexts \mbox{\cite[\S{9}]{bourke2021}}. In general, a \emph{2-theory} is a
small 2-category with finite weighted limits, for some choice of weights, and a
\emph{model} of a 2-theory is a 2-functor out of it preserving those weighted
limits. Since horizontally trivial cartesian double categories are the same as
2-categories with strict finite 2-products, horizontally trivial cartesian
double theories can be identified with \emph{finite product 2-theories}.
Effectively, then, we give a double category of models to each finite product
2-theory.

The result that our models form a double category when the theory is purely
2-categorical is proved in stages. A first step is repackaging lax double
functors valued in a double category $\dbl{E}$ as \emph{normal}, or actually
\emph{unitary} given a choice of units, lax functors valued in
$\Module{\dbl{E}}$, the double category of category objects and profunctors in
$\dbl{E}$, under mild conditions ensuring that $\Module{\dbl{E}}$ is itself a
double category. In the special case when $\dbl{E} = \Span$ is the double
category of spans, $\Module{\dbl{E}} = \Prof$ is just the double category of
profunctors. Thus, our \emph{normalization} result
(\cref{prop:unitalizationisoof1categories}) shows that span-valued models of
double theories are equivalently described as normal profunctor-valued models.
This equivalence is a useful tool in confirming that the proposed theories do in
fact have the intended categorical structures as their models. But,
additionally, it is a tool to establish
(\cref{cor:2-categories-of-models-non-cartesian,cor:2-categories-of-models-cartesian})
that when the theory $\dbl{T}$ is a 2-category, its models are the objects of a
2-category with lax transformations and special modulations as morphisms and
2-cells. The 2-cells turn out to be essentially ordinary modifications of
suitably structured transformations. Such 2-categories of models thus underlie
the double categories of models when $\dbl{T}$ is a 2-category. That this
assumption on $\dbl{T}$ enables modules to be composed in the simple and
cartesian cases is proved directly (\cref{prop:lax-func-form-double-category}).
The proof is phrased in terms of unitary lax functors to simplify the arguments
and the normalization results are then applied to show that the correct models
are recovered in several cases of interest
(\cref{cor:double-categories-of-models-examples}).

Most double categories that would be semantics for double theories, including
double categories of spans and of matrices, are \emph{equipments}
\cite{wood1982,shulman2008}, a structure in which proarrows have universal
restrictions along pairs of incoming arrows. We show that when the target double
category $\dbl{E}$ is a cartesian equipment, the virtual double category
$\vCartLaxFun_\pseudo(\dbl{D},\dbl{E})$ of cartesian lax functors, cartesian
\emph{pseudo} transformations, cartesian modules, and multimodulations is a
virtual equipment \cite{cruttwell2010}. Hence, it is an equipment when it is
representable as a double category. For example, the theory of pseudomonoids
yields an equipment of monoidal categories, strong monoidal functors, monoidal
profunctors, and monoidal natural transformations. Equipments and related
structures have been identified as ideal environments for formal category
theory, enabling many parts of ordinary category theory to be reproduced
abstractly. That any double theory has an equipment or virtual equipment of
models is further evidence that double categories can serve as a foundation for
categorical doctrines.

There has been at least one previous attempt to use double categories in
categorical logic. In unpublished work reported in talks
\cite{pare2009,pare2010}, Paré extends the 1-categorical notion of a coherent
theory to a double theory, in which proarrows represent formulas or relations
and cells represent logical entailments or implications. Although analogous in
some respects, the aims and technical approach of that work are notably
different from our own. Paré studies a one-dimensional logic (namely, coherent
logic) from a double-categorical point of view, with models of a double theory
being pseudo double functors into $\Rel$. We study two-dimensional categorical
logic, with models of a double theory being lax double functors into $\Span$ or
equivalently normal lax functors into $\Prof$. Although Paré does not specify a
definite double category of models, he seems to have in mind pseudo or lax
horizontal transformations (``protransformations'') as proarrows between models.
In contrast, we take proarrows between models to be modules between lax double
functors, as these generalize profunctors between categories.

\paragraph{Background}

We have endeavored to write a largely self-contained paper. We take for granted
the definitions of a (pseudo) double category, a (pseudo) double functor, and a
natural transformation of double functors, which can be found in many sources,
including the textbook by Grandis \cite{grandis2019}. The final section of the
paper also assumes knowledge of virtual double categories
\cite{leinster2004,cruttwell2010}. We begin in \cref{sec:lax-functors} by
reviewing the concept of a lax double functor, which plays such a central role
in the development. In \cref{sec:cartesian-equipments}, we review cartesian
double categories, equipments, and their conjunction as cartesisan equipments.
We also give detailed definitions of the higher morphisms involving double
categories, such as lax natural transformations, modules, and modulations, in
part because we have needed to generalize these notions beyond the definitions
available in the literature.

\paragraph{Conventions}

Unless otherwise stated, double categories and double functors are assumed to be
pseudo. We take our double categories to be strict in the vertical direction and
weak in the horizontal direction, although we mostly avoid the terminology of
``vertical'' and ``horizontal'' morphisms, speaking instead of ``arrows'' and
``proarrows.''

We write categories $\cat{C}, \cat{D}, \dots$ in sans-serif font; 2-categories
and bicategories $\bicat{B}, \bicat{C}, \dots$ in bold font; and double
categories $\dbl{D}, \dbl{E}, \dots$ in blackboard bold font. Composites of
morphisms $x \xto{f} y \xto{g} z$ in a category are written variously in
diagrammatic order as $f \cdot g$ or applicative order as $g \circ f$.
Composites of proarrows $x \xproto{m} y \xproto{n} z$ in a double category are
always written in diagrammatic order as $m \odot n$. Identity arrows are written
as $1_x: x \to x$ and identity proarrows as $\id_x: x \proto x$.

\section{Lax functors}
\label{sec:lax-functors}

Lax double functors are a natural foundation for categorical doctrines because,
as we sketched in \cref{sec:introduction} and will elaborate in this section,
the concept of a lax double functor contains within it all of the most
fundamental definitions of category theory: categories, functors, and natural
transformations, as well as profunctors and maps of profunctors.

Lax functors between double categories are defined in many sources, including
the textbook by Grandis \cite[Definition 3.5.1]{grandis2019}. For ease of
reference, we recall the complete definition.

\begin{definition}[Lax functor] \label{def:lax-functor}
  A \define{lax double functor} $F: \dbl{D} \to \dbl{E}$ between pseudo double
  categories $\dbl{D}$ and $\dbl{E}$ consists of
  \begin{itemize}
    \item a pair of functors $F_0: \dbl{D}_0 \to \dbl{E}_0$ and
      $F_1: \dbl{D}_1 \to \dbl{E}_1$ between the underlying categories of
      objects and morphisms, which preserve the external source and target:
      \begin{equation*}
        \begin{tikzcd}
          {\dbl{D}_1} & {\dbl{E}_1} \\
          {\dbl{D}_0} & {\dbl{E}_0}
          \arrow["{F_1}", from=1-1, to=1-2]
          \arrow["s"', from=1-1, to=2-1]
          \arrow["{F_0}"', from=2-1, to=2-2]
          \arrow["s", from=1-2, to=2-2]
        \end{tikzcd}
        \qquad\text{and}\qquad
        \begin{tikzcd}
          {\dbl{D}_1} & {\dbl{E}_1} \\
          {\dbl{D}_0} & {\dbl{E}_0}
          \arrow["{F_1}", from=1-1, to=1-2]
          \arrow["t"', from=1-1, to=2-1]
          \arrow["{F_0}"', from=2-1, to=2-2]
          \arrow["t", from=1-2, to=2-2]
        \end{tikzcd};
      \end{equation*}
    \item for every consecutive pair of proarrows $x \xproto{m} y \xproto{n} z$
      in $\dbl{D}$, a globular cell in $\dbl{E}$
      \begin{equation*}
        \begin{tikzcd}
          Fx & Fy & Fz \\
          Fx && Fz
          \arrow[""{name=0, anchor=center, inner sep=0}, "{F(m \odot n)}"', "\shortmid"{marking}, from=2-1, to=2-3]
          \arrow[Rightarrow, no head, from=1-1, to=2-1]
          \arrow[Rightarrow, no head, from=1-3, to=2-3]
          \arrow["Fm", "\shortmid"{marking}, from=1-1, to=1-2]
          \arrow["Fn", "\shortmid"{marking}, from=1-2, to=1-3]
          \arrow["{F_{m,n}}"{description}, draw=none, from=1-2, to=0]
        \end{tikzcd},
      \end{equation*}
      the \define{laxator} or \define{composition comparison} at $m$ and $n$;
    \item for every object $x \in \dbl{D}$, a globular cell in $\dbl{E}$
    \begin{equation*}
      \begin{tikzcd}
        Fx & Fx \\
        Fx & Fx
        \arrow[""{name=0, anchor=center, inner sep=0}, "{\mathrm{id}_{Fx}}", "\shortmid"{marking}, from=1-1, to=1-2]
        \arrow[""{name=1, anchor=center, inner sep=0}, "{F \mathrm{id}_x}"', "\shortmid"{marking}, from=2-1, to=2-2]
        \arrow[Rightarrow, no head, from=1-1, to=2-1]
        \arrow[Rightarrow, no head, from=1-2, to=2-2]
        \arrow["{F_x}"{description}, draw=none, from=0, to=1]
      \end{tikzcd},
    \end{equation*}
    the \define{unitor} or \define{identity comparison} at $x$.
  \end{itemize}
  The following axioms must be satisfied.
  \begin{itemize}
    \item Naturality of laxators: for any cells
      $\inlineCell{x}{y}{x'}{y'}{m}{m'}{f}{g}{\alpha}$ and
      $\inlineCell{y}{z}{y'}{z'}{n}{n'}{g}{h}{\beta}$ in $\dbl{D}$,
      \begin{equation} \label{eq:naturality-laxators}
        \begin{tikzcd}
          Fx & Fy & Fz \\
          {Fx'} & {Fy'} & {Fz'} \\
          {Fx'} && {Fz'}
          \arrow[""{name=0, anchor=center, inner sep=0}, "{F(m' \odot n')}"', "\shortmid"{marking}, from=3-1, to=3-3]
          \arrow[Rightarrow, no head, from=2-1, to=3-1]
          \arrow[Rightarrow, no head, from=2-3, to=3-3]
          \arrow[""{name=1, anchor=center, inner sep=0}, "Fm", "\shortmid"{marking}, from=1-1, to=1-2]
          \arrow[""{name=2, anchor=center, inner sep=0}, "Fn", "\shortmid"{marking}, from=1-2, to=1-3]
          \arrow[""{name=3, anchor=center, inner sep=0}, "{Fm'}"', "\shortmid"{marking}, from=2-1, to=2-2]
          \arrow[""{name=4, anchor=center, inner sep=0}, "{Fn'}"', "\shortmid"{marking}, from=2-2, to=2-3]
          \arrow["Ff"', from=1-1, to=2-1]
          \arrow["Fg"{description}, from=1-2, to=2-2]
          \arrow["Fh", from=1-3, to=2-3]
          \arrow["F\alpha"{description}, draw=none, from=1, to=3]
          \arrow["F\beta"{description}, draw=none, from=2, to=4]
          \arrow["{F_{m',n'}}"{description}, draw=none, from=2-2, to=0]
        \end{tikzcd}
        \quad=\quad
        \begin{tikzcd}
          Fx & Fy & Fz \\
          Fx && Fz \\
          {Fx'} && {Fz'}
          \arrow[""{name=0, anchor=center, inner sep=0}, "{F(m' \odot n')}"', "\shortmid"{marking}, from=3-1, to=3-3]
          \arrow["Fm", "\shortmid"{marking}, from=1-1, to=1-2]
          \arrow["Fn", "\shortmid"{marking}, from=1-2, to=1-3]
          \arrow[""{name=1, anchor=center, inner sep=0}, "{F(m\odot n)}"', "\shortmid"{marking}, from=2-1, to=2-3]
          \arrow["Ff"', from=2-1, to=3-1]
          \arrow["Fh", from=2-3, to=3-3]
          \arrow[Rightarrow, no head, from=1-1, to=2-1]
          \arrow[Rightarrow, no head, from=1-3, to=2-3]
          \arrow["{F(\alpha \odot \beta)}"{description}, draw=none, from=1, to=0]
          \arrow["{F_{m,n}}"{description}, draw=none, from=1-2, to=1]
        \end{tikzcd}.
      \end{equation}
    \item Naturality of unitors: for every arrow $f: x \to y$ in $\dbl{D}$,
      \begin{equation} \label{eq:naturality-unitors}
        \begin{tikzcd}
          Fx & Fx \\
          Fy & Fy \\
          Fy & Fy
          \arrow[""{name=0, anchor=center, inner sep=0}, "{F \mathrm{id}_y}"', "\shortmid"{marking}, from=3-1, to=3-2]
          \arrow[Rightarrow, no head, from=2-2, to=3-2]
          \arrow[Rightarrow, no head, from=2-1, to=3-1]
          \arrow[""{name=1, anchor=center, inner sep=0}, "{\mathrm{id}_{Fx}}", "\shortmid"{marking}, from=1-1, to=1-2]
          \arrow[""{name=2, anchor=center, inner sep=0}, "{\mathrm{id}_{Fy}}"', "\shortmid"{marking}, from=2-1, to=2-2]
          \arrow["Ff"', from=1-1, to=2-1]
          \arrow["Ff", from=1-2, to=2-2]
          \arrow["{\mathrm{id}_{Ff}}"{description}, draw=none, from=1, to=2]
          \arrow["{F_y}"{description}, draw=none, from=2, to=0]
        \end{tikzcd}
        \quad=\quad
        \begin{tikzcd}
          Fx & Fx \\
          Fx & Fx \\
          Fy & Fy
          \arrow[""{name=0, anchor=center, inner sep=0}, "{\mathrm{id}_{Fx}}", "\shortmid"{marking}, from=1-1, to=1-2]
          \arrow[""{name=1, anchor=center, inner sep=0}, "{F \mathrm{id}_x}"', "\shortmid"{marking}, from=2-1, to=2-2]
          \arrow[Rightarrow, no head, from=1-1, to=2-1]
          \arrow[Rightarrow, no head, from=1-2, to=2-2]
          \arrow["Ff"', from=2-1, to=3-1]
          \arrow["Ff", from=2-2, to=3-2]
          \arrow[""{name=2, anchor=center, inner sep=0}, "{F \mathrm{id}_y}"', "\shortmid"{marking}, from=3-1, to=3-2]
          \arrow["{F_x}"{description}, draw=none, from=0, to=1]
          \arrow["{F \mathrm{id}_f}"{description}, draw=none, from=1, to=2]
        \end{tikzcd}.
      \end{equation}
    \item Associativity: for every triple of consecutive proarrows
      $w \xproto{m} x \xproto{n} y \xproto{p} z$ in $\dbl{D}$, the diagram in
      $\dbl{E}_1$ commutes:
      \begin{equation*}
        \begin{tikzcd}
          {(Fm \odot Fn) \odot Fp} & {Fm \odot (Fn \odot Fp)} \\
          {F(m \odot n) \odot Fp} & {Fm \odot F(n \odot p)} \\
          {F((m \odot n) \odot p)} & {F(m \odot (n \odot p))}
          \arrow["{1_{Fm} \odot F_{n,p}}", from=1-2, to=2-2]
          \arrow["{F_{m, n \odot p}}", from=2-2, to=3-2]
          \arrow["{F_{m,n} \odot 1_{Fp}}"', from=1-1, to=2-1]
          \arrow["{F_{m \odot n, p}}"', from=2-1, to=3-1]
          \arrow["\cong", from=1-1, to=1-2]
          \arrow["\cong"', from=3-1, to=3-2]
        \end{tikzcd}.
      \end{equation*}
    \item Unitality: for every proarrow $m: x \proto y$ in $\dbl{D}$, the
      diagrams in $\dbl{E}_1$ commute:
      \begin{equation*}
        \begin{tikzcd}[column sep=large]
          {\mathrm{id}_{Fx} \odot Fm} & {F\mathrm{id}_x \odot Fm} \\
          Fm & {F(\mathrm{id}_x \odot m)}
          \arrow["{F_x \odot 1_{Fm}}", from=1-1, to=1-2]
          \arrow["{F_{x,m}}", from=1-2, to=2-2]
          \arrow["\cong"', from=1-1, to=2-1]
          \arrow["\cong", from=2-2, to=2-1]
        \end{tikzcd}
        \qquad
        \begin{tikzcd}[column sep=large]
          {Fm \odot \mathrm{id}_{Fy}} & {Fm \odot F\mathrm{id}_y} \\
          Fm & {F(m \odot \mathrm{id}_y)}
          \arrow["{1_{Fm} \odot F_y}", from=1-1, to=1-2]
          \arrow["{F_{m,y}}", from=1-2, to=2-2]
          \arrow["\cong"', from=1-1, to=2-1]
          \arrow["\cong", from=2-2, to=2-1]
        \end{tikzcd}.
      \end{equation*}
  \end{itemize}
  If the laxators and unitors are isomorphisms in $\dbl{E}_1$, the double
  functor is called \define{pseudo}; if they are identities, the double functor
  is \define{strict}. If just the unitors are invertible, then $F$ is said to be
  \define{normal}; if the unitors are strict identities, then $F$ is
  \define{unitary}.
\end{definition}

In addition to the 2-category $\Dbl$ of double categories, double functors, and
natural transformations, there is a 2-category $\DblLax$ of double categories,
\emph{lax} functors, and natural transformations.

Lax functors abound. Here are a few naturally occurring examples.

\begin{example}[Hom functor]
  Given a double category $\dbl{D}$, the Hom double functor
  \begin{equation*}
    \dbl{D}(-,-) \coloneqq \Hom_{\dbl{D}}: \dbl{D}^\op \times \dbl{D} \to \Span,
  \end{equation*}
  as well as the representable double functors $\dbl{D}(x,-): \dbl{D} \to \Span$
  and $\dbl{D}(-,y): \dbl{D}^\op \to \Span$, are all in general lax
  \cite[\S{2.1}]{pare2011}.
\end{example}

\begin{example}[Ob functor] \label{ex:ob-functor}
  The forgetful functor $\Ob: \cat{Cat} \to \Set$ that extracts a category's set
  of objects upgrades to a lax double functor $\Ob\colon \Prof \to \Span$ that
  sends a profunctor to the span having the profunctor's set of heteromorphisms
  as its apex and the heteromorphism source and target maps as its legs
  \cite[\S{1.2}]{pare2011}.
\end{example}

Before we see how lax functors give rise to categories, functors, and natural
transformations, we need to know how to interpret these concepts inside any
double category. Category objects in a double category have been variously
called \emph{monoids} \cite{shulman2008} and \emph{monads} \cite{fiore2011} in
the double category. We prefer to call them simply \emph{categories},
emphasizing that double categories are an appropriate categorified structure in
which to interpret categories, just as monoidal categories and symmetric
monoidal categories are for monoids and commutative monoids. We also wish to
avoid confusion with the double theories of monoids and monads introduced later,
which give different ways to interpret monoids and monads inside a double
category.

\begin{definition}[Category object] \label{def:category-object}
  Let $\dbl{D}$ be a double category.
  \begin{enumerate}[(i)]
    \item A \define{category object}, or simply a \define{category}, in
      $\dbl{D}$ consists of an object $x \in \dbl{D}$, a proarrow
      $r: x \proto x$, and cells
      \begin{equation*}
        \begin{tikzcd}
          x & x & x \\
          x && x
          \arrow["r", "\shortmid"{marking}, from=1-1, to=1-2]
          \arrow["r", "\shortmid"{marking}, from=1-2, to=1-3]
          \arrow[""{name=0, anchor=center, inner sep=0}, "r"', "\shortmid"{marking}, from=2-1, to=2-3]
          \arrow[Rightarrow, no head, from=1-1, to=2-1]
          \arrow[Rightarrow, no head, from=1-3, to=2-3]
          \arrow["\mu"{description}, draw=none, from=1-2, to=0]
        \end{tikzcd}
        \qquad\text{and}\qquad
        \begin{tikzcd}
          x & x \\
          x & x
          \arrow[""{name=0, anchor=center, inner sep=0}, "r"', "\shortmid"{marking}, from=2-1, to=2-2]
          \arrow[""{name=1, anchor=center, inner sep=0}, "{\mathrm{id}_x}", "\shortmid"{marking}, from=1-1, to=1-2]
          \arrow[Rightarrow, no head, from=1-1, to=2-1]
          \arrow[Rightarrow, no head, from=1-2, to=2-2]
          \arrow["\eta"{description}, draw=none, from=1, to=0]
        \end{tikzcd}
      \end{equation*}
      satisfying the usual associativity and unitality axioms.
    \item A \define{functor} in $\dbl{D}$ from one category object
      $(x,r,\mu,\eta)$ to another $(y,s,\nu,\theta)$ consists of an arrow
      $f: x \to y$ along with a cell $\inlineCell{x}{x}{y}{y}{r}{s}{f}{f}{\phi}$
      that preserves composition and units:
      \begin{equation*}
        \begin{tikzcd}
          x & x & x \\
          y & y & y \\
          y && y
          \arrow[""{name=0, anchor=center, inner sep=0}, "r", "\shortmid"{marking}, from=1-1, to=1-2]
          \arrow[""{name=1, anchor=center, inner sep=0}, "r", "\shortmid"{marking}, from=1-2, to=1-3]
          \arrow["f"', from=1-1, to=2-1]
          \arrow["f", from=1-3, to=2-3]
          \arrow["f"{description}, from=1-2, to=2-2]
          \arrow[""{name=2, anchor=center, inner sep=0}, "s"', "\shortmid"{marking}, from=2-1, to=2-2]
          \arrow[""{name=3, anchor=center, inner sep=0}, "s"', "\shortmid"{marking}, from=2-2, to=2-3]
          \arrow[Rightarrow, no head, from=2-1, to=3-1]
          \arrow[Rightarrow, no head, from=2-3, to=3-3]
          \arrow[""{name=4, anchor=center, inner sep=0}, "s"', "\shortmid"{marking}, from=3-1, to=3-3]
          \arrow["\phi"{description}, draw=none, from=0, to=2]
          \arrow["\phi"{description}, draw=none, from=1, to=3]
          \arrow["\nu"{description}, draw=none, from=2-2, to=4]
        \end{tikzcd}
        =
        \begin{tikzcd}
          x & x & x \\
          x && x \\
          y && y
          \arrow["r", "\shortmid"{marking}, from=1-1, to=1-2]
          \arrow["r", "\shortmid"{marking}, from=1-2, to=1-3]
          \arrow[Rightarrow, no head, from=1-1, to=2-1]
          \arrow[Rightarrow, no head, from=1-3, to=2-3]
          \arrow[""{name=0, anchor=center, inner sep=0}, "s"', "\shortmid"{marking}, from=3-1, to=3-3]
          \arrow[""{name=1, anchor=center, inner sep=0}, "r"', "\shortmid"{marking}, from=2-1, to=2-3]
          \arrow["f"', from=2-1, to=3-1]
          \arrow["f", from=2-3, to=3-3]
          \arrow["\phi"{description}, draw=none, from=1, to=0]
          \arrow["\mu"{description}, draw=none, from=1-2, to=1]
        \end{tikzcd}
      \end{equation*}
      and
      \begin{equation*}
        \begin{tikzcd}
          x & x \\
          y & y \\
          y & y
          \arrow[""{name=0, anchor=center, inner sep=0}, "s"', "\shortmid"{marking}, from=3-1, to=3-2]
          \arrow[""{name=1, anchor=center, inner sep=0}, "{\mathrm{id}_y}"', "\shortmid"{marking}, from=2-1, to=2-2]
          \arrow[Rightarrow, no head, from=2-1, to=3-1]
          \arrow[Rightarrow, no head, from=2-2, to=3-2]
          \arrow[""{name=2, anchor=center, inner sep=0}, "{\mathrm{id}_x}", "\shortmid"{marking}, from=1-1, to=1-2]
          \arrow["f"', from=1-1, to=2-1]
          \arrow["f", from=1-2, to=2-2]
          \arrow["\theta"{description, pos=0.6}, draw=none, from=1, to=0]
          \arrow["{\mathrm{id}_f}"{description}, draw=none, from=2, to=1]
        \end{tikzcd}
        =
        \begin{tikzcd}
          x & x \\
          x & x \\
          y & y
          \arrow["f"', from=2-1, to=3-1]
          \arrow["f", from=2-2, to=3-2]
          \arrow[""{name=0, anchor=center, inner sep=0}, "s"', "\shortmid"{marking}, from=3-1, to=3-2]
          \arrow[""{name=1, anchor=center, inner sep=0}, "r"', "\shortmid"{marking}, from=2-1, to=2-2]
          \arrow[""{name=2, anchor=center, inner sep=0}, "{\mathrm{id}_x}", "\shortmid"{marking}, from=1-1, to=1-2]
          \arrow[Rightarrow, no head, from=1-1, to=2-1]
          \arrow[Rightarrow, no head, from=1-2, to=2-2]
          \arrow["\phi"{description}, draw=none, from=1, to=0]
          \arrow["\eta"{description}, draw=none, from=2, to=1]
        \end{tikzcd}.
      \end{equation*}
    \item A \define{natural transformation} in $\dbl{D}$ from one functor
      $(f,\phi)$ to another $(g,\psi)$ with the same domain and codomain is a
      cell $\inlineCell{x}{x}{y}{y}{\id_x}{s}{f}{g}{\alpha}$ satisfying the
      naturality axiom:
      \begin{equation} \label{eq:naturality-in-dbl}
        \begin{tikzcd}
          x & x & x \\
          y & y & y \\
          y && y
          \arrow["f"', from=1-1, to=2-1]
          \arrow["g"{description}, from=1-2, to=2-2]
          \arrow[""{name=0, anchor=center, inner sep=0}, "s"', "\shortmid"{marking}, from=2-1, to=2-2]
          \arrow["g", from=1-3, to=2-3]
          \arrow[""{name=1, anchor=center, inner sep=0}, "s"', "\shortmid"{marking}, from=2-2, to=2-3]
          \arrow[""{name=2, anchor=center, inner sep=0}, "s"', "\shortmid"{marking}, from=3-1, to=3-3]
          \arrow[Rightarrow, no head, from=2-3, to=3-3]
          \arrow[Rightarrow, no head, from=2-1, to=3-1]
          \arrow[""{name=3, anchor=center, inner sep=0}, "{\mathrm{id}_x}", "\shortmid"{marking}, from=1-1, to=1-2]
          \arrow[""{name=4, anchor=center, inner sep=0}, "r", "\shortmid"{marking}, from=1-2, to=1-3]
          \arrow["\alpha"{description}, draw=none, from=3, to=0]
          \arrow["\psi"{description}, draw=none, from=4, to=1]
          \arrow["\nu"{description}, draw=none, from=2-2, to=2]
        \end{tikzcd}
        \quad=\quad
        \begin{tikzcd}
          x & x & x \\
          y & y & y \\
          y && y
          \arrow["f"', from=1-1, to=2-1]
          \arrow["f"{description}, from=1-2, to=2-2]
          \arrow[""{name=0, anchor=center, inner sep=0}, "s"', "\shortmid"{marking}, from=2-1, to=2-2]
          \arrow["g", from=1-3, to=2-3]
          \arrow[""{name=1, anchor=center, inner sep=0}, "s"', "\shortmid"{marking}, from=2-2, to=2-3]
          \arrow[""{name=2, anchor=center, inner sep=0}, "s"', "\shortmid"{marking}, from=3-1, to=3-3]
          \arrow[Rightarrow, no head, from=2-3, to=3-3]
          \arrow[Rightarrow, no head, from=2-1, to=3-1]
          \arrow[""{name=3, anchor=center, inner sep=0}, "r", "\shortmid"{marking}, from=1-1, to=1-2]
          \arrow[""{name=4, anchor=center, inner sep=0}, "{\mathrm{id}_x}", "\shortmid"{marking}, from=1-2, to=1-3]
          \arrow["\alpha"{description}, draw=none, from=4, to=1]
          \arrow["\phi"{description}, draw=none, from=3, to=0]
          \arrow["\nu"{description}, draw=none, from=2-2, to=2]
        \end{tikzcd}.
      \end{equation}
  \end{enumerate}
  Categories, functors, and natural transformations in $\dbl{D}$ form a
  2-category, denoted $\Cat(\dbl{D})$, with vertical and horizontal composition
  defined in the obvious way.
\end{definition}

\begin{remark}[Orientations]
  It may seem backward that vertical and horizontal composition in the
  2-category $\Cat(\dbl{D})$ correspond to horizontal and vertical composition
  in the double category $\dbl{D}$, but this is consistent with our orientation
  convention for double categories, in which 2-categories embed in double
  categories vertically and bicategories embed horizontally. Under this
  convention, a 2-cell $\alpha: f \To g: x \to y$ in a 2-category should be
  depicted as
  $%
  \begin{tikzcd}
    x \\
    y
    \arrow[""{name=0, anchor=center, inner sep=0}, "f"', curve={height=12pt}, from=1-1, to=2-1]
    \arrow[""{name=1, anchor=center, inner sep=0}, "g", curve={height=-12pt}, from=1-1, to=2-1]
    \arrow["\alpha", shorten <=5pt, shorten >=5pt, Rightarrow, from=0, to=1]
  \end{tikzcd}$.
\end{remark}

The two most important instances of category objects are internal and enriched
categories, obtained as category objects in the following two double categories.

\begin{example}[Spans] \label{ex:spans}
  Let $\cat{S}$ be a category with pullbacks. The double category
  $\Span{\cat{S}}$ has as objects and arrows, objects and morphisms in
  $\cat{S}$; as proarrows, spans in $\cat{S}$; and as cells, maps of spans in
  $\cat{S}$. Composition of proarrows is by pullback in $\cat{S}$. For details,
  see \cite[\S{3.4.1}]{grandis2019}.
\end{example}

\begin{example}[Matrices] \label{ex:matrices}
  Let $\catV$ be an (infinitary) distributive monoidal category: a monoidal
  category with coproducts over which the monoidal product distributes. The
  double category $\Mat{\catV}$ has as objects, sets; as arrows, functions; as
  proarrows $X \proto Y$, functions $M: X \times Y \to \catV_0$, the $X$-by-$Y$
  \define{$\catV$-matrices}; and as cells
  $\inlineCell{X}{Y}{W}{Z}{M}{N}{f}{g}{\alpha}$, families of morphisms in
  $\catV$
  \begin{equation*}
    \alpha_{x,y}: M(x,y) \to N(f(x), g(y)), \qquad x \in X,\ y \in Y.
  \end{equation*}
  For composition of $\catV$-matrices and properties of $\Mat{\catV}$, see
  \cite[\S{4}]{vasilakopoulou2019}.
\end{example}

The double categories of spans and of matrices are equivalent in the fundamental
case of spans of sets and set-valued matrices:
\begin{equation*}
  \Span{\Set} \eqqcolon \Span \simeq \Mat \coloneqq \Mat{\Set}.
\end{equation*}
In terms of the following proposition, this equivalence means that category
objects in both $\Span$ and $\Mat$ can be identified with ordinary categories.

\begin{proposition}[Internal and enriched categories]
  For any category $\cat{S}$ with pullbacks, we have
  \begin{equation*}
   \Cat(\Span{\cat{S}}) = \Cat(\cat{S}),
  \end{equation*}
  the 2-category of internal categories, internal functors, and internal natural
  transformations in $\cat{S}$.

  Similarly, for any distributive monoidal category $\catV$,
  \begin{equation*}
    \Cat(\Mat{\catV}) = \EnCat{\catV},
  \end{equation*}
  the 2-category of $\catV$-categories, $\catV$-functors, and $\catV$-natural
  transformations.
\end{proposition}
\begin{proof}
  We prove the statement for internal categories; the case of enriched
  categories is similar. For the main definitions of internal category theory,
  see \cite[Vol 1, Chapter 8]{borceux1994}.
  \begin{enumerate}[(i)]
    \item A category object in $\Span{\cat{S}}$ consists of a span
      $(X_0 \xproto{R} X_0) = (X_0 \xleftarrow{s} X_1 \xrightarrow{t} X_0)$
      together with maps of spans
      \begin{equation*}
        \begin{tikzcd}
          {X_0} & {X_0} & {X_0} \\
          {X_0} && {X_0}
          \arrow["R", "\shortmid"{marking}, from=1-1, to=1-2]
          \arrow["R", "\shortmid"{marking}, from=1-2, to=1-3]
          \arrow[""{name=0, anchor=center, inner sep=0}, "R"', "\shortmid"{marking}, from=2-1, to=2-3]
          \arrow[Rightarrow, no head, from=1-1, to=2-1]
          \arrow[Rightarrow, no head, from=1-3, to=2-3]
          \arrow["\mu"{description}, draw=none, from=1-2, to=0]
        \end{tikzcd}
        \quad=\quad
        \begin{tikzcd}
          {X_0} & {X_1 \times_{X_0} X_1} & {X_0} \\
          {X_0} & {X_1} & {X_0}
          \arrow["{s \circ \pi_1}"', from=1-2, to=1-1]
          \arrow["{t \circ \pi_2}", from=1-2, to=1-3]
          \arrow["c"', from=1-2, to=2-2]
          \arrow[Rightarrow, no head, from=1-1, to=2-1]
          \arrow[Rightarrow, no head, from=1-3, to=2-3]
          \arrow["s", from=2-2, to=2-1]
          \arrow["t"', from=2-2, to=2-3]
        \end{tikzcd}
      \end{equation*}
      and
      \begin{equation*}
        \begin{tikzcd}
          {X_0} & {X_0} \\
          {X_0} & {X_0}
          \arrow[""{name=0, anchor=center, inner sep=0}, "R"', "\shortmid"{marking}, from=2-1, to=2-2]
          \arrow[Rightarrow, no head, from=1-1, to=2-1]
          \arrow[Rightarrow, no head, from=1-2, to=2-2]
          \arrow[""{name=1, anchor=center, inner sep=0}, "{\mathrm{id}_{X_0}}", "\shortmid"{marking}, from=1-1, to=1-2]
          \arrow["\eta"{description}, draw=none, from=1, to=0]
        \end{tikzcd}
        \quad=\quad
        \begin{tikzcd}
          {X_0} & {X_0} & {X_0} \\
          {X_0} & {X_1} & {X_0}
          \arrow["s", from=2-2, to=2-1]
          \arrow["t"', from=2-2, to=2-3]
          \arrow[Rightarrow, no head, from=1-2, to=1-1]
          \arrow[Rightarrow, no head, from=1-1, to=2-1]
          \arrow["i"', from=1-2, to=2-2]
          \arrow[Rightarrow, no head, from=1-2, to=1-3]
          \arrow[Rightarrow, no head, from=1-3, to=2-3]
        \end{tikzcd}
      \end{equation*}
      satisfying associativity and unitality axioms, which is precisely a
      category internal to $\cat{S}$.
    \item A functor in $\Span{\cat{S}}$, or equivalently a functor internal to
      $\cat{S}$, consists of a map $f: X_0 \to Y_0$ in $\cat{S}$, the object
      map, and a map of spans
      \begin{equation*}
        \begin{tikzcd}
          {X_0} & {X_1} & {X_0} \\
          {Y_0} & {Y_1} & {Y_0}
          \arrow["{f_0}"', from=1-1, to=2-1]
          \arrow["s", from=2-2, to=2-1]
          \arrow["s"', from=1-2, to=1-1]
          \arrow["t"', from=2-2, to=2-3]
          \arrow["t", from=1-2, to=1-3]
          \arrow["{f_0}", from=1-3, to=2-3]
          \arrow["{f_1}"', from=1-2, to=2-2]
        \end{tikzcd},
      \end{equation*}
      the morphism map, satisfying the functor axioms.
    \item A natural transformation in $\Span{\cat{S}}$ consists of a map of
      spans
      \begin{equation*}
        \begin{tikzcd}
          {X_0} & {X_0} \\
          {Y_0} & {Y_0}
          \arrow["{f_0}"', from=1-1, to=2-1]
          \arrow["{g_0}", from=1-2, to=2-2]
          \arrow[""{name=0, anchor=center, inner sep=0}, "{\mathrm{id}_{X_0}}", "\shortmid"{marking}, from=1-1, to=1-2]
          \arrow[""{name=1, anchor=center, inner sep=0}, "S"', "\shortmid"{marking}, from=2-1, to=2-2]
          \arrow["\alpha"{description}, draw=none, from=0, to=1]
        \end{tikzcd}
        \quad=\quad
        \begin{tikzcd}
          {X_0} & {X_0} & {X_0} \\
          {Y_0} & {Y_1} & {Y_0}
          \arrow["s", from=2-2, to=2-1]
          \arrow["t"', from=2-2, to=2-3]
          \arrow["{f_0}"', from=1-1, to=2-1]
          \arrow["\alpha"', from=1-2, to=2-2]
          \arrow[Rightarrow, no head, from=1-2, to=1-1]
          \arrow["{g_0}", from=1-3, to=2-3]
          \arrow[Rightarrow, no head, from=1-2, to=1-3]
        \end{tikzcd}
      \end{equation*}
      satisfying the naturality axiom, which is precisely a natural
      transformation $\alpha: X_0 \to Y_1$ internal to $\cat{S}$. \qedhere
  \end{enumerate}
\end{proof}

The above definition of a natural transformation inside a double category is the
obvious generalization of the usual concept, but we will see that a different
formulation arises more immediately from lax double functors.

\begin{proposition} \label{prop:natural-transformation-in-dbl}
  A natural transformation $\inlineCell{x}{x}{y}{y}{\id_x}{s}{f}{g}{\alpha}$ in
  a double category $\dbl{D}$ (\cref{def:category-object}) is a equivalent to a
  cell $\inlineCell{x}{x}{y}{y}{r}{s}{f}{g}{\bar\alpha}$ in $\dbl{D}$ satisfying
  the equations:
  \begin{equation} \label{eq:naturality-in-dbl-alt}
    \begin{tikzcd}
      x & x & x \\
      y & y & y \\
      y && y
      \arrow["f"', from=1-1, to=2-1]
      \arrow["g"{description}, from=1-2, to=2-2]
      \arrow[""{name=0, anchor=center, inner sep=0}, "s"', "\shortmid"{marking}, from=2-1, to=2-2]
      \arrow["g", from=1-3, to=2-3]
      \arrow[""{name=1, anchor=center, inner sep=0}, "s"', "\shortmid"{marking}, from=2-2, to=2-3]
      \arrow[""{name=2, anchor=center, inner sep=0}, "s"', "\shortmid"{marking}, from=3-1, to=3-3]
      \arrow[Rightarrow, no head, from=2-3, to=3-3]
      \arrow[Rightarrow, no head, from=2-1, to=3-1]
      \arrow[""{name=3, anchor=center, inner sep=0}, "r", "\shortmid"{marking}, from=1-1, to=1-2]
      \arrow[""{name=4, anchor=center, inner sep=0}, "r", "\shortmid"{marking}, from=1-2, to=1-3]
      \arrow["\bar\alpha"{description}, draw=none, from=3, to=0]
      \arrow["\psi"{description}, draw=none, from=4, to=1]
      \arrow["\nu"{description}, draw=none, from=2-2, to=2]
    \end{tikzcd}
    \quad=\quad
    \begin{tikzcd}
      x & x & x \\
      x && x \\
      y && y
      \arrow["r", "\shortmid"{marking}, from=1-1, to=1-2]
      \arrow["r", "\shortmid"{marking}, from=1-2, to=1-3]
      \arrow[""{name=0, anchor=center, inner sep=0}, "r"', "\shortmid"{marking}, from=2-1, to=2-3]
      \arrow[Rightarrow, no head, from=1-3, to=2-3]
      \arrow[Rightarrow, no head, from=1-1, to=2-1]
      \arrow[""{name=1, anchor=center, inner sep=0}, "s"', "\shortmid"{marking}, from=3-1, to=3-3]
      \arrow["f"', from=2-1, to=3-1]
      \arrow["g", from=2-3, to=3-3]
      \arrow["\bar\alpha"{description}, draw=none, from=0, to=1]
      \arrow["\mu"{description}, draw=none, from=1-2, to=0]
    \end{tikzcd}
    \quad=\quad
    \begin{tikzcd}
      x & x & x \\
      y & y & y \\
      y && y
      \arrow["f"', from=1-1, to=2-1]
      \arrow["f"{description}, from=1-2, to=2-2]
      \arrow[""{name=0, anchor=center, inner sep=0}, "s"', "\shortmid"{marking}, from=2-1, to=2-2]
      \arrow["g", from=1-3, to=2-3]
      \arrow[""{name=1, anchor=center, inner sep=0}, "s"', "\shortmid"{marking}, from=2-2, to=2-3]
      \arrow[""{name=2, anchor=center, inner sep=0}, "s"', "\shortmid"{marking}, from=3-1, to=3-3]
      \arrow[Rightarrow, no head, from=2-3, to=3-3]
      \arrow[Rightarrow, no head, from=2-1, to=3-1]
      \arrow[""{name=3, anchor=center, inner sep=0}, "r", "\shortmid"{marking}, from=1-1, to=1-2]
      \arrow[""{name=4, anchor=center, inner sep=0}, "r", "\shortmid"{marking}, from=1-2, to=1-3]
      \arrow["\bar\alpha"{description}, draw=none, from=4, to=1]
      \arrow["\phi"{description}, draw=none, from=3, to=0]
      \arrow["\nu"{description}, draw=none, from=2-2, to=2]
    \end{tikzcd}.
  \end{equation}
\end{proposition}
\begin{proof}
  Given a natural transformation $\alpha$ in $\dbl{D}$, define the cell
  $\bar\alpha$ by either side of \cref{eq:naturality-in-dbl}. To prove, for
  example, the left equality in \cref{eq:naturality-in-dbl-alt}, we calculate
  \begin{equation*}
    \begin{dblArray}{cc}
      \bar\alpha & \psi \\
      \Block{1-2}{\nu} &
    \end{dblArray} =
    \begin{dblArray}{ccc}
      \alpha & \psi & \psi \\
      \Block{1-2}{\nu} & & 1_s \\
      \Block{1-3}{\nu} & &
    \end{dblArray} =
    \begin{dblArray}{ccc}
      \alpha & \psi & \psi \\
      1_s & \Block{1-2}{\nu} & \\
      \Block{1-3}{\nu} & &
    \end{dblArray} =
    \begin{dblArray}{cc}
      1_{\id_x} & \mu \\
      \alpha & \psi \\
      \Block{1-2}{\nu} &
    \end{dblArray} =
    \begin{dblArray}{cc}
      \id_{1_x} & \mu \\
      \Block{1-2}{\bar\alpha} &
    \end{dblArray} =
    \begin{dblArray}{c}
      \mu \\ \bar\alpha
    \end{dblArray}.
  \end{equation*}
  Conversely, given a cell $\bar\alpha$, define the natural transformation $\alpha$ by
  \begin{equation*}
    \begin{tikzcd}
      x & x \\
      y & y
      \arrow[""{name=0, anchor=center, inner sep=0}, "{\mathrm{id}_x}", "\shortmid"{marking}, from=1-1, to=1-2]
      \arrow["f"', from=1-1, to=2-1]
      \arrow["g", from=1-2, to=2-2]
      \arrow[""{name=1, anchor=center, inner sep=0}, "s"', "\shortmid"{marking}, from=2-1, to=2-2]
      \arrow["\alpha"{description}, draw=none, from=0, to=1]
    \end{tikzcd}
    \quad\coloneqq\quad
    \begin{tikzcd}
      x & x \\
      x & x \\
      y & y
      \arrow[""{name=0, anchor=center, inner sep=0}, "r", "\shortmid"{marking}, from=2-1, to=2-2]
      \arrow["f"', from=2-1, to=3-1]
      \arrow["g", from=2-2, to=3-2]
      \arrow[""{name=1, anchor=center, inner sep=0}, "s"', "\shortmid"{marking}, from=3-1, to=3-2]
      \arrow[""{name=2, anchor=center, inner sep=0}, "{\mathrm{id}_x}", "\shortmid"{marking}, from=1-1, to=1-2]
      \arrow[Rightarrow, no head, from=1-1, to=2-1]
      \arrow[Rightarrow, no head, from=1-2, to=2-2]
      \arrow["\bar\alpha"{description}, draw=none, from=0, to=1]
      \arrow["\eta"{description}, draw=none, from=2, to=0]
    \end{tikzcd}.
  \end{equation*}
  Composing on the top of \cref{eq:naturality-in-dbl-alt} with $\eta \odot 1_r$
  and $1_r \odot \eta$ yields \cref{eq:naturality-in-dbl}.

  These operations put the two types of cells $\alpha$ and $\bar\alpha$ in
  bijective correspondence. On the one hand, since $\psi$ and $\nu$ are
  unit-preserving, we have
  \begin{equation*}
    \begin{dblArray}{cc}
      \Block{1-2}{\eta} & \\
      \Block{1-2}{\cong} & \\
      \alpha & \psi \\
      \Block{1-2}{\nu} &
    \end{dblArray} =
    \begin{dblArray}{cc}
      \Block{1-2}{\cong} & \\
      1 & \eta \\
      \alpha & \psi \\
      \Block{1-2}{\nu} &
    \end{dblArray} =
    \begin{dblArray}{cc}
      \Block{1-2}{\cong} & \\
      \alpha & \theta \\
      \Block{1-2}{\nu} &
    \end{dblArray} =
    \begin{dblArray}{c}
      \alpha
    \end{dblArray}.
  \end{equation*}
  The isomorphisms in the computation are the canonical ones. On the other hand, we have
  \begin{equation*}
    \begin{dblArray}{cc}
      \eta & 1 \\
      \bar\alpha & \psi \\
      \Block{1-2}{\nu} &
    \end{dblArray} =
    \begin{dblArray}{cc}
      \eta & 1 \\
      \Block{1-2}{\mu} & \\
      \Block{1-2}{\bar\alpha} &
    \end{dblArray} =
    \begin{dblArray}{c}
      \bar\alpha
    \end{dblArray}
  \end{equation*}
  by \cref{eq:naturality-in-dbl-alt} and the fact that $\mu$ is unit-preserving.
\end{proof}

In the fundamental situation that $\dbl{D} = \Mat \simeq \Span$, the result says
that a natural transformation $\alpha: F \To G: \cat{C} \to \cat{D}$ is
equivalent to a family of maps $\bar\alpha_f: Fx \to Gy$ in $\cat{D}$ indexed by
maps $f: x \to y$ in $\cat{C}$, such that
$Ff \cdot \bar\alpha_g = \bar\alpha_{fg} = \bar\alpha_f \cdot Gg$ for all
consecutive maps $x \xto{f} y \xto{g} z$, i.e., to a natural transformation
\begin{equation*}
  \bar\alpha: \Hom_{\cat{C}} \To \Hom_{\cat{D}} \circ (F^\op \times G):
    \cat{C}^\op \times \cat{C} \to \Set.
\end{equation*}

We can now state our first important fact about lax functors.

\begin{proposition}[Lax functors give category objects]
  \label{prop:lax-functor-category-objects}
  Let $\dbl{D}$ be a strict double category and let $F: \dbl{D} \to \dbl{E}$ be
  a lax double functor.
  \begin{enumerate}[(i)]
    \item For every object $x \in \dbl{D}$, the data
      $(Fx, F\id_x, F_{x,x}, F_x)$, where
      $F_{x,x}: F\id_x \odot F\id_x \to F\id_x$ is the laxator and
      $F_x: \id_{Fx} \to F\id_x$ is the unitor, is a category object in
      $\dbl{E}$.
    \item For every arrow $f: x \to y$ in $\dbl{D}$, the arrow $Ff: Fx \to Fy$
      together with the cell $F\id_f: F\id_x \to F\id_y$ is a functor in
      $\dbl{E}$.
    \item For every cell of the form
      $\inlineCell{x}{x}{y}{y}{\id_x}{\id_y}{f}{g}{\alpha}$ in $\dbl{D}$, the
      cell $F\alpha: F\id_x \to F\id_y$ is a natural transformation in
      $\dbl{E}$.
  \end{enumerate}
  Moreover, these assignments are functorial with respect to composition in
  $\dbl{D}_0$ and $\dbl{D}_1$.
\end{proposition}
\begin{proof}
  \begin{enumerate}[(i)]
    \item Associativity and unitality of the category object are precisely the
      associativity and unitality axioms for the laxator $F_{x,x}$ and unitor
      $F_x$ of the lax functor $F$.
    \item Functorality of the data $(Ff, F\id_f)$ is precisely the naturality
      axioms $(F\id_f \odot F\id_f) \cdot F_{y,y} = F_{x,x} \cdot F\id_f$ and
      $\id_{Ff} \cdot F_y = F_x \cdot F\id_f$ of the lax functor $F$.
    \item Given a cell $\inlineCell{x}{x}{y}{y}{\id_x}{\id_y}{f}{g}{\alpha}$ in
      $\dbl{D}$, the equation $\id_f \odot \alpha = \alpha = \alpha \odot \id_g$ implies that
      \begin{align*}
        \begin{tikzcd}[ampersand replacement=\&]
          Fx \& Fx \& Fx \\
          Fy \& Fy \& Fy \\
          Fy \&\& Fy
          \arrow[""{name=0, anchor=center, inner sep=0}, "{F \mathrm{id}_x}", "\shortmid"{marking}, from=1-1, to=1-2]
          \arrow[""{name=1, anchor=center, inner sep=0}, "{F \mathrm{id}_x}", "\shortmid"{marking}, from=1-2, to=1-3]
          \arrow["Ff"{description}, from=1-2, to=2-2]
          \arrow["Ff"', from=1-1, to=2-1]
          \arrow["Fg", from=1-3, to=2-3]
          \arrow[""{name=2, anchor=center, inner sep=0}, "{F \mathrm{id}_y}"', "\shortmid"{marking}, from=2-2, to=2-3]
          \arrow[""{name=3, anchor=center, inner sep=0}, "{F \mathrm{id}_y}"', "\shortmid"{marking}, from=2-1, to=2-2]
          \arrow[""{name=4, anchor=center, inner sep=0}, "{F \mathrm{id}_y}"', "\shortmid"{marking}, from=3-1, to=3-3]
          \arrow[Rightarrow, no head, from=2-1, to=3-1]
          \arrow[Rightarrow, no head, from=2-3, to=3-3]
          \arrow["F\alpha"{description}, draw=none, from=1, to=2]
          \arrow["{F_{y,y}}"{description}, draw=none, from=2-2, to=4]
          \arrow["{F \mathrm{id}_f}"{description}, draw=none, from=0, to=3]
        \end{tikzcd}
        &=
        \begin{tikzcd}[ampersand replacement=\&]
          Fx \& Fx \& Fx \\
          Fx \&\& Fx \\
          Fy \&\& Fy
          \arrow[""{name=0, anchor=center, inner sep=0}, "{F \mathrm{id}_y}"', "\shortmid"{marking}, from=3-1, to=3-3]
          \arrow[""{name=1, anchor=center, inner sep=0}, "{F \mathrm{id}_x}"', "\shortmid"{marking}, from=2-1, to=2-3]
          \arrow["Ff"', from=2-1, to=3-1]
          \arrow["Fg", from=2-3, to=3-3]
          \arrow["{F \mathrm{id}_x}", "\shortmid"{marking}, from=1-1, to=1-2]
          \arrow["{F \mathrm{id}_x}", "\shortmid"{marking}, from=1-2, to=1-3]
          \arrow[Rightarrow, no head, from=1-1, to=2-1]
          \arrow[Rightarrow, no head, from=1-3, to=2-3]
          \arrow["F\alpha"{description}, draw=none, from=1, to=0]
          \arrow["{F_{x,x}}"{description}, draw=none, from=1-2, to=1]
        \end{tikzcd}
        \\
        &=
        \begin{tikzcd}[ampersand replacement=\&]
          Fx \& Fx \& Fx \\
          Fy \& Fy \& Fy \\
          Fy \&\& Fy
          \arrow[""{name=0, anchor=center, inner sep=0}, "{F \mathrm{id}_x}", "\shortmid"{marking}, from=1-1, to=1-2]
          \arrow[""{name=1, anchor=center, inner sep=0}, "{F \mathrm{id}_x}", "\shortmid"{marking}, from=1-2, to=1-3]
          \arrow["Fg"{description}, from=1-2, to=2-2]
          \arrow["Ff"', from=1-1, to=2-1]
          \arrow["Fg", from=1-3, to=2-3]
          \arrow[""{name=2, anchor=center, inner sep=0}, "{F \mathrm{id}_y}"', "\shortmid"{marking}, from=2-2, to=2-3]
          \arrow[""{name=3, anchor=center, inner sep=0}, "{F \mathrm{id}_y}"', "\shortmid"{marking}, from=2-1, to=2-2]
          \arrow[""{name=4, anchor=center, inner sep=0}, "{F \mathrm{id}_y}"', "\shortmid"{marking}, from=3-1, to=3-3]
          \arrow[Rightarrow, no head, from=2-1, to=3-1]
          \arrow[Rightarrow, no head, from=2-3, to=3-3]
          \arrow["{F \mathrm{id}_g}"{description}, draw=none, from=1, to=2]
          \arrow["{F_{y,y}}"{description}, draw=none, from=2-2, to=4]
          \arrow["{F \alpha}"{description}, draw=none, from=0, to=3]
        \end{tikzcd}.
      \end{align*}
      This is precisely \cref{eq:naturality-in-dbl-alt} characterizing natural
      transformations, hence we obtain a natural transformation in $\dbl{E}$ via
      the composite $F_x \cdot F\alpha: \id_x \to F\id_y$. \qedhere
  \end{enumerate}
\end{proof}

There is a forgetful 2-functor $\VerTwoCat: \Dbl \to \TwoCat$ that sends a
double category to its underlying or vertical 2-category, whose cells are the
cells of the double category bounded by identity proarrows. In the other
direction, there is a 2-functor $\VerDbl: \TwoCat \to \Dbl$ sends a 2-category
$\bicat{B}$ to the double category $\VerDbl(\bicat{B})$ with underlying
2-category $\bicat{B}$ and only trivial proarrows.

\begin{corollary} \label{cor:lax-double-fs-corresp-to-2-fs}
  For any 2-category $\bicat{B}$ and double category $\dbl{D}$, a lax double
  functor $\VerDbl(\bicat{B}) \to \dbl{D}$ is the same thing as a 2-functor
  $\bicat{B} \to \Cat(\dbl{D})$. Likewise, a \emph{unitary} lax double functor 
  $\VerDbl(\bicat{B}) \to \dbl{D}$ is the same as a 2-functor $\bicat{B} \to 
  \bicat{V}(\dbl{D})$.
\end{corollary}

Expanding on this corollary, \emph{normal} lax double functors
$\VerDbl(\bicat{B})\to\dbl{D}$ are in pseudo-inverse correspondence with
2-functors $\bicat{B}\to\Cat(\dbl{D})$. This observation is the basis of the
equivalence that we will examine in
\cref{cor:characterize-unitary-models-on-veritical-dbl-cat-I}.

\begin{corollary} \label{cor:lax-double-fs-corresp-to-2-fs-upto-iso}
  For any 2-category $\bicat{B}$ and double category $\dbl{D}$, a \emph{normal}
  lax double functor $\VerDbl(\bicat{B}) \to \dbl{D}$ is up to invertible
  transformation the same thing as a 2-functor $\bicat{B} \to \bicat{V}(\dbl{D})$.
\end{corollary}

We have yet to consider how lax functors act on non-identity proarrows. It turns
out that they send proarrows to \emph{profunctors} or \emph{bimodules}, which
can be defined inside any double category \cite[\S{11}]{shulman2008}.

\begin{definition}[Profunctor object] \label{def:profunctor-object}
  Let $\dbl{D}$ be a double category.
  \begin{enumerate}[(i)]
    \item A \define{profunctor} or \define{bimodule} between categories
      $(x,r,\mu,\eta)$ and $(y,s,\nu,\theta)$ in $\dbl{D}$ consists of a
      proarrow $m: x \proto y$ and cells
      \begin{equation*}
        \begin{tikzcd}
          x & x & y \\
          x && y
          \arrow["r", "\shortmid"{marking}, from=1-1, to=1-2]
          \arrow["m", "\shortmid"{marking}, from=1-2, to=1-3]
          \arrow[""{name=0, anchor=center, inner sep=0}, "m"', "\shortmid"{marking}, from=2-1, to=2-3]
          \arrow[Rightarrow, no head, from=1-1, to=2-1]
          \arrow[Rightarrow, no head, from=1-3, to=2-3]
          \arrow["\lambda"{description}, draw=none, from=1-2, to=0]
        \end{tikzcd}
        \qquad\text{and}\qquad
        \begin{tikzcd}
          x & y & y \\
          x && y
          \arrow[""{name=0, anchor=center, inner sep=0}, "m"', "\shortmid"{marking}, from=2-1, to=2-3]
          \arrow[Rightarrow, no head, from=1-1, to=2-1]
          \arrow[Rightarrow, no head, from=1-3, to=2-3]
          \arrow["m", "\shortmid"{marking}, from=1-1, to=1-2]
          \arrow["s", "\shortmid"{marking}, from=1-2, to=1-3]
          \arrow["\rho"{description}, draw=none, from=1-2, to=0]
        \end{tikzcd}
      \end{equation*}
      in $\dbl{D}$, the \define{left} and \define{right actions}, satisfying the
      associativity and unitality axioms
      \begin{equation*}
        \begin{tikzcd}[sep=small]
          x & x & x & y \\
          x && x & y \\
          x &&& y
          \arrow[""{name=0, anchor=center, inner sep=0}, "r"', "\shortmid"{marking}, from=2-1, to=2-3]
          \arrow[""{name=1, anchor=center, inner sep=0}, "m"', "\shortmid"{marking}, from=2-3, to=2-4]
          \arrow[""{name=2, anchor=center, inner sep=0}, "m"', "\shortmid"{marking}, from=3-1, to=3-4]
          \arrow[Rightarrow, no head, from=2-1, to=3-1]
          \arrow[Rightarrow, no head, from=2-4, to=3-4]
          \arrow["r", "\shortmid"{marking}, from=1-1, to=1-2]
          \arrow[""{name=3, anchor=center, inner sep=0}, "r", "\shortmid"{marking}, from=1-2, to=1-3]
          \arrow[Rightarrow, no head, from=1-1, to=2-1]
          \arrow[Rightarrow, no head, from=1-3, to=2-3]
          \arrow[""{name=4, anchor=center, inner sep=0}, "m", "\shortmid"{marking}, from=1-3, to=1-4]
          \arrow[Rightarrow, no head, from=1-4, to=2-4]
          \arrow["{1_m}"{description}, draw=none, from=4, to=1]
          \arrow["\mu"{description}, draw=none, from=1-2, to=0]
          \arrow["\lambda"{description, pos=0.8}, draw=none, from=3, to=2]
        \end{tikzcd}
        =
        \begin{tikzcd}[sep=small]
          x & x & x & y \\
          x & x && y \\
          x &&& y
          \arrow[""{name=0, anchor=center, inner sep=0}, "r"', "\shortmid"{marking}, from=2-1, to=2-2]
          \arrow[""{name=1, anchor=center, inner sep=0}, "m"', "\shortmid"{marking}, from=2-2, to=2-4]
          \arrow[""{name=2, anchor=center, inner sep=0}, "m"', "\shortmid"{marking}, from=3-1, to=3-4]
          \arrow[Rightarrow, no head, from=2-1, to=3-1]
          \arrow[Rightarrow, no head, from=2-4, to=3-4]
          \arrow[Rightarrow, no head, from=1-1, to=2-1]
          \arrow[Rightarrow, no head, from=1-2, to=2-2]
          \arrow[""{name=3, anchor=center, inner sep=0}, "r", "\shortmid"{marking}, from=1-2, to=1-3]
          \arrow[""{name=4, anchor=center, inner sep=0}, "r", "\shortmid"{marking}, from=1-1, to=1-2]
          \arrow["m", "\shortmid"{marking}, from=1-3, to=1-4]
          \arrow[Rightarrow, no head, from=1-4, to=2-4]
          \arrow["{1_r}"{description}, draw=none, from=4, to=0]
          \arrow["\lambda"{description, pos=0.8}, draw=none, from=3, to=2]
          \arrow["\lambda"{description}, draw=none, from=1-3, to=1]
        \end{tikzcd}
        \qquad\text{and}\qquad
        \begin{tikzcd}[sep=small]
          x & x & y \\
          x & x & y \\
          x && y
          \arrow[""{name=0, anchor=center, inner sep=0}, "r"', "\shortmid"{marking}, from=2-1, to=2-2]
          \arrow[""{name=1, anchor=center, inner sep=0}, "m"', "\shortmid"{marking}, from=2-2, to=2-3]
          \arrow[""{name=2, anchor=center, inner sep=0}, "m"', "\shortmid"{marking}, from=3-1, to=3-3]
          \arrow[Rightarrow, no head, from=2-1, to=3-1]
          \arrow[Rightarrow, no head, from=2-3, to=3-3]
          \arrow[Rightarrow, no head, from=1-1, to=2-1]
          \arrow[Rightarrow, no head, from=1-2, to=2-2]
          \arrow[""{name=3, anchor=center, inner sep=0}, "{\id_x}", from=1-1, to=1-2]
          \arrow[Rightarrow, no head, from=1-3, to=2-3]
          \arrow[""{name=4, anchor=center, inner sep=0}, "m", "\shortmid"{marking}, from=1-2, to=1-3]
          \arrow["\lambda"{description}, draw=none, from=2-2, to=2]
          \arrow["\eta"{description}, draw=none, from=3, to=0]
          \arrow["{1_m}"{description}, draw=none, from=4, to=1]
        \end{tikzcd}
        =
        \begin{tikzcd}[sep=small]
          x & y \\
          x & y
          \arrow[""{name=0, anchor=center, inner sep=0}, "m"', "\shortmid"{marking}, from=2-1, to=2-2]
          \arrow[Rightarrow, no head, from=1-1, to=2-1]
          \arrow[Rightarrow, no head, from=1-2, to=2-2]
          \arrow[""{name=1, anchor=center, inner sep=0}, "m", "\shortmid"{marking}, from=1-1, to=1-2]
          \arrow["{1_m}"{description}, draw=none, from=1, to=0]
        \end{tikzcd}
      \end{equation*}
      for the left action and similarly for the right action, as well as the
      compatibility axiom
      \begin{equation*}
        \begin{tikzcd}[sep=small]
          x & x & y & y \\
          x && y & y \\
          x &&& y
          \arrow["r", "\shortmid"{marking}, from=1-1, to=1-2]
          \arrow[""{name=0, anchor=center, inner sep=0}, "m", "\shortmid"{marking}, from=1-2, to=1-3]
          \arrow[""{name=1, anchor=center, inner sep=0}, "m"', "\shortmid"{marking}, from=2-1, to=2-3]
          \arrow[Rightarrow, no head, from=1-1, to=2-1]
          \arrow[Rightarrow, no head, from=1-3, to=2-3]
          \arrow[""{name=2, anchor=center, inner sep=0}, "s", "\shortmid"{marking}, from=1-3, to=1-4]
          \arrow[Rightarrow, no head, from=1-4, to=2-4]
          \arrow[""{name=3, anchor=center, inner sep=0}, "s"', "\shortmid"{marking}, from=2-3, to=2-4]
          \arrow[""{name=4, anchor=center, inner sep=0}, "m"', "\shortmid"{marking}, from=3-1, to=3-4]
          \arrow[Rightarrow, no head, from=2-4, to=3-4]
          \arrow[Rightarrow, no head, from=2-1, to=3-1]
          \arrow["\lambda"{description}, draw=none, from=1-2, to=1]
          \arrow["{1_s}"{description}, draw=none, from=2, to=3]
          \arrow["\rho"{description, pos=0.8}, draw=none, from=0, to=4]
        \end{tikzcd}
        =
        \begin{tikzcd}[sep=small]
          x & x & y & y \\
          x & x && y \\
          x &&& y
          \arrow[""{name=0, anchor=center, inner sep=0}, "m"', "\shortmid"{marking}, from=2-2, to=2-4]
          \arrow[Rightarrow, no head, from=1-2, to=2-2]
          \arrow[Rightarrow, no head, from=1-4, to=2-4]
          \arrow[""{name=1, anchor=center, inner sep=0}, "m", "\shortmid"{marking}, from=1-2, to=1-3]
          \arrow["s", "\shortmid"{marking}, from=1-3, to=1-4]
          \arrow[""{name=2, anchor=center, inner sep=0}, "r", "\shortmid"{marking}, from=1-1, to=1-2]
          \arrow[""{name=3, anchor=center, inner sep=0}, "r"', "\shortmid"{marking}, from=2-1, to=2-2]
          \arrow[Rightarrow, no head, from=1-1, to=2-1]
          \arrow[""{name=4, anchor=center, inner sep=0}, "m"', "\shortmid"{marking}, from=3-1, to=3-4]
          \arrow[Rightarrow, no head, from=2-4, to=3-4]
          \arrow[Rightarrow, no head, from=2-1, to=3-1]
          \arrow["\rho"{description}, draw=none, from=1-3, to=0]
          \arrow["{1_r}"{description}, draw=none, from=2, to=3]
          \arrow["\lambda"{description, pos=0.8}, draw=none, from=1, to=4]
        \end{tikzcd}.
      \end{equation*}
      \item Given functors $(f,\phi): (x,r) \to (w,t)$ and
      $(g,\psi): (y,s) \to (z,u)$ in $\dbl{D}$, a \define{natural
      transformation} or \define{bimodule homomorphism} from a profunctor
      $(m,\lambda^m,\rho^m): (x,r) \proto (y,s)$ to another profunctor
      $(n,\lambda^n,\rho^n): (w,t) \proto (z,u)$, having source $(f,\phi)$ and
      target $(g,\psi)$, is a cell
      \begin{equation*}
        \begin{tikzcd}
          x & y \\
          w & z
          \arrow["f"', from=1-1, to=2-1]
          \arrow["g", from=1-2, to=2-2]
          \arrow[""{name=0, anchor=center, inner sep=0}, "m", "\shortmid"{marking}, from=1-1, to=1-2]
          \arrow[""{name=1, anchor=center, inner sep=0}, "n"', "\shortmid"{marking}, from=2-1, to=2-2]
          \arrow["\alpha"{description}, draw=none, from=0, to=1]
        \end{tikzcd}
      \end{equation*}
      in $\dbl{D}$ satisfying the equivariance axioms
      \begin{equation*}
        \begin{tikzcd}[sep=small]
          x & x & y \\
          w & w & z \\
          w && z
          \arrow["f"', from=1-2, to=2-2]
          \arrow["g", from=1-3, to=2-3]
          \arrow[""{name=0, anchor=center, inner sep=0}, "m", "\shortmid"{marking}, from=1-2, to=1-3]
          \arrow[""{name=1, anchor=center, inner sep=0}, "n"', "\shortmid"{marking}, from=2-2, to=2-3]
          \arrow[""{name=2, anchor=center, inner sep=0}, "r", "\shortmid"{marking}, from=1-1, to=1-2]
          \arrow[""{name=3, anchor=center, inner sep=0}, "t"', "\shortmid"{marking}, from=2-1, to=2-2]
          \arrow["f"', from=1-1, to=2-1]
          \arrow[""{name=4, anchor=center, inner sep=0}, "n"', "\shortmid"{marking}, from=3-1, to=3-3]
          \arrow[Rightarrow, no head, from=2-3, to=3-3]
          \arrow[Rightarrow, no head, from=2-1, to=3-1]
          \arrow["\alpha"{description}, draw=none, from=0, to=1]
          \arrow["\phi"{description}, draw=none, from=2, to=3]
          \arrow["{\lambda^n}"{description}, draw=none, from=2-2, to=4]
        \end{tikzcd}
        =
        \begin{tikzcd}[sep=small]
          x & x & y \\
          x && y \\
          w && z
          \arrow["f"', from=2-1, to=3-1]
          \arrow["g", from=2-3, to=3-3]
          \arrow[""{name=0, anchor=center, inner sep=0}, "m"', "\shortmid"{marking}, from=2-1, to=2-3]
          \arrow[""{name=1, anchor=center, inner sep=0}, "n"', "\shortmid"{marking}, from=3-1, to=3-3]
          \arrow["r", "\shortmid"{marking}, from=1-1, to=1-2]
          \arrow["m", "\shortmid"{marking}, from=1-2, to=1-3]
          \arrow[Rightarrow, no head, from=1-3, to=2-3]
          \arrow[Rightarrow, no head, from=1-1, to=2-1]
          \arrow["\alpha"{description}, draw=none, from=0, to=1]
          \arrow["{\lambda^m}"{description}, draw=none, from=1-2, to=0]
        \end{tikzcd}
        \qquad\text{and}\qquad
        \begin{tikzcd}[sep=small]
          x & y & y \\
          w & z & z \\
          w && z
          \arrow["f"', from=1-1, to=2-1]
          \arrow["g", from=1-2, to=2-2]
          \arrow[""{name=0, anchor=center, inner sep=0}, "m", "\shortmid"{marking}, from=1-1, to=1-2]
          \arrow[""{name=1, anchor=center, inner sep=0}, "n"', "\shortmid"{marking}, from=2-1, to=2-2]
          \arrow["g", from=1-3, to=2-3]
          \arrow[""{name=2, anchor=center, inner sep=0}, "s", "\shortmid"{marking}, from=1-2, to=1-3]
          \arrow[""{name=3, anchor=center, inner sep=0}, "u"', "\shortmid"{marking}, from=2-2, to=2-3]
          \arrow[""{name=4, anchor=center, inner sep=0}, "n"', "\shortmid"{marking}, from=3-1, to=3-3]
          \arrow[Rightarrow, no head, from=2-1, to=3-1]
          \arrow[Rightarrow, no head, from=2-3, to=3-3]
          \arrow["\alpha"{description}, draw=none, from=0, to=1]
          \arrow["\psi"{description}, draw=none, from=2, to=3]
          \arrow["{\rho^n}"{description}, draw=none, from=2-2, to=4]
        \end{tikzcd}
        =
        \begin{tikzcd}[sep=small]
          x & y & y \\
          x && y \\
          w && z
          \arrow["f"', from=2-1, to=3-1]
          \arrow["g", from=2-3, to=3-3]
          \arrow[""{name=0, anchor=center, inner sep=0}, "m"', "\shortmid"{marking}, from=2-1, to=2-3]
          \arrow[""{name=1, anchor=center, inner sep=0}, "n"', "\shortmid"{marking}, from=3-1, to=3-3]
          \arrow["m", "\shortmid"{marking}, from=1-1, to=1-2]
          \arrow["s", "\shortmid"{marking}, from=1-2, to=1-3]
          \arrow[Rightarrow, no head, from=1-3, to=2-3]
          \arrow[Rightarrow, no head, from=1-1, to=2-1]
          \arrow["\alpha"{description}, draw=none, from=0, to=1]
          \arrow["{\rho^m}"{description}, draw=none, from=1-2, to=0]
        \end{tikzcd}.
      \end{equation*}
  \end{enumerate}
  There is a category of profunctors and natural transformations in $\dbl{D}$,
  with composition inherited from $\dbl{D}_1$.
\end{definition}

Every category object $(x,r,\mu,\eta)$ in a double category has an associated
\define{hom-profunctor}, the endoprofunctor $(r,\mu,\mu)$ with both left and
right action given by the composition $\mu$. For any two functors
$(f,\phi), (g,\psi): (x,r,\mu,\eta) \to (y,s,\nu,\theta)$, a natural
transformation between the hom-profunctors $(r,\mu,\mu)$ and $(s,\nu,\nu)$ with
source $(f,\phi)$ and target $(g,\psi)$ is precisely a natural transformation
from $(f,\phi)$ to $(g,\psi)$, as shown by
\cref{prop:natural-transformation-in-dbl}. Thus, the definition of a natural
transformation just given generalizes the previous one in
\cref{def:category-object}.

With these definitions, one might expect that the 2-category $\Cat(\dbl{D})$ of
category objects in a double category $\dbl{D}$ upgrades to a double category,
traditionally denoted $\Module{\dbl{D}}$, whose proarrows are
profunctors/bimodules in $\dbl{D}$. In general, $\Module{\dbl{D}}$ is only a
\emph{virtual} double category \cite[\S 5.3]{leinster2004}, \cite[\S
2.8]{cruttwell2010}. However, in the important cases of spans and matrices, we
can sketch a direct proof that there is a well-defined double category of
bimodules.

\begin{proposition}[Internal and enriched profunctors]
  For any category $\cat{S}$ with pullbacks, profunctor objects in
  $\Span{\cat{S}}$ are internal profunctors in $\cat{S}$ and natural
  transformations between profunctor objects are internal natural
  transformations in $\cat{S}$. Morever, if $\cat{S}$ has coequalizers that are
  preserved by pullbacks, then
  \begin{equation*}
    \Module{\Span{\cat{S}}} = \Prof{\cat{S}}
  \end{equation*}
  is the double category of internal categories, internal functors, internal
  profunctors, and internal natural transformations in $\cat{S}$.

  Similarly, for any distributive monoidal category $\catV$, profunctor objects
  in $\Mat{\catV}$ are profunctors enriched in $\catV$ and likewise for natural
  transformations between profunctors. Moreover, if $\catV$ is a cocomplete
  closed monoidal category, then
  \begin{equation*}
    \Module{\Mat{\catV}} = \EnProf{\catV}
  \end{equation*}
  is the double category of $\catV$-categories, $\catV$-functors,
  $\catV$-profunctors, and $\catV$-natural transformations.
\end{proposition}
\begin{proof}[Proof sketch]
  We prove that profunctors in $\Span{\cat{S}}$ are internal profunctors in
  $\cat{S}$; the case of enriched profunctors is similar. For the definitions of
  internal $\cat{S}$-valued functors, including internal profunctors, and of
  morphisms between these, see \cite[Vol 1, \S 8.2]{borceux1994} and \cite[\S
  B2.7]{johnstone2002}.
  \begin{enumerate}[(i)]
    \item A profunctor object in $\Span{\cat{S}}$ consists of a span
      $(X_0 \xproto{M} Y_0) = (X_0 \xleftarrow{\ell} M_0 \xrightarrow{r} Y_0)$
      in $\cat{S}$ and two maps of spans in $\cat{S}$,
      \begin{equation*}
        \begin{tikzcd}
          {X_0} & {X_1 \times_{X_0} M_0} & {Y_0} \\
          {X_0} & {M_0} & {Y_0}
          \arrow["\ell", from=2-2, to=2-1]
          \arrow[Rightarrow, no head, from=1-1, to=2-1]
          \arrow["r"', from=2-2, to=2-3]
          \arrow["\lambda"', from=1-2, to=2-2]
          \arrow[Rightarrow, no head, from=1-3, to=2-3]
          \arrow["{\pi_2 \circ r}", from=1-2, to=1-3]
          \arrow["{s \circ \pi_1}"', from=1-2, to=1-1]
        \end{tikzcd}
        \qquad\text{and}\qquad
        \begin{tikzcd}
          {X_0} & {M_0 \times_{Y_0} Y_1} & {Y_0} \\
          {X_0} & {M_0} & {Y_0}
          \arrow["\ell", from=2-2, to=2-1]
          \arrow[Rightarrow, no head, from=1-1, to=2-1]
          \arrow["r"', from=2-2, to=2-3]
          \arrow["\rho"', from=1-2, to=2-2]
          \arrow[Rightarrow, no head, from=1-3, to=2-3]
          \arrow["{t \circ \pi_2}", from=1-2, to=1-3]
          \arrow["{\ell \circ \pi_1}"', from=1-2, to=1-1]
        \end{tikzcd},
      \end{equation*}
      obeying the associativity, unitality, and compatibility axioms. This is
      precisely an internal profunctor in $\cat{S}$.
    \item A natural transformation between profunctors in $\Span{\cat{S}}$
      consists of a map of spans in $\cat{S}$
      \begin{equation*}
        \begin{tikzcd}
          {X_0} & {M_0} & {Y_0} \\
          {W_0} & {N_0} & {Z_0}
          \arrow["{f_0}"', from=1-1, to=2-1]
          \arrow[from=1-2, to=1-1]
          \arrow[from=1-2, to=1-3]
          \arrow[from=2-2, to=2-1]
          \arrow[from=2-2, to=2-3]
          \arrow["\alpha"', from=1-2, to=2-2]
          \arrow["{g_0}", from=1-3, to=2-3]
        \end{tikzcd}
      \end{equation*}
      satisfying the equivariance axioms, which is precisely a natural
      transformation $\alpha: M_0 \to N_0$ between internal profunctors in
      $\cat{S}$.
  \end{enumerate}
  The statements about obtaining double categories of profunctors are special
  cases of a result by Shulman on bimodules in equipments \cite[Theorem
  11.5]{shulman2008}; see, in particular, \cite[Examples 11.7 and
  11.8]{shulman2008}.
\end{proof}

The fact that lax functors send objects to categories
(\cref{prop:lax-functor-category-objects}) now extends to:

\begin{proposition}[Lax functors give profunctor objects]
  \label{prop:lax-functors-give-profunctors}
  Let $\dbl{D}$ be a strict double category and let $F: \dbl{D} \to \dbl{E}$ be
  a lax double functor.
  \begin{enumerate}[(i)]
    \item For every proarrow $m: x \proto y$ in $\dbl{D}$, the proarrow
      $Fm: Fx \proto Fy$ together with the laxators
      $F_{x,m}: F\id_x \odot Fm \to Fm$ and $F_{m,y}: Fm \odot F\id_y \to Fm$
      are a profunctor in $\dbl{E}$ from the category
      $(Fx, F\id_x, F_{x,x}, F_x)$ to the category $(Fy, F\id_y, F_{y,y}, F_y)$.
    \item For every cell $\stdInlineCell{\alpha}$ in $\dbl{D}$, the cell
      $\inlineCell{Fx}{Fy}{Fw}{Fz}{Fm}{Fn}{Ff}{Fg}{\alpha}$ in $\dbl{E}$ is a
      natural transformation from the profunctor $(Fm,F_{x,m},F_{m,y})$ to
      $(Fn,F_{w,n},F_{n,z})$ with source $(Ff, F\id_f)$ and target
      $(Fg,F\id_g)$.
  \end{enumerate}
  Moreover, these assignments are functorial with respect to composition in
  $\dbl{D}_1$.
\end{proposition}
\begin{proof}
  \begin{enumerate}[(i)]
    \item Associativity, unitality, and compatiblity for the left and right
      profunctor actions follow directly from the associativity and unitality of
      the laxators of the lax functor $F$.
    \item By the naturality of the laxators of $F$, the equation
      $\id_f \odot \alpha = \alpha = \alpha \odot \id_g$ in $\dbl{D}$ implies
      that
      \begin{equation*}
        \begin{tikzcd}[row sep=scriptsize]
          Fx & Fx & Fx \\
          Fw & Fw & Fz \\
          Fw && Fz
          \arrow[""{name=0, anchor=center, inner sep=0}, "{F \id_x}", "\shortmid"{marking}, from=1-1, to=1-2]
          \arrow[""{name=1, anchor=center, inner sep=0}, "Fm", "\shortmid"{marking}, from=1-2, to=1-3]
          \arrow["Ff"{description}, from=1-2, to=2-2]
          \arrow["Ff"', from=1-1, to=2-1]
          \arrow["Fg", from=1-3, to=2-3]
          \arrow[""{name=2, anchor=center, inner sep=0}, "Fn"', "\shortmid"{marking}, from=2-2, to=2-3]
          \arrow[""{name=3, anchor=center, inner sep=0}, "{F \id_w}"', "\shortmid"{marking}, from=2-1, to=2-2]
          \arrow[""{name=4, anchor=center, inner sep=0}, "Fn"', "\shortmid"{marking}, from=3-1, to=3-3]
          \arrow[Rightarrow, no head, from=2-1, to=3-1]
          \arrow[Rightarrow, no head, from=2-3, to=3-3]
          \arrow["{F_{w,n}}"{description}, draw=none, from=2-2, to=4]
          \arrow["{F \id_f}"{description}, draw=none, from=0, to=3]
          \arrow["F\alpha"{description}, draw=none, from=1, to=2]
        \end{tikzcd}
        \quad=\quad
        \begin{tikzcd}[row sep=scriptsize]
          Fx & Fx & Fx \\
          Fx && Fy \\
          Fw && Fz
          \arrow[""{name=0, anchor=center, inner sep=0}, "Fn"', "\shortmid"{marking}, from=3-1, to=3-3]
          \arrow[""{name=1, anchor=center, inner sep=0}, "Fm"', "\shortmid"{marking}, from=2-1, to=2-3]
          \arrow["Ff"', from=2-1, to=3-1]
          \arrow["Fg", from=2-3, to=3-3]
          \arrow["{F \id_x}", "\shortmid"{marking}, from=1-1, to=1-2]
          \arrow["Fm", "\shortmid"{marking}, from=1-2, to=1-3]
          \arrow[Rightarrow, no head, from=1-1, to=2-1]
          \arrow[Rightarrow, no head, from=1-3, to=2-3]
          \arrow["F\alpha"{description}, draw=none, from=1, to=0]
          \arrow["{F_{x,m}}"{description}, draw=none, from=1-2, to=1]
        \end{tikzcd}
      \end{equation*}
      and
      \begin{equation*}
        \begin{tikzcd}[row sep=scriptsize]
          Fx & Fy & Fy \\
          Fw & Fz & Fz \\
          Fw && Fz
          \arrow[""{name=0, anchor=center, inner sep=0}, "Fm", "\shortmid"{marking}, from=1-1, to=1-2]
          \arrow[""{name=1, anchor=center, inner sep=0}, "{F \id_y}", "\shortmid"{marking}, from=1-2, to=1-3]
          \arrow["Fg"{description}, from=1-2, to=2-2]
          \arrow["Ff"', from=1-1, to=2-1]
          \arrow["Fg", from=1-3, to=2-3]
          \arrow[""{name=2, anchor=center, inner sep=0}, "{F \id_z}"', "\shortmid"{marking}, from=2-2, to=2-3]
          \arrow[""{name=3, anchor=center, inner sep=0}, "Fn"', "\shortmid"{marking}, from=2-1, to=2-2]
          \arrow[""{name=4, anchor=center, inner sep=0}, "Fn"', "\shortmid"{marking}, from=3-1, to=3-3]
          \arrow[Rightarrow, no head, from=2-1, to=3-1]
          \arrow[Rightarrow, no head, from=2-3, to=3-3]
          \arrow["{F\id_g}"{description}, draw=none, from=1, to=2]
          \arrow["{F_{n,z}}"{description}, draw=none, from=2-2, to=4]
          \arrow["{F \alpha}"{description}, draw=none, from=0, to=3]
        \end{tikzcd}
        \quad=\quad
        \begin{tikzcd}[row sep=scriptsize]
          Fx & Fy & Fy \\
          Fx && Fy \\
          Fw && Fz
          \arrow[""{name=0, anchor=center, inner sep=0}, "Fn"', "\shortmid"{marking}, from=3-1, to=3-3]
          \arrow[""{name=1, anchor=center, inner sep=0}, "Fm"', "\shortmid"{marking}, from=2-1, to=2-3]
          \arrow["Ff"', from=2-1, to=3-1]
          \arrow["Fg", from=2-3, to=3-3]
          \arrow["Fm", "\shortmid"{marking}, from=1-1, to=1-2]
          \arrow["{F \id_y}", "\shortmid"{marking}, from=1-2, to=1-3]
          \arrow[Rightarrow, no head, from=1-1, to=2-1]
          \arrow[Rightarrow, no head, from=1-3, to=2-3]
          \arrow["F\alpha"{description}, draw=none, from=1, to=0]
          \arrow["{F_{m,y}}"{description}, draw=none, from=1-2, to=1]
        \end{tikzcd},
      \end{equation*}
      which are precisely the equivariance axioms of a natural transformation
      between profunctors. \qedhere
  \end{enumerate}
\end{proof}

\begin{remark}[External functorality] \label{rmk:unitalization-preview}
  \cref{prop:lax-functor-category-objects,prop:lax-functors-give-profunctors}
  show that a lax functor $F\colon \dbl{D}\to\dbl{E}$ sends objects and
  proarrows to category and profunctor objects in $\dbl{E}$. A natural question
  arises as to whether these assignments are coherent in the sense of being
  laxly functorial. The answer, as we will see, is not only an affirmative one,
  but moreover the passage from $\dbl{E}$-valued lax functors to
  $\Module{\dbl{E}}$-valued lax functors exhibits a universal property. The most
  direct way to state and prove this result involves knowing more about the
  structure of the ordinarily merely virtual double category $\Module{\dbl{E}}$,
  so we will defer this development until
  \cref{sec:cartesian-equipments,sec:lax-functors-into-cartesian-equipments}.
  However, since this development will also be helpful in understanding models
  of (simple) double theories, we will preview it now. The default semantics of
  double theories will be $\Span$. As span-valued lax functors give categories
  and profunctors, the embryonic form of the result is that this assignment of
  categories and profunctors is coherent in the sense that span-valued lax
  functors are in one-to-one correspondence with certain \emph{normal} lax
  functors. In fact, these lax functors in actually \emph{unitary} provided we 
  make a choice of units in $\Prof$ and $\Span$. The following corollary can in 
  fact be proved directly by hand, but is probably best deduced as a result of
  \cref{prop:unitalizationoflaxdoublefunctor} which we prove later.
\end{remark}

\begin{corollary}[One-dimensional universal property of profunctors]
  \label{corollary:universalpropertyofPROF}
  For any lax double functor $F\colon \dbl{D}\to \Span$ there is a unique,
  profunctor-valued, unitary lax double functor $\bar F: \dbl{D} \to \Prof$
  making the following triangle commute
    \begin{equation*}
      \begin{tikzcd}
        & \Prof \\
        {\dbl{D}} & \Span
        \arrow["{\Ob}", from=1-2, to=2-2]
        \arrow["F"', from=2-1, to=2-2]
        \arrow["{\bar F}", dashed, from=2-1, to=1-2]
      \end{tikzcd},
    \end{equation*}
  where $\Ob\colon \Prof\to\Span$ is the lax functor of \cref{ex:ob-functor}.
\end{corollary}

\section{Simple double theories and models}
\label{sec:simple-double-theories}

Having seen how lax functors internalize the main concepts of category theory,
we already have enough machinery to define a primitive kind of double theory and
its models. While relatively inexpressive, such theories still include a number
of interesting examples. They also serve to motivate the general approach. We
call these double theories ``simple'' to distinguish them from the cartesian
double theories developed later.

\begin{definition}[Simple double theory]
  A \define{simple double theory} is a small, strict double category. A
  \define{morphism} between double theories $\dbl{T}$ and $\dbl{T}'$ is a strict
  double functor $\dbl{T} \to \dbl{T}'$.

  A \define{model} of a simple double theory $\dbl{T}$ in a double category
  $\dbl{S}$ is a lax double functor $\dbl{T} \to \dbl{S}$. The receiving double
  category $\dbl{S}$ is called the \define{semantics}.
\end{definition}

When not explicitly stated, the semantics $\dbl{S}$ is assumed to be $\Span$,
the double category of spans, or equivalently $\Mat$, the double category of
set-valued matrices. We will also occasionally speak of \define{strict} or
\define{pseudo} models of a theory $\dbl{T}$, meaning strict or pseudo double
functors $\dbl{T} \to \dbl{S}$, mainly as a contrast to the default lax notion.

\begin{remark}[Presenting simple double theories]
  \label{rem:presenting-simple-theories}
  In what follows we will present simple double theories by generators and
  relations. The formal justification for these informal presentations lies in
  the fact that the category of simple double theories, i.e., the category of
  small, strict double categories and strict double functors, is the category of
  models of a finite limit theory, or equivalently of a finite limit sketch. We
  can thus rely on classical results about finite limit sketches. In particular,
  a strict double category presented by generators and relations can be
  constructed as the free model of the finite limit sketch for double categories
  augmented with global elements for each generator and suitable equations for
  the relations. The existence of free models of finite limit sketches is
  provided by \cite[Theorem 4.4.1]{barr1985}.
\end{remark}

The validity of the first several examples is immediate from
\cref{prop:lax-functor-category-objects}.

\begin{theory}[Categories] \label{th:unit-theory}
  The \define{unit theory} is the terminal double category $\dbl{1}$. A model of
  the unit theory is a category, restating the famous fact that a lax double
  functor from $\dbl{1}$ to $\Span$ or $\Mat$ is equivalent to a category.
  Equivalently, by \cref{corollary:universalpropertyofPROF}, models of the unit
  theory are unitary lax functors into $\Prof$. A model valued in $\Rel$ is a
  preorder. A model in $\Span(\cat{Top})$ is a \emph{2-space}, that is, a
  category internal to topological spaces.

  A strict model of the unit theory is merely a set.
\end{theory}

\begin{theory}[Functors] \label{th:walking-arrow}
  The \define{walking arrow theory} is $\VerDbl(\cat{2})$, the walking arrow
  $\cat{2} \coloneqq \{0 \to 1\}$ regarded as a double category with trivial
  proarrows and cells. A model is a pair of categories along with a functor
  between them. A model valued in $\Rel \cong \Mat(\cat{2})$ is a monotone
  function between preorders.

  A strict model is a pair of sets and a function between them.
\end{theory}

\begin{theory}[Transformations] \label{th:special-cell}
  The \define{theory of special cells} is the double category freely generated
  by two objects $x$ and $y$, two arrows $f,g: x \to y$, and one cell of the
  form
  \begin{equation*}
    \begin{tikzcd}
      x & x \\
      y & y
      \arrow["f"', from=1-1, to=2-1]
      \arrow["g", from=1-2, to=2-2]
      \arrow[""{name=0, anchor=center, inner sep=0}, "{\id_x}", "\shortmid"{marking}, from=1-1, to=1-2]
      \arrow[""{name=1, anchor=center, inner sep=0}, "{\id_y}"', "\shortmid"{marking}, from=2-1, to=2-2]
      \arrow["\alpha"{description}, draw=none, from=0, to=1]
    \end{tikzcd}.
  \end{equation*}
  A model is a parallel pair of functors and a natural transformation between them.
\end{theory}

\begin{theory}[Adjunctions] \label{th:adjunctions}
  The \define{theory of adjunctions} is generated by two objects $x$ and $y$,
  two arrows $f\colon x\to y$ and $g\colon y\to x$ and two cells
  \begin{equation*}
    \begin{tikzcd}
      x & x \\
      & y \\
      x & x
      \arrow[""{name=0, anchor=center, inner sep=0}, "{\id_x}", "\shortmid"{marking}, from=1-1, to=1-2]
      \arrow["f", from=1-2, to=2-2]
      \arrow[Rightarrow, no head, from=1-1, to=3-1]
      \arrow[""{name=1, anchor=center, inner sep=0}, "{\id_x}"', "\shortmid"{marking}, from=3-1, to=3-2]
      \arrow["g", from=2-2, to=3-2]
      \arrow["\eta"{description}, draw=none, from=0, to=1]
    \end{tikzcd}
    \qquad\text{and}\qquad
    \begin{tikzcd}
      y & y \\
      x \\
      y & y
      \arrow["g"', from=1-1, to=2-1]
      \arrow["f"', from=2-1, to=3-1]
      \arrow[""{name=0, anchor=center, inner sep=0}, "{\id_y}"', "\shortmid"{marking}, from=3-1, to=3-2]
      \arrow[""{name=1, anchor=center, inner sep=0}, "{\id_y}", "\shortmid"{marking}, from=1-1, to=1-2]
      \arrow[Rightarrow, no head, from=1-2, to=3-2]
      \arrow["\epsilon"{description}, draw=none, from=1, to=0]
    \end{tikzcd},
  \end{equation*}
  the \define{unit} and \define{counit}, satisfying the \define{triangle equations}
  \begin{equation*}
    \begin{tikzcd}[row sep=scriptsize]
      x & x & x &&&& y & y & y \\
      & y & y & x & x && x & x && y & y \\
      x & x && y & y &&& y & y & x & x \\
      y & y & y &&&& x & x & x
      \arrow[""{name=0, anchor=center, inner sep=0}, "{\id_x}", "\shortmid"{marking}, from=1-1, to=1-2]
      \arrow["f"', from=1-2, to=2-2]
      \arrow[Rightarrow, no head, from=1-1, to=3-1]
      \arrow[""{name=1, anchor=center, inner sep=0}, "{\id_x}", "\shortmid"{marking}, from=3-1, to=3-2]
      \arrow[""{name=2, anchor=center, inner sep=0}, "{\id_y}"', "\shortmid"{marking}, from=4-2, to=4-3]
      \arrow[""{name=3, anchor=center, inner sep=0}, Rightarrow, no head, from=2-3, to=4-3]
      \arrow["g", from=2-2, to=3-2]
      \arrow[""{name=4, anchor=center, inner sep=0}, "{\id_y}"', "\shortmid"{marking}, from=2-2, to=2-3]
      \arrow["f", from=3-2, to=4-2]
      \arrow[""{name=5, anchor=center, inner sep=0}, "{\id_x}", "\shortmid"{marking}, from=1-2, to=1-3]
      \arrow["f", from=1-3, to=2-3]
      \arrow["f"', from=3-1, to=4-1]
      \arrow[""{name=6, anchor=center, inner sep=0}, "{\id_y}"', "\shortmid"{marking}, from=4-1, to=4-2]
      \arrow[""{name=7, anchor=center, inner sep=0}, "f"', from=2-4, to=3-4]
      \arrow[""{name=8, anchor=center, inner sep=0}, "{\id_x}", "\shortmid"{marking}, from=2-4, to=2-5]
      \arrow["f", from=2-5, to=3-5]
      \arrow[""{name=9, anchor=center, inner sep=0}, "{\id_y}"', "\shortmid"{marking}, from=3-4, to=3-5]
      \arrow[""{name=10, anchor=center, inner sep=0}, "{\id_y}", "\shortmid"{marking}, from=1-7, to=1-8]
      \arrow[""{name=11, anchor=center, inner sep=0}, "{\id_y}", "\shortmid"{marking}, from=1-8, to=1-9]
      \arrow[""{name=12, anchor=center, inner sep=0}, Rightarrow, no head, from=1-9, to=3-9]
      \arrow["g", from=1-8, to=2-8]
      \arrow["f", from=2-8, to=3-8]
      \arrow[""{name=13, anchor=center, inner sep=0}, "\shortmid"{marking}, from=3-8, to=3-9]
      \arrow[""{name=14, anchor=center, inner sep=0}, "\shortmid"{marking}, from=2-7, to=2-8]
      \arrow[Rightarrow, no head, from=2-7, to=4-7]
      \arrow[""{name=15, anchor=center, inner sep=0}, "{\id_x}"', "\shortmid"{marking}, from=4-7, to=4-8]
      \arrow["g"', from=3-8, to=4-8]
      \arrow[""{name=16, anchor=center, inner sep=0}, "{\id_x}"', "\shortmid"{marking}, from=4-8, to=4-9]
      \arrow["g", from=3-9, to=4-9]
      \arrow["g"', from=1-7, to=2-7]
      \arrow[""{name=17, anchor=center, inner sep=0}, "g"', from=2-10, to=3-10]
      \arrow[""{name=18, anchor=center, inner sep=0}, "{\id_x}"', "\shortmid"{marking}, from=3-10, to=3-11]
      \arrow[""{name=19, anchor=center, inner sep=0}, "{\id_y}", "\shortmid"{marking}, from=2-10, to=2-11]
      \arrow["g", from=2-11, to=3-11]
      \arrow["\eta"{description}, draw=none, from=0, to=1]
      \arrow["\epsilon"{description}, draw=none, from=4, to=2]
      \arrow["{\id_f}"{description}, draw=none, from=5, to=4]
      \arrow["{\id_f}"{description}, draw=none, from=1, to=6]
      \arrow["{\id_f}"{description}, draw=none, from=8, to=9]
      \arrow["{=}"{description, pos=0.6}, shift right=3, draw=none, from=12, to=17]
      \arrow["{=}"{description, pos=0.6}, shift left=4, draw=none, from=3, to=7]
      \arrow["{\id_g}"{description}, draw=none, from=13, to=16]
      \arrow["\eta"{description}, draw=none, from=14, to=15]
      \arrow["{\id_g}"{description}, draw=none, from=10, to=14]
      \arrow["\epsilon"{description}, draw=none, from=11, to=13]
      \arrow["{\id_g}"{description}, draw=none, from=19, to=18]
    \end{tikzcd}.
  \end{equation*}
  A model is a pair of functors $F: \cat{C} \to \cat{D}$ and
  $G: \cat{D} \to \cat{C}$ equipped with an adjunction $F \dashv G$. To see
  this, notice that a model, or $\Span$-valued lax functor, is by
  \cref{corollary:universalpropertyofPROF} equivalent to a $\Prof$-valued
  \emph{normal} or better \emph{unitary} lax functor, which is evidently 
  an adjunction specified by unit and counit cells.
\end{theory}

\begin{theory}[Dual pairs] \label{th:dualpair}
  The \define{theory of dual pairs} consists of two objects $x$ and $y$ along
  with two proarrows $u\colon x \proto y$ and $v\colon y\proto x$ and two cells
  \begin{equation*}
    \begin{tikzcd}
      x && x && y & x & y \\
      x & y & x && y && y
      \arrow["v", "\shortmid"{marking}, from=1-5, to=1-6]
      \arrow["u", "\shortmid"{marking}, from=1-6, to=1-7]
      \arrow[Rightarrow, no head, from=1-5, to=2-5]
      \arrow[""{name=0, anchor=center, inner sep=0}, "{\id_y}"', "\shortmid"{marking}, from=2-5, to=2-7]
      \arrow[Rightarrow, no head, from=1-7, to=2-7]
      \arrow[""{name=1, anchor=center, inner sep=0}, "{\id_x}", "\shortmid"{marking}, from=1-1, to=1-3]
      \arrow["u"', "\shortmid"{marking}, from=2-1, to=2-2]
      \arrow["v"', "\shortmid"{marking}, from=2-2, to=2-3]
      \arrow[Rightarrow, no head, from=1-3, to=2-3]
      \arrow[Rightarrow, no head, from=1-1, to=2-1]
      \arrow["\eta"{description, pos=0.6}, draw=none, from=1, to=2-2]
      \arrow["\epsilon"{description, pos=0.4}, draw=none, from=1-6, to=0]
    \end{tikzcd}
  \end{equation*}
  satisfying the bicategorical analogue of the equations from \cref{th:adjunctions},
  namely,
  \begin{equation*}
    \begin{tikzcd}[sep=scriptsize]
      x && x & y & x & y & y & x && x & y & x \\
      x & y & x & y &&& y & x & y & x \\
      x & y && y & x & y & y && y & x & y & x
      \arrow[""{name=0, anchor=center, inner sep=0}, "{\id_y}"', "\shortmid"{marking}, from=3-2, to=3-4]
      \arrow[Rightarrow, no head, from=2-4, to=3-4]
      \arrow[""{name=1, anchor=center, inner sep=0}, "{\id_x}", "\shortmid"{marking}, from=1-1, to=1-3]
      \arrow[""{name=2, anchor=center, inner sep=0}, "u", "\shortmid"{marking}, from=2-1, to=2-2]
      \arrow["v"', "\shortmid"{marking}, from=2-2, to=2-3]
      \arrow[Rightarrow, no head, from=1-3, to=2-3]
      \arrow[Rightarrow, no head, from=1-1, to=2-1]
      \arrow[Rightarrow, no head, from=2-2, to=3-2]
      \arrow[""{name=3, anchor=center, inner sep=0}, "u"', "\shortmid"{marking}, from=2-3, to=2-4]
      \arrow[""{name=4, anchor=center, inner sep=0}, "u", "\shortmid"{marking}, from=1-3, to=1-4]
      \arrow[Rightarrow, no head, from=1-4, to=2-4]
      \arrow[Rightarrow, no head, from=2-1, to=3-1]
      \arrow[""{name=5, anchor=center, inner sep=0}, "u"', "\shortmid"{marking}, from=3-1, to=3-2]
      \arrow[""{name=6, anchor=center, inner sep=0}, Rightarrow, no head, from=1-5, to=3-5]
      \arrow[""{name=7, anchor=center, inner sep=0}, "u", "\shortmid"{marking}, from=1-5, to=1-6]
      \arrow[Rightarrow, no head, from=1-6, to=3-6]
      \arrow[""{name=8, anchor=center, inner sep=0}, "u"', "\shortmid"{marking}, from=3-5, to=3-6]
      \arrow[""{name=9, anchor=center, inner sep=0}, "v"', "\shortmid"{marking}, from=2-7, to=2-8]
      \arrow["u", "\shortmid"{marking}, from=2-8, to=2-9]
      \arrow[Rightarrow, no head, from=2-9, to=3-9]
      \arrow[Rightarrow, no head, from=2-7, to=3-7]
      \arrow[""{name=10, anchor=center, inner sep=0}, "{\id_y}"', "\shortmid"{marking}, from=3-7, to=3-9]
      \arrow[Rightarrow, no head, from=1-8, to=2-8]
      \arrow[""{name=11, anchor=center, inner sep=0}, "{\id_x}", "\shortmid"{marking}, from=1-8, to=1-10]
      \arrow[Rightarrow, no head, from=1-10, to=2-10]
      \arrow[""{name=12, anchor=center, inner sep=0}, "v", "\shortmid"{marking}, from=2-9, to=2-10]
      \arrow[""{name=13, anchor=center, inner sep=0}, "v"', "\shortmid"{marking}, from=3-9, to=3-10]
      \arrow[Rightarrow, no head, from=2-10, to=3-10]
      \arrow[Rightarrow, no head, from=1-7, to=2-7]
      \arrow[""{name=14, anchor=center, inner sep=0}, "v", "\shortmid"{marking}, from=1-7, to=1-8]
      \arrow[""{name=15, anchor=center, inner sep=0}, Rightarrow, no head, from=1-11, to=3-11]
      \arrow[""{name=16, anchor=center, inner sep=0}, "v"', "\shortmid"{marking}, from=3-11, to=3-12]
      \arrow[""{name=17, anchor=center, inner sep=0}, "v", "\shortmid"{marking}, from=1-11, to=1-12]
      \arrow[Rightarrow, no head, from=1-12, to=3-12]
      \arrow["\eta"{description, pos=0.6}, draw=none, from=1, to=2-2]
      \arrow["\epsilon"{description, pos=0.4}, draw=none, from=2-3, to=0]
      \arrow["{1_u}"{description}, draw=none, from=4, to=3]
      \arrow["{1_u}"{description}, draw=none, from=7, to=8]
      \arrow["{1_u}"{description}, draw=none, from=2, to=5]
      \arrow["{=}"{description, pos=0.4}, draw=none, from=2-4, to=6]
      \arrow["{=}"{description, pos=0.4}, draw=none, from=2-10, to=15]
      \arrow["{1_v}"{description}, draw=none, from=14, to=9]
      \arrow["\epsilon"{description, pos=0.4}, draw=none, from=2-8, to=10]
      \arrow["{1_v}"{description}, draw=none, from=12, to=13]
      \arrow["\eta"{description, pos=0.6}, draw=none, from=11, to=2-9]
      \arrow["{1_v}"{description}, draw=none, from=17, to=16]
    \end{tikzcd}.
  \end{equation*}
  A strict model in a double category $\dbl{S}$ is a \define{dual pair}
  \cite[\S{5}]{shulman2008}, that is, a pair of proarrows in $\dbl{S}$ that are
  internally adjoint in the horizontal bicategory of $\dbl{S}$. Dual pairs
  formally resemble Street's \emph{biexact pairings}; see
  \cref{th:biexactpairing} below. The image of $u$ is called the \define{left
    dual} and that of $v$ the \define{right dual}. In several cases, dual pairs
  are well-understood \cite[Example 5.6]{shulman2008}. For example, when
  $\dbl{S} = \Span$, they are companions and conjoints.

  The theory of dual pairs is atypical from the perspective of double theories.
  Of the theories presented in this section, it is the only one where strict,
  rather lax, models are of primary interest. It is also the only theory that
  requires a cell whose codomain is a nontrivial composite of proarrows.
\end{theory}

\begin{theory}[Monads] \label{th:monad}
  The \define{theory of monads} is generated by
  \begin{itemize}[noitemsep]
    \item an object $x$,
    \item an arrow $t: x \to x$, and
    \item \define{multiplication} and \define{unit} cells
    \begin{equation*}
      \begin{tikzcd}[row sep=scriptsize]
        x & x \\
        x \\
        x & x
        \arrow["t"', from=1-1, to=2-1]
        \arrow["t"', from=2-1, to=3-1]
        \arrow[""{name=0, anchor=center, inner sep=0}, "{\id_x}", "\shortmid"{marking}, from=1-1, to=1-2]
        \arrow[""{name=1, anchor=center, inner sep=0}, "{\id_x}"', "\shortmid"{marking}, from=3-1, to=3-2]
        \arrow["t", from=1-2, to=3-2]
        \arrow["\mu"{description}, draw=none, from=0, to=1]
      \end{tikzcd}
      \qquad\text{and}\qquad
      \begin{tikzcd}[row sep=scriptsize]
        x & x \\
        x & x
        \arrow[Rightarrow, no head, from=1-1, to=2-1]
        \arrow[""{name=0, anchor=center, inner sep=0}, "{\id_x}", "\shortmid"{marking}, from=1-1, to=1-2]
        \arrow["t", from=1-2, to=2-2]
        \arrow[""{name=1, anchor=center, inner sep=0}, "{\id_x}"', "\shortmid"{marking}, from=2-1, to=2-2]
        \arrow["\eta"{description}, draw=none, from=0, to=1]
      \end{tikzcd}
    \end{equation*}
  \end{itemize}
  subject to the equations of associativity
  \begin{equation*}
    \begin{tikzcd}[row sep=small]
      x & x & x \\
      x \\
      x & x \\
      x & x & x
      \arrow["t"', from=1-1, to=2-1]
      \arrow["t"', from=2-1, to=3-1]
      \arrow[""{name=0, anchor=center, inner sep=0}, "{\id_x}", "\shortmid"{marking}, from=1-1, to=1-2]
      \arrow[""{name=1, anchor=center, inner sep=0}, "{\id_x}", "\shortmid"{marking}, from=3-1, to=3-2]
      \arrow["t", from=1-2, to=3-2]
      \arrow["t"', from=3-1, to=4-1]
      \arrow["t", from=3-2, to=4-2]
      \arrow[""{name=2, anchor=center, inner sep=0}, "{\id_x}"', "\shortmid"{marking}, from=4-1, to=4-2]
      \arrow[""{name=3, anchor=center, inner sep=0}, "{\id_x}", "\shortmid"{marking}, from=1-2, to=1-3]
      \arrow[""{name=4, anchor=center, inner sep=0}, "{\id_x}"', "\shortmid"{marking}, from=4-2, to=4-3]
      \arrow["t", from=1-3, to=4-3]
      \arrow["\mu"{description}, draw=none, from=0, to=1]
      \arrow["{\id_t}"{description}, draw=none, from=1, to=2]
      \arrow["\mu"{description}, draw=none, from=3, to=4]
    \end{tikzcd}
    \quad=\quad
    \begin{tikzcd}[row sep=small]
      x & x & x \\
      x & x \\
      x \\
      x & x & x
      \arrow[""{name=0, anchor=center, inner sep=0}, "{\id_x}", "\shortmid"{marking}, from=1-1, to=1-2]
      \arrow[""{name=1, anchor=center, inner sep=0}, "{\id_x}"', "\shortmid"{marking}, from=2-1, to=2-2]
      \arrow["t", from=1-2, to=2-2]
      \arrow["t", from=2-2, to=4-2]
      \arrow[""{name=2, anchor=center, inner sep=0}, "{\id_x}"', "\shortmid"{marking}, from=4-1, to=4-2]
      \arrow[""{name=3, anchor=center, inner sep=0}, "{\id_x}", "\shortmid"{marking}, from=1-2, to=1-3]
      \arrow[""{name=4, anchor=center, inner sep=0}, "{\id_x}"', "\shortmid"{marking}, from=4-2, to=4-3]
      \arrow["t", from=1-3, to=4-3]
      \arrow["t"', from=1-1, to=2-1]
      \arrow["t"', from=2-1, to=3-1]
      \arrow["t"', from=3-1, to=4-1]
      \arrow["{\id_t}"{description}, draw=none, from=0, to=1]
      \arrow["\mu"{description}, draw=none, from=3, to=4]
      \arrow["\mu"{description}, draw=none, from=1, to=2]
    \end{tikzcd}
  \end{equation*}
  and unitality
  \begin{equation*}
    \begin{tikzcd}[row sep=scriptsize]
      x & x & x \\
      x & x \\
      x & x & x
      \arrow["t"', from=1-2, to=2-2]
      \arrow["t"', from=2-2, to=3-2]
      \arrow[""{name=0, anchor=center, inner sep=0}, "{\id_x}", "\shortmid"{marking}, from=1-2, to=1-3]
      \arrow[""{name=1, anchor=center, inner sep=0}, "{\id_x}"', "\shortmid"{marking}, from=3-2, to=3-3]
      \arrow["t", from=1-3, to=3-3]
      \arrow[Rightarrow, no head, from=1-1, to=2-1]
      \arrow[""{name=2, anchor=center, inner sep=0}, "{\id_x}", "\shortmid"{marking}, from=1-1, to=1-2]
      \arrow[""{name=3, anchor=center, inner sep=0}, "{\id_x}", "\shortmid"{marking}, from=2-1, to=2-2]
      \arrow["t"', from=2-1, to=3-1]
      \arrow[""{name=4, anchor=center, inner sep=0}, "{\id_x}"', "\shortmid"{marking}, from=3-1, to=3-2]
      \arrow["\mu"{description}, draw=none, from=0, to=1]
      \arrow["{\id_t}"{description}, draw=none, from=3, to=4]
      \arrow["\eta"{description, pos=0.4}, draw=none, from=2, to=3]
    \end{tikzcd}
    \quad=\quad
    \id_t
    \quad=\quad
    \begin{tikzcd}[row sep=scriptsize]
      x & x & x \\
      x & x \\
      x & x & x
      \arrow["t"', from=1-2, to=2-2]
      \arrow["t"', from=2-2, to=3-2]
      \arrow[""{name=0, anchor=center, inner sep=0}, "{\id_x}", "\shortmid"{marking}, from=1-2, to=1-3]
      \arrow[""{name=1, anchor=center, inner sep=0}, "{\id_x}"', "\shortmid"{marking}, from=3-2, to=3-3]
      \arrow["t", from=1-3, to=3-3]
      \arrow["t"', from=1-1, to=2-1]
      \arrow[""{name=2, anchor=center, inner sep=0}, "{\id_x}", "\shortmid"{marking}, from=1-1, to=1-2]
      \arrow[""{name=3, anchor=center, inner sep=0}, "{\id_x}", "\shortmid"{marking}, from=2-1, to=2-2]
      \arrow[Rightarrow, no head, from=2-1, to=3-1]
      \arrow[""{name=4, anchor=center, inner sep=0}, "{\id_x}"', "\shortmid"{marking}, from=3-1, to=3-2]
      \arrow["\mu"{description}, draw=none, from=0, to=1]
      \arrow["\eta"{description}, draw=none, from=3, to=4]
      \arrow["{\id_t}"{description, pos=0.4}, draw=none, from=2, to=3]
    \end{tikzcd}.
  \end{equation*}
  A model of the theory of monads is a category $\cat{C}$ along with a monad
  $(T, \mu, \eta)$ on $\cat{C}$. If $\dbl{T}$ denotes the double theory of
  monads, then $\dbl{T}^\co$ is the theory of comonads, whose models are
  categories equipped with comonads.
\end{theory}

\begin{theory}[Frobenius monads] \label{th:frob-monad}
  The \define{theory of Frobenius monads} is the theory of monads
  (\cref{th:monad}) augmented with further cells
  \begin{equation*}
    \begin{tikzcd}
      x & x \\
      & x \\
      x & x
      \arrow[Rightarrow, no head, from=1-1, to=3-1]
      \arrow[""{name=0, anchor=center, inner sep=0}, "{\id_x}", "\shortmid"{marking}, from=1-1, to=1-2]
      \arrow["t", from=1-2, to=2-2]
      \arrow["t", from=2-2, to=3-2]
      \arrow[""{name=1, anchor=center, inner sep=0}, "{\id_x}"', "\shortmid"{marking}, from=3-1, to=3-2]
      \arrow["\rho"{description}, draw=none, from=0, to=1]
    \end{tikzcd}
    \qquad\text{and}\qquad
    \begin{tikzcd}
      x & x \\
      x & x
      \arrow[""{name=0, anchor=center, inner sep=0}, "{\id_x}", "\shortmid"{marking}, from=1-1, to=1-2]
      \arrow[Rightarrow, no head, from=1-2, to=2-2]
      \arrow["t"', from=1-1, to=2-1]
      \arrow[""{name=1, anchor=center, inner sep=0}, "{\id_x}"', "\shortmid"{marking}, from=2-1, to=2-2]
      \arrow["\epsilon"{description}, draw=none, from=0, to=1]
    \end{tikzcd}
  \end{equation*}
  satisfying the additional equations
  \begin{equation*}
    \begin{tikzcd}[row sep=scriptsize]
      x & x & x && x & x & x \\
      & x & x && x & x \\
      x & x &&&& x & x \\
      x & x & x && x & x & x
      \arrow[Rightarrow, no head, from=1-1, to=3-1]
      \arrow[""{name=0, anchor=center, inner sep=0}, "{\id_x}", "\shortmid"{marking}, from=1-1, to=1-2]
      \arrow["t", from=1-2, to=2-2]
      \arrow[""{name=1, anchor=center, inner sep=0}, "t", from=2-2, to=3-2]
      \arrow[""{name=2, anchor=center, inner sep=0}, "{\id_x}"', "\shortmid"{marking}, from=2-2, to=2-3]
      \arrow["t", from=2-3, to=4-3]
      \arrow["t", from=3-2, to=4-2]
      \arrow[""{name=3, anchor=center, inner sep=0}, "{\id_x}"', "\shortmid"{marking}, from=4-2, to=4-3]
      \arrow[""{name=4, anchor=center, inner sep=0}, "{\id_x}", "\shortmid"{marking}, from=1-2, to=1-3]
      \arrow["t", from=1-3, to=2-3]
      \arrow["t"', from=3-1, to=4-1]
      \arrow[""{name=5, anchor=center, inner sep=0}, "{\id_x}"', "\shortmid"{marking}, from=4-1, to=4-2]
      \arrow[Rightarrow, no head, from=2-5, to=4-5]
      \arrow[""{name=6, anchor=center, inner sep=0}, "{\id_x}"', "\shortmid"{marking}, from=4-5, to=4-6]
      \arrow[""{name=7, anchor=center, inner sep=0}, "{\id_x}"', "\shortmid"{marking}, from=2-5, to=2-6]
      \arrow[""{name=8, anchor=center, inner sep=0}, "t", from=2-6, to=3-6]
      \arrow["t", from=3-6, to=4-6]
      \arrow["t", from=1-6, to=2-6]
      \arrow[""{name=9, anchor=center, inner sep=0}, "{\id_x}", "\shortmid"{marking}, from=1-6, to=1-7]
      \arrow["t", from=1-7, to=3-7]
      \arrow[""{name=10, anchor=center, inner sep=0}, "{\id_x}", "\shortmid"{marking}, from=3-6, to=3-7]
      \arrow["t"', from=1-5, to=2-5]
      \arrow[""{name=11, anchor=center, inner sep=0}, "{\id_x}", "\shortmid"{marking}, from=1-5, to=1-6]
      \arrow[""{name=12, anchor=center, inner sep=0}, "{\id_x}"', "\shortmid"{marking}, from=4-6, to=4-7]
      \arrow["t", from=3-7, to=4-7]
      \arrow[""{name=13, anchor=center, inner sep=0}, "{\id_x}", "\shortmid"{marking}, from=3-1, to=3-2]
      \arrow["\mu"{description}, draw=none, from=2, to=3]
      \arrow["\rho"{description}, draw=none, from=7, to=6]
      \arrow["\mu"{description}, draw=none, from=9, to=10]
      \arrow["{\id_t}"{description}, draw=none, from=10, to=12]
      \arrow["{\id_t}"{description}, draw=none, from=4, to=2]
      \arrow["\rho"{description}, draw=none, from=0, to=13]
      \arrow["{\id_t}"{description}, draw=none, from=13, to=5]
      \arrow["{\id_t}"{description}, draw=none, from=11, to=7]
      \arrow["{=}"{description}, draw=none, from=1, to=8]
    \end{tikzcd}
  \end{equation*}
  and
  \begin{equation*}
    \begin{tikzcd}[row sep=scriptsize]
      x & x & x & x & x & x & x & x \\
      & x & x &&&& x & x \\
      x & x & x & x & x & x & x & x
      \arrow[Rightarrow, no head, from=1-1, to=3-1]
      \arrow[""{name=0, anchor=center, inner sep=0}, "{\id_x}", "\shortmid"{marking}, from=1-1, to=1-2]
      \arrow["t"', from=1-2, to=2-2]
      \arrow["t"', from=2-2, to=3-2]
      \arrow[""{name=1, anchor=center, inner sep=0}, "{\id_x}"', "\shortmid"{marking}, from=3-1, to=3-2]
      \arrow[""{name=2, anchor=center, inner sep=0}, "{\id_x}", "\shortmid"{marking}, from=1-2, to=1-3]
      \arrow[Rightarrow, no head, from=1-3, to=2-3]
      \arrow[""{name=3, anchor=center, inner sep=0}, "{\id_x}", "\shortmid"{marking}, from=2-2, to=2-3]
      \arrow[""{name=4, anchor=center, inner sep=0}, "{\id_x}"', "\shortmid"{marking}, from=3-2, to=3-3]
      \arrow["t", from=2-3, to=3-3]
      \arrow[""{name=5, anchor=center, inner sep=0}, Rightarrow, no head, from=1-4, to=3-4]
      \arrow[""{name=6, anchor=center, inner sep=0}, "{\id_x}", "\shortmid"{marking}, from=1-4, to=1-5]
      \arrow[""{name=7, anchor=center, inner sep=0}, "t", from=1-5, to=3-5]
      \arrow[""{name=8, anchor=center, inner sep=0}, "{\id_x}"', "\shortmid"{marking}, from=3-4, to=3-5]
      \arrow[""{name=9, anchor=center, inner sep=0}, Rightarrow, no head, from=1-6, to=3-6]
      \arrow[""{name=10, anchor=center, inner sep=0}, "{\id_x}"', "\shortmid"{marking}, from=3-6, to=3-7]
      \arrow[""{name=11, anchor=center, inner sep=0}, "{\id_x}", "\shortmid"{marking}, from=1-6, to=1-7]
      \arrow["t"', from=1-7, to=2-7]
      \arrow["t"', from=2-7, to=3-7]
      \arrow[""{name=12, anchor=center, inner sep=0}, "{\id_x}", "\shortmid"{marking}, from=1-7, to=1-8]
      \arrow["t", from=1-8, to=2-8]
      \arrow[""{name=13, anchor=center, inner sep=0}, "{\id_x}"', "\shortmid"{marking}, from=2-7, to=2-8]
      \arrow[Rightarrow, no head, from=2-8, to=3-8]
      \arrow[""{name=14, anchor=center, inner sep=0}, "{\id_x}"', "\shortmid"{marking}, from=3-7, to=3-8]
      \arrow["\epsilon"{description}, draw=none, from=13, to=14]
      \arrow["{\id_t}"{description}, draw=none, from=12, to=13]
      \arrow["\rho"{description}, draw=none, from=11, to=10]
      \arrow["\eta"{description}, draw=none, from=6, to=8]
      \arrow["{\id_t}"{description}, draw=none, from=3, to=4]
      \arrow["\epsilon"{description}, draw=none, from=2, to=3]
      \arrow["\rho"{description}, draw=none, from=0, to=1]
      \arrow["{=}"{description, pos=0.4}, draw=none, from=2-3, to=5]
      \arrow["{=}"{description}, draw=none, from=7, to=9]
    \end{tikzcd}.
  \end{equation*}
  A model is a Frobenius monad in the sense of \cite[Definition
  1.1]{street2004}. A Frobenius monad amounts to a self-adjoint endofunctor,
  called \emph{ambidextrous} in \cite{lauda2006}. Examples include the
  oppositization functor $(-)^\op\colon \cat{Cat} \to \cat{Cat}$ and, for any
  category $\cat{C}$ with biproducts, the composite functor
  $\oplus \circ \Delta: \cat{C} \to \cat{C}$. The power-object functor
  $P\colon \cat{E}^\op \to \cat{E}$ is contravariantly self-adjoint for any
  topos $\cat{E}$.
\end{theory}

\begin{theory}[Promonads]
  The \define{theory of promonads} is generated by
  \begin{itemize}[noitemsep]
    \item an object $x$,
    \item a proarrow $p: x \proto x$, and
    \item a globular cell $\eta: \id_x \to p$, the \define{unit},
  \end{itemize}
  subject to the axioms of idempotency $p \odot p = p$ and unitality
  $\eta \odot 1_p = 1_p = 1_p \odot \eta$.

  A model of the theory is a \define{promonad} \cite{daystreet2007}. It consists
  of a category $\cat{C}$ and a profunctor $P: \cat{C} \proto \cat{C}$ along
  with natural transformations $\mu: P \odot P \To P$ and
  $\eta: \Hom_{\cat{C}} \To P$. The multiplication $\mu$, given by the laxators
  for $p$, is associative. By the naturality of the laxators, the multiplication
  is also unital in the sense that
  \begin{equation*}
    \mu_{w,y}(\eta_{w,x}(f), u) = f \cdot u
    \qquad\text{and}\qquad
    \mu_{x,z}(u, \eta_{y,z}(g)) = u \cdot g
  \end{equation*}
  for all morphisms $f: w \to x$ and $g: y \to z$ in $\cat{C}$ and all
  heteromorphisms $u \in P(x,y)$. A promonad on a category $\cat{C}$ can thus be
  regarded as giving an extension of the morphisms in $\cat{C}$, generalizing
  the Kleisli category of a monad on $\cat{C}$.
\end{theory}

\section{Cartesian equipments}
\label{sec:cartesian-equipments}

In the passage from concrete to formal category theory, categories with finite
products are abstracted as \emph{cartesian objects} in a 2-category with finite
2-categorical products \cite[\S 5.1]{carboni1991}. Cartesian double categories
are succinctly defined using this notion.

\begin{definition}[Cartesian double category]
  A \define{precartesian double category} is a cartesian object in $\DblLax$.
  Similarly, a \define{cartesian double category} is a cartesian object in
  $\Dbl$.

  In other words, a double category $\dbl{D}$ is \define{precartesian} when the
  diagonal and terminal double functors
  \begin{equation*}
    \Delta: \dbl{D} \to \dbl{D} \times \dbl{D}
    \quad\text{and}\quad
    !: \dbl{D} \to \dbl{1}
  \end{equation*}
  have lax right adjoints, denoted
  \begin{equation*}
    \times: \dbl{D} \times \dbl{D} \to \dbl{D}
    \quad\text{and}\quad
    I: \dbl{1} \to \dbl{D}.
  \end{equation*}
  If the right adjoints are pseudo, then the double category $\dbl{D}$ is
  \define{cartesian}.
\end{definition}

The main reference for cartesian double categories is Aleiferi's PhD thesis
\cite{aleiferi2018}. The short, conceptual definition of a cartesian double
category should be contrasted with the far more complicated definition of a
\emph{cartesian bicategory} as a monoidal bicategory possessing extra structure
\cite{carboni1987,carboni2008}. In particular, unlike the situation for
cartesian bicategories, it is immediate from the definition that being
(pre)cartesian is a property of, not a structure on, a double category, as it
merely asserts the existence of certain (lax) right adjoints.

It is useful to have a more explicit description of a cartesian double category.
Unpacking the definition, a precartesian double category is seen to be a double
category $\dbl{D}$ such that
\begin{itemize}
  \item the diagonal and terminal functors on the underlying categories
    $\dbl{D}_0$ and $\dbl{D}_1$ have right adjoints
    \begin{equation*}
      \times_i: \dbl{D}_i \times \dbl{D}_i \to \dbl{D}_i
      \quad\text{and}\quad
      I_i: \cat{1} \to \dbl{D}_i,
      \qquad i = 0,1,
    \end{equation*}
    which are preserved by the source and target functors:
    \begin{equation*}
      \begin{tikzcd}
        {\dbl{D}_1 \times \dbl{D}_1} & {\dbl{D}_1} \\
        {\dbl{D}_0 \times \dbl{D}_0} & {\dbl{D}_0}
        \arrow["{s \times s}"', from=1-1, to=2-1]
        \arrow["{\times_0}"', from=2-1, to=2-2]
        \arrow["{\times_1}", from=1-1, to=1-2]
        \arrow["s", from=1-2, to=2-2]
      \end{tikzcd}
      \qquad
      \begin{tikzcd}
        {\dbl{D}_1 \times \dbl{D}_1} & {\dbl{D}_1} \\
        {\dbl{D}_0 \times \dbl{D}_0} & {\dbl{D}_0}
        \arrow["{t \times t}"', from=1-1, to=2-1]
        \arrow["{\times_0}"', from=2-1, to=2-2]
        \arrow["{\times_1}", from=1-1, to=1-2]
        \arrow["t", from=1-2, to=2-2]
      \end{tikzcd}
      \qquad
      \begin{tikzcd}
        {\mathsf{1}} & {\dbl{D}_1} \\
        & {\dbl{D}_0}
        \arrow["{I_0}"', from=1-1, to=2-2]
        \arrow["s", from=1-2, to=2-2]
        \arrow["{I_1}", from=1-1, to=1-2]
      \end{tikzcd}
      \qquad
      \begin{tikzcd}
        {\mathsf{1}} & {\dbl{D}_1} \\
        & {\dbl{D}_0}
        \arrow["{I_0}"', from=1-1, to=2-2]
        \arrow["t", from=1-2, to=2-2]
        \arrow["{I_1}", from=1-1, to=1-2]
      \end{tikzcd};
    \end{equation*}
  \item diagonals and projections in $\dbl{D}_0$ and $\dbl{D}_1$, given by the
    unit and counit of the adjunctions $\Delta_0 \dashv \times_0$ and
    $\Delta_1 \dashv \times_1$, also respect source and target, in that for any
    proarrows $m: x \proto y$ and $m': x' \proto y'$ in $\dbl{D}$, the unit and
    counit cells have form
    \begin{equation*}
      \begin{tikzcd}
        x & y \\
        {x \times x} & {y \times y}
        \arrow[""{name=0, anchor=center, inner sep=0}, "m", "\shortmid"{marking}, from=1-1, to=1-2]
        \arrow[""{name=1, anchor=center, inner sep=0}, "{m \times m}"', "\shortmid"{marking}, from=2-1, to=2-2]
        \arrow["{\Delta_x}"', from=1-1, to=2-1]
        \arrow["{\Delta_y}", from=1-2, to=2-2]
        \arrow["{\Delta_m}"{description}, draw=none, from=0, to=1]
      \end{tikzcd}
      \qquad
      \begin{tikzcd}
        {x \times x'} & {y \times y'} \\
        x & y
        \arrow[""{name=0, anchor=center, inner sep=0}, "{m \times m'}", "\shortmid"{marking}, from=1-1, to=1-2]
        \arrow[""{name=1, anchor=center, inner sep=0}, "m"', "\shortmid"{marking}, from=2-1, to=2-2]
        \arrow["{\pi_{x,x'}}"', from=1-1, to=2-1]
        \arrow["{\pi_{y,y'}}", from=1-2, to=2-2]
        \arrow["{\pi_{m,m'}}"{description}, draw=none, from=0, to=1]
      \end{tikzcd}
      \qquad
      \begin{tikzcd}
        {x \times x'} & {y \times y'} \\
        {x'} & {y'}
        \arrow[""{name=0, anchor=center, inner sep=0}, "{m \times m'}", "\shortmid"{marking}, from=1-1, to=1-2]
        \arrow[""{name=1, anchor=center, inner sep=0}, "{m'}"', "\shortmid"{marking}, from=2-1, to=2-2]
        \arrow["{\pi_{x,x'}'}"', from=1-1, to=2-1]
        \arrow["{\pi_{y,y'}'}", from=1-2, to=2-2]
        \arrow["{\pi_{m,m'}'}"{description}, draw=none, from=0, to=1]
      \end{tikzcd};
    \end{equation*}
  \item for each pair of consecutive proarrows $x \xproto{m} y \xproto{n} z$ and
    $x' \xproto{m'} y' \xproto{n'} z'$ and of objects $x, x'$ in $\dbl{D}$,
    there are comparison cells
    \begin{equation} \label{eq:dbl-product-comparisons}
      \begin{tikzcd}
        {x \times x'} & {y \times y'} & {z \times z'} \\
        {x \times x'} && {z \times z'}
        \arrow[Rightarrow, no head, from=1-3, to=2-3]
        \arrow[Rightarrow, no head, from=1-1, to=2-1]
        \arrow[""{name=0, anchor=center, inner sep=0}, "{(m \odot n) \times (m' \odot n')}"', "\shortmid"{marking}, from=2-1, to=2-3]
        \arrow["{m \times m'}", "\shortmid"{marking}, from=1-1, to=1-2]
        \arrow["{n \times n'}", "\shortmid"{marking}, from=1-2, to=1-3]
        \arrow["{\times_{(m,m'),(n,n')}}"{description}, draw=none, from=1-2, to=0]
      \end{tikzcd}
      \qquad\text{and}\qquad
      \begin{tikzcd}
        {x \times x'} & {x \times x'} \\
        {x \times x'} & {x \times x'}
        \arrow[""{name=0, anchor=center, inner sep=0}, "{\mathrm{id}_{x \times x'}}", "\shortmid"{marking}, from=1-1, to=1-2]
        \arrow[""{name=1, anchor=center, inner sep=0}, "{\mathrm{id}_x \times \mathrm{id}_{x'}}"', "\shortmid"{marking}, from=2-1, to=2-2]
        \arrow[Rightarrow, no head, from=1-1, to=2-1]
        \arrow[Rightarrow, no head, from=1-2, to=2-2]
        \arrow["{\times_{(x,x')}}"{description}, draw=none, from=0, to=1]
      \end{tikzcd}
    \end{equation}
    and also comparison cells
    \begin{equation} \label{eq:dlb-terminal-comparisons}
      \begin{tikzcd}
        {I_0} & {I_0} & {I_0} \\
        {I_0} && {I_0}
        \arrow["{I_1}", "\shortmid"{marking}, from=1-1, to=1-2]
        \arrow["{I_1}", "\shortmid"{marking}, from=1-2, to=1-3]
        \arrow[""{name=0, anchor=center, inner sep=0}, "{I_1}"', "\shortmid"{marking}, from=2-1, to=2-3]
        \arrow[Rightarrow, no head, from=1-1, to=2-1]
        \arrow[Rightarrow, no head, from=1-3, to=2-3]
        \arrow["{\mu_I}"{description}, draw=none, from=1-2, to=0]
      \end{tikzcd}
      \qquad\text{and}\qquad
      \begin{tikzcd}
        {I_0} & {I_0} \\
        {I_0} & {I_0}
        \arrow[""{name=0, anchor=center, inner sep=0}, "{\mathrm{id}_{I_0}}", "\shortmid"{marking}, from=1-1, to=1-2]
        \arrow[""{name=1, anchor=center, inner sep=0}, "{I_1}"', "\shortmid"{marking}, from=2-1, to=2-2]
        \arrow[Rightarrow, no head, from=1-1, to=2-1]
        \arrow[Rightarrow, no head, from=1-2, to=2-2]
        \arrow["{\eta_I}"{description}, draw=none, from=0, to=1]
      \end{tikzcd},
    \end{equation}
    obeying the axioms of a lax double functor;
  \item diagonals and projections are externally functorial, meaning that the
    diagonal cells satisfy
    \begin{equation*}
      \begin{tikzcd}
        x && z \\
        {x \times x} && {z \times z}
        \arrow[""{name=0, anchor=center, inner sep=0}, "{m \odot n}", "\shortmid"{marking}, from=1-1, to=1-3]
        \arrow[""{name=1, anchor=center, inner sep=0}, "{(m \odot n) \times (m \odot n)}"', "\shortmid"{marking}, from=2-1, to=2-3]
        \arrow["{\Delta_x}"', from=1-1, to=2-1]
        \arrow["{\Delta_z}", from=1-3, to=2-3]
        \arrow["{\Delta_{m \odot n}}"{description}, draw=none, from=0, to=1]
      \end{tikzcd}
      \quad=\quad
      \begin{tikzcd}
        x & y & z \\
        {x \times x} & {y \times y} & {z \times z} \\
        {x \times x} && {z \times z}
        \arrow["{\Delta_x}"', from=1-1, to=2-1]
        \arrow["{\Delta_z}", from=1-3, to=2-3]
        \arrow["{\Delta_y}"{description}, from=1-2, to=2-2]
        \arrow[""{name=0, anchor=center, inner sep=0}, "{m \times m}"', "\shortmid"{marking}, from=2-1, to=2-2]
        \arrow[""{name=1, anchor=center, inner sep=0}, "m", "\shortmid"{marking}, from=1-1, to=1-2]
        \arrow[""{name=2, anchor=center, inner sep=0}, "n", "\shortmid"{marking}, from=1-2, to=1-3]
        \arrow[""{name=3, anchor=center, inner sep=0}, "{n \times n}"', "\shortmid"{marking}, from=2-2, to=2-3]
        \arrow[""{name=4, anchor=center, inner sep=0}, "{(m \odot n) \times (m \odot n)}"', "\shortmid"{marking}, from=3-1, to=3-3]
        \arrow[Rightarrow, no head, from=2-1, to=3-1]
        \arrow[Rightarrow, no head, from=2-3, to=3-3]
        \arrow["{\Delta_m}"{description}, draw=none, from=1, to=0]
        \arrow["{\Delta_n}"{description}, draw=none, from=2, to=3]
        \arrow["{\times_{(m,m),(n,n)}}"{description}, draw=none, from=2-2, to=4]
      \end{tikzcd}
    \end{equation*}
    \begin{equation*}
      \begin{tikzcd}
        x & x \\
        {x \times x} & {x \times x}
        \arrow["{\Delta_x}"', from=1-1, to=2-1]
        \arrow[""{name=0, anchor=center, inner sep=0}, "{\mathrm{id}_x}", "\shortmid"{marking}, from=1-1, to=1-2]
        \arrow["{\Delta_x}", from=1-2, to=2-2]
        \arrow[""{name=1, anchor=center, inner sep=0}, "{\mathrm{id}_x \times \mathrm{id}_x}"', "\shortmid"{marking}, from=2-1, to=2-2]
        \arrow["{\Delta_{\mathrm{id}_x}}"{description}, draw=none, from=0, to=1]
      \end{tikzcd}
      \quad=\quad
      \begin{tikzcd}
        x & x \\
        {x \times x} & {x \times x} \\
        {x \times x} & {x \times x}
        \arrow[""{name=0, anchor=center, inner sep=0}, "{\mathrm{id}_{x \times x}}"', "\shortmid"{marking}, from=2-1, to=2-2]
        \arrow["{\Delta_x}"', from=1-1, to=2-1]
        \arrow["{\Delta_x}", from=1-2, to=2-2]
        \arrow[""{name=1, anchor=center, inner sep=0}, "{\mathrm{id}_x}", "\shortmid"{marking}, from=1-1, to=1-2]
        \arrow[""{name=2, anchor=center, inner sep=0}, "{\mathrm{id}_x \times \mathrm{id}_x}"', "\shortmid"{marking}, from=3-1, to=3-2]
        \arrow[Rightarrow, no head, from=2-2, to=3-2]
        \arrow[Rightarrow, no head, from=2-1, to=3-1]
        \arrow["{\mathrm{id}_{\Delta_x}}"{description}, draw=none, from=1, to=0]
        \arrow["{\times_{(x,x)}}"{description, pos=0.6}, draw=none, from=0, to=2]
      \end{tikzcd},
    \end{equation*}
    and the projection cells satisfy
    \begin{equation} \label{eq:dbl-proj-external-compose}
      \begin{tikzcd}
        {x \times x'} & {y \times y'} & {z \times z'} \\
        {x \times x'} && {z \times z'} \\
        x && z
        \arrow[Rightarrow, no head, from=1-3, to=2-3]
        \arrow[Rightarrow, no head, from=1-1, to=2-1]
        \arrow[""{name=0, anchor=center, inner sep=0}, "{(m \odot n) \times (m' \odot n')}"', "\shortmid"{marking}, from=2-1, to=2-3]
        \arrow["{m \times m'}", "\shortmid"{marking}, from=1-1, to=1-2]
        \arrow["{n \times n'}", "\shortmid"{marking}, from=1-2, to=1-3]
        \arrow[""{name=1, anchor=center, inner sep=0}, "{m \odot n}"', "\shortmid"{marking}, from=3-1, to=3-3]
        \arrow["{\pi_{x,x'}}"', from=2-1, to=3-1]
        \arrow["{\pi_{z,z'}}", from=2-3, to=3-3]
        \arrow["{\times_{(m,m'),(n,n')}}"{description}, draw=none, from=1-2, to=0]
        \arrow["{\pi_{m \odot n, m' \odot n'}}"{description, pos=0.6}, draw=none, from=0, to=1]
      \end{tikzcd}
      \quad=\quad
      \begin{tikzcd}
        {x \times x'} & {y \times y'} & {z \times z'} \\
        x & y & z
        \arrow["{\pi_{x,x'}}"', from=1-1, to=2-1]
        \arrow[""{name=0, anchor=center, inner sep=0}, "{m \times m'}", "\shortmid"{marking}, from=1-1, to=1-2]
        \arrow[""{name=1, anchor=center, inner sep=0}, "{n \times n'}", "\shortmid"{marking}, from=1-2, to=1-3]
        \arrow[""{name=2, anchor=center, inner sep=0}, "m"', "\shortmid"{marking}, from=2-1, to=2-2]
        \arrow[""{name=3, anchor=center, inner sep=0}, "n"', "\shortmid"{marking}, from=2-2, to=2-3]
        \arrow["{\pi_{z,z'}}", from=1-3, to=2-3]
        \arrow["{\pi_{y,y'}}"{description}, from=1-2, to=2-2]
        \arrow["{\pi_{m,m'}}"{description}, draw=none, from=0, to=2]
        \arrow["{\pi_{n,n'}}"{description}, draw=none, from=1, to=3]
      \end{tikzcd}
    \end{equation}
    \begin{equation} \label{eq:dbl-proj-external-id}
      \begin{tikzcd}
        {x \times x'} & {x \times x'} \\
        {x \times x'} & {x \times x'} \\
        x & x
        \arrow["{\pi_{x,x'}}"', from=2-1, to=3-1]
        \arrow[""{name=0, anchor=center, inner sep=0}, "{\mathrm{id}_x}"', "\shortmid"{marking}, from=3-1, to=3-2]
        \arrow["{\pi_{x,x'}}", from=2-2, to=3-2]
        \arrow[""{name=1, anchor=center, inner sep=0}, "{\mathrm{id}_x \times \mathrm{id}_{x'}}"', "\shortmid"{marking}, from=2-1, to=2-2]
        \arrow[Rightarrow, no head, from=1-1, to=2-1]
        \arrow[Rightarrow, no head, from=1-2, to=2-2]
        \arrow[""{name=2, anchor=center, inner sep=0}, "{\mathrm{id}_{x \times x'}}", "\shortmid"{marking}, from=1-1, to=1-2]
        \arrow["{\pi_{\mathrm{id}_x, \mathrm{id}_{x'}}}"{description, pos=0.6}, draw=none, from=1, to=0]
        \arrow["{\times_{(x,x')}}"{description}, draw=none, from=2, to=1]
      \end{tikzcd}
      \quad=\quad
      \begin{tikzcd}
        {x \times x'} & {x \times x'} \\
        x & x
        \arrow[""{name=0, anchor=center, inner sep=0}, "{\mathrm{id}_x}"', "\shortmid"{marking}, from=2-1, to=2-2]
        \arrow["{\pi_{x,x'}}"', from=1-1, to=2-1]
        \arrow["{\pi_{x,x'}}", from=1-2, to=2-2]
        \arrow[""{name=1, anchor=center, inner sep=0}, "{\mathrm{id}_{x \times x'}}", "\shortmid"{marking}, from=1-1, to=1-2]
        \arrow["{\mathrm{id}_{\pi_{x,x'}}}"{description}, draw=none, from=1, to=0]
      \end{tikzcd}
    \end{equation}
    and similarly for the projections onto the second component.
\end{itemize}
When the comparison cells for the products and unit are all isomorphisms, the
double category is cartesian. Even in this case, it can be important to keep
track of the comparisons since the double functors
$\times: \dbl{D} \times \dbl{D} \to \dbl{D}$ and $I: \dbl{1} \to \dbl{D}$ are
generally not strict. This happens, for example, in the prototypical cartesian
double category $\Span$, where external composition and products exchange only
up to isomorphism, since limits commute with limits (particularly, pullbacks
commute with products) only up to canonical isomorphism.

Pairing in a precartesian double category $\dbl{D}$ is defined as usual by the
universal properties of the products in $\dbl{D}_0$ and $\dbl{D}_1$.
Alternatively, the \define{pairing} of two cells
$\inlineCell{x}{y}{w}{z}{m}{n}{f}{g}{\alpha}$ and
$\inlineCell{x}{y}{w'}{z'}{m}{n'}{f'}{g'}{\alpha'}$ with common domain is given
by the formula
\begin{equation*}
  \begin{tikzcd}
    x & y \\
    {w \times w'} & {z \times z'}
    \arrow[""{name=0, anchor=center, inner sep=0}, "{n \times n'}"', "\shortmid"{marking}, from=2-1, to=2-2]
    \arrow["{\langle f, f' \rangle}"', from=1-1, to=2-1]
    \arrow[""{name=1, anchor=center, inner sep=0}, "m", "\shortmid"{marking}, from=1-1, to=1-2]
    \arrow["{\langle g, g' \rangle}", from=1-2, to=2-2]
    \arrow["{\langle \alpha, \alpha' \rangle}"{description}, draw=none, from=1, to=0]
  \end{tikzcd}
  \quad=\quad
  \begin{tikzcd}
    x & y \\
    {x \times x} & {y \times y} \\
    {w \times w'} & {z \times z'}
    \arrow[""{name=0, anchor=center, inner sep=0}, "{n \times n'}"', "\shortmid"{marking}, from=3-1, to=3-2]
    \arrow["{f \times f'}"', from=2-1, to=3-1]
    \arrow["{g \times g'}", from=2-2, to=3-2]
    \arrow[""{name=1, anchor=center, inner sep=0}, "{m \times m}"', "\shortmid"{marking}, from=2-1, to=2-2]
    \arrow["{\Delta_x}"', from=1-1, to=2-1]
    \arrow["{\Delta_y}", from=1-2, to=2-2]
    \arrow[""{name=2, anchor=center, inner sep=0}, "m", "\shortmid"{marking}, from=1-1, to=1-2]
    \arrow["{\alpha \times \alpha'}"{description, pos=0.6}, draw=none, from=1, to=0]
    \arrow["{\Delta_m}"{description}, draw=none, from=2, to=1]
  \end{tikzcd}.
\end{equation*}

The definition of a cartesian double category is evidently motivated by the
famous theorem that a category $\cat{C}$ has finite products if and only if the
diagonal and terminal functors
$\Delta_{\cat{C}}: \cat{C} \to \cat{C} \times \cat{C}$ and
$!_{\cat{C}}: \cat{C} \to \cat{1}$ have right adjoints \cite[\mbox{Propositions
  8.2.1-2}]{reyes2004}. Cartesian double categories too can be characterized
using universal properties. This result significantly simplifies checking that a
double category is cartesian since the choices apparently involved in defining
the right adjoint double functors are not choices at all, but are uniquely
determined by the universal properties combined with the double functor axioms.
While not explicitly stated in Aleiferi's thesis, a series of similar but weaker
results are obtained \cite[Propositions 3.4.13, 3.4.16, and
4.1.2]{aleiferi2018}.

\begin{proposition}[Cartesian double categories via universal properties]
  A double category $\dbl{D}$ is precartesian if and only if
  \begin{itemize}[noitemsep]
    \item the categories $\dbl{D}_0$ and $\dbl{D}_1$ have finite products, and
    \item the source and target functors $s, t: \dbl{D}_1 \to \dbl{D}_0$
      preserve finite products.
  \end{itemize}
  In this case, $\dbl{D}$ is cartesian if and only if the external composition
  $\odot: \dbl{D}_1 \times_{\dbl{D}_0} \dbl{D}_1 \to \dbl{D}_1$ and identity
  $\id: \dbl{D}_0 \to \dbl{D}_1$ also preserve finite products, meaning that the
  canonical comparison cells
  \begin{equation} \label{eq:dbl-product-compose-comparison}
    \begin{tikzcd}[column sep=huge]
      {x \times x'} & {z \times z'} \\
      {x \times x'} & {z \times z'}
      \arrow[""{name=0, anchor=center, inner sep=0}, "{(m \times m') \odot (n \times n')}", "\shortmid"{marking}, from=1-1, to=1-2]
      \arrow[Rightarrow, no head, from=1-1, to=2-1]
      \arrow[Rightarrow, no head, from=1-2, to=2-2]
      \arrow[""{name=1, anchor=center, inner sep=0}, "{(m \odot n) \times (m' \odot n')}"', "\shortmid"{marking}, from=2-1, to=2-2]
      \arrow["{\times_{(m,m'),(n,n')}}"{description}, draw=none, from=0, to=1]
    \end{tikzcd}
    \quad\coloneqq\quad
    \begin{tikzcd}
      {x \times x'} & {y \times y'} & {z \times z'} \\
      {x \times x'} && {z \times z'}
      \arrow["{m \times m'}", "\shortmid"{marking}, from=1-1, to=1-2]
      \arrow["{n \times n'}", "\shortmid"{marking}, from=1-2, to=1-3]
      \arrow[""{name=0, anchor=center, inner sep=0}, "{(m \odot n) \times (m' \odot n')}"', "\shortmid"{marking}, from=2-1, to=2-3]
      \arrow[Rightarrow, no head, from=1-1, to=2-1]
      \arrow[Rightarrow, no head, from=1-3, to=2-3]
      \arrow["{\langle\pi_{m,m'} \odot \pi_{n,n'},\ \pi_{m,m'}' \odot \pi_{n,n'}'\rangle}"{description, pos=0.4}, draw=none, from=1-2, to=0]
    \end{tikzcd}
  \end{equation}
  and
  \begin{equation} \label{eq:dbl-product-id-comparison}
    \begin{tikzcd}[column sep=large]
      {x \times x'} & {x \times x'} \\
      {x \times x'} & {x \times x'}
      \arrow[Rightarrow, no head, from=1-2, to=2-2]
      \arrow[Rightarrow, no head, from=1-1, to=2-1]
      \arrow[""{name=0, anchor=center, inner sep=0}, "{\id_x \times \id_{x'}}"', "\shortmid"{marking}, from=2-1, to=2-2]
      \arrow[""{name=1, anchor=center, inner sep=0}, "{\id_{x \times x'}}", "\shortmid"{marking}, from=1-1, to=1-2]
      \arrow["{\times_{(x,x')}}"{description, pos=0.6}, draw=none, from=1, to=0]
    \end{tikzcd}
    \quad\coloneqq
    \begin{tikzcd}[column sep=large]
      {x \times x'} & {x \times x'} \\
      {x \times x'} & {x \times x'}
      \arrow[Rightarrow, no head, from=1-2, to=2-2]
      \arrow[Rightarrow, no head, from=1-1, to=2-1]
      \arrow[""{name=0, anchor=center, inner sep=0}, "{\id_x \times \id_{x'}}"', "\shortmid"{marking}, from=2-1, to=2-2]
      \arrow[""{name=1, anchor=center, inner sep=0}, "{\id_{x \times x'}}", "\shortmid"{marking}, from=1-1, to=1-2]
      \arrow["{\langle \id_{\pi_{x,x'}},\ \id_{\pi_{x,x'}'} \rangle}"{description, pos=0.6}, draw=none, from=1, to=0]
    \end{tikzcd}
  \end{equation}
  given by the universal property of products, as well as the comparisons
  $I_1 \odot I_1 \xto{!} I_1$ and $\id_{I_0} \xto{!} I_1$ given by the universal
  property of terminal objects, are all isomorphisms in $\dbl{D}_1$.
\end{proposition}
\begin{proof}
  This characterization of (pre)cartesian double categories follows from a
  general result about double adjunctions \cite[Corollary 4.3.7]{grandis2019}.
  To illustrate, we give a direct proof in one direction.

  If $\dbl{D}$ is precartesian, then the existence of adjunctions
  $\Delta_i \dashv \times_i$ and $!_i \dashv I_i$, $i=0,1$, in $\Cat$ implies
  that $\dbl{D}_0$ and $\dbl{D}_1$ have finite products. Moreover, since the
  functors $\times_i$ and $I_i$ assemble into double functors
  $\times = (\times_0, \times_1)$ and $I = (I_0, I_1)$, the source and target
  functors $s,t$ preserve finite products. By the universal property of products
  in $\dbl{D}_1$, the comparison cells labeled $\times_{(m,m'),(n,n')}$ and
  $\times_{(x,x')}$ in \cref{eq:dbl-product-comparisons} are uniquely determined
  by \cref{eq:dbl-proj-external-compose,eq:dbl-proj-external-id} and coincide
  with those in
  \cref{eq:dbl-product-compose-comparison,eq:dbl-product-id-comparison}.
  Similarly for the comparison cells labeled $\mu_I$ and $\eta_I$ in
  \cref{eq:dlb-terminal-comparisons}, by the universal property of terminal
  objects. In particular, if $\dbl{D}$ is cartesian, so that the right adjoints
  are pseudo, then all of these cells are isomorphisms.
\end{proof}

\begin{example}[Cartesian 2-categories]
  Any 2-category with finite (strict) 2-products, viewed as a double category
  with trivial proarrow structure, is a cartesian double category. Conversely,
  any cartesian double category has an underlying 2-category with finite
  products.
\end{example}

An \emph{equipment} is a double category in which any proarrow can be restricted
along a pair of incoming arrows or, dually, extended along a pair of outgoing
arrows, in universal ways. A great many commonly occurring double categories,
such as those of relations, spans, cospans, modules, and profunctors, are
equipments. Equipments appear in many guises and have also been called
\emph{proarrow equipments} \cite{wood1982}, \emph{framed bicategories}
\cite{shulman2008}, \emph{fibrant double categories} \cite{aleiferi2018}, and
\emph{gregarious double categories} \cite{dawson2010}. As shown by Shulman
\cite[Theorem 4.1]{shulman2008}, equipments can be defined in several equivalent
ways.

\begin{theorem}[Equipment] \label{thm:equipment}
  A double category $\dbl{D}$ is an \define{equipment} if any of the following
  equivalent statements hold:
  \begin{enumerate}[(i)]
    \item The source-target projection
      $\langle s, t \rangle: \dbl{D}_1 \to \dbl{D}_0 \times \dbl{D}_0$ is a
      fibration.
    \item The source-target projection
      $\langle s, t \rangle: \dbl{D}_1 \to \dbl{D}_0 \times \dbl{D}_0$ is an
      opfibration.
    \item Each arrow $f: x \to y$ in $\dbl{D}$ can be associated with a
      \define{companion} proarrow $f_!: x \proto y$ and a \define{conjoint}
      proarrow $f^*: y \proto x$, along with cells
      \begin{equation*}
        \begin{tikzcd}
          x & y \\
          y & y
          \arrow[""{name=0, anchor=center, inner sep=0}, "{\mathrm{id}_y}"', "\shortmid"{marking}, from=2-1, to=2-2]
          \arrow["f"', from=1-1, to=2-1]
          \arrow[""{name=1, anchor=center, inner sep=0}, "{f_!}", "\shortmid"{marking}, from=1-1, to=1-2]
          \arrow[Rightarrow, no head, from=1-2, to=2-2]
          \arrow["{\mathrm{res}}"{description}, draw=none, from=1, to=0]
        \end{tikzcd}
        \qquad
        \begin{tikzcd}
          x & x \\
          x & y
          \arrow[""{name=0, anchor=center, inner sep=0}, "{\mathrm{id}_x}", "\shortmid"{marking}, from=1-1, to=1-2]
          \arrow["f", from=1-2, to=2-2]
          \arrow[""{name=1, anchor=center, inner sep=0}, "{f_!}"', "\shortmid"{marking}, from=2-1, to=2-2]
          \arrow[Rightarrow, no head, from=1-1, to=2-1]
          \arrow["{\mathrm{ext}}"{description}, draw=none, from=0, to=1]
        \end{tikzcd}
        \qquad
        \begin{tikzcd}
          y & x \\
          y & y
          \arrow[""{name=0, anchor=center, inner sep=0}, "{\mathrm{id}_y}"', "\shortmid"{marking}, from=2-1, to=2-2]
          \arrow["f", from=1-2, to=2-2]
          \arrow[Rightarrow, no head, from=1-1, to=2-1]
          \arrow[""{name=1, anchor=center, inner sep=0}, "{f^*}", "\shortmid"{marking}, from=1-1, to=1-2]
          \arrow["{\mathrm{res}}"{description}, draw=none, from=1, to=0]
        \end{tikzcd}
        \qquad
        \begin{tikzcd}
          x & x \\
          y & x
          \arrow[""{name=0, anchor=center, inner sep=0}, "{\mathrm{id}_x}", "\shortmid"{marking}, from=1-1, to=1-2]
          \arrow[Rightarrow, no head, from=1-2, to=2-2]
          \arrow[""{name=1, anchor=center, inner sep=0}, "{f^*}"', "\shortmid"{marking}, from=2-1, to=2-2]
          \arrow["f"', from=1-1, to=2-1]
          \arrow["{\mathrm{ext}}"{description}, draw=none, from=0, to=1]
        \end{tikzcd}
      \end{equation*}
      satisfying the following equations.
      \begin{gather*}
        \begin{tikzcd}[ampersand replacement=\&]
          x \& x \\
          x \& y \\
          y \& y
          \arrow[""{name=0, anchor=center, inner sep=0}, "{\mathrm{id}_y}"', "\shortmid"{marking}, from=3-1, to=3-2]
          \arrow["f"', from=2-1, to=3-1]
          \arrow[""{name=1, anchor=center, inner sep=0}, "{f_!}", "\shortmid"{marking}, from=2-1, to=2-2]
          \arrow[Rightarrow, no head, from=2-2, to=3-2]
          \arrow[""{name=2, anchor=center, inner sep=0}, "{\mathrm{id}_x}", "\shortmid"{marking}, from=1-1, to=1-2]
          \arrow["f", from=1-2, to=2-2]
          \arrow[Rightarrow, no head, from=1-1, to=2-1]
          \arrow["{\mathrm{res}}"{description}, draw=none, from=1, to=0]
          \arrow["{\mathrm{ext}}"{description}, draw=none, from=2, to=1]
        \end{tikzcd}
        =
        \begin{tikzcd}[ampersand replacement=\&]
          x \& x \\
          y \& y
          \arrow["f"', from=1-1, to=2-1]
          \arrow["f", from=1-2, to=2-2]
          \arrow[""{name=0, anchor=center, inner sep=0}, "{\mathrm{id}_x}", "\shortmid"{marking}, from=1-1, to=1-2]
          \arrow[""{name=1, anchor=center, inner sep=0}, "{\mathrm{id}_y}"', "\shortmid"{marking}, from=2-1, to=2-2]
          \arrow["{\mathrm{id}_f}"{description}, draw=none, from=0, to=1]
        \end{tikzcd}
        \qquad\qquad
        \begin{tikzcd}[ampersand replacement=\&]
          x \& x \& y \\
          x \& y \& y
          \arrow[""{name=0, anchor=center, inner sep=0}, "{\mathrm{id}_y}"', "\shortmid"{marking}, from=2-2, to=2-3]
          \arrow["f"', from=1-2, to=2-2]
          \arrow[""{name=1, anchor=center, inner sep=0}, "{f_!}", "\shortmid"{marking}, from=1-2, to=1-3]
          \arrow[Rightarrow, no head, from=1-3, to=2-3]
          \arrow[""{name=2, anchor=center, inner sep=0}, "{\mathrm{id}_x}", "\shortmid"{marking}, from=1-1, to=1-2]
          \arrow[""{name=3, anchor=center, inner sep=0}, "{f_!}"', "\shortmid"{marking}, from=2-1, to=2-2]
          \arrow[Rightarrow, no head, from=1-1, to=2-1]
          \arrow["{\mathrm{res}}"{description}, draw=none, from=1, to=0]
          \arrow["{\mathrm{ext}}"{description}, draw=none, from=2, to=3]
        \end{tikzcd}
        =
        \begin{tikzcd}[ampersand replacement=\&]
          x \& y \\
          x \& y
          \arrow[""{name=0, anchor=center, inner sep=0}, "{f_!}"', "\shortmid"{marking}, from=2-1, to=2-2]
          \arrow[Rightarrow, no head, from=1-1, to=2-1]
          \arrow[""{name=1, anchor=center, inner sep=0}, "{f_!}", "\shortmid"{marking}, from=1-1, to=1-2]
          \arrow[Rightarrow, no head, from=1-2, to=2-2]
          \arrow["{1_{f_!}}"{description}, draw=none, from=1, to=0]
        \end{tikzcd}
        \\
        \begin{tikzcd}[ampersand replacement=\&]
          y \& x \& x \\
          y \& y \& x
          \arrow[""{name=0, anchor=center, inner sep=0}, "{\mathrm{id}_y}"', "\shortmid"{marking}, from=2-1, to=2-2]
          \arrow["f", from=1-2, to=2-2]
          \arrow[Rightarrow, no head, from=1-1, to=2-1]
          \arrow[""{name=1, anchor=center, inner sep=0}, "{f^*}", "\shortmid"{marking}, from=1-1, to=1-2]
          \arrow[""{name=2, anchor=center, inner sep=0}, "{f^*}"', "\shortmid"{marking}, from=2-2, to=2-3]
          \arrow[Rightarrow, no head, from=1-3, to=2-3]
          \arrow[""{name=3, anchor=center, inner sep=0}, "{\mathrm{id}_x}", "\shortmid"{marking}, from=1-2, to=1-3]
          \arrow["{\mathrm{res}}"{description}, draw=none, from=1, to=0]
          \arrow["{\mathrm{ext}}"{description}, draw=none, from=3, to=2]
        \end{tikzcd}
        =
        \begin{tikzcd}[ampersand replacement=\&]
          y \& x \\
          y \& x
          \arrow[""{name=0, anchor=center, inner sep=0}, "{f^*}"', "\shortmid"{marking}, from=2-1, to=2-2]
          \arrow[Rightarrow, no head, from=1-2, to=2-2]
          \arrow[Rightarrow, no head, from=1-1, to=2-1]
          \arrow[""{name=1, anchor=center, inner sep=0}, "{f^*}", "\shortmid"{marking}, from=1-1, to=1-2]
          \arrow["{1_{f^*}}"{description}, draw=none, from=1, to=0]
        \end{tikzcd}
        \qquad\quad
        \begin{tikzcd}[ampersand replacement=\&]
          x \& x \\
          y \& x \\
          y \& y
          \arrow[""{name=0, anchor=center, inner sep=0}, "{\mathrm{id}_x}", "\shortmid"{marking}, from=1-1, to=1-2]
          \arrow[Rightarrow, no head, from=1-2, to=2-2]
          \arrow[""{name=1, anchor=center, inner sep=0}, "{f^*}"', "\shortmid"{marking}, from=2-1, to=2-2]
          \arrow["f"', from=1-1, to=2-1]
          \arrow["f", from=2-2, to=3-2]
          \arrow[""{name=2, anchor=center, inner sep=0}, "{\mathrm{id}_y}"', "\shortmid"{marking}, from=3-1, to=3-2]
          \arrow[Rightarrow, no head, from=2-1, to=3-1]
          \arrow["{\mathrm{ext}}"{description}, draw=none, from=0, to=1]
          \arrow["{\mathrm{res}}"{description}, draw=none, from=1, to=2]
        \end{tikzcd}
        =
        \begin{tikzcd}[ampersand replacement=\&]
          x \& x \\
          y \& y
          \arrow["f"', from=1-1, to=2-1]
          \arrow["f", from=1-2, to=2-2]
          \arrow[""{name=0, anchor=center, inner sep=0}, "{\mathrm{id}_x}", "\shortmid"{marking}, from=1-1, to=1-2]
          \arrow[""{name=1, anchor=center, inner sep=0}, "{\mathrm{id}_y}"', "\shortmid"{marking}, from=2-1, to=2-2]
          \arrow["{\mathrm{id}_f}"{description}, draw=none, from=0, to=1]
        \end{tikzcd}
      \end{gather*}
  \end{enumerate}
\end{theorem}

Explicitly, condition (i) that the functor $\langle s,t \rangle$ be a fibration
means that every niche in $\dbl{D}$ of the form on the left can be completed to
a cell as on the right
\begin{equation*}
  \begin{tikzcd}
    x & y \\
    w & z
    \arrow["f"', from=1-1, to=2-1]
    \arrow["g", from=1-2, to=2-2]
    \arrow["n"', "\shortmid"{marking}, from=2-1, to=2-2]
  \end{tikzcd}
  \qquad\leadsto\qquad
  \begin{tikzcd}
    x & y \\
    w & z
    \arrow["f"', from=1-1, to=2-1]
    \arrow["g", from=1-2, to=2-2]
    \arrow[""{name=0, anchor=center, inner sep=0}, "n"', "\shortmid"{marking}, from=2-1, to=2-2]
    \arrow[""{name=1, anchor=center, inner sep=0}, "{n(f,g)}", "\shortmid"{marking}, from=1-1, to=1-2]
    \arrow["{\mathrm{res}}"{description}, draw=none, from=1, to=0]
  \end{tikzcd},
\end{equation*}
called a \define{cartesian} or \define{restriction} cell, with the universal
property that for any arrows $h: x' \to x$ and $k: y' \to y$ in $\dbl{D}$, each
cell of form $\inlineCell{x'}{y'}{w}{z}{m'}{n}{f \circ h}{g \circ k}{\alpha}$
factors uniquely through it:
\begin{equation*}
  \begin{tikzcd}
    {x'} & {y'} \\
    x & y \\
    w & z
    \arrow["f"', from=2-1, to=3-1]
    \arrow["g", from=2-2, to=3-2]
    \arrow[""{name=0, anchor=center, inner sep=0}, "n"', "\shortmid"{marking}, from=3-1, to=3-2]
    \arrow["h"', from=1-1, to=2-1]
    \arrow["k", from=1-2, to=2-2]
    \arrow[""{name=1, anchor=center, inner sep=0}, "{m'}", "\shortmid"{marking}, from=1-1, to=1-2]
    \arrow["\alpha"{description}, draw=none, from=1, to=0]
  \end{tikzcd}
  \quad=\quad
  \begin{tikzcd}
    {x'} & {y'} \\
    x & y \\
    w & z
    \arrow["f"', from=2-1, to=3-1]
    \arrow["g", from=2-2, to=3-2]
    \arrow[""{name=0, anchor=center, inner sep=0}, "n"', "\shortmid"{marking}, from=3-1, to=3-2]
    \arrow[""{name=1, anchor=center, inner sep=0}, "{n(f,g)}", "\shortmid"{marking}, from=2-1, to=2-2]
    \arrow["h"', from=1-1, to=2-1]
    \arrow["k", from=1-2, to=2-2]
    \arrow[""{name=2, anchor=center, inner sep=0}, "{m'}", "\shortmid"{marking}, from=1-1, to=1-2]
    \arrow["{\mathrm{res}}"{description}, draw=none, from=1, to=0]
    \arrow["{\exists !}"{description, pos=0.4}, draw=none, from=2, to=1]
  \end{tikzcd}.
\end{equation*}
The restricted proarrow $n(f,g)$ is unique up to unique globular isomorphism, as
can be seen by applying the universal property with $h = 1_x$ and $k = 1_y$.
Dually, condition (ii) that the functor $\langle s,t \rangle$ be an opfibration
means that every co-niche in $\dbl{D}$ of the form on the left can be completed
to a cell as on the right
\begin{equation*}
  \begin{tikzcd}
    x & y \\
    w & z
    \arrow["m", "\shortmid"{marking}, from=1-1, to=1-2]
    \arrow["f"', from=1-1, to=2-1]
    \arrow["g", from=1-2, to=2-2]
  \end{tikzcd}
  \qquad\leadsto\qquad
  \begin{tikzcd}
    x & y \\
    w & z
    \arrow[""{name=0, anchor=center, inner sep=0}, "m", "\shortmid"{marking}, from=1-1, to=1-2]
    \arrow["f"', from=1-1, to=2-1]
    \arrow["g", from=1-2, to=2-2]
    \arrow[""{name=1, anchor=center, inner sep=0}, "{(f,g)m}"', "\shortmid"{marking}, from=2-1, to=2-2]
    \arrow["{\mathrm{ext}}"{description}, draw=none, from=0, to=1]
  \end{tikzcd},
\end{equation*}
called an \define{opcartesian} or \define{extension} cell, with a universal
property dual to cartesian cells.

As suggested by the defining characterization of an equipment, restriction and
extension cells can be generated from companions and conjoints. Specifically,
restrictions and extensions are given by the composites
\begin{equation} \label{eq:restrictions-via-companions-conjoints}
  \begin{tikzcd}
    x & y \\
    w & z
    \arrow["f"', from=1-1, to=2-1]
    \arrow["g", from=1-2, to=2-2]
    \arrow[""{name=0, anchor=center, inner sep=0}, "n"', "\shortmid"{marking}, from=2-1, to=2-2]
    \arrow[""{name=1, anchor=center, inner sep=0}, "{n(f,g)}", "\shortmid"{marking}, from=1-1, to=1-2]
    \arrow["{\mathrm{res}}"{description}, draw=none, from=1, to=0]
  \end{tikzcd}
  \quad=\quad
  \begin{tikzcd}
    x & w & z & y \\
    w & w & z & z
    \arrow[""{name=0, anchor=center, inner sep=0}, "{\mathrm{id}_w}"', "\shortmid"{marking}, from=2-1, to=2-2]
    \arrow["f"', from=1-1, to=2-1]
    \arrow[""{name=1, anchor=center, inner sep=0}, "{f_!}", "\shortmid"{marking}, from=1-1, to=1-2]
    \arrow[Rightarrow, no head, from=1-2, to=2-2]
    \arrow[""{name=2, anchor=center, inner sep=0}, "n", "\shortmid"{marking}, from=1-2, to=1-3]
    \arrow[""{name=3, anchor=center, inner sep=0}, "n"', "\shortmid"{marking}, from=2-2, to=2-3]
    \arrow[Rightarrow, no head, from=1-3, to=2-3]
    \arrow["g", from=1-4, to=2-4]
    \arrow[""{name=4, anchor=center, inner sep=0}, "{\mathrm{id}_z}"', "\shortmid"{marking}, from=2-3, to=2-4]
    \arrow[""{name=5, anchor=center, inner sep=0}, "{g^*}", "\shortmid"{marking}, from=1-3, to=1-4]
    \arrow["{\mathrm{res}}"{description}, draw=none, from=1, to=0]
    \arrow["{1_n}"{description}, draw=none, from=2, to=3]
    \arrow["{\mathrm{res}}"{description}, draw=none, from=5, to=4]
  \end{tikzcd}
\end{equation}
and
\begin{equation} \label{eq:extensions-via-companions-conjoints}
  \begin{tikzcd}
    x & y \\
    w & z
    \arrow[""{name=0, anchor=center, inner sep=0}, "m", "\shortmid"{marking}, from=1-1, to=1-2]
    \arrow["f"', from=1-1, to=2-1]
    \arrow["g", from=1-2, to=2-2]
    \arrow[""{name=1, anchor=center, inner sep=0}, "{(f,g)m}"', "\shortmid"{marking}, from=2-1, to=2-2]
    \arrow["{\mathrm{ext}}"{description}, draw=none, from=0, to=1]
  \end{tikzcd}
  \quad=\quad
  \begin{tikzcd}
    x & x & y & y \\
    w & x & y & z
    \arrow[""{name=0, anchor=center, inner sep=0}, "m", "\shortmid"{marking}, from=1-2, to=1-3]
    \arrow[Rightarrow, no head, from=1-2, to=2-2]
    \arrow[Rightarrow, no head, from=1-3, to=2-3]
    \arrow[""{name=1, anchor=center, inner sep=0}, "m"', "\shortmid"{marking}, from=2-2, to=2-3]
    \arrow["f"', from=1-1, to=2-1]
    \arrow[""{name=2, anchor=center, inner sep=0}, "{\mathrm{id}_x}", "\shortmid"{marking}, no head, from=1-1, to=1-2]
    \arrow[""{name=3, anchor=center, inner sep=0}, "{\mathrm{id}_y}", "\shortmid"{marking}, from=1-3, to=1-4]
    \arrow["g", from=1-4, to=2-4]
    \arrow[""{name=4, anchor=center, inner sep=0}, "{f^*}"', "\shortmid"{marking}, from=2-1, to=2-2]
    \arrow[""{name=5, anchor=center, inner sep=0}, "{g_!}"', "\shortmid"{marking}, from=2-3, to=2-4]
    \arrow["{1_m}"{description}, draw=none, from=0, to=1]
    \arrow["{\mathrm{ext}}"{description}, draw=none, from=2, to=4]
    \arrow["{\mathrm{ext}}"{description}, draw=none, from=3, to=5]
  \end{tikzcd}.
\end{equation}
These formulas, part of the proof of \cref{thm:equipment}, are important in
their own right. The cases $\id_w(f,1_w) = f_!$ and $\id_z(1_z,g) = g^*$ are
especially useful.

Equipments, double functors, and natural transformations form a 2-category
$\Eqp$. Similarly, there are 2-categories $\EqpLax$ and $\EqpLaxNormal$ having
lax functors and normal lax functors, respectively, as morphisms. The absence of
extra conditions on double functors and transformations between equipments will
be explained in \cref{sec:lax-functors}.

Having defined cartesian double categories and equipments, there are no
surprises in the definition of a cartesian equipment.

\begin{definition}[Cartesian equipment]
  A \define{precartesian equipment} is a cartesian object in $\EqpLax$.
  Similarly, a \define{cartesian equipment} is a cartesian object in $\Eqp$.

  In other words, a (pre)cartesian equipment is a double category that is both
  (pre)cartesian and an equipment.
\end{definition}

Our two main semantics for double theories, the double categories of spans and
matrices, are cartesian equipments under additional assumptions.

\begin{example}[Spans]
  When $\cat{S}$ is a category with finite limits, the double category
  $\Span{\cat{S}}$ of spans in $\cat{S}$ from \cref{ex:spans} is a cartesian
  equipment.
\end{example}

\begin{example}[Matrices]
  For any (infinitary) distributive monoidal category $\catV$, the double
  category $\Mat{\catV}$ of $\catV$-matrices from \cref{ex:matrices} is an
  equipment \cite[Proposition 4.1]{vasilakopoulou2019}. When $\catV$ is an
  (infinitary) distributive category, i.e., its monoidal product is cartesian,
  $\Mat{\catV}$ is also a cartesian double category \cite[Proposition
  4.2.5]{aleiferi2018}, hence is a cartesian equipment.
\end{example}

\begin{example}[Relations]
  Another example of a cartesian equipment is $\Rel$, the double category of
  relations. Relations are the special case of $\catV$-matrices where
  $\catV = \{\bot \to \top\}$ is the poset of booleans. More generally,
  $\Rel{\cat{S}}$, the double category of relations in a regular category
  $\cat{S}$, is a cartesian equipment \cite{lambert2022}.
\end{example}

The next example of a cartesian double category is less significant as a
semantics for double theories but is occasionally useful. It is usually not an
equipment as it lacks conjoints.

\begin{example}[Quintets] \label{ex:quintets}
  For any 2-category $\bicat{C}$, the \define{quintet construction}
  $\Quintet(\bicat{C})$ is the strict double category whose objects are the
  objects of $\bicat{C}$ and whose arrows and proarrows are morphisms of
  $\bicat{C}$. A cell in $\Quintet(\bicat{C})$ as on the left
  \begin{equation*}
    \begin{tikzcd}
      x & y \\
      w & z
      \arrow["f"', from=1-1, to=2-1]
      \arrow[""{name=0, anchor=center, inner sep=0}, "h", "\shortmid"{marking}, from=1-1, to=1-2]
      \arrow["g", from=1-2, to=2-2]
      \arrow[""{name=1, anchor=center, inner sep=0}, "k"', "\shortmid"{marking}, from=2-1, to=2-2]
      \arrow["\alpha"{description}, draw=none, from=0, to=1]
    \end{tikzcd}
    \qquad\leftrightsquigarrow\qquad
    \begin{tikzcd}
      x & y \\
      w & z
      \arrow["f"', from=1-1, to=2-1]
      \arrow["h", from=1-1, to=1-2]
      \arrow["g", from=1-2, to=2-2]
      \arrow["k"', from=2-1, to=2-2]
      \arrow["\alpha"', shorten <=4pt, shorten >=4pt, Rightarrow, from=1-2, to=2-1]
    \end{tikzcd}
  \end{equation*}
  is a 2-morphism in $\bicat{C}$ as on the right. Then $\Quintet(\bicat{C})$ is
  a cartesian double category precisely when $\bicat{C}$ is a cartesian
  2-category, i.e., $\bicat{C}$ has finite 2-products; see
  \cite[\S{6.1}]{grandis1999} or \cite[\S{C5.11}]{grandis2019}. The orientation
  of the 2-cells in quintets is a matter of convention as the 2-cell dual
  $\bicat{C}^\co$ is cartesian whenever $\bicat{C}$ is.
\end{example}

Under reasonable conditions on a cartesian equipment $\dbl{E}$, the double
category $\Module{\dbl{E}}$ of bimodules in $\dbl{E}$ is again a cartesian
equipment. It first needs to be seen that $\Module{\dbl{E}}$ is a double
category at all. For composites of bimodules to exist, the base double category
$\dbl{E}$ must have local coequalizers \cite[Definition 11.4]{shulman2008}. We
review external composition of bimodules in some detail as it will help with
later calculations.

\begin{definition}[Local coequalizers] \label{def:localcoequalizersinadoublecategory}
  A double category $\dbl{D}$ has \define{local coequalizers} if each
  hom-category $\dbl{D}(x,y)$ has coequalizers that are preserved by external
  composition in each argument.

  Let $\EqpLax^q$ and $\EqpLaxNormal^q$ denote the 2-categories of equipments
  with local coequalizers, (normal) lax functors, and natural transformations.
\end{definition}

Suppose that $\dbl{E}$ is a double category with local coequalizers. Given two
composable bimodules $m\colon a\proto b$ and $n\colon b\proto c$ between
categories $a\colon x\proto x$, $b\colon y\proto y$, and $c\colon z \proto z$ in
$\dbl{E}$, their external composite is defined as the coequalizer in
$\dbl{E}(x,z)$ of the action cells
\begin{equation*}
  \begin{tikzcd}
    x & y & y & z \\
    x && y & z
    \arrow["m", "\shortmid"{marking}, from=1-1, to=1-2]
    \arrow["b", "\shortmid"{marking}, from=1-2, to=1-3]
    \arrow[Rightarrow, no head, from=1-3, to=2-3]
    \arrow[Rightarrow, no head, from=1-1, to=2-1]
    \arrow[""{name=0, anchor=center, inner sep=0}, "m"', "\shortmid"{marking}, from=2-1, to=2-3]
    \arrow[""{name=1, anchor=center, inner sep=0}, "n", "\shortmid"{marking}, from=1-3, to=1-4]
    \arrow[Rightarrow, no head, from=1-4, to=2-4]
    \arrow[""{name=2, anchor=center, inner sep=0}, "n"', "\shortmid"{marking}, from=2-3, to=2-4]
    \arrow["\rho"{description, pos=0.4}, draw=none, from=1-2, to=0]
    \arrow["{1_n}"{description}, draw=none, from=1, to=2]
  \end{tikzcd}
  \qquad\text{and}\qquad
  \begin{tikzcd}
    x & y & y & z \\
    x & y && z
    \arrow[""{name=0, anchor=center, inner sep=0}, "m", "\shortmid"{marking}, from=1-1, to=1-2]
    \arrow["b", "\shortmid"{marking}, from=1-2, to=1-3]
    \arrow["n", "\shortmid"{marking}, from=1-3, to=1-4]
    \arrow[Rightarrow, no head, from=1-1, to=2-1]
    \arrow[""{name=1, anchor=center, inner sep=0}, "m"', "\shortmid"{marking}, from=2-1, to=2-2]
    \arrow[Rightarrow, no head, from=1-2, to=2-2]
    \arrow[""{name=2, anchor=center, inner sep=0}, "n"', "\shortmid"{marking}, from=2-2, to=2-4]
    \arrow[Rightarrow, no head, from=1-4, to=2-4]
    \arrow["{1_m}"{description}, draw=none, from=0, to=1]
    \arrow["\lambda"{description, pos=0.4}, draw=none, from=1-3, to=2]
  \end{tikzcd}
\end{equation*}
as displayed by
\begin{equation} \label{equation:modulecompositionascoequalizer}
  \begin{tikzcd}
    {m\odot b\odot n} & {m\odot n} & {m \otimes n}
    \arrow["1\odot\lambda"', shift right, Rightarrow, from=1-1, to=1-2]
    \arrow["{\rho\odot 1}", shift left, Rightarrow, from=1-1, to=1-2]
    \arrow["{\text{coeq}}", Rightarrow, from=1-2, to=1-3]
  \end{tikzcd}.
\end{equation}
To avoid confusion with external composition in $\dbl{E}$, we are using the
tensor symbol ``$\otimes$'' to denote external composition in
$\Module{\dbl{E}}$.

Now suppose that the double category $\dbl{E}$ is also an equipment, which will
be used to define external composition of cells in $\Module{\dbl{E}}$. We review
this construction here as we will use it in proofs below, although we are just
expanding the proof of \cite[Proposition 11.10]{shulman2008}. Given cells
$\alpha$ and $\beta$, first form the restriction along the given external source
and target; then there is induced a unique globular cell (I) as in the diagram
\begin{equation*}
  \begin{tikzcd}
    \cdot & \cdot & \cdot \\
    \cdot & \cdot & \cdot \\
    \cdot && \cdot
    \arrow[""{name=0, anchor=center, inner sep=0}, "m", "\shortmid"{marking}, from=1-1, to=1-2]
    \arrow["g"', from=1-2, to=2-2]
    \arrow["f"', from=1-1, to=2-1]
    \arrow[""{name=1, anchor=center, inner sep=0}, "n", "\shortmid"{marking}, from=1-2, to=1-3]
    \arrow["h", from=1-3, to=2-3]
    \arrow[""{name=2, anchor=center, inner sep=0}, "p"', "\shortmid"{marking}, from=2-1, to=2-2]
    \arrow[""{name=3, anchor=center, inner sep=0}, "q"', "\shortmid"{marking}, from=2-2, to=2-3]
    \arrow[Rightarrow, no head, from=2-1, to=3-1]
    \arrow[""{name=4, anchor=center, inner sep=0}, "{p\otimes q}"', "\shortmid"{marking}, from=3-1, to=3-3]
    \arrow[Rightarrow, no head, from=2-3, to=3-3]
    \arrow["{\mathrm{coeq}}"{description}, draw=none, from=2-2, to=4]
    \arrow["\alpha"{description}, draw=none, from=0, to=2]
    \arrow["\beta"{description}, draw=none, from=1, to=3]
  \end{tikzcd}
  \quad=\quad
  \begin{tikzcd}
    \cdot & \cdot \\
    \cdot & \cdot \\
    \cdot & \cdot
    \arrow[Rightarrow, no head, from=1-1, to=2-1]
    \arrow[""{name=0, anchor=center, inner sep=0}, "\shortmid"{marking}, from=2-1, to=2-2]
    \arrow[""{name=1, anchor=center, inner sep=0}, "{m\odot n}", "\shortmid"{marking}, from=1-1, to=1-2]
    \arrow[Rightarrow, no head, from=1-2, to=2-2]
    \arrow["f"', from=2-1, to=3-1]
    \arrow[""{name=2, anchor=center, inner sep=0}, "{p\otimes q}"', "\shortmid"{marking}, from=3-1, to=3-2]
    \arrow["h", from=2-2, to=3-2]
    \arrow["{\text{(I)}}"{description}, draw=none, from=1, to=0]
    \arrow["\res"{description}, draw=none, from=0, to=2]
  \end{tikzcd}.
\end{equation*}
But, by the equivariance of the cells $\alpha$ and $\beta$, the cell (I)
coequalizes the actions forming $m\otimes n$. So, there exists a further unique
globular cell (II) factoring (I) through this coequalizer. Take the external
composite $\alpha\otimes \beta$ to be the internal composite of (II) with the
restriction cell:
\begin{equation*}
  \begin{tikzcd}
    \cdot & \cdot \\
    \cdot & \cdot
    \arrow[""{name=0, anchor=center, inner sep=0}, "{m\otimes n}", "\shortmid"{marking}, from=1-1, to=1-2]
    \arrow[""{name=1, anchor=center, inner sep=0}, "{p\otimes q}"', "\shortmid"{marking}, from=2-1, to=2-2]
    \arrow["f"', from=1-1, to=2-1]
    \arrow["h", from=1-2, to=2-2]
    \arrow["{\alpha \otimes \beta}"{description}, draw=none, from=0, to=1]
  \end{tikzcd}
  \quad\coloneqq\quad
  \begin{tikzcd}
    \cdot & \cdot \\
    \cdot & \cdot \\
    \cdot & \cdot
    \arrow[""{name=0, anchor=center, inner sep=0}, "{m\otimes n}", "\shortmid"{marking}, from=1-1, to=1-2]
    \arrow[Rightarrow, no head, from=1-1, to=2-1]
    \arrow[""{name=1, anchor=center, inner sep=0}, "\shortmid"{marking}, from=2-1, to=2-2]
    \arrow["f"', from=2-1, to=3-1]
    \arrow[""{name=2, anchor=center, inner sep=0}, "{p\otimes q}"', "\shortmid"{marking}, from=3-1, to=3-2]
    \arrow["h", from=2-2, to=3-2]
    \arrow[Rightarrow, no head, from=1-2, to=2-2]
    \arrow["\res"{description}, draw=none, from=1, to=2]
    \arrow["{\text{(II)}}"{description}, draw=none, from=0, to=1]
  \end{tikzcd}.
\end{equation*}
Now, since (II) factors (I) through the coequalizer giving $m \otimes n$, this
means that the external composite $\alpha\otimes \beta$ satisfies the equation
\begin{equation} \label{equation:auxiliaryequationforexternalcompositeofmodulations}
  \begin{tikzcd}
    \cdot & \cdot \\
    \cdot & \cdot \\
    \cdot & \cdot
    \arrow[""{name=0, anchor=center, inner sep=0}, "{m\odot n}", "\shortmid"{marking}, from=1-1, to=1-2]
    \arrow[Rightarrow, no head, from=1-1, to=2-1]
    \arrow[""{name=1, anchor=center, inner sep=0}, "{m\otimes n}"', "\shortmid"{marking}, from=2-1, to=2-2]
    \arrow["f"', from=2-1, to=3-1]
    \arrow[""{name=2, anchor=center, inner sep=0}, "{p\otimes q}"', "\shortmid"{marking}, from=3-1, to=3-2]
    \arrow["h", from=2-2, to=3-2]
    \arrow[Rightarrow, no head, from=1-2, to=2-2]
    \arrow["\alpha\otimes\beta"{description, pos=0.6}, draw=none, from=1, to=2]
    \arrow["{\mathrm{coeq}}"{description}, draw=none, from=0, to=1]
  \end{tikzcd}
  \quad=\quad
  \begin{tikzcd}
    \cdot & \cdot & \cdot \\
    \cdot & \cdot & \cdot \\
    \cdot && \cdot
    \arrow[""{name=0, anchor=center, inner sep=0}, "m", "\shortmid"{marking}, from=1-1, to=1-2]
    \arrow[from=1-2, to=2-2]
    \arrow["f"', from=1-1, to=2-1]
    \arrow[""{name=1, anchor=center, inner sep=0}, "p"', "\shortmid"{marking}, from=2-1, to=2-2]
    \arrow[Rightarrow, no head, from=2-1, to=3-1]
    \arrow[""{name=2, anchor=center, inner sep=0}, "{p\otimes q}"', "\shortmid"{marking}, from=3-1, to=3-3]
    \arrow[""{name=3, anchor=center, inner sep=0}, "n", "\shortmid"{marking}, from=1-2, to=1-3]
    \arrow["h", from=1-3, to=2-3]
    \arrow[""{name=4, anchor=center, inner sep=0}, "q"', "\shortmid"{marking}, from=2-2, to=2-3]
    \arrow[Rightarrow, no head, from=2-3, to=3-3]
    \arrow["\alpha"{description}, draw=none, from=0, to=1]
    \arrow["\beta"{description}, draw=none, from=3, to=4]
    \arrow["{\mathrm{coeq}}"{description}, draw=none, from=2-2, to=2]
  \end{tikzcd}.
\end{equation}
This equation will be used in one of the proofs below. It is thus worth noting
that even if we are unconcerned with the fact that $\Module{\dbl{E}}$ is an
equipment, the fact that $\dbl{E}$ is one is used in showing that bimodules form
at least a double category.

We have described the construction behind the following result. For further
details, see \cite[Theorem 11.5 and Proposition 11.10]{shulman2008}.

\begin{lemma} \label{lemma:MODisanequipmentwithlocalcoequalizers}
  When $\dbl{E}$ is an equipment with local coequalizers, $\Module{\dbl{E}}$ is
  an equipment and has local coequalizers too.
\end{lemma}

Under these hypotheses, the property of being cartesian also carries over from
the base double category to the double category of bimodules.

\begin{proposition} \label{prop:MODiscartesian}
  If $\dbl{E}$ is a cartesian equipment with local coequalizers, then
  $\Module{\dbl{E}}$ is a cartesian equipment with local coequalizers.
\end{proposition}
\begin{proof}
  For an abstract proof, taking modules defines a 2-functor
  $\Module\colon \EqpLax^{q} \to \EqpLaxNormal^{q}$ from the 2-category of
  equipments with local coequalizers to the 2-category of equipments with local
  coequalizers and \emph{normal} lax functors between them \cite[Proposition
  11.11]{shulman2008}. Now, since cartesian double categories are defined by a
  2-adjunction and any 2-functor preserves such 2-adjunctions, the result
  follows. However, it is also possible to give a concrete \emph{ground-level}
  proof, constructing the required products by hand and showing that they have
  the right universal properties and that they are double-categorically
  coherent. This proof is actually the preferred one, since it tells us how to
  compute such products. The product of category objects is the product in
  $\dbl{E}$ of the underlying objects; likewise, the product of bimodules is the
  product in $\dbl{E}$ of their underlying proarrows. That these work as
  products in $\Module{\dbl{E}}_0$ and $\Module{\dbl{E}}_1$ is tedious but
  straightforward to check using the existing product structure in $\dbl{E}$.
  The only tricky part is well-definition, which in the case of bimodules
  requires defining actions over restrictions of product categories along the
  diagonals.
\end{proof}

\begin{lemma}[Closure properties of restrictions] \label{lem:closure-restrictions}
  Restriction cells in a double category $\dbl{D}$ (possibly, but not
  necessarily, an equipment) satisfy the following closure properties.
  \begin{enumerate}[(i),noitemsep]
    \item Isomorphisms: all isomorphisms in $\dbl{D}_1$ are restriction cells.
    \item Internal composites: restriction cells are closed under composition in
    $\dbl{D}_1$.
    \item External composites of companions and conjoints: for any
      $x \xto{f} y \xto{g} z$ in $\dbl{D}$, cells of the form below are
      restrictions.
      \begin{equation*}
        \begin{tikzcd}
          x & y & z \\
          y & y \\
          z & z & z
          \arrow["f"', from=1-1, to=2-1]
          \arrow[""{name=0, anchor=center, inner sep=0}, "{f_!}", "\shortmid"{marking}, from=1-1, to=1-2]
          \arrow[Rightarrow, no head, from=1-2, to=2-2]
          \arrow["g", from=2-2, to=3-2]
          \arrow["g"', from=2-1, to=3-1]
          \arrow[""{name=1, anchor=center, inner sep=0}, "{\id_z}"', "\shortmid"{marking}, from=3-1, to=3-2]
          \arrow[""{name=2, anchor=center, inner sep=0}, "{\id_y}", "\shortmid"{marking}, from=2-1, to=2-2]
          \arrow[""{name=3, anchor=center, inner sep=0}, "{g_!}", "\shortmid"{marking}, from=1-2, to=1-3]
          \arrow[Rightarrow, no head, from=1-3, to=3-3]
          \arrow[""{name=4, anchor=center, inner sep=0}, "{\id_z}"', "\shortmid"{marking}, from=3-2, to=3-3]
          \arrow["{\id_g}"{description}, draw=none, from=2, to=1]
          \arrow["\res"{description, pos=0.4}, draw=none, from=0, to=2]
          \arrow["\res"{description}, draw=none, from=3, to=4]
        \end{tikzcd}
        \qquad\qquad
        \begin{tikzcd}
          z & y & x \\
          & y & y \\
          z & z & z
          \arrow[""{name=0, anchor=center, inner sep=0}, "{g^*}", "\shortmid"{marking}, from=1-1, to=1-2]
          \arrow[""{name=1, anchor=center, inner sep=0}, "{f^*}", "\shortmid"{marking}, from=1-2, to=1-3]
          \arrow[""{name=2, anchor=center, inner sep=0}, "{\id_z}"', "\shortmid"{marking}, from=3-1, to=3-2]
          \arrow[""{name=3, anchor=center, inner sep=0}, "{\id_z}"', "\shortmid"{marking}, from=3-2, to=3-3]
          \arrow["f", from=1-3, to=2-3]
          \arrow[Rightarrow, no head, from=1-2, to=2-2]
          \arrow[""{name=4, anchor=center, inner sep=0}, "{\id_y}", "\shortmid"{marking}, from=2-2, to=2-3]
          \arrow["g", from=2-3, to=3-3]
          \arrow["g"', from=2-2, to=3-2]
          \arrow[Rightarrow, no head, from=1-1, to=3-1]
          \arrow["\res"{description, pos=0.4}, draw=none, from=1, to=4]
          \arrow["{\id_g}"{description}, draw=none, from=4, to=3]
          \arrow["\res"{description}, draw=none, from=0, to=2]
        \end{tikzcd}
      \end{equation*}
    \item General restrictions via companions and conjoints: external composites
      of the form \eqref{eq:restrictions-via-companions-conjoints} are
      restriction cells.
    \item When $\dbl{D}$ is a precartesian double category, restriction cells
      are closed under finite products.
  \end{enumerate}
\end{lemma}
\begin{proof}
  Statements (i) and (ii) are general facts about cartesian morphisms with
  respect to a functor \cite[Proposition 9.1.4]{johnson2021}, applied to the
  functor $\langle s,t \rangle: \dbl{D}_1 \to \dbl{D}_0 \times \dbl{D}_0$.
  Statement (iii) is \cite[Lemma 3.13]{shulman2010} and (iv) is proved in the
  course of \cite[Theorem 4.1]{shulman2008}. Statement (v) is a general fact
  about cartesian morphisms with respect to a cartesian functor; alternatively,
  under the assumption that $\dbl{D}$ is a precartesian equipment, it is
  \cite[Lemma 4.3.1]{aleiferi2018}.
\end{proof}

\section{Lax functors into cartesian equipments}
\label{sec:lax-functors-into-cartesian-equipments}

Bringing together the threads of
\cref{sec:lax-functors,sec:cartesian-equipments}, we turn to how lax double
functors interact with the extra structure present in cartesian double
categories and equipments.

Just as cartesian double categories can be defined as cartesian objects in the
2-category $\Dbl$, the concept of a \emph{cartesian} lax functor can be
extracted from the general notion of a cartesian morphism between cartesian
objects \cite[\S 5.2]{carboni1991}.

\begin{definition}[Cartesian lax functor] \label{def:cartesian-lax-functor}
  Let $\dbl{D}$ and $\dbl{E}$ be precartesian double categories (which could be
  cartesian). A lax double functor $F: \dbl{D} \to \dbl{E}$ is
  \define{cartesian} or \define{preserves finite products} if it is a cartesian
  morphism between $\dbl{D}$ and $\dbl{E}$, viewing them as cartesian objects in
  $\DblLax$.
\end{definition}

Spelling this out, if we take the mates of the identity transformations
$\Delta_{\dbl{E}} \circ F \To (F \times F) \circ \Delta_{\dbl{D}}$ and
$!_{\dbl{E}} \circ F \To !_{\dbl{D}}$, we obtain natural transformations
\begin{equation*}
  \Phi: F \circ \times_{\dbl{D}} \To \times_{\dbl{E}} \circ (F \times F):
  \dbl{D} \times \dbl{D} \to \dbl{E}
  \qquad\text{and}\qquad
  \phi: F \circ I_{\dbl{D}} \To I_{\dbl{E}}: \dbl{1} \to \dbl{E}
\end{equation*}
with components
\begin{equation*}
  \begin{tikzcd}
    {F(x \times x')} & {F(y \times y')} \\
    {Fx \times Fx'} & {Fy \times Fy'}
    \arrow["{\Phi_{x,x'}}"', from=1-1, to=2-1]
    \arrow[""{name=0, anchor=center, inner sep=0}, "{F(m \times m')}", "\shortmid"{marking}, from=1-1, to=1-2]
    \arrow[""{name=1, anchor=center, inner sep=0}, "{Fm \times Fm'}"', "\shortmid"{marking}, from=2-1, to=2-2]
    \arrow["{\Phi_{y,y'}}", from=1-2, to=2-2]
    \arrow["{\Phi_{m,m'}}"{description}, draw=none, from=0, to=1]
  \end{tikzcd}
  \quad\coloneqq\quad
  \begin{tikzcd}
    {F(x \times x')} && {F(y \times y')} \\
    {Fx \times Fx'} && {Fy \times Fy'}
    \arrow["{\langle F\pi_{x,x'}, F\pi_{x,x'}' \rangle}"', from=1-1, to=2-1]
    \arrow[""{name=0, anchor=center, inner sep=0}, "{F(m \times m')}", "\shortmid"{marking}, from=1-1, to=1-3]
    \arrow[""{name=1, anchor=center, inner sep=0}, "{Fm \times Fm'}"', "\shortmid"{marking}, from=2-1, to=2-3]
    \arrow["{\langle F\pi_{y,y'}, F\pi_{y,y'}' \rangle}", from=1-3, to=2-3]
    \arrow["{\langle F\pi_{m,m'}, F\pi_{m,m'}' \rangle}"{description}, draw=none, from=0, to=1]
  \end{tikzcd}
\end{equation*}
for proarrows $m: x \proto y$ and $m': x' \proto y'$ in $\dbl{D}$, and
\begin{equation*}
  \begin{tikzcd}
    {F(I_0)} & {F(I_0)} \\
    {I_0} & {I_0}
    \arrow["{\phi_0}"', from=1-1, to=2-1]
    \arrow["{\phi_0}", from=1-2, to=2-2]
    \arrow[""{name=0, anchor=center, inner sep=0}, "{F(I_1)}", "\shortmid"{marking}, from=1-1, to=1-2]
    \arrow[""{name=1, anchor=center, inner sep=0}, "{I_1}"', "\shortmid"{marking}, from=2-1, to=2-2]
    \arrow["{\phi_{\mathrm{id}_0}}"{description}, draw=none, from=0, to=1]
  \end{tikzcd}
  \quad\coloneqq\quad
  \begin{tikzcd}
    {F(I_0)} & {F(I_0)} \\
    {I_0} & {I_0}
    \arrow["{!}"', from=1-1, to=2-1]
    \arrow["{!}", from=1-2, to=2-2]
    \arrow[""{name=0, anchor=center, inner sep=0}, "{F(I_1)}", "\shortmid"{marking}, from=1-1, to=1-2]
    \arrow[""{name=1, anchor=center, inner sep=0}, "{I_1}"', "\shortmid"{marking}, from=2-1, to=2-2]
    \arrow["{!}"{description}, draw=none, from=0, to=1]
  \end{tikzcd}.
\end{equation*}
The lax double functor $F: \dbl{D} \to \dbl{E}$ is \define{cartesian} if both
natural transformations $\Phi$ and $\phi$ are natural isomorphisms, i.e., their
components are isomorphisms in $\dbl{E}$. This, in turn, is equivalent to both
underlying functors $F_0: \dbl{D}_0 \to \dbl{E}_0$ and
$F_1: \dbl{D}_1 \to \dbl{E}_1$ preserving finite products in the ordinary sense.
The property of a lax functor being cartesian thus reduces to a simple criterion
that is easily checked in examples.

The laxators and unitors of a cartesian lax functor preserve products in the
sense that they commute with products up to the product comparison cells. They
also preserve terminal objects, hence all finite products, although for the sake
of brevity we will not spell that out.

\begin{lemma}[Laxators and unitors for products] \label{lem:cartesian-laxators-unitors}
  Let $F: \dbl{D} \to \dbl{E}$ be a lax double functor between precartesian
  double categories $\dbl{D}$ and $\dbl{E}$. Then for any proarrows
  \mbox{$x \xproto{m} y \xproto{n} z$} and
  \mbox{$x' \xproto{m'} y' \xproto{n'} z'$} in $\dbl{D}$, we have
  \begin{equation*}
    \begin{tikzcd}
      {F(x \times x')} & {F(y \times y')} & {F(z \times z')} \\
      {F(x \times x')} && {F(z \times z')} \\
      {F(x \times x')} && {F(z \times z')} \\
      {Fx \times Fx'} && {Fz \times Fz'}
      \arrow["{F(m \times m')}", "\shortmid"{marking}, from=1-1, to=1-2]
      \arrow["{F(n \times n')}", "\shortmid"{marking}, from=1-2, to=1-3]
      \arrow[""{name=0, anchor=center, inner sep=0}, "{F((m \times m') \odot (n \times n'))}"', "\shortmid"{marking}, from=2-1, to=2-3]
      \arrow[Rightarrow, no head, from=1-1, to=2-1]
      \arrow[Rightarrow, no head, from=1-3, to=2-3]
      \arrow[""{name=1, anchor=center, inner sep=0}, "{F((m \odot n) \times (m' \odot n'))}"', "\shortmid"{marking}, from=3-1, to=3-3]
      \arrow[Rightarrow, no head, from=2-1, to=3-1]
      \arrow[Rightarrow, no head, from=2-3, to=3-3]
      \arrow["{\Phi_{x,x'}}"', from=3-1, to=4-1]
      \arrow["{\Phi_{z,z'}}", from=3-3, to=4-3]
      \arrow[""{name=2, anchor=center, inner sep=0}, "{F(m \odot n) \times F(m' \odot n')}"', "\shortmid"{marking}, from=4-1, to=4-3]
      \arrow["{F_{m \times m', n \times n'}}"{description}, draw=none, from=1-2, to=0]
      \arrow["{F \times_{(m,m'),(n,n')}}"{description, pos=0.6}, draw=none, from=0, to=1]
      \arrow["{\Phi_{m \odot n, m' \odot n'}}"{description, pos=0.6}, draw=none, from=1, to=2]
    \end{tikzcd}
    \quad=\quad
    \begin{tikzcd}
      {F(x \times x')} & {F(y \times y')} & {F(z \times z')} \\
      {Fx \times Fx'} & {Fy \times Fy'} & {Fz \times Fz'} \\
      {Fx \times Fx'} && {Fz \times Fz'} \\
      {Fx \times Fx'} && {Fz \times Fz'}
      \arrow[""{name=0, anchor=center, inner sep=0}, "{F(m \times m')}", "\shortmid"{marking}, from=1-1, to=1-2]
      \arrow[""{name=1, anchor=center, inner sep=0}, "{F(n \times n')}", "\shortmid"{marking}, from=1-2, to=1-3]
      \arrow[""{name=2, anchor=center, inner sep=0}, "{Fm \times Fm'}"', "\shortmid"{marking}, from=2-1, to=2-2]
      \arrow[""{name=3, anchor=center, inner sep=0}, "{Fn \times Fn'}"', "\shortmid"{marking}, from=2-2, to=2-3]
      \arrow[Rightarrow, no head, from=2-1, to=3-1]
      \arrow[Rightarrow, no head, from=2-3, to=3-3]
      \arrow[""{name=4, anchor=center, inner sep=0}, "{(Fm \odot Fn) \times (Fm' \odot Fn')}"', "\shortmid"{marking}, from=3-1, to=3-3]
      \arrow["{\Phi_{x,x'}}"', from=1-1, to=2-1]
      \arrow["{\Phi_{y,y'}}"{description}, from=1-2, to=2-2]
      \arrow["{\Phi_{z,z'}}", from=1-3, to=2-3]
      \arrow[""{name=5, anchor=center, inner sep=0}, "{F(m \odot n) \times F(m' \odot n')}"', "\shortmid"{marking}, from=4-1, to=4-3]
      \arrow[Rightarrow, no head, from=3-1, to=4-1]
      \arrow[Rightarrow, no head, from=3-3, to=4-3]
      \arrow["{\times_{(Fm,Fm'),(Fn,Fn')}}"{description}, draw=none, from=2-2, to=4]
      \arrow["{\Phi_{m,m'}}"{description}, draw=none, from=0, to=2]
      \arrow["{\Phi_{n,n'}}"{description}, draw=none, from=1, to=3]
      \arrow["{F_{m,n} \times F_{m',n'}}"{description, pos=0.6}, draw=none, from=4, to=5]
    \end{tikzcd}.
  \end{equation*}
  Also, for any objects $x$ and $x'$ in $\dbl{D}$, we have
  \begin{equation*}
    \begin{tikzcd}[column sep=large]
      {F(x \times x')} & {F(x \times x')} \\
      {F(x \times x')} & {F(x \times x')} \\
      {F(x \times x')} & {F(x \times x')} \\
      {Fx \times Fx'} & {Fx \times Fx'}
      \arrow[""{name=0, anchor=center, inner sep=0}, "{\id_{F(x \times x')}}", "\shortmid"{marking}, from=1-1, to=1-2]
      \arrow[Rightarrow, no head, from=1-1, to=2-1]
      \arrow[Rightarrow, no head, from=1-2, to=2-2]
      \arrow[""{name=1, anchor=center, inner sep=0}, "{F \id_{x \times x'}}", "\shortmid"{marking}, from=2-1, to=2-2]
      \arrow[Rightarrow, no head, from=2-1, to=3-1]
      \arrow[Rightarrow, no head, from=2-2, to=3-2]
      \arrow[""{name=2, anchor=center, inner sep=0}, "{F(\id_x \times \id_{x'})}", "\shortmid"{marking}, from=3-1, to=3-2]
      \arrow["{\Phi_{x,x'}}"', from=3-1, to=4-1]
      \arrow["{\Phi_{x,x'}}", from=3-2, to=4-2]
      \arrow[""{name=3, anchor=center, inner sep=0}, "{F\id_x \times F\id_{x'}}"', "\shortmid"{marking}, from=4-1, to=4-2]
      \arrow["{F_{x \times x'}}"{description, pos=0.4}, draw=none, from=0, to=1]
      \arrow["{\Phi_{\id_x, \id_{x'}}}"{description}, draw=none, from=2, to=3]
      \arrow["{F \times_{(x,x')}}"{description, pos=0.4}, draw=none, from=1, to=2]
    \end{tikzcd}
    \quad=\quad
    \begin{tikzcd}[column sep=large]
      {F(x \times x')} & {F(x \times x')} \\
      {Fx \times Fx'} & {Fx \times Fx'} \\
      {Fx \times Fx'} & {Fx \times Fx'} \\
      {Fx \times Fx'} & {Fx \times Fx'}
      \arrow[""{name=0, anchor=center, inner sep=0}, "{\mathrm{id}_{F(x \times x')}}", "\shortmid"{marking}, from=1-1, to=1-2]
      \arrow["{\Phi_{x,x'}}"', from=1-1, to=2-1]
      \arrow["{\Phi_{x,x'}}", from=1-2, to=2-2]
      \arrow[""{name=1, anchor=center, inner sep=0}, "{\mathrm{id}_{Fx \times Fx'}}", "\shortmid"{marking}, from=2-1, to=2-2]
      \arrow[Rightarrow, no head, from=2-2, to=3-2]
      \arrow[Rightarrow, no head, from=2-1, to=3-1]
      \arrow[""{name=2, anchor=center, inner sep=0}, "{\mathrm{id}_{Fx} \times \mathrm{id}_{Fx'}}", "\shortmid"{marking}, from=3-1, to=3-2]
      \arrow[Rightarrow, no head, from=3-1, to=4-1]
      \arrow[Rightarrow, no head, from=3-2, to=4-2]
      \arrow[""{name=3, anchor=center, inner sep=0}, "{F\id_x \times F\id_{x'}}"', "\shortmid"{marking}, from=4-1, to=4-2]
      \arrow["{\mathrm{id}_{\Phi_{x,x'}}}"{description, pos=0.4}, draw=none, from=0, to=1]
      \arrow["{F_x \times F_{x'}}"{description}, draw=none, from=2, to=3]
      \arrow["{\times_{(Fx, Fx')}}"{description, pos=0.4}, draw=none, from=1, to=2]
    \end{tikzcd}.
  \end{equation*}
  In particular, when the double category $\dbl{D}$ and the lax functor $F$ are
  both cartesian, then the laxator $F_{m \times m', n \times n'}$ is completely
  determined by the product of the laxators $F_{m,n}$ and $F_{m',n'}$ and,
  similarly, the unitor $F_{x \times x'}$ is determined by the product of the
  unitors $F_x$ and $F_{x'}$.
\end{lemma}
\begin{proof}
  By the naturality of the laxators for the pair of projection cells
  $\pi_{m,m'}: m \times m' \to m$ and $\pi_{n,n'}: n \times n' \to n$, we have
  \begin{equation*}
    \begin{tikzcd}
      {F(x \times x')} & {F(y \times y')} & {F(z \times z')} \\
      {F(x \times x')} && {F(z \times z')} \\
      Fx && Fz
      \arrow["{F(m \times m')}", "\shortmid"{marking}, from=1-1, to=1-2]
      \arrow["{F(n \times n')}", "\shortmid"{marking}, from=1-2, to=1-3]
      \arrow[Rightarrow, no head, from=1-1, to=2-1]
      \arrow[Rightarrow, no head, from=1-3, to=2-3]
      \arrow[""{name=0, anchor=center, inner sep=0}, "{F((m \times m') \odot (n \times n'))}"', "\shortmid"{marking}, from=2-1, to=2-3]
      \arrow["{F\pi_{x,x'}}"', from=2-1, to=3-1]
      \arrow["{F\pi_{z,z'}}", from=2-3, to=3-3]
      \arrow[""{name=1, anchor=center, inner sep=0}, "{F(m \odot n)}"', "\shortmid"{marking}, from=3-1, to=3-3]
      \arrow["{F_{m \times m', n \times n'}}"{description}, draw=none, from=1-2, to=0]
      \arrow["{F(\pi_{m,m'} \odot \pi_{n,n'})}"{description, pos=0.6}, draw=none, from=0, to=1]
    \end{tikzcd}
    \quad=\quad
    \begin{tikzcd}
      {F(x \times x')} & {F(y \times y')} & {F(z \times z')} \\
      Fx & Fy & Fz \\
      Fx && Fz
      \arrow[""{name=0, anchor=center, inner sep=0}, "{F(m \times m')}", "\shortmid"{marking}, from=1-1, to=1-2]
      \arrow[""{name=1, anchor=center, inner sep=0}, "{F(n \times n')}", "\shortmid"{marking}, from=1-2, to=1-3]
      \arrow["{F\pi_{x,x'}}"', from=1-1, to=2-1]
      \arrow["{F\pi_{y,y'}}"{description}, from=1-2, to=2-2]
      \arrow["{F\pi_{z,z'}}", from=1-3, to=2-3]
      \arrow[""{name=2, anchor=center, inner sep=0}, "Fm"', "\shortmid"{marking}, from=2-1, to=2-2]
      \arrow[""{name=3, anchor=center, inner sep=0}, "Fn"', "\shortmid"{marking}, from=2-2, to=2-3]
      \arrow[""{name=4, anchor=center, inner sep=0}, "{F(m \odot n)}"', "\shortmid"{marking}, from=3-1, to=3-3]
      \arrow[Rightarrow, no head, from=2-1, to=3-1]
      \arrow[Rightarrow, no head, from=2-3, to=3-3]
      \arrow["{F\pi_{m,m'}}"{description}, draw=none, from=0, to=2]
      \arrow["{F\pi_{n,n'}}"{description}, draw=none, from=1, to=3]
      \arrow["{F_{m,n}}"{description}, draw=none, from=2-2, to=4]
    \end{tikzcd}.
  \end{equation*}
  Applying the naturality of the laxators for the other pair of projection cells
  $\pi_{m,m'}': m \times m' \to m'$ and $\pi_{n,n'}': n \times n' \to n'$ yields
  a similar equation. The pairing of these two equations is
  \begin{equation*}
    \begin{tikzcd}
      {F(x \times x')} & {F(y \times y')} & {F(z \times z')} \\
      {F(x \times x')} && {F(z \times z')} \\
      {Fx \times Fx'} && {Fz \times Fz'}
      \arrow["{F(m \times m')}", "\shortmid"{marking}, from=1-1, to=1-2]
      \arrow["{F(n \times n')}", "\shortmid"{marking}, from=1-2, to=1-3]
      \arrow[Rightarrow, no head, from=1-1, to=2-1]
      \arrow[Rightarrow, no head, from=1-3, to=2-3]
      \arrow[""{name=0, anchor=center, inner sep=0}, "{F((m \times m') \odot (n \times n'))}"', "\shortmid"{marking}, from=2-1, to=2-3]
      \arrow["{\Phi_{x,x'}}"', from=2-1, to=3-1]
      \arrow["{\Phi_{z,z'}}", from=2-3, to=3-3]
      \arrow[""{name=1, anchor=center, inner sep=0}, "{F(m \odot n) \times F(m' \odot n')}"', "\shortmid"{marking}, from=3-1, to=3-3]
      \arrow["{F_{m \times m', n \times n'}}"{description}, draw=none, from=1-2, to=0]
      \arrow["{\langle F(\pi_{m,m'} \odot \pi_{n,n'}), F(\pi_{m,m'}' \odot \pi_{n,n'}') \rangle}"{description, pos=0.6}, draw=none, from=0, to=1]
    \end{tikzcd}
    \quad=\quad
    \begin{tikzcd}
      {F(x \times x')} & {F(y \times y')} & {F(z \times z')} \\
      {Fx \times Fx'} && {Fz \times Fz'} \\
      {Fx \times Fx'} && {Fz \times Fz'}
      \arrow["{F(m \times m')}", "\shortmid"{marking}, from=1-1, to=1-2]
      \arrow["{F(n \times n')}", "\shortmid"{marking}, from=1-2, to=1-3]
      \arrow["{\Phi_{x,x'}}"', from=1-1, to=2-1]
      \arrow["{\Phi_{z,z'}}", from=1-3, to=2-3]
      \arrow[""{name=0, anchor=center, inner sep=0}, "{F(m \odot n) \times F(m' \odot n')}"', "\shortmid"{marking}, from=3-1, to=3-3]
      \arrow[Rightarrow, no head, from=2-1, to=3-1]
      \arrow[Rightarrow, no head, from=2-3, to=3-3]
      \arrow[""{name=1, anchor=center, inner sep=0}, "{(Fm \odot Fn) \times (Fm' \odot Fn')}"', "\shortmid"{marking}, from=2-1, to=2-3]
      \arrow["{F_{m,n} \times F_{m',n'}}"{description, pos=0.6}, draw=none, from=1, to=0]
      \arrow["{\langle F\pi_{m,m'} \odot F\pi_{n,n'}, F\pi_{m,m'}' \odot F\pi_{n,n'}' \rangle}"{description}, draw=none, from=1-2, to=1]
    \end{tikzcd}.
  \end{equation*}
  Thus, the first statement of the lemma is established provided that
  \begin{equation} \label{eq:functor-product-compose-comparison-1}
    \begin{tikzcd}[column sep=huge]
      {F(x \times x')} && {F(z \times z')} \\
      {Fx \times Fx'} && {Fz \times Fz'}
      \arrow[""{name=0, anchor=center, inner sep=0}, "{F((m \times m') \odot (n \times n'))}", "\shortmid"{marking}, from=1-1, to=1-3]
      \arrow["{\Phi_{x,x'}}"', from=1-1, to=2-1]
      \arrow["{\Phi_{z,z'}}", from=1-3, to=2-3]
      \arrow[""{name=1, anchor=center, inner sep=0}, "{F(m \odot n) \times F(m' \odot n')}"', "\shortmid"{marking}, from=2-1, to=2-3]
      \arrow["{\langle F(\pi_{m,m'} \odot \pi_{n,n'}), F(\pi_{m,m'}' \odot \pi_{n,n'}') \rangle}"{description}, draw=none, from=0, to=1]
    \end{tikzcd}
    =
    \begin{tikzcd}[row sep=scriptsize, column sep=large]
      {F(x \times x')} && {F(z \times z')} \\
      {F(x \times x')} && {F(z \times z')} \\
      {Fx \times Fx'} && {Fz \times Fz'}
      \arrow[""{name=0, anchor=center, inner sep=0}, "{F((m \times m') \odot (n \times n'))}", "\shortmid"{marking}, from=1-1, to=1-3]
      \arrow[""{name=1, anchor=center, inner sep=0}, "{F((m \odot n) \times (m' \odot n'))}"', "\shortmid"{marking}, from=2-1, to=2-3]
      \arrow[Rightarrow, no head, from=1-1, to=2-1]
      \arrow[Rightarrow, no head, from=1-3, to=2-3]
      \arrow["{\Phi_{x,x'}}"', from=2-1, to=3-1]
      \arrow["{\Phi_{z,z'}}", from=2-3, to=3-3]
      \arrow[""{name=2, anchor=center, inner sep=0}, "{F(m \odot n) \times F(m' \odot n')}"', "\shortmid"{marking}, from=3-1, to=3-3]
      \arrow["{F \times_{(m,m'),(n,n')}}"{description}, draw=none, from=0, to=1]
      \arrow["{\Phi_{m \odot n, m' \odot n'}}"{description, pos=0.6}, draw=none, from=1, to=2]
    \end{tikzcd}
  \end{equation}
  and
  \begin{equation} \label{eq:functor-product-compose-comparison-2}
    \begin{tikzcd}
      {F(x \times x')} & {F(y \times y')} & {F(z \times z')} \\
      {Fx \times Fx'} && {Fz \times Fz'}
      \arrow["{F(m \times m')}", "\shortmid"{marking}, from=1-1, to=1-2]
      \arrow["{F(n \times n')}", "\shortmid"{marking}, from=1-2, to=1-3]
      \arrow["{\Phi_{x,x'}}"', from=1-1, to=2-1]
      \arrow["{\Phi_{z,z'}}", from=1-3, to=2-3]
      \arrow[""{name=0, anchor=center, inner sep=0}, "{(Fm \odot Fn) \times (Fm' \odot Fn')}"', "\shortmid"{marking}, from=2-1, to=2-3]
      \arrow["{\langle F\pi_{m,m'} \odot F\pi_{n,n'}, F\pi_{m,m'}' \odot F\pi_{n,n'}' \rangle}"{description}, draw=none, from=1-2, to=0]
    \end{tikzcd}
    =
    \begin{tikzcd}[row sep=scriptsize]
      {F(x \times x')} & {F(y \times y')} & {F(z \times z')} \\
      {Fx \times Fx'} & {Fy \times Fy'} & {Fz \times Fz'} \\
      {Fx \times Fx'} && {Fz \times Fz'}
      \arrow[""{name=0, anchor=center, inner sep=0}, "{F(m \times m')}", "\shortmid"{marking}, from=1-1, to=1-2]
      \arrow[""{name=1, anchor=center, inner sep=0}, "{F(n \times n')}", "\shortmid"{marking}, from=1-2, to=1-3]
      \arrow[""{name=2, anchor=center, inner sep=0}, "{Fm \times Fm'}"', "\shortmid"{marking}, from=2-1, to=2-2]
      \arrow[""{name=3, anchor=center, inner sep=0}, "{Fn \times Fn'}"', "\shortmid"{marking}, from=2-2, to=2-3]
      \arrow[Rightarrow, no head, from=2-1, to=3-1]
      \arrow[Rightarrow, no head, from=2-3, to=3-3]
      \arrow[""{name=4, anchor=center, inner sep=0}, "{(Fm \odot Fn) \times (Fm' \odot Fn')}"', "\shortmid"{marking}, from=3-1, to=3-3]
      \arrow["{\Phi_{x,x'}}"', from=1-1, to=2-1]
      \arrow["{\Phi_{y,y'}}"{description}, from=1-2, to=2-2]
      \arrow["{\Phi_{z,z'}}", from=1-3, to=2-3]
      \arrow["{\times_{(Fm,Fm'),(Fn,Fn')}}"{description}, draw=none, from=2-2, to=4]
      \arrow["{\Phi_{m,m'}}"{description}, draw=none, from=0, to=2]
      \arrow["{\Phi_{n,n'}}"{description}, draw=none, from=1, to=3]
    \end{tikzcd}.
  \end{equation}

  To prove \cref{eq:functor-product-compose-comparison-1}, observe that
  postcomposing with the projection $\pi_{F(m \odot n), F(m' \odot n')}$ on left
  yields $F(\pi_{m,m'} \odot \pi_{n,n'})$ by definition and on the right yields
  the same thing, by \cref{eq:dbl-proj-external-compose} or
  \eqref{eq:dbl-product-compose-comparison}. Similarly, postcomposing with the
  projection $\pi_{F(m \odot n), F(m' \odot n')}'$ yields
  $F(\pi_{m,m'}' \odot \pi_{n,n'}')$ on both sides. Thus,
  \cref{eq:functor-product-compose-comparison-1} holds by the universal property
  of products. It remains to prove
  \cref{eq:functor-product-compose-comparison-2}. Postcomposing with the
  projection $\pi_{Fm \odot Fn, Fm' \odot Fn'}$ on the left gives
  $F\pi_{m,m'} \odot F \pi_{n,n'}$ by definition, and on the right gives the
  same thing:
  \begin{equation*}
    \begin{tikzcd}
      {F(x \times x')} & {F(y \times y')} & {F(z \times z')} \\
      {Fx \times Fx'} & {Fy \times Fy'} & {Fz \times Fz'} \\
      Fx & Fy & Fz
      \arrow[""{name=0, anchor=center, inner sep=0}, "{F(m \times m')}", "\shortmid"{marking}, from=1-1, to=1-2]
      \arrow[""{name=1, anchor=center, inner sep=0}, "{F(n \times n')}", "\shortmid"{marking}, from=1-2, to=1-3]
      \arrow[""{name=2, anchor=center, inner sep=0}, "{Fm \times Fm'}"', "\shortmid"{marking}, from=2-1, to=2-2]
      \arrow[""{name=3, anchor=center, inner sep=0}, "{Fn \times Fn'}"', "\shortmid"{marking}, from=2-2, to=2-3]
      \arrow["{\pi_{Fx,Fx'}}"', from=2-1, to=3-1]
      \arrow["{\Phi_{x,x'}}"', from=1-1, to=2-1]
      \arrow["{\Phi_{y,y'}}"{description}, from=1-2, to=2-2]
      \arrow["{\Phi_{z,z'}}", from=1-3, to=2-3]
      \arrow["{\pi_{Fy,Fy'}}"{description}, from=2-2, to=3-2]
      \arrow["{\pi_{Fz,Fz'}}", from=2-3, to=3-3]
      \arrow[""{name=4, anchor=center, inner sep=0}, "Fm"', "\shortmid"{marking}, from=3-1, to=3-2]
      \arrow[""{name=5, anchor=center, inner sep=0}, "Fn"', "\shortmid"{marking}, from=3-2, to=3-3]
      \arrow["{\Phi_{n,n'}}"{description}, draw=none, from=1, to=3]
      \arrow["{\Phi_{m,m'}}"{description}, draw=none, from=0, to=2]
      \arrow["{\pi_{Fm,Fm'}}"{description, pos=0.6}, draw=none, from=2, to=4]
      \arrow["{\pi_{Fn,Fn'}}"{description, pos=0.6}, draw=none, from=3, to=5]
    \end{tikzcd}
    =
    \begin{tikzcd}
      {F(x \times x')} & {F(y \times y')} & {F(z \times z')} \\
      Fx & Fy & Fz
      \arrow[""{name=0, anchor=center, inner sep=0}, "{F(m \times m')}", "\shortmid"{marking}, from=1-1, to=1-2]
      \arrow[""{name=1, anchor=center, inner sep=0}, "{F(n \times n')}", "\shortmid"{marking}, from=1-2, to=1-3]
      \arrow[""{name=2, anchor=center, inner sep=0}, "Fm"', "\shortmid"{marking}, from=2-1, to=2-2]
      \arrow[""{name=3, anchor=center, inner sep=0}, "Fn"', "\shortmid"{marking}, from=2-2, to=2-3]
      \arrow["{F\pi_{x,x'}}"', from=1-1, to=2-1]
      \arrow["{F\pi_{y,y'}}"{description}, from=1-2, to=2-2]
      \arrow["{F\pi_{z,z'}}", from=1-3, to=2-3]
      \arrow["{F\pi_{m,m'}}"{description}, draw=none, from=0, to=2]
      \arrow["{F\pi_{n,n'}}"{description}, draw=none, from=1, to=3]
    \end{tikzcd}.
  \end{equation*}
  Similarly, postcomposing with the other projection gives
  $F\pi_{m,m'}' \odot F\pi_{n,n'}'$ on both sides. This completes the proof of
  \cref{eq:functor-product-compose-comparison-2} and the first part of the
  lemma.

  Next, for any objects $x$ and $x'$ in $\dbl{D}$, we have
  \begin{equation*}
    \begin{tikzcd}[row sep=scriptsize]
      {F(x \times x')} & {F(x \times x')} \\
      {F(x \times x')} & {F(x \times x')} \\
      Fx & Fx
      \arrow[""{name=0, anchor=center, inner sep=0}, "{\id_{F(x \times x')}}", "\shortmid"{marking}, from=1-1, to=1-2]
      \arrow[""{name=1, anchor=center, inner sep=0}, "{F(\id_{x \times x'})}", "\shortmid"{marking}, from=2-1, to=2-2]
      \arrow[""{name=2, anchor=center, inner sep=0}, "{F \id_x}"', "\shortmid"{marking}, from=3-1, to=3-2]
      \arrow["{F\pi_{x,x'}}"', from=2-1, to=3-1]
      \arrow["{F\pi_{x,x'}}", from=2-2, to=3-2]
      \arrow[Rightarrow, no head, from=1-2, to=2-2]
      \arrow[Rightarrow, no head, from=1-1, to=2-1]
      \arrow["{F\id_{\pi_{x,x'}}}"{description}, draw=none, from=1, to=2]
      \arrow["{F_{x \times x'}}"{description, pos=0.4}, draw=none, from=0, to=1]
    \end{tikzcd}
    \quad=\quad
    \begin{tikzcd}[row sep=scriptsize]
      {F(x \times x')} & {F(x \times x')} \\
      Fx & {Fx'} \\
      Fx & {Fx'}
      \arrow[""{name=0, anchor=center, inner sep=0}, "{\id_{F(x \times x')}}", "\shortmid"{marking}, from=1-1, to=1-2]
      \arrow["{F\pi_{x,x'}}"', from=1-1, to=2-1]
      \arrow["{F\pi_{x,x'}}", from=1-2, to=2-2]
      \arrow[""{name=1, anchor=center, inner sep=0}, "{\id_{Fx}}", "\shortmid"{marking}, from=2-1, to=2-2]
      \arrow[Rightarrow, no head, from=2-1, to=3-1]
      \arrow[Rightarrow, no head, from=2-2, to=3-2]
      \arrow[""{name=2, anchor=center, inner sep=0}, "{F \id_x}"', from=3-1, to=3-2]
      \arrow["{\id_{F\pi_{x,x'}}}"{description, pos=0.4}, draw=none, from=0, to=1]
      \arrow["{F_x}"{description}, draw=none, from=1, to=2]
    \end{tikzcd}
  \end{equation*}
  by the naturality of the unitors for the projection arrow
  $\pi_{x,x'}: x \times x' \to x$. The pairing of this equation with the
  analogous equation for the other projection $\pi_{x,x'}': x \times x' \to x'$
  is
  \begin{equation*}
    \begin{tikzcd}[row sep=scriptsize, column sep=large]
      {F(x \times x')} & {F(x \times x')} \\
      {F(x \times x')} & {F(x \times x')} \\
      {Fx \times Fx'} & {Fx \times Fx'}
      \arrow[""{name=0, anchor=center, inner sep=0}, "{\id_{F(x \times x')}}", "\shortmid"{marking}, from=1-1, to=1-2]
      \arrow[""{name=1, anchor=center, inner sep=0}, "{F\id_x \times F\id_{x'}}"', "\shortmid"{marking}, from=3-1, to=3-2]
      \arrow[Rightarrow, no head, from=1-2, to=2-2]
      \arrow[Rightarrow, no head, from=1-1, to=2-1]
      \arrow[""{name=2, anchor=center, inner sep=0}, "{F \id_{x \times x'}}", "\shortmid"{marking}, from=2-1, to=2-2]
      \arrow["{\Phi_{x,x'}}"', from=2-1, to=3-1]
      \arrow["{\Phi_{x,x'}}", from=2-2, to=3-2]
      \arrow["{F_{x \times x'}}"{description, pos=0.4}, draw=none, from=0, to=2]
      \arrow["{\langle F \id_{\pi_{x,x'}}, F\id_{\pi_{x,x'}'} \rangle}"{description}, draw=none, from=2, to=1]
    \end{tikzcd}
    \quad=\quad
    \begin{tikzcd}[row sep=scriptsize, column sep=large]
      {F(x \times x')} & {F(x \times x')} \\
      {Fx \times Fx'} & {Fx \times Fx'} \\
      {Fx \times Fx'} & {Fx \times Fx'}
      \arrow[""{name=0, anchor=center, inner sep=0}, "{\mathrm{id}_{F(x \times x')}}", "\shortmid"{marking}, from=1-1, to=1-2]
      \arrow[""{name=1, anchor=center, inner sep=0}, "{\mathrm{id}_{Fx} \times \mathrm{id}_{Fx'}}"', "\shortmid"{marking}, from=2-1, to=2-2]
      \arrow[Rightarrow, no head, from=2-1, to=3-1]
      \arrow[Rightarrow, no head, from=2-2, to=3-2]
      \arrow[""{name=2, anchor=center, inner sep=0}, "{F\id_x \times F\id_{x'}}"', "\shortmid"{marking}, from=3-1, to=3-2]
      \arrow["{\Phi_{x,x'}}"', from=1-1, to=2-1]
      \arrow["{\Phi_{x,x'}}", from=1-2, to=2-2]
      \arrow["{F_x \times F_{x'}}"{description, pos=0.6}, draw=none, from=1, to=2]
      \arrow["{\langle \id_{F\pi_{x,x'}}, \id_{F\pi_{x,x'}'} \rangle}"{description}, draw=none, from=0, to=1]
    \end{tikzcd}.
  \end{equation*}
  So the second statement in the lemma is proved once we know that
  \begin{equation} \label{eq:functor-product-id-comparison-1}
    \begin{tikzcd}[column sep=large]
      {F(x \times x')} & {F(x \times x')} \\
      {Fy \times Fy'} & {Fx \times Fx'}
      \arrow[""{name=0, anchor=center, inner sep=0}, "{F \id_{x \times x'}}", "\shortmid"{marking}, from=1-1, to=1-2]
      \arrow[""{name=1, anchor=center, inner sep=0}, "{F \id_x \times F \id_{x'}}"', "\shortmid"{marking}, from=2-1, to=2-2]
      \arrow["{\Phi_{x,x'}}"', from=1-1, to=2-1]
      \arrow["{\Phi_{x,x'}}", from=1-2, to=2-2]
      \arrow["{\langle F\id_{\pi_{x,x'}}, F\id_{\pi_{x,x'}'}\rangle}"{description}, draw=none, from=0, to=1]
    \end{tikzcd}
    \quad=\quad
    \begin{tikzcd}[row sep=scriptsize, column sep=large]
      {F(x \times x')} & {F(x \times x')} \\
      {F(x \times x')} & {F(x \times x')} \\
      {Fx \times Fx'} & {Fx \times Fx'}
      \arrow[""{name=0, anchor=center, inner sep=0}, "{F \id_{x \times x'}}", "\shortmid"{marking}, from=1-1, to=1-2]
      \arrow[""{name=1, anchor=center, inner sep=0}, "{F(\id_x \times \id_{x'})}"', "\shortmid"{marking}, from=2-1, to=2-2]
      \arrow[Rightarrow, no head, from=1-1, to=2-1]
      \arrow[Rightarrow, no head, from=1-2, to=2-2]
      \arrow["{\Phi_{x,x'}}"', from=2-1, to=3-1]
      \arrow["{\Phi_{x,x'}}", from=2-2, to=3-2]
      \arrow[""{name=2, anchor=center, inner sep=0}, "{F \id_x \times F \id_{x'}}"', "\shortmid"{marking}, from=3-1, to=3-2]
      \arrow["{F \times_{(x,x')}}"{description}, draw=none, from=0, to=1]
      \arrow["{\Phi_{\id_x, \id_{x'}}}"{description, pos=0.6}, draw=none, from=1, to=2]
    \end{tikzcd}
  \end{equation}
  and
  \begin{equation} \label{eq:functor-product-id-comparison-2}
    \begin{tikzcd}[column sep=large]
      {F(x \times x')} & {F(x \times x')} \\
      {Fx \times Fx'} & {Fx \times Fx'}
      \arrow[""{name=0, anchor=center, inner sep=0}, "{\id_{Fx} \times \id_{Fx'}}"', "\shortmid"{marking}, from=2-1, to=2-2]
      \arrow["{\Phi_{x,x'}}"', from=1-1, to=2-1]
      \arrow["{\Phi_{x,x'}}", from=1-2, to=2-2]
      \arrow[""{name=1, anchor=center, inner sep=0}, "{\id_{F(x \times x')}}", "\shortmid"{marking}, from=1-1, to=1-2]
      \arrow["{\langle \id_{F\pi_{x,x'}}, \id_{F\pi_{x,x'}'} \rangle}"{description}, draw=none, from=1, to=0]
    \end{tikzcd}
    \quad=\quad
    \begin{tikzcd}[row sep=scriptsize, column sep=large]
      {F(x \times x')} & {F(x \times x')} \\
      {Fx \times Fx'} & {Fx \times Fx'} \\
      {Fx \times Fx'} & {Fx \times Fx'}
      \arrow[""{name=0, anchor=center, inner sep=0}, "{\id_{F(x \times x')}}", "\shortmid"{marking}, from=1-1, to=1-2]
      \arrow["{\Phi_{x,x'}}"', from=1-1, to=2-1]
      \arrow["{\Phi_{x,x'}}", from=1-2, to=2-2]
      \arrow[""{name=1, anchor=center, inner sep=0}, "{\id_{Fx \times Fx'}}"', "\shortmid"{marking}, from=2-1, to=2-2]
      \arrow[Rightarrow, no head, from=2-2, to=3-2]
      \arrow[Rightarrow, no head, from=2-1, to=3-1]
      \arrow[""{name=2, anchor=center, inner sep=0}, "{\id_{Fx} \times \id_{Fx'}}"', "\shortmid"{marking}, from=3-1, to=3-2]
      \arrow["{\id_{\Phi_{x,x'}}}"{description}, draw=none, from=0, to=1]
      \arrow["{\times_{(Fx, Fx')}}"{description, pos=0.6}, draw=none, from=1, to=2]
    \end{tikzcd}.
  \end{equation}
  These last equations are proved similarly as before, using
  \cref{eq:dbl-proj-external-id} or \eqref{eq:dbl-product-id-comparison}.
\end{proof}

The essential fact about the interaction between lax functors and equipments,
due to Shulman \cite[Proposition 6.4]{shulman2008}, is:

\begin{proposition}[Lax functors preserve restrictions]
  \label{prop:lax-functor-restrictions}
  Lax double functors between equipments preserve restriction cells. That is, if
  $\dbl{D}$ and $\dbl{E}$ are equipments and $F: \dbl{D} \to \dbl{E}$ is a lax
  double functor, then for any cartesian cell
  $\inlineCell[normal]{x}{y}{w}{z}{n(f,g)}{n}{f}{g}{\res}$ in $\dbl{D}$, its
  image $\inlineCell[normal]{Fx}{Fy}{Fw}{Fz}{F(n(f,g))}{Fn}{Ff}{Fg}{F(\res)}$
  under $F$ is again a cartesian cell in $\dbl{E}$.
\end{proposition}

This fact justifies the lack of extra conditions on lax double functors in the
2-category $\EqpLax$ (\cref{sec:cartesian-equipments}). But it is also important
to understand what the result does not say. Extension cells are generally
preserved by \emph{op}lax functors, not lax ones. Meanwhile, companions and
conjoints are preserved by a lax or oplax functor only when it is \emph{normal}
\cite[Proposition 3.8]{dawson2010}. Thus, when working with lax double functors,
the ``legitimate'' structure available in an equipment appears to be the
restriction cells, in contrast to the extra structure afforded by other
characterizations of an equipment (\cref{thm:equipment}). Due to this excess of
structure we do not take our double theories to be equipments; instead, we will
introduce a notion of ``restriction sketch'' to incorporate restrictions into
theories.

When using restriction sketches, it will be important to know when restricted
proarrows have uniquely determined laxators.

\begin{lemma}[Laxators for restricted proarrows] \label{lem:laxators-restrictions}
  Let $F: \dbl{D} \to \dbl{E}$ be a lax double functor.
  \begin{enumerate}[(i),nosep]
    \item For any cells $\inlineCell{x}{y}{x'}{y'}{m}{m'}{f}{g}{\alpha}$ and
      $\inlineCell{y}{z}{y'}{z'}{n}{n'}{g}{h}{\beta}$ in $\dbl{D}$, if the image
      $F(\alpha \odot \beta)$ of their composite is a restriction cell in
      $\dbl{E}$, then the laxator $F_{m,n}$ is uniquely determined by the
      laxator $F_{m',n'}$.
    \item When $\dbl{D}$ is a unitary double category, for any cell
      $\stdInlineCell{\alpha}$ in $\dbl{D}$, if its image $F(\alpha)$ is a
      restriction cell in $\dbl{E}$, then the laxators $F_{x,m}$ and $F_{m,y}$
      are uniquely determined by the laxators $F_{w,n}$ and $F_{n,z}$.
  \end{enumerate}
\end{lemma}
\begin{proof}
  Using the naturality of the laxators \eqref{eq:naturality-laxators} and then
  the universal property of the restriction cell $F(\alpha \odot \beta)$ proves
  the first statement. The second statement then follows from the first using
  the equations $\id_f \odot \alpha = \alpha = \alpha \odot \id_g$.
\end{proof}

In addition to preserving restrictions, lax functors preserve category and
profunctor objects. This is partly the observation behind the fact that the
bimodule construction defines a 2-functor
$\Module: \EqpLax^q \to \EqpLaxNormal^q$ \cite[Proposition 11.11]{shulman2008}.
Having studied the receiving structure of cartesian equipments such as
$\Module{\dbl{E}}$ when $\dbl{E}$ is a cartesian equipment with local
coequalizers, we can now make good on the promise of
\cref{rmk:unitalization-preview}. There we indicated that a lax functor's
assignment of objects to categories and proarrows to profunctors is coherent. To
be more precise, we now show that lax functors $\dbl{D} \to\dbl{E}$ are in
one-to-one correspondence with certain bimodule-valued lax functors
$\dbl{D}\to\Module{\dbl{E}}$. The correspondence is based on that of
\cite[Proposition 5.14]{cruttwell2010} which shows that in the virtual setting
$\Module: \vDbl \to \vDblNormal$ is a right adjoint. We will specialize this
result to the case where $\Module{\dbl{E}}$ is a genuine equipment and then also
cartesian, under the assumption that $\dbl{E}$ is a (cartesian) equipment with
local coequalizers. Our elaboration of the quoted results of \cite{shulman2008}
and \cite{cruttwell2010} makes one adjustment, namely, to note that under the 
hypothesis that $\Module{\dbl{E}}$ has \emph{chosen} units the 
ordinarily merely normal lax functor $\Module{F}$ is actually \emph{unitary}. We
are happy to make this assumption in light of the simplifications it allows in
many of the subsequent arguments.

In more detail, as described in the reference, the unit of the adjunction
$\eta\colon \dbl{D} \to \Module{\dbl{D}}$ makes the assignments
\begin{enumerate}[noitemsep]
  \item $x\mapsto \id_x$, regarded as the trivial category with carrier $x$;
  \item $f\mapsto \id_f$, regarded as a functor between trivial categories;
  \item $m\mapsto m$, regarded as profunctor between the trivial categories on
    its source and target;
  \item $\theta\mapsto \theta$, regarded as a map of profunctors between trivial categories.
\end{enumerate}
In particular, the left and right actions of trivial categories on a proarrow
$m: x \proto y$ are given by the external unit isomorphisms present in the
double category $\dbl{D}$:
\begin{equation*}
  \begin{tikzcd}
    x & x & y \\
    x && x
    \arrow["{\id_x}", "\shortmid"{marking}, from=1-1, to=1-2]
    \arrow["m", "\shortmid"{marking}, from=1-2, to=1-3]
    \arrow[Rightarrow, no head, from=1-1, to=2-1]
    \arrow[""{name=0, anchor=center, inner sep=0}, "m"', "\shortmid"{marking}, from=2-1, to=2-3]
    \arrow[Rightarrow, no head, from=1-3, to=2-3]
    \arrow["\cong"{description, pos=0.4}, draw=none, from=1-2, to=0]
  \end{tikzcd}
  \qquad\qquad
  \begin{tikzcd}
    x & y & y \\
    x && y
    \arrow["m", "\shortmid"{marking}, from=1-1, to=1-2]
    \arrow["{\id_y}", "\shortmid"{marking}, from=1-2, to=1-3]
    \arrow[Rightarrow, no head, from=1-3, to=2-3]
    \arrow[Rightarrow, no head, from=1-1, to=2-1]
    \arrow[""{name=0, anchor=center, inner sep=0}, "m"', "\shortmid"{marking}, from=2-1, to=2-3]
    \arrow["\cong"{description, pos=0.4}, draw=none, from=1-2, to=0]
  \end{tikzcd}.
\end{equation*}
By construction, the functor $\eta: \dbl{D} \to \Module{\dbl{D}}$ is unitary,
assuming a choice of units, hence is normal in particular. Now let
$F\colon \dbl{D}\to\dbl{E}$ denote a lax functor into an equipment $\dbl{E}$
with local coequalizers. The composite functor
$\Module{F}\eta\colon\dbl{D} \to\Module{\dbl{E}}$ thus makes the assignments
\begin{enumerate}[noitemsep]
  \item $x\mapsto F\id_x$, regarded as the image of the trivial category;
  \item $f\mapsto F\id_f$, regarded as a functor between images of trivial
    categories;
  \item $m\mapsto Fm$, regarded as profunctor between images of trivial categories;
  \item $\theta\mapsto F\theta$, regarded as a map of profunctors between images
    of trivial categories.
\end{enumerate}
These, however, need to be padded by appropriate structural morphisms. For example, 
the left and right actions on $Fm$ are the composites
\begin{equation} \label{equation:inducedactionoftrivialmonoidunderlaxfunctor}
  \begin{tikzcd}
    Fx & Fx & Fy \\
    Fx && Fy \\
    Fx && Fy
    \arrow[Rightarrow, no head, from=2-1, to=3-1]
    \arrow[""{name=0, anchor=center, inner sep=0}, "Fm"', "\shortmid"{marking}, from=3-1, to=3-3]
    \arrow[Rightarrow, no head, from=2-3, to=3-3]
    \arrow["{F\id_x}", "\shortmid"{marking}, from=1-1, to=1-2]
    \arrow["Fm", "\shortmid"{marking}, from=1-2, to=1-3]
    \arrow[""{name=1, anchor=center, inner sep=0}, "{F(\id_x\odot m)}"', "\shortmid"{marking}, from=2-1, to=2-3]
    \arrow[from=1-1, to=2-1]
    \arrow[from=1-3, to=2-3]
    \arrow["{F(\cong)}"{description, pos=0.6}, draw=none, from=1, to=0]
    \arrow["{F_{\id_x,m}}"{description, pos=0.4}, draw=none, from=1-2, to=1]
  \end{tikzcd}
  \qquad\text{and}\qquad
  \begin{tikzcd}
    Fx & Fy & Fy \\
    Fx && Fy \\
    Fx && Fy
    \arrow[Rightarrow, no head, from=2-3, to=3-3]
    \arrow[Rightarrow, no head, from=2-1, to=3-1]
    \arrow[""{name=0, anchor=center, inner sep=0}, "Fm"', "\shortmid"{marking}, from=3-1, to=3-3]
    \arrow[from=1-1, to=2-1]
    \arrow["Fm", "\shortmid"{marking}, from=1-1, to=1-2]
    \arrow["{F\id_y}", "\shortmid"{marking}, from=1-2, to=1-3]
    \arrow[from=1-3, to=2-3]
    \arrow[""{name=1, anchor=center, inner sep=0}, "{F(m\odot \id_y)}"', "\shortmid"{marking}, from=2-1, to=2-3]
    \arrow["{F(\cong)}"{description, pos=0.6}, draw=none, from=1, to=0]
    \arrow["{F_{m,\id_y}}"{description, pos=0.4}, draw=none, from=1-2, to=1]
  \end{tikzcd},
\end{equation}
where $F(\cong)$ in each case denotes the image under $F$ of the canonical unit
isomorphism. When $\dbl{D}$ is a double theory, strict unitality will hold and
so these isomorphisms and their images under $F$ will be identities. Note that
we do not yet know that the composite
$\Module{F}\eta: \dbl{D} \to \Module{\dbl{E}}$ is a lax functor between double
categories since the intermediate object $\Module{\dbl{D}}$ is merely virtual.

The counit $\epsilon\colon \Module{\dbl{E}} \to \dbl{E}$ of the adjunction is
the functor of virtual double categories taking the underlying object and
underlying arrow of an internal category and internal functor, and the
underlying proarrow and cell of any internal profunctor and any cell. For any
normal or unitary lax functor $H\colon \dbl{D}\to\Module{\dbl{E}}$, pushforward 
by the counit $\epsilon\colon \Module{\dbl{E}}\to\dbl{E}$ defines a lax functor
$\epsilon H\colon \dbl{D}\to\dbl{E}$ that is not necessarily normal since
$\epsilon$ is not normal.

\begin{example} \label{ex:unitandcounitMODadjunctionforSPANandPROF}
  In the case where $\dbl{E} = \Span$, the unit $\eta: \Span\to\Prof$ is the
  inclusion of sets as discrete categories, whereas the counit $\epsilon$ is the
  forgetful lax functor $\Ob\colon \Prof\to\Span$ taking the underlying set of
  objects of a given category (\cref{ex:ob-functor}). This example generalizes
  to the double category $\dbl{E} = \Span(\cat{S})$ for any finitely complete
  category $\cat{S}$.
\end{example}

The assignments $F\mapsto \Module{F}\eta$ and $H\mapsto \epsilon H$ are easily
seen to be mutually inverse provided they are well-defined. The crucial point is
that the assignment $F\mapsto \Module{F}\eta$ results in a genuine unitary lax
functor of double categories. This would be the case if both $\eta$ and
$\Module{F}$ were unitary lax functors. If both $\dbl{D}$ and $\dbl{E}$ are
equipments with local coequalizers, then $\Module{F}$ \emph{is} a unitary lax
functor is proved in in \cite[Proposition 11.11]{shulman2008}. However, we are
not working under the assumption that $\dbl{D}$ is even an equipment. Thus, the
result needs to be proved.

\begin{proposition}[Unitalization of lax functors, zero-dimensional]
  \label{prop:unitalizationoflaxdoublefunctor}
  If $\dbl{D}$ is a double category and $\dbl{E}$ is an equipment with local
  coequalizers, then the composite
  \begin{equation*}
    \dbl{D} \xrightarrow{\eta} \Module{\dbl{D}} \xrightarrow{\Module{F}} \Module{\dbl{E}}
  \end{equation*}
  amounts canonically to a unitary lax double functor. The assignment
  $H\mapsto \epsilon H$ thus induces a one-to-one correspondence between unitary
  lax functors $H\colon \dbl{D}\to\Module{\dbl{E}}$ and lax functors
  $F\colon \dbl{D}\to\dbl{E}$ whose inverse correspondence is
  $F\mapsto \Module{F}\eta$.
\end{proposition}
\begin{proof}
  We exhibit the laxators for the proposed composite $\Module{F}\eta$. Suppose
  $m\colon x\proto y$ and $n\colon y\proto z$ are composable proarrows in
  $\dbl{D}$. The given laxator
  $F_{m,n}\colon Fm\odot Fn \Rightarrow F(m\odot n)$ for $F$ then coequalizes
  the actions that define bimodule composition in $\Module{\dbl{E}}$ as above in
  \cref{equation:modulecompositionascoequalizer}. Therefore, there exists a
  unique cell $F_{m,n}^\otimes$ making the diagram
  \begin{equation*}
    \begin{tikzcd}
      {Fm\odot F\id_y\odot Fn} & {Fm\odot Fn} & {Fm \otimes Fn} \\
      && {F(m\odot n)}
      \arrow["{\rho\odot 1}", shift left=2, Rightarrow, from=1-1, to=1-2]
      \arrow["1\odot\lambda"', shift right=2, Rightarrow, from=1-1, to=1-2]
      \arrow["{\text{coeq}}", Rightarrow, from=1-2, to=1-3]
      \arrow["{F_{m,n}}"', Rightarrow, from=1-2, to=2-3]
      \arrow["{F_{m,n}^\otimes}", Rightarrow, dashed, from=1-3, to=2-3]
    \end{tikzcd}
  \end{equation*}
  commute. The computation showing that $F_{m,n}$ does coequalize the left and
  right actions \eqref{equation:inducedactionoftrivialmonoidunderlaxfunctor} is
  as follows:
  \begin{equation*}
    \begin{dblArray}{cccc}
      \Block{1-3}{F_{m,\id_y}} &&& 1_{Fn} \\
      \Block{1-3}{F(\cong)} &&& 1_{Fn} \\
      \Block{1-4}{F_{m,n}} &&&
    \end{dblArray} =
    \begin{dblArray}{cccc}
      \Block{1-3}{F_{m,\id_y}} &&& 1_{Fn} \\
      \Block{1-4}{F_{m\odot\id_y,n}} &&&  \\
      \Block{1-4}{F(\cong\odot 1)} &&&
    \end{dblArray} =
    \begin{dblArray}{cccc}
      1_{Fm} & \Block{1-3}{F_{\id_y,n}} && \\
      \Block{1-4}{F_{m,\id_y\odot n}} &&&  \\
      \Block{1-4}{F(1\odot \cong)} &&&
    \end{dblArray} =
    \begin{dblArray}{cccc}
      1_{Fm} & \Block{1-3}{F_{\id_y,n}} && \\
      1_{Fm} & \Block{1-3}{F(\cong)} && \\
      \Block{1-4}{F_{m,n}} &&&
    \end{dblArray}.
  \end{equation*}
  Note that all the external compositions in the computation involve the functor
  $\odot$ included in $\dbl{E}$, not the new external composition $\otimes$ for
  bimodules in $\dbl{E}$. The first equality is the laxator naturality
  condition; the second is the associativity condition for laxators; the final
  is again naturality.

  It needs to be seen that we have a unitary lax functor. To this end,
  there are two conditions to verify, namely, the naturality and associativity
  of the laxators as induced above in \cref{equation:modulecompositionascoequalizer}.
  Both computations employ the same strategy, namely, showing that the purportedly
  equal diagrams are equalized by the same coequalizers and therefore must be
  equal by uniqueness. We shall leave naturality to the reader and show the
  computations for the associativity condition since these are slighly more
  involved. We need to show that the equality
  \begin{equation} \label{equation:associativityconditiontobeproved}
    F_{m\odot n,p}^\otimes (F_{m,n}^\otimes \otimes 1_{Fp}) =
    F_{m,n\odot p}^\otimes (1_{Fm} \otimes F_{n,p}^\otimes)
  \end{equation}
  holds in $\Module{\dbl{E}}$. Some care is required in forming the external
  cell composites on each side. For example, $F_{m,n}^\otimes \otimes 1_{Fp}$ is
  induced from a local coequalizer and thus by construction satisfies
  \begin{equation*}
    \begin{tikzcd}
      \cdot && \cdot \\
      \cdot && \cdot \\
      \cdot && \cdot
      \arrow[Rightarrow, no head, from=1-1, to=2-1]
      \arrow[Rightarrow, no head, from=1-3, to=2-3]
      \arrow[Rightarrow, no head, from=2-1, to=3-1]
      \arrow[""{name=0, anchor=center, inner sep=0}, "{F(m\odot n)\otimes Fp}"', "\shortmid"{marking}, from=3-1, to=3-3]
      \arrow[Rightarrow, no head, from=2-3, to=3-3]
      \arrow[""{name=1, anchor=center, inner sep=0}, "{(Fm\otimes Fn)\otimes Fp}"', "\shortmid"{marking}, from=2-1, to=2-3]
      \arrow[""{name=2, anchor=center, inner sep=0}, "{(Fm\otimes Fn)\odot Fp}", "\shortmid"{marking}, from=1-1, to=1-3]
      \arrow["{F_{m,n}^\otimes \otimes 1_{Fp}}"{description, pos=0.6}, draw=none, from=1, to=0]
      \arrow["{\mathrm{coeq}}"{description}, draw=none, from=2, to=1]
    \end{tikzcd}
    \quad=\quad
    \begin{tikzcd}
      \cdot & \cdot & \cdot \\
      \cdot & \cdot & \cdot \\
      \cdot && \cdot
      \arrow[""{name=0, anchor=center, inner sep=0}, "{Fm\otimes Fn}", "\shortmid"{marking}, from=1-1, to=1-2]
      \arrow[""{name=1, anchor=center, inner sep=0}, "Fp", "\shortmid"{marking}, from=1-2, to=1-3]
      \arrow[Rightarrow, no head, from=1-3, to=2-3]
      \arrow[Rightarrow, no head, from=1-1, to=2-1]
      \arrow[Rightarrow, no head, from=2-1, to=3-1]
      \arrow[""{name=2, anchor=center, inner sep=0}, "{F(m\odot n)\otimes Fp}"', "\shortmid"{marking}, from=3-1, to=3-3]
      \arrow[Rightarrow, no head, from=2-3, to=3-3]
      \arrow[""{name=3, anchor=center, inner sep=0}, "{F(m\odot n)}"', "\shortmid"{marking}, from=2-1, to=2-2]
      \arrow[""{name=4, anchor=center, inner sep=0}, "Fp"', "\shortmid"{marking}, from=2-2, to=2-3]
      \arrow[Rightarrow, no head, from=1-2, to=2-2]
      \arrow["{\mathrm{coeq}}"{description}, draw=none, from=2-2, to=2]
      \arrow["{F_{m,n}^\otimes}"{description}, draw=none, from=0, to=3]
      \arrow["{1_{Fp}}"{description}, draw=none, from=1, to=4]
    \end{tikzcd},
  \end{equation*}
  which is just instantiating the observation made in
  \cref{equation:auxiliaryequationforexternalcompositeofmodulations}. There is an
  analogous equation for $1_{Fm} \otimes F_{n,p}^\otimes$. Now, as a result,
  the required computation for the associativity condition is:
  \begin{equation*}
    \begin{dblArray}{cc}
      \mathrm{coeq} & 1_{Fp} \\
      \Block{1-2}{\mathrm{coeq}} \\
      \Block{1-2}{F_{m,n}^\otimes \otimes 1_{Fp}} \\
      \Block{1-2}{F_{m\odot n,p}^\otimes}
    \end{dblArray} =
    \begin{dblArray}{cc}
      \mathrm{coeq} & 1_{Fp} \\
      F_{m,n}^\otimes & 1_{Fp} \\
      \Block{1-2}{\mathrm{coeq}} \\
      \Block{1-2}{F_{m\odot n,p}^\otimes}
    \end{dblArray} =
    \begin{dblArray}{cc}
      F_{m,n} & 1_{Fp} \\
      \Block{1-2}{F_{m\odot n,p}}
    \end{dblArray} =
    \begin{dblArray}{cc}
      1_{Fm} & F_{n,p} \\
      \Block{1-2}{F_{m,n\odot p}}
    \end{dblArray} =
    \begin{dblArray}{cc}
      1_{Fm} & \mathrm{coeq} \\
      1_{Fm} & F_{n,p}^\otimes \\
      \Block{1-2}{\mathrm{coeq}} \\
      \Block{1-2}{F_{m,n\odot p}^\otimes}
    \end{dblArray} =
    \begin{dblArray}{cc}
      1_{Fm} & \mathrm{coeq} \\
      \Block{1-2}{\mathrm{coeq}} \\
      \Block{1-2}{1_{Fm}\otimes F_{n,p}^\otimes} \\
      \Block{1-2}{F_{m,n\odot p}^\otimes}
    \end{dblArray}.
  \end{equation*}
  The middle equality is the assumed associativity condition for the given
  laxators coming with $F$. The leftmost and rightmost follow by the
  construction of the external composite of cells as noted in the penultimate
  display. The middle-left and middle-right follow by construction of the
  induced laxators $F^\otimes$ in
  \cref{equation:modulecompositionascoequalizer}. This computation proves the
  required associativity condition by uniqueness of morphisms induced by
  coequalizers. For the composite coequalizers on the far sides of the
  computation are each coequalizers of the same diagram by an application of the
  \emph{3x3 lemma} \cite[Lemma 0.17]{johnstone1977}.
\end{proof}

\begin{remark}
  As noted in the proof of \cite[Proposition 5.14]{cruttwell2010} for the
  virtual case of the above result, the one-to-one correspondence depends upon a
  choice of units in our (virtual) double categories of the form
  $\Module{\dbl{D}}$, which we are happy to assume. Without such a choice, the
  correspondence is only that of pseudo-inverses.
\end{remark}

The lemma finally proves \cref{corollary:universalpropertyofPROF} since
$\Module{\Span} = \Prof$. Now we specialize the correspondence of
\cref{prop:unitalizationoflaxdoublefunctor} to \emph{cartesian} lax functors.

\begin{corollary}[Unitalization of cartesian lax functors, zero-dimensional]
  \label{cor:unitalizationofCartesianLaxFunctors}
  Suppose $\dbl{D}$ is a cartesian double category, $\dbl{E}$ is a cartesian
  equipment with local coequalizers, and $F\colon \dbl{D}\to\dbl{E}$ is a
  cartesian lax double functor. Then the composite
  \begin{equation*}
    \dbl{D} \xrightarrow{\eta} \Module{\dbl{D}} \xrightarrow{\Module{F}} \Module{\dbl{E}}
  \end{equation*}
  is a cartesian unitary lax double functor. The assignment
  $H\mapsto \epsilon H$ thus induces a one-to-one correspondence between
  cartesian unitary lax double functors $H\colon \dbl{D}\to\Module{\dbl{E}}$
  and cartesian lax double functors $F\colon \dbl{D}\to\dbl{E}$ whose inverse
  correspondence is $F\mapsto \Module{F}\eta$.
\end{corollary}
\begin{proof}
  The lax functor $\Module{F} \eta$ is cartesian when $F$ is, owing to the
  construction of products in $\Module{\dbl{E}}$. Similarly, the counit lax
  functor $\epsilon\colon \Module{\dbl{E}}\to\dbl{E}$ is cartesian by
  construction of products in $\Module{\dbl{E}}$. Since cartesian lax double
  functors, normal or not, are closed under composition, it follows that the
  composite $\epsilon H$ is cartesian when $H$ is.
\end{proof}

In particular, taking $\dbl{E} = \Span$, we deduce the analogue of
\cref{corollary:universalpropertyofPROF} for cartesian lax double functors.

\begin{corollary} \label{cor:spanvaluedprofvaluedcartesianfunctorsonetooneocorrespondence}
  For any cartesian lax functor $F\colon \dbl{D}\to \Span$, there is a unique
  cartesian \emph{unitary} lax double functor $\bar F: \dbl{D} \to \Prof$
  making the triangle commute:
  \begin{equation*}
    \begin{tikzcd}
      & \Prof \\
      {\dbl{D}} & \Span
      \arrow["{\Ob}", from=1-2, to=2-2]
      \arrow["F"', from=2-1, to=2-2]
      \arrow["{\bar F}", dashed, from=2-1, to=1-2]
    \end{tikzcd}.
  \end{equation*}
\end{corollary}

\section{Cartesian double theories and models}
\label{sec:double-theories}

The expressivity of double theories is significantly increased by introducing
finite products, in the sense of cartesian double categories. Models of
cartesian double theories encompass most kinds of monoidal categories, among
many other things.

\begin{definition}[Cartesian double theory]
  A \define{cartesian double theory} is a small, cartesian, strict double
  category $\dbl{T}$. A \define{morphism} between cartesian double theories
  $\dbl{T}$ and $\dbl{T}'$ is a cartesian strict double functor
  $\dbl{T} \to \dbl{T}'$.

  A \define{model} of a cartesian double theory $\dbl{T}$ in a cartesian double
  category $\dbl{S}$ is a cartesian lax double functor $\dbl{T} \to \dbl{S}$, in
  which case $\dbl{S}$ is called the \define{semantics}.
\end{definition}

As before, we will occasionally speak of \define{strict} or \define{pseudo}
models of a cartesian double theory, meaning cartesian strict or pseudo double
functors out of the theory.

\begin{remark}[Presenting cartesian double theories]
  Our remarks on presenting simple double theories
  (\cref{rem:presenting-simple-theories}) carry over to cartesian double
  theories because the category of small, strict double categories with chosen
  finite products and strict double functors strictly preserving those products
  is again the category of models of a finite limit sketch. Moreover, chosen
  finite products can be freely added to a strict double category, using the
  result that pullback functors between categories of models of finite limit
  sketches have left adjoints \cite[Theorem 4.4.1]{barr1985}.
\end{remark}

Since any simple double theory can be turned into a cartesian double theory by
taking its free finite product completion, all of the examples of simple double
theories in \cref{sec:simple-double-theories} are also examples of cartesian
double theories. The examples in this section will make nontrivial use of the
cartesian structure. The first example uses only the terminal object.

\begin{theory}[Copresheaves]
  The \define{theory of families} $\Th{\cat{Fam}}$ is freely generated by an
  object $x$ and a proarrow $p: 1 \proto x$. A strict model in $\Mat$ of the
  theory of families is a set $X$ together a family of sets indexed by $X$,
  denoted $P: X \to \Set$. A model in $\Mat$ is a category $\cat{C}$ together
  with a copresheaf $P: \cat{C} \to \Set$.

  Dually, the theory freely generated by a proarrow $p: x \proto 1$ has
  presheaves as models.
\end{theory}

We now begin to use nonempty products in cartesian double theories.

\begin{theory}[Strict monoidal categories] \label{th:monoid}
  The \define{theory of monoids} $\Th{\cat{Mon}}$ is generated by
  \begin{itemize}[noitemsep]
    \item an object $x$, and
    \item arrows $\otimes: x^2 \to x$ and $I: 1 \to x$,
  \end{itemize}
  subject to the usual associativity and unitality equations:
  \begin{equation*}
    \begin{tikzcd}
      {x^3} & {x^2} \\
      {x^2} & x
      \arrow["{1_x \times \otimes}", from=1-1, to=1-2]
      \arrow["\otimes", from=1-2, to=2-2]
      \arrow["{\otimes \times 1_x}"', from=1-1, to=2-1]
      \arrow["\otimes"', from=2-1, to=2-2]
    \end{tikzcd}
    \hspace{1in}
    \begin{tikzcd}
      x & {x^2} & x \\
      & x
      \arrow[Rightarrow, no head, from=1-1, to=2-2]
      \arrow["{I \times 1_x}", from=1-1, to=1-2]
      \arrow["{1_x \times I}"', from=1-3, to=1-2]
      \arrow[Rightarrow, no head, from=1-3, to=2-2]
      \arrow["\otimes"', from=1-2, to=2-2]
    \end{tikzcd}.
  \end{equation*}
  In other words, $\Th{\cat{Mon}}$ is the double theory $\dbl{T}$ whose
  cartesian category of objects $\dbl{T}_0$ is the usual Lawvere theory of
  monoids and whose category of morphisms $\dbl{T}_1$ is trivial.

  A model of the theory of monoids is a strict monoidal category. A strict model
  of the theory is merely a monoid.
\end{theory}

\begin{theory}[Strict 2-groups]
  The \define{theory of groups} is the theory of monoids augmented with an arrow
  $i: x\to x$ satisfying two further equations,
  \begin{equation*}
    \begin{tikzcd}
      x & 1 & x \\
      {x^2} && {x^2}
      \arrow["{!}", from=1-1, to=1-2]
      \arrow["I", from=1-2, to=1-3]
      \arrow["\Delta"', from=1-1, to=2-1]
      \arrow["{i \times 1}"', from=2-1, to=2-3]
      \arrow["\otimes"', from=2-3, to=1-3]
    \end{tikzcd}
    \qquad\text{and}\qquad
    \begin{tikzcd}
      x & 1 & x \\
      {x^2} && {x^2}
      \arrow["\Delta"', from=1-1, to=2-1]
      \arrow["{1 \times i}"', from=2-1, to=2-3]
      \arrow["{!}", from=1-1, to=1-2]
      \arrow["I", from=1-2, to=1-3]
      \arrow["\otimes"', from=2-3, to=1-3]
    \end{tikzcd},
  \end{equation*}
  saying that the arrow $i$ picks out an inverse for each element of the
  underlying object $x$.

  A strict model in $\Span$ is an ordinary group, while a strict model in
  $\Span(\mathbf{Top})$ is a topological group. Generally, for a category
  $\cat{C}$ with finite limits, a strict model in $\Span(\cat{C})$ is a group
  object in $\cat{C}$. Thus, a strict model in $\Span(\Cat)$ is a strict 2-group
  \cite{baezlauda2004}. Equivalently, a strict 2-group is a (lax) model in
  $\Span$, or by
  \cref{cor:spanvaluedprofvaluedcartesianfunctorsonetooneocorrespondence} a
  unitary lax model in $\Prof$.
\end{theory}

\begin{theory}[Commutative monoidal categories] \label{th:commutative-monoid}
  The \define{theory of commutative monoids} $\Th{\cat{CMon}}$ is the theory of
  monoids augmented with the commutativity equation
  \begin{equation*}
    \begin{tikzcd}
      {x^2} & {x^2} \\
      & x
      \arrow["{\sigma_{x,x}}", from=1-1, to=1-2]
      \arrow["\otimes"', from=1-1, to=2-2]
      \arrow["\otimes", from=1-2, to=2-2]
    \end{tikzcd}.
  \end{equation*}
  A model of $\Th{\cat{CMon}}$ is a \define{commutative monoidal category}: a
  symmetric monoidal category whose associators, unitors, and braidings are all
  identities. A strict model of the theory is merely a commutative monoid.
\end{theory}

\begin{theory}[Monoidal categories] \label{th:pseudomonoid}
  The \define{theory of pseudomonoids} $\Th{\cat{PsMon}}$ is generated by
  \begin{itemize}[noitemsep]
    \item an object $x$,
    \item arrows $\otimes: x^2 \to x$ and $I: 1 \to x$, and
    \item \define{associator} and \define{unitor} cells
      \begin{equation*}
        \begin{tikzcd}
          {x^3} & {x^3} \\
          {x^2} & {x^2} \\
          x & x
          \arrow[""{name=0, anchor=center, inner sep=0}, "{\mathrm{id}_x^3}", "\shortmid"{marking}, from=1-1, to=1-2]
          \arrow["{\otimes \times 1_x}"', from=1-1, to=2-1]
          \arrow["\otimes"', from=2-1, to=3-1]
          \arrow["{1_x \times \otimes}", from=1-2, to=2-2]
          \arrow[""{name=1, anchor=center, inner sep=0}, "{\mathrm{id}_x}"', "\shortmid"{marking}, from=3-1, to=3-2]
          \arrow["\otimes", from=2-2, to=3-2]
          \arrow["\alpha"{description}, draw=none, from=0, to=1]
        \end{tikzcd}
        \qquad
        \begin{tikzcd}
          x & x \\
          {x^2} \\
          x & x
          \arrow["{I \times 1_x}"', from=1-1, to=2-1]
          \arrow["\otimes"', from=2-1, to=3-1]
          \arrow[""{name=0, anchor=center, inner sep=0}, "{\mathrm{id}_x}"', "\shortmid"{marking}, from=3-1, to=3-2]
          \arrow[""{name=1, anchor=center, inner sep=0}, "{\mathrm{id}_x}", "\shortmid"{marking}, from=1-1, to=1-2]
          \arrow[Rightarrow, no head, from=1-2, to=3-2]
          \arrow["\lambda"{description}, draw=none, from=1, to=0]
        \end{tikzcd}
        \qquad
        \begin{tikzcd}
          x & x \\
          {x^2} \\
          x & x
          \arrow["{1_x \times I}"', from=1-1, to=2-1]
          \arrow["\otimes"', from=2-1, to=3-1]
          \arrow[""{name=0, anchor=center, inner sep=0}, "{\mathrm{id}_x}"', "\shortmid"{marking}, from=3-1, to=3-2]
          \arrow[""{name=1, anchor=center, inner sep=0}, "{\mathrm{id}_x}", "\shortmid"{marking}, from=1-1, to=1-2]
          \arrow[Rightarrow, no head, from=1-2, to=3-2]
          \arrow["\rho"{description}, draw=none, from=1, to=0]
        \end{tikzcd}
      \end{equation*}
      along with their inverses $\alpha^{-1}: \id_x^3 \to \id_x$ and
      $\lambda^{-1}, \rho^{-1}: \id_x \to \id_x$,
  \end{itemize}
  subject to the \define{pentagon identity}
  \begin{equation*}
    \begin{tikzcd}
      {x^4} & {x^4} && {x^4} & {x^4} \\
      {x^3} & {x^3} && {x^3} & {x^3} \\
      {x^2} & {x^2} && {x^2} & {x^2} \\
      x & x && x & x
      \arrow[""{name=0, anchor=center, inner sep=0}, "\shortmid"{marking}, Rightarrow, no head, from=3-1, to=3-2]
      \arrow["\otimes"', from=3-1, to=4-1]
      \arrow["\otimes", from=3-2, to=4-2]
      \arrow[""{name=1, anchor=center, inner sep=0}, "\shortmid"{marking}, Rightarrow, no head, from=4-1, to=4-2]
      \arrow[""{name=2, anchor=center, inner sep=0}, "\shortmid"{marking}, Rightarrow, no head, from=4-2, to=4-4]
      \arrow["{\otimes 1_x^2}"', from=1-1, to=2-1]
      \arrow["{\otimes 1_x}"', from=2-1, to=3-1]
      \arrow["{1_x \otimes 1_x}"{description}, from=1-2, to=2-2]
      \arrow["{\otimes 1_x}", from=2-2, to=3-2]
      \arrow["{1_x \otimes 1_x}"{description}, from=1-4, to=2-4]
      \arrow[""{name=3, anchor=center, inner sep=0}, "\shortmid"{marking}, Rightarrow, no head, from=1-2, to=1-4]
      \arrow[""{name=4, anchor=center, inner sep=0}, "\shortmid"{marking}, Rightarrow, no head, from=1-1, to=1-2]
      \arrow[""{name=5, anchor=center, inner sep=0}, "\shortmid"{marking}, Rightarrow, no head, from=2-2, to=2-4]
      \arrow["{1_x \otimes}"', from=2-4, to=3-4]
      \arrow["\otimes"', from=3-4, to=4-4]
      \arrow[""{name=6, anchor=center, inner sep=0}, "\shortmid"{marking}, Rightarrow, no head, from=3-4, to=3-5]
      \arrow[""{name=7, anchor=center, inner sep=0}, "\shortmid"{marking}, Rightarrow, no head, from=4-4, to=4-5]
      \arrow["\otimes", from=3-5, to=4-5]
      \arrow[""{name=8, anchor=center, inner sep=0}, "\shortmid"{marking}, Rightarrow, no head, from=1-4, to=1-5]
      \arrow["{1_x^2 \otimes}", from=1-5, to=2-5]
      \arrow["{1_x \otimes}", from=2-5, to=3-5]
      \arrow["{\mathrm{id}_\otimes}"{description}, draw=none, from=0, to=1]
      \arrow["{\alpha\;\mathrm{id}_{1_x}}"{description}, draw=none, from=4, to=0]
      \arrow["{\mathrm{id}_{1_x \otimes 1_x}}"{description}, draw=none, from=3, to=5]
      \arrow["\alpha"{description}, draw=none, from=5, to=2]
      \arrow["{\mathrm{id}_\otimes}"{description}, draw=none, from=6, to=7]
      \arrow["{\mathrm{id}_{1_x}\;\alpha}"{description}, draw=none, from=8, to=6]
    \end{tikzcd}
    \quad=\quad
    \begin{tikzcd}
      {x^4} & {x^4} & {x^4} & {x^4} \\
      {x^3} & {x^3} & {x^3} & {x^3} \\
      {x^2} & {x^2} & {x^2} & {x^2} \\
      x & x & x & x
      \arrow[""{name=0, anchor=center, inner sep=0}, "\shortmid"{marking}, Rightarrow, no head, from=1-1, to=1-2]
      \arrow["{\otimes 1_x^2}"', from=1-1, to=2-1]
      \arrow["{\otimes 1_x^2}", from=1-2, to=2-2]
      \arrow[""{name=1, anchor=center, inner sep=0}, "\shortmid"{marking}, Rightarrow, no head, from=2-1, to=2-2]
      \arrow["{\otimes 1_x}"', from=2-1, to=3-1]
      \arrow["\otimes"', from=3-1, to=4-1]
      \arrow["{1_x \otimes}", from=2-2, to=3-2]
      \arrow["\otimes", from=3-2, to=4-2]
      \arrow[""{name=2, anchor=center, inner sep=0}, "\shortmid"{marking}, Rightarrow, no head, from=4-1, to=4-2]
      \arrow[""{name=3, anchor=center, inner sep=0}, "\shortmid"{marking}, Rightarrow, no head, from=1-2, to=1-3]
      \arrow["{1_x^2 \otimes}"', from=1-3, to=2-3]
      \arrow["{\otimes 1_x}"', from=2-3, to=3-3]
      \arrow["\otimes"', from=3-3, to=4-3]
      \arrow[""{name=4, anchor=center, inner sep=0}, "\shortmid"{marking}, Rightarrow, no head, from=3-2, to=3-3]
      \arrow[""{name=5, anchor=center, inner sep=0}, "\shortmid"{marking}, Rightarrow, no head, from=4-2, to=4-3]
      \arrow[""{name=6, anchor=center, inner sep=0}, "\shortmid"{marking}, Rightarrow, no head, from=1-3, to=1-4]
      \arrow[""{name=7, anchor=center, inner sep=0}, "\shortmid"{marking}, Rightarrow, no head, from=2-3, to=2-4]
      \arrow["{1_x^2 \otimes}", from=1-4, to=2-4]
      \arrow["{1_x \otimes}", from=2-4, to=3-4]
      \arrow["\otimes", from=3-4, to=4-4]
      \arrow[""{name=8, anchor=center, inner sep=0}, "\shortmid"{marking}, Rightarrow, no head, from=4-3, to=4-4]
      \arrow["{\mathrm{id}_{\otimes 1_x^2}}"{description}, draw=none, from=0, to=1]
      \arrow["\alpha"{description}, draw=none, from=1, to=2]
      \arrow["{\mathrm{id}_{\otimes^2}}"{description}, draw=none, from=3, to=4]
      \arrow["{\mathrm{id}_\otimes}"{description}, draw=none, from=4, to=5]
      \arrow["{\mathrm{id}_{1_x^2 \otimes}}"{description}, draw=none, from=6, to=7]
      \arrow["\alpha"{description}, draw=none, from=7, to=8]
    \end{tikzcd}
  \end{equation*}
  and the \define{triangle identity}
  \begin{equation*}
    \begin{tikzcd}
      {x^2} & {x^2} & {x^2} \\
      {x^3} & {x^3} \\
      {x^2} & {x^2} & {x^2} \\
      x & {x } & x
      \arrow["{1_x I 1_x}"', from=1-1, to=2-1]
      \arrow["{1_x I 1_x}", from=1-2, to=2-2]
      \arrow[""{name=0, anchor=center, inner sep=0}, "\shortmid"{marking}, Rightarrow, no head, from=1-1, to=1-2]
      \arrow[""{name=1, anchor=center, inner sep=0}, "\shortmid"{marking}, Rightarrow, no head, from=2-1, to=2-2]
      \arrow["{\otimes 1_x}"', from=2-1, to=3-1]
      \arrow["\otimes"', from=3-1, to=4-1]
      \arrow["{1_x \otimes}"', from=2-2, to=3-2]
      \arrow["\otimes"', from=3-2, to=4-2]
      \arrow[""{name=2, anchor=center, inner sep=0}, "\shortmid"{marking}, Rightarrow, no head, from=4-1, to=4-2]
      \arrow[Rightarrow, no head, from=1-3, to=3-3]
      \arrow[""{name=3, anchor=center, inner sep=0}, "\shortmid"{marking}, Rightarrow, no head, from=3-2, to=3-3]
      \arrow["\otimes", from=3-3, to=4-3]
      \arrow[""{name=4, anchor=center, inner sep=0}, "\shortmid"{marking}, Rightarrow, no head, from=4-2, to=4-3]
      \arrow[""{name=5, anchor=center, inner sep=0}, "\shortmid"{marking}, Rightarrow, no head, from=1-2, to=1-3]
      \arrow["{\mathrm{id}_{1_x I 1_x}}"{description}, draw=none, from=0, to=1]
      \arrow["\alpha"{description}, draw=none, from=1, to=2]
      \arrow["{\mathrm{id}_\otimes}"{description}, draw=none, from=3, to=4]
      \arrow["{\mathrm{id}_{1_x}\; \lambda}"{description}, draw=none, from=5, to=3]
    \end{tikzcd}
    \quad=\quad
    \begin{tikzcd}
      {x^2} & {x^2} \\
      {x^3} \\
      {x^2} & {x^2} \\
      {x } & x
      \arrow["{1_x I 1_x}"', from=1-1, to=2-1]
      \arrow["{\otimes 1_x}"', from=2-1, to=3-1]
      \arrow["\otimes"', from=3-1, to=4-1]
      \arrow[Rightarrow, no head, from=1-2, to=3-2]
      \arrow[""{name=0, anchor=center, inner sep=0}, "\shortmid"{marking}, Rightarrow, no head, from=3-1, to=3-2]
      \arrow["\otimes", from=3-2, to=4-2]
      \arrow[""{name=1, anchor=center, inner sep=0}, "\shortmid"{marking}, Rightarrow, no head, from=4-1, to=4-2]
      \arrow[""{name=2, anchor=center, inner sep=0}, "\shortmid"{marking}, Rightarrow, no head, from=1-1, to=1-2]
      \arrow["{\mathrm{id}_\otimes}"{description}, draw=none, from=0, to=1]
      \arrow["{\rho\; \mathrm{id}_{1_x}}"{description}, draw=none, from=2, to=0]
    \end{tikzcd}.
  \end{equation*}
  A model of the theory of pseudomonoids is a (weak) monoidal category, whereas
  a strict model is again just a monoid.
\end{theory}

\begin{theory}[Symmetric monoidal categories] \label{th:symmetric-pseudomonoid}
  The \define{theory of symmetric pseudomonoids} $\Th{\cat{SPsMon}}$ is the
  theory of pseudomonoids augmented with a \define{braiding cell}
  \begin{equation*}
    \begin{tikzcd}
      {x^2} & {x^2} \\
      & {x^2} \\
      x & x
      \arrow[""{name=0, anchor=center, inner sep=0}, "{\mathrm{id}_x^2}", "\shortmid"{marking}, from=1-1, to=1-2]
      \arrow["\otimes"', from=1-1, to=3-1]
      \arrow["{\sigma_{x,x}}", from=1-2, to=2-2]
      \arrow["\otimes", from=2-2, to=3-2]
      \arrow[""{name=1, anchor=center, inner sep=0}, "{\mathrm{id}_x}"', "\shortmid"{marking}, from=3-1, to=3-2]
      \arrow["\sigma"{description}, draw=none, from=0, to=1]
    \end{tikzcd}
  \end{equation*}
  subject to the involutivity equation
  \begin{equation*}
    \begin{tikzcd}[row sep=small]
      {x^2} & {x^2} & {x^2} \\
      & {x^2} & {x^2} \\
      && {x^2} \\
      x & x & x
      \arrow[""{name=0, anchor=center, inner sep=0}, "\shortmid"{marking}, Rightarrow, no head, from=2-2, to=2-3]
      \arrow["\otimes"', from=2-2, to=4-2]
      \arrow["{\sigma_{x,x}}", from=2-3, to=3-3]
      \arrow["\otimes", from=3-3, to=4-3]
      \arrow[""{name=1, anchor=center, inner sep=0}, "\shortmid"{marking}, Rightarrow, no head, from=4-2, to=4-3]
      \arrow["{\sigma_{x,x}}"', from=1-2, to=2-2]
      \arrow["\otimes"', from=1-1, to=4-1]
      \arrow[""{name=2, anchor=center, inner sep=0}, "\shortmid"{marking}, Rightarrow, no head, from=4-1, to=4-2]
      \arrow[""{name=3, anchor=center, inner sep=0}, "\shortmid"{marking}, Rightarrow, no head, from=1-1, to=1-2]
      \arrow["{\sigma_{x,x}}", from=1-3, to=2-3]
      \arrow[""{name=4, anchor=center, inner sep=0}, "\shortmid"{marking}, Rightarrow, no head, from=1-2, to=1-3]
      \arrow["\sigma"{description}, draw=none, from=0, to=1]
      \arrow["{\mathrm{id}_{\sigma_{x,x}}}"{description}, draw=none, from=4, to=0]
      \arrow["\sigma"{description}, draw=none, from=3, to=2]
    \end{tikzcd}
    \quad=\quad
    \begin{tikzcd}
      {x^2} & {x^2} \\
      x & x
      \arrow[""{name=0, anchor=center, inner sep=0}, "\shortmid"{marking}, Rightarrow, no head, from=1-1, to=1-2]
      \arrow[""{name=1, anchor=center, inner sep=0}, "\shortmid"{marking}, Rightarrow, no head, from=2-1, to=2-2]
      \arrow["\otimes"', from=1-1, to=2-1]
      \arrow["\otimes", from=1-2, to=2-2]
      \arrow["{\mathrm{id}_\otimes}"{description}, draw=none, from=0, to=1]
    \end{tikzcd}
  \end{equation*}
  and the first \define{hexagon identity}
  \begin{equation*}
    \begin{tikzcd}[row sep=small]
      {x^3} & {x^3} & {x^3} & {x^3} \\
      & {x^3} & {x^3} & {x^3} \\
      &&& {x^3} \\
      {x^2} & {x^2} & {x^2} & {x^2} \\
      x & x & x & x
      \arrow["{\otimes 1_x}"', from=1-1, to=4-1]
      \arrow["{\sigma_{x,x} 1_x}"', from=1-2, to=2-2]
      \arrow["{\otimes 1_x}", from=2-2, to=4-2]
      \arrow[""{name=0, anchor=center, inner sep=0}, "\shortmid"{marking}, Rightarrow, no head, from=4-1, to=4-2]
      \arrow[""{name=1, anchor=center, inner sep=0}, "\shortmid"{marking}, Rightarrow, no head, from=1-1, to=1-2]
      \arrow["\otimes"', from=4-1, to=5-1]
      \arrow["\otimes", from=4-2, to=5-2]
      \arrow[""{name=2, anchor=center, inner sep=0}, "\shortmid"{marking}, Rightarrow, no head, from=5-1, to=5-2]
      \arrow[""{name=3, anchor=center, inner sep=0}, "\shortmid"{marking}, Rightarrow, no head, from=2-2, to=2-3]
      \arrow["{1_x \otimes}"', from=2-3, to=4-3]
      \arrow["\otimes"', from=4-3, to=5-3]
      \arrow["\shortmid"{marking}, Rightarrow, no head, from=1-2, to=1-3]
      \arrow[""{name=4, anchor=center, inner sep=0}, "\shortmid"{marking}, Rightarrow, no head, from=5-2, to=5-3]
      \arrow["\otimes", from=4-4, to=5-4]
      \arrow[""{name=5, anchor=center, inner sep=0}, "\shortmid"{marking}, Rightarrow, no head, from=5-3, to=5-4]
      \arrow[""{name=6, anchor=center, inner sep=0}, "\shortmid"{marking}, Rightarrow, no head, from=4-3, to=4-4]
      \arrow["{1_x \otimes}", from=3-4, to=4-4]
      \arrow["{1_x \sigma_{x,x}}", from=2-4, to=3-4]
      \arrow[""{name=7, anchor=center, inner sep=0}, "\shortmid"{marking}, Rightarrow, no head, from=2-3, to=2-4]
      \arrow["{\sigma_{x,x} 1_x}", from=1-4, to=2-4]
      \arrow["\shortmid"{marking}, Rightarrow, no head, from=1-3, to=1-4]
      \arrow["{\mathrm{id}_{\sigma_{x,x} 1_x}}"{description}, draw=none, from=1-3, to=2-3]
      \arrow["\alpha"{description}, draw=none, from=3, to=4]
      \arrow["{\sigma\; \mathrm{id}_{1_x}}"{description}, draw=none, from=1, to=0]
      \arrow["{\mathrm{id}_\otimes}"{description}, draw=none, from=0, to=2]
      \arrow["{\mathrm{id}_{1_x}\; \sigma}"{description}, draw=none, from=7, to=6]
      \arrow["{\mathrm{id}_\otimes}"{description}, draw=none, from=6, to=5]
    \end{tikzcd}
    =
    \begin{tikzcd}[row sep=small]
      {x^3} & {x^3} & {x^3} & {x^3} & {x^3} \\
      &&& {x^3} & {x^3} \\
      {x^2} & {x^2} & {x^2} & {x^3} & {x^3} \\
      && {x^2} & {x^2} & {x^2} \\
      x & x & x & x & x
      \arrow[""{name=0, anchor=center, inner sep=0}, "\shortmid"{marking}, Rightarrow, no head, from=1-2, to=1-3]
      \arrow["{1_x \otimes}"', from=1-2, to=3-2]
      \arrow["{1_x \otimes}", from=1-3, to=3-3]
      \arrow[""{name=1, anchor=center, inner sep=0}, "\shortmid"{marking}, Rightarrow, no head, from=3-2, to=3-3]
      \arrow["\otimes"', from=3-2, to=5-2]
      \arrow["{\sigma_{x,x}}", from=3-3, to=4-3]
      \arrow["\otimes"', from=4-3, to=5-3]
      \arrow[""{name=2, anchor=center, inner sep=0}, "\shortmid"{marking}, Rightarrow, no head, from=5-2, to=5-3]
      \arrow["{\otimes 1_x}"', from=1-1, to=3-1]
      \arrow["{\sigma_{x,x} 1_x}", from=1-4, to=2-4]
      \arrow["{1_x \sigma_{xx}}", from=2-4, to=3-4]
      \arrow["{\otimes 1_x}", from=3-4, to=4-4]
      \arrow["\otimes", from=4-4, to=5-4]
      \arrow["{1_x \otimes}", from=3-5, to=4-5]
      \arrow["\otimes", from=4-5, to=5-5]
      \arrow[""{name=3, anchor=center, inner sep=0}, "\shortmid"{marking}, Rightarrow, no head, from=4-3, to=4-4]
      \arrow[""{name=4, anchor=center, inner sep=0}, "\shortmid"{marking}, Rightarrow, no head, from=5-3, to=5-4]
      \arrow[""{name=5, anchor=center, inner sep=0}, "\shortmid"{marking}, Rightarrow, no head, from=3-4, to=3-5]
      \arrow[""{name=6, anchor=center, inner sep=0}, "\shortmid"{marking}, Rightarrow, no head, from=5-4, to=5-5]
      \arrow["\otimes"', from=3-1, to=5-1]
      \arrow[""{name=7, anchor=center, inner sep=0}, "\shortmid"{marking}, Rightarrow, no head, from=1-3, to=1-4]
      \arrow["{\sigma_{x,x} 1_x}", from=1-5, to=2-5]
      \arrow["{1_x \sigma_{x,x}}", from=2-5, to=3-5]
      \arrow[""{name=8, anchor=center, inner sep=0}, "\shortmid"{marking}, Rightarrow, no head, from=1-4, to=1-5]
      \arrow[""{name=9, anchor=center, inner sep=0}, "\shortmid"{marking}, Rightarrow, no head, from=1-1, to=1-2]
      \arrow[""{name=10, anchor=center, inner sep=0}, "\shortmid"{marking}, Rightarrow, no head, from=5-1, to=5-2]
      \arrow["\sigma"{description}, draw=none, from=1, to=2]
      \arrow["{\mathrm{id}_{1_x \otimes}}"{description}, draw=none, from=0, to=1]
      \arrow["\alpha"{description}, draw=none, from=5, to=6]
      \arrow["{\mathrm{id}}"{description}, draw=none, from=7, to=3]
      \arrow["{\mathrm{id}_\otimes}"{description}, draw=none, from=3, to=4]
      \arrow["{\mathrm{id}}"{description}, draw=none, from=8, to=5]
      \arrow["\alpha"{description}, draw=none, from=9, to=10]
    \end{tikzcd}.
  \end{equation*}
  Note that the identity cell on the right-hand side is well-defined due to the
  naturality of the braiding inside the cartesian category of objects
  $\dbl{T}_0$ underlying $\dbl{T} = \Th{\cat{SPsMon}}$, an equation that is most
  clearly expressed using string diagrams.

  A model of the theory of symmetric pseudomonoids is a (weak) symmetric
  monoidal category, whereas a strict model is again just a commutative monoid.
\end{theory}

A theory of braided pseudomonoids, having braided model categories as models,
can be defined similarly. It replaces the involutivity axiom with a distinct
inverse to the braiding and adds a second hexagon identity.

\begin{theory}[Multicategories] \label{th:promonoid}
  The \define{theory of promonoids} is generated by
  \begin{itemize}[noitemsep]
    \item an object $x$ and
    \item proarrows $p: x^2 \proto x$ and $j: 1 \proto x$
  \end{itemize}
  subject to the axioms of associativity
  $(p \times \id_x) \odot p = (\id_x \times p) \odot p$ and unitality
  \begin{equation*}
    (j \times \id_x) \odot p = \id_x = (\id_x \times j) \odot p.
  \end{equation*}
  This is a finite presentation of a cartesian double theory such that for each
  arity $n \geq 0$, there is a \emph{unique} proarrow $p_n: x^n \proto x$ with
  source $x^n$ and target $x$. For example, we have $p_0 = j$, $p_1 = \id_x$,
  $p_2 = p$, and $p_4 = (p \times p) \odot p$.

  A model of the theory of promonoids is a \define{multicategory}. It consists
  of a set of objects $\cat{C}_0$; for each arity $n \geq 0$, a family of sets
  of $n$-ary multimorphisms
  \begin{equation*}
    \cat{C}(c_1, \dots, c_n; c) \in \Set, \qquad
    c_1, \dots, c_n, c \in \cat{C}_0;
  \end{equation*}
  and identity morphisms $1_c \in \cat{C}(c,c)$ for each $c \in \cat{C}_0$.
  Multimorphisms compose associatively and unitally, due to the associativity
  and unitality of the laxators. This example can be seen as a reformulation of
  Day and Street's observation that lax monoids in the monoidal bicategory of
  spans are multicategories \cite[\S{1}]{daystreet2003}.

  The theory of promonoids is so named first, because it is the transpose of the
  theory of monoids (\cref{th:monoid}) and second, because when the laxators of
  a model (in $\Span$) involving $p: x^2 \proto x$ or $j: 1 \proto x$ are
  required to be isomorphisms, the model is a \define{promonoidal category}
  \mbox{\cite[\S{3}]{day1970}}. Relatedly, a \emph{pseudo} model of the theory
  of promonoids in a cartesian double category $\dbl{D}$ is a pseudomonoid in
  $\HorBicat(\dbl{D})$, the horizontal monoidal bicategory of $\dbl{D}$
  \cite[\S{3}]{daystreet1997}. So a pseudo model in $\Prof$ is again a
  promonoidal category, whereas a pseudo model in $\Quintet(\Cat)$, the
  cartesian double category of quintets in $\Cat$ (\cref{ex:quintets}), is a
  monoidal category.
\end{theory}

In contrast with symmetric monoidal categories, there seems to be no obvious
presentation by generators and relations of a cartesian double theory whose
models are \emph{symmetric} multicategories. Yet, remarkably, it is possible to
directly construct a cartesian double theory whose models are symmetric
multicategories, as established by Pisani in his previous study of ``operads as
double functors'' \mbox{\cite[\S{3}]{pisani2022}}. At the end of this section,
we generalize the notion of cartesian double theory in order to present a theory
of symmetric promonoids that mirrors the standard axioms for symmetric
multicategories.

The following two double theories are inspired by concepts from Day and Street
\cite{daystreet2004,street2004}. They are the only examples in this section
where we are primarily interested in pseudo, rather than lax, models.

\begin{theory}[Biexact pairings] \label{th:biexactpairing}
  The \define{theory of biexact pairings} is generated by two objects $x$ and
  $y$ and two proarrows $\eta\colon 1\proto y\times x$ and
  $\epsilon\colon x\times y \proto 1$ satisfying the two identities
  \begin{equation*}
    \begin{tikzcd}
      {x\times 1} & {x\times y\times x} & {1\times x} \\
      x && x
      \arrow["{\id_x\times\eta}", "\shortmid"{marking}, from=1-1, to=1-2]
      \arrow["{\epsilon\times \id_x}", "\shortmid"{marking}, from=1-2, to=1-3]
      \arrow["\cong"', "\shortmid"{marking}, from=1-1, to=2-1]
      \arrow["\cong", "\shortmid"{marking}, from=1-3, to=2-3]
      \arrow["{\id_x}"', "\shortmid"{marking}, from=2-1, to=2-3]
    \end{tikzcd}
    \qquad\text{and}\qquad
    \begin{tikzcd}
      {1\times y} & {y\times x\times y} & {y\times 1} \\
      y && y
      \arrow["{\eta\times\id_y}", "\shortmid"{marking}, from=1-1, to=1-2]
      \arrow["{\id_y\times \epsilon}", "\shortmid"{marking}, from=1-2, to=1-3]
      \arrow["\cong"', "\shortmid"{marking}, from=1-1, to=2-1]
      \arrow["\cong", "\shortmid"{marking}, from=1-3, to=2-3]
      \arrow["{\id_y}"', "\shortmid"{marking}, from=2-1, to=2-3]
    \end{tikzcd}.
  \end{equation*}
  In any model, the image of $x$ is the \define{left bidual} and the image of
  $y$ is the \define{right bidual}.

  A pseudo model in a cartesian double category is a biexact pairing in its
  horizontal monoidal bicategory, as defined by Street
  \cite[\S{III}]{street2004}. It turns out that a biexact pairing in $\Prof$
  consists of a category and its opposite.
\end{theory}

\begin{theory}[Frobenius pseudomonoids]
  The \define{theory of a form} on a promonoid is the theory of promonoids
  (\cref{th:promonoid}) augmented with a proarrow $s\colon x^2 \proto 1$
  satisfying the equation
  \begin{equation*}
    \begin{tikzcd}
      {x^3} & {x^2} \\
      {x^2} & 1
      \arrow["{p \times \id_x}"', "\shortmid"{marking}, from=1-1, to=2-1]
      \arrow["s"', "\shortmid"{marking}, from=2-1, to=2-2]
      \arrow["{\id_x \times p}", "\shortmid"{marking}, from=1-1, to=1-2]
      \arrow["s", "\shortmid"{marking}, from=1-2, to=2-2]
    \end{tikzcd}.
  \end{equation*}
  A pseudo model in a cartesian double category is a pseudomonoid in the
  horizontal monoidal bicategory along with a \emph{form} for the pseudomonoid
  \cite[\S{9}]{daystreet2004}.

  The \define{theory of *-autonomous promonoids} requires additionally that
  $s: x^2 \proto 1$ be the counit of a biexact pairing as in
  \cref{th:biexactpairing}. A pseudo model in a cartesian double category is a
  \emph{*-autonomous pseudomonoid} in the horizontal monoidal bicategory
  \cite[\S{9}]{daystreet2004}, \cite[\S{IV}]{street2004}. A different but
  equivalent axiomatization of a *-autonomous pseudomonoid has been called a
  \emph{Frobenius pseudomonoid} \cite[Proposition 3.2]{street2004}. Intuitively,
  a Frobenius pseudomonoid categorifies a Frobenius algebra and so should
  involve a pseudomonoid and a pseudocomonoid interacting with each other. In
  these axiomatizations, the pseudocomonoid is derived rather than primitive
  structure \cite[Proposition 3.1]{street2004}.
\end{theory}

We now consider a series of double theories whose models are monoidal categories
equipped with extra structure.

\begin{theory}[Monoidal copresheaves] \label{th:pseudomonoid-action}
  The \define{theory of pseudomonoid actions} is the theory of pseudomonoids
  augmented with a generating proarrow $p: 1 \proto x$ and generating cells
  \begin{equation*}
    \begin{tikzcd}
      1 & {x^2} \\
      1 & x
      \arrow[""{name=0, anchor=center, inner sep=0}, "{p^2}", "\shortmid"{marking}, from=1-1, to=1-2]
      \arrow[""{name=1, anchor=center, inner sep=0}, "p"', "\shortmid"{marking}, from=2-1, to=2-2]
      \arrow["\otimes", from=1-2, to=2-2]
      \arrow[Rightarrow, no head, from=1-1, to=2-1]
      \arrow["\mu"{description}, draw=none, from=0, to=1]
    \end{tikzcd}
    \qquad\text{and}\qquad
    \begin{tikzcd}
      1 & 1 \\
      1 & x
      \arrow[""{name=0, anchor=center, inner sep=0}, "{\mathrm{id}_1}", "\shortmid"{marking}, from=1-1, to=1-2]
      \arrow["I", from=1-2, to=2-2]
      \arrow[""{name=1, anchor=center, inner sep=0}, "p"', "\shortmid"{marking}, from=2-1, to=2-2]
      \arrow[Rightarrow, no head, from=1-1, to=2-1]
      \arrow["\eta"{description}, draw=none, from=0, to=1]
    \end{tikzcd}
  \end{equation*}
  subject to the associativity axiom
  \begin{equation*}
    \begin{tikzcd}
      1 & {x^3} & {x^3} \\
      1 & {x^2} & {x^2} \\
      1 & x & x
      \arrow[""{name=0, anchor=center, inner sep=0}, "{\mathrm{id}_x^3}", "\shortmid"{marking}, from=1-2, to=1-3]
      \arrow["{\otimes 1_x}", from=1-2, to=2-2]
      \arrow["\otimes", from=2-2, to=3-2]
      \arrow["{1_x \otimes}", from=1-3, to=2-3]
      \arrow[""{name=1, anchor=center, inner sep=0}, "{\mathrm{id}_x}"', "\shortmid"{marking}, from=3-2, to=3-3]
      \arrow["\otimes", from=2-3, to=3-3]
      \arrow[Rightarrow, no head, from=1-1, to=2-1]
      \arrow[""{name=2, anchor=center, inner sep=0}, "{p^3}", "\shortmid"{marking}, from=1-1, to=1-2]
      \arrow[""{name=3, anchor=center, inner sep=0}, "{p^2}", "\shortmid"{marking}, from=2-1, to=2-2]
      \arrow[Rightarrow, no head, from=2-1, to=3-1]
      \arrow[""{name=4, anchor=center, inner sep=0}, "p"', "\shortmid"{marking}, from=3-1, to=3-2]
      \arrow["\alpha"{description}, draw=none, from=0, to=1]
      \arrow["{\mu\; 1_p}"{description}, draw=none, from=2, to=3]
      \arrow["\mu"{description}, draw=none, from=3, to=4]
    \end{tikzcd}
    \quad=\quad
    \begin{tikzcd}
      1 & {x^3} \\
      1 & {x^2} \\
      1 & x
      \arrow["{1_x \otimes}", from=1-2, to=2-2]
      \arrow["\otimes", from=2-2, to=3-2]
      \arrow[Rightarrow, no head, from=1-1, to=2-1]
      \arrow[""{name=0, anchor=center, inner sep=0}, "{p^3}", "\shortmid"{marking}, from=1-1, to=1-2]
      \arrow[""{name=1, anchor=center, inner sep=0}, "{p^2}", "\shortmid"{marking}, from=2-1, to=2-2]
      \arrow[Rightarrow, no head, from=2-1, to=3-1]
      \arrow[""{name=2, anchor=center, inner sep=0}, "p"', "\shortmid"{marking}, from=3-1, to=3-2]
      \arrow["{1_p\;\mu}"{description}, draw=none, from=0, to=1]
      \arrow["\mu"{description}, draw=none, from=1, to=2]
    \end{tikzcd}
  \end{equation*}
  and the unitality axioms
  \begin{equation*}
    \begin{tikzcd}
      1 & x & x \\
      1 & {x^2} \\
      1 & x & x
      \arrow["{I\; 1_x}", from=1-2, to=2-2]
      \arrow["\otimes", from=2-2, to=3-2]
      \arrow[""{name=0, anchor=center, inner sep=0}, "{\mathrm{id}_x}"', "\shortmid"{marking}, from=3-2, to=3-3]
      \arrow[""{name=1, anchor=center, inner sep=0}, "{\mathrm{id}_x}", "\shortmid"{marking}, from=1-2, to=1-3]
      \arrow[Rightarrow, no head, from=1-3, to=3-3]
      \arrow[""{name=2, anchor=center, inner sep=0}, "p", "\shortmid"{marking}, from=1-1, to=1-2]
      \arrow[""{name=3, anchor=center, inner sep=0}, "{p^2}", "\shortmid"{marking}, from=2-1, to=2-2]
      \arrow[""{name=4, anchor=center, inner sep=0}, "p"', "\shortmid"{marking}, from=3-1, to=3-2]
      \arrow[Rightarrow, no head, from=2-1, to=3-1]
      \arrow[Rightarrow, no head, from=1-1, to=2-1]
      \arrow["\lambda"{description}, draw=none, from=1, to=0]
      \arrow["\mu"{description}, draw=none, from=3, to=4]
      \arrow["{\eta\; 1_p}"{description}, draw=none, from=2, to=3]
    \end{tikzcd}
    \quad=\quad
    1_p
    \quad=\quad
    \begin{tikzcd}
      1 & x & x \\
      1 & {x^2} \\
      1 & x & x
      \arrow["{1_x\; I}", from=1-2, to=2-2]
      \arrow["\otimes", from=2-2, to=3-2]
      \arrow[""{name=0, anchor=center, inner sep=0}, "{\mathrm{id}_x}"', "\shortmid"{marking}, from=3-2, to=3-3]
      \arrow[""{name=1, anchor=center, inner sep=0}, "{\mathrm{id}_x}", "\shortmid"{marking}, from=1-2, to=1-3]
      \arrow[Rightarrow, no head, from=1-3, to=3-3]
      \arrow[""{name=2, anchor=center, inner sep=0}, "p", "\shortmid"{marking}, from=1-1, to=1-2]
      \arrow[""{name=3, anchor=center, inner sep=0}, "{p^2}", "\shortmid"{marking}, from=2-1, to=2-2]
      \arrow[""{name=4, anchor=center, inner sep=0}, "p"', "\shortmid"{marking}, from=3-1, to=3-2]
      \arrow[Rightarrow, no head, from=2-1, to=3-1]
      \arrow[Rightarrow, no head, from=1-1, to=2-1]
      \arrow["\rho"{description}, draw=none, from=1, to=0]
      \arrow["\mu"{description}, draw=none, from=3, to=4]
      \arrow["{1_p\; \eta}"{description}, draw=none, from=2, to=3]
    \end{tikzcd}.
  \end{equation*}
  A model of this theory is a monoidal category $(\cat{C}, \otimes, I)$ together
  with a \define{monoidal copresheaf} on $(\cat{C}, \otimes, I)$, i.e., a lax
  monoidal functor
  $(P, \mu, \eta): (\cat{C}, \otimes, I) \to (\Set, \times, 1)$. The laxators
  and unitors are natural transformations of form
  \begin{equation*}
    \mu_{c,c'}: P(c) \times P(c') \to P(c \otimes c'), \quad c,c' \in \cat{C}
    \qquad\text{and}\qquad
    \eta: 1 \to P(I). \qedhere
  \end{equation*}
\end{theory}

\begin{theory}[Symmetric monoidal copresheaves]
  \label{th:symmetric-pseudomonoid-action}
  The \define{theory of symmetric pseudomonoid actions} is the theory of
  symmetric pseudomonoids augmented as in the theory of pseudomonoid actions,
  along with the further equation
  \begin{equation*}
    \begin{tikzcd}
      1 & {x^2} \\
      1 & {x^2} \\
      1 & x
      \arrow[""{name=0, anchor=center, inner sep=0}, "{p^2}", "\shortmid"{marking}, from=2-1, to=2-2]
      \arrow[""{name=1, anchor=center, inner sep=0}, "p"', "\shortmid"{marking}, from=3-1, to=3-2]
      \arrow["\otimes", from=2-2, to=3-2]
      \arrow[Rightarrow, no head, from=2-1, to=3-1]
      \arrow["{\sigma_{x,x}}", from=1-2, to=2-2]
      \arrow[""{name=2, anchor=center, inner sep=0}, "{p^2}", "\shortmid"{marking}, from=1-1, to=1-2]
      \arrow[Rightarrow, no head, from=1-1, to=2-1]
      \arrow["\mu"{description}, draw=none, from=0, to=1]
      \arrow["{\sigma_{p,p}}"{description, pos=0.4}, draw=none, from=2, to=0]
    \end{tikzcd}
    \quad=\quad
    \begin{tikzcd}
      1 & {x^2} & {x^2} \\
      && {x^2} \\
      1 & x & x
      \arrow[""{name=0, anchor=center, inner sep=0}, "{p^2}", "\shortmid"{marking}, from=1-1, to=1-2]
      \arrow[""{name=1, anchor=center, inner sep=0}, "p"', "\shortmid"{marking}, from=3-1, to=3-2]
      \arrow["\otimes", from=1-2, to=3-2]
      \arrow[Rightarrow, no head, from=1-1, to=3-1]
      \arrow["\otimes", from=2-3, to=3-3]
      \arrow["{\sigma_{x,x}}", from=1-3, to=2-3]
      \arrow[""{name=2, anchor=center, inner sep=0}, "{\id_x^2}", "\shortmid"{marking}, from=1-2, to=1-3]
      \arrow[""{name=3, anchor=center, inner sep=0}, "{\id_x}"', "\shortmid"{marking}, from=3-2, to=3-3]
      \arrow["\mu"{description}, draw=none, from=0, to=1]
      \arrow["\sigma"{description}, draw=none, from=2, to=3]
    \end{tikzcd}.
  \end{equation*}
  A model of the theory is a symmetric monoidal category $(\cat{C},\otimes,I)$
  together with a \define{symmetric monoidal copresheaf} on it, i.e., a lax
  symmetric monoidal functor
  $(P,\mu,\eta): (\cat{C},\otimes,I) \to (\Set,\times,1)$. The laxators satisfy
  the symmetry axiom
  \begin{equation*}
    \begin{tikzcd}
      {P(c) \times P(c')} & {P(c') \times P(c)} \\
      {P(c \otimes c')} & {P(c' \otimes c)}
      \arrow["{\sigma_{Pc,Pc'}}", from=1-1, to=1-2]
      \arrow["{P(\sigma_{c,c'})}"', from=2-1, to=2-2]
      \arrow["{\mu_{c,c'}}"', from=1-1, to=2-1]
      \arrow["{\mu_{c',c}}", from=1-2, to=2-2]
    \end{tikzcd},
    \qquad c,c' \in \cat{C}. \qedhere
  \end{equation*}
\end{theory}

\begin{theory}[Cartesian monoidal categories] \label{th:cart-mon-cat-I}
  A first pass on cartesian categories unites the theory of adjunctions
  (\cref{th:adjunctions}) with the hypothesis that the double theory itself is
  cartesian. The \define{theory of cartesian objects} is generated by a single
  object $x$, two arrows $\otimes\colon x^2\to x$ and $I: 1 \to x$, and four
  cells
  \begin{equation*}
    \begin{tikzcd}[row sep=scriptsize]
      x & x \\
      & {x^2} \\
      x & x
      \arrow[""{name=0, anchor=center, inner sep=0}, "{\id_x}", "\shortmid"{marking}, from=1-1, to=1-2]
      \arrow["{\Delta_x}", from=1-2, to=2-2]
      \arrow["\otimes", from=2-2, to=3-2]
      \arrow[Rightarrow, no head, from=1-1, to=3-1]
      \arrow[""{name=1, anchor=center, inner sep=0}, "{\id_x}"', "\shortmid"{marking}, from=3-1, to=3-2]
      \arrow["\delta"{description}, draw=none, from=0, to=1]
    \end{tikzcd}
    \qquad\qquad
    \begin{tikzcd}[row sep=scriptsize]
      {x^2} & {x^2} \\
      x \\
      {x^2} & {x^2}
      \arrow["\otimes"', from=1-1, to=2-1]
      \arrow["{\Delta_x}"', from=2-1, to=3-1]
      \arrow[""{name=0, anchor=center, inner sep=0}, "{\id_{x^2}}"', "\shortmid"{marking}, from=3-1, to=3-2]
      \arrow[""{name=1, anchor=center, inner sep=0}, "{\id_{x^2}}", "\shortmid"{marking}, from=1-1, to=1-2]
      \arrow[Rightarrow, no head, from=1-2, to=3-2]
      \arrow["\pi"{description}, draw=none, from=1, to=0]
    \end{tikzcd}
    \qquad\qquad
    \begin{tikzcd}[row sep=scriptsize]
      x & x \\
      & 1 \\
      x & x
      \arrow[Rightarrow, no head, from=1-1, to=3-1]
      \arrow[""{name=0, anchor=center, inner sep=0}, "{\id_x}"', "\shortmid"{marking}, from=3-1, to=3-2]
      \arrow[""{name=1, anchor=center, inner sep=0}, "{\id_x}", "\shortmid"{marking}, from=1-1, to=1-2]
      \arrow["{!}", from=1-2, to=2-2]
      \arrow["I", from=2-2, to=3-2]
      \arrow["\epsilon"{description}, draw=none, from=1, to=0]
    \end{tikzcd}
    \qquad\qquad
    \begin{tikzcd}[row sep=scriptsize]
      1 & 1 \\
      x \\
      1 & 1
      \arrow["I"', from=1-1, to=2-1]
      \arrow["{!}"', from=2-1, to=3-1]
      \arrow[""{name=0, anchor=center, inner sep=0}, "{\id_1}"', "\shortmid"{marking}, from=3-1, to=3-2]
      \arrow[""{name=1, anchor=center, inner sep=0}, "{\id_1}", "\shortmid"{marking}, from=1-1, to=1-2]
      \arrow[Rightarrow, no head, from=1-2, to=3-2]
      \arrow["\eta"{description}, draw=none, from=1, to=0]
    \end{tikzcd}
  \end{equation*}
  meant to represent the internal diagonal and projections. The cells should be
  the units and counits of two adjunctions, so we also require the equations
  \begin{equation*}
    \begin{tikzcd}[row sep=scriptsize]
      x & x & x \\
      & {x^2} & {x^2} \\
      x & x \\
      {x^2} & {x^2} & {x^2}
      \arrow[""{name=0, anchor=center, inner sep=0}, "{\id_x}", "\shortmid"{marking}, from=1-1, to=1-2]
      \arrow["\Delta"', from=1-2, to=2-2]
      \arrow[Rightarrow, no head, from=1-1, to=3-1]
      \arrow[""{name=1, anchor=center, inner sep=0}, "{\id_x}", "\shortmid"{marking}, from=3-1, to=3-2]
      \arrow[""{name=2, anchor=center, inner sep=0}, "{\id_x^2}"', "\shortmid"{marking}, from=4-2, to=4-3]
      \arrow[Rightarrow, no head, from=2-3, to=4-3]
      \arrow["\otimes", from=2-2, to=3-2]
      \arrow[""{name=3, anchor=center, inner sep=0}, "{\id_x^2}"', "\shortmid"{marking}, from=2-2, to=2-3]
      \arrow["\Delta", from=3-2, to=4-2]
      \arrow[""{name=4, anchor=center, inner sep=0}, "{\id_x}", "\shortmid"{marking}, from=1-2, to=1-3]
      \arrow["\Delta", from=1-3, to=2-3]
      \arrow["\Delta"', from=3-1, to=4-1]
      \arrow[""{name=5, anchor=center, inner sep=0}, "{\id_x^2}"', "\shortmid"{marking}, from=4-1, to=4-2]
      \arrow["\delta"{description}, draw=none, from=0, to=1]
      \arrow["{\id_\Delta}"{description}, draw=none, from=4, to=3]
      \arrow["{\id_\Delta}"{description}, draw=none, from=1, to=5]
      \arrow["\pi"{description}, draw=none, from=3, to=2]
    \end{tikzcd}
    \quad=\quad
    \begin{tikzcd}
      x & x \\
      {x^2} & {x^2}
      \arrow["\Delta"', from=1-1, to=2-1]
      \arrow[""{name=0, anchor=center, inner sep=0}, "{\id_x}", "\shortmid"{marking}, from=1-1, to=1-2]
      \arrow["\Delta", from=1-2, to=2-2]
      \arrow[""{name=1, anchor=center, inner sep=0}, "{\id_x^2}"', "\shortmid"{marking}, from=2-1, to=2-2]
      \arrow["{\id_\Delta}"{description}, draw=none, from=0, to=1]
    \end{tikzcd}
  \end{equation*}
  and
  \begin{equation*}
    \begin{tikzcd}[row sep=scriptsize]
      {x^2} & {x^2} & {x^2} \\
      x & x \\
      & {x^2} & {x^2} \\
      x & x & x
      \arrow[""{name=0, anchor=center, inner sep=0}, "{\id_x^2}", "\shortmid"{marking}, from=1-1, to=1-2]
      \arrow[""{name=1, anchor=center, inner sep=0}, "{\id_x^2}", "\shortmid"{marking}, from=1-2, to=1-3]
      \arrow[Rightarrow, no head, from=1-3, to=3-3]
      \arrow["\otimes", from=1-2, to=2-2]
      \arrow["\Delta", from=2-2, to=3-2]
      \arrow[""{name=2, anchor=center, inner sep=0}, "\shortmid"{marking}, from=3-2, to=3-3]
      \arrow[""{name=3, anchor=center, inner sep=0}, "\shortmid"{marking}, from=2-1, to=2-2]
      \arrow[Rightarrow, no head, from=2-1, to=4-1]
      \arrow[""{name=4, anchor=center, inner sep=0}, "{\id_x}"', "\shortmid"{marking}, from=4-1, to=4-2]
      \arrow["\otimes"', from=3-2, to=4-2]
      \arrow[""{name=5, anchor=center, inner sep=0}, "{\id_x}"', "\shortmid"{marking}, from=4-2, to=4-3]
      \arrow["\otimes", from=3-3, to=4-3]
      \arrow["\otimes"', from=1-1, to=2-1]
      \arrow["{\id_\otimes}"{description}, draw=none, from=2, to=5]
      \arrow["\delta"{description}, draw=none, from=3, to=4]
      \arrow["{\id_\otimes}"{description}, draw=none, from=0, to=3]
      \arrow["\pi"{description}, draw=none, from=1, to=2]
    \end{tikzcd}
    \quad=\quad
    \begin{tikzcd}
      {x^2} & {x^2} \\
      x & x
      \arrow["\otimes"', from=1-1, to=2-1]
      \arrow[""{name=0, anchor=center, inner sep=0}, "{\id_x}"', "\shortmid"{marking}, from=2-1, to=2-2]
      \arrow[""{name=1, anchor=center, inner sep=0}, "{\id_x^2}", "\shortmid"{marking}, from=1-1, to=1-2]
      \arrow["\otimes", from=1-2, to=2-2]
      \arrow["{\id_\otimes}"{description}, draw=none, from=1, to=0]
    \end{tikzcd}
  \end{equation*}
  and similarly for $\epsilon$ and $\eta$.

  By \cref{cor:unitalizationofCartesianLaxFunctors}, a model of the theory of
  cartesian objects in $\Span$ is the same as a unitary lax functor into $\Prof$.
  This is a cartesian monoidal category since the theory has trivial proarrow
  structure, hence is equivalent to a 2-category, and thus any such model is a
  2-functor, which must preserve the adjunction. More generally, for any
  cartesian double category $\dbl{S}$, a model in $\dbl{S}$ is a \emph{cartesian
    object} in the cartesian 2-category underlying $\dbl{S}$.
\end{theory}

\begin{theory}[Cartesian monoidal categories, again]
  \label{th:internal-commutative-comonoid}
  The \define{theory of internal commutative comonoids} is the theory of
  symmetric pseudomonoids augmented with generators
  \begin{equation*}
    \begin{tikzcd}[row sep=scriptsize]
      x & x \\
      & {x^2} \\
      x & x
      \arrow[""{name=0, anchor=center, inner sep=0}, "{\mathrm{id}_x}", "\shortmid"{marking}, from=1-1, to=1-2]
      \arrow["{\Delta_x}", from=1-2, to=2-2]
      \arrow["\otimes", from=2-2, to=3-2]
      \arrow[Rightarrow, no head, from=1-1, to=3-1]
      \arrow[""{name=1, anchor=center, inner sep=0}, "{\mathrm{id}_x}"', "\shortmid"{marking}, from=3-1, to=3-2]
      \arrow["\delta"{description}, draw=none, from=0, to=1]
    \end{tikzcd}
    \qquad\text{and}\qquad
    \begin{tikzcd}[row sep=scriptsize]
      x & x \\
      & 1 \\
      x & x
      \arrow["{!}", from=1-2, to=2-2]
      \arrow["I", from=2-2, to=3-2]
      \arrow[""{name=0, anchor=center, inner sep=0}, "{\mathrm{id}_x}", "\shortmid"{marking}, from=1-1, to=1-2]
      \arrow[Rightarrow, no head, from=1-1, to=3-1]
      \arrow[""{name=1, anchor=center, inner sep=0}, "{\mathrm{id}_x}"', "\shortmid"{marking}, from=3-1, to=3-2]
      \arrow["\varepsilon"{description}, draw=none, from=0, to=1]
    \end{tikzcd},
  \end{equation*}
  the \define{comultiplication} and \define{counit} cells, subject to laws of
  associativity, unitality, and commutativity, as well as four coherence axioms.
  For example, the associativity law is
  \begin{equation*}
    \begin{tikzcd}[row sep=small]
      x & x & x & x \\
      & {x^2} & {x^2} & {x^2} \\
      && {x^3} & {x^3} \\
      & {x^2} & {x^2} & {x^2} \\
      x & x & x & x
      \arrow[""{name=0, anchor=center, inner sep=0}, "\shortmid"{marking}, Rightarrow, no head, from=2-2, to=2-3]
      \arrow["{\Delta_x 1_x}", from=2-3, to=3-3]
      \arrow["{\otimes 1_x}", from=3-3, to=4-3]
      \arrow[Rightarrow, no head, from=2-2, to=4-2]
      \arrow[""{name=1, anchor=center, inner sep=0}, "\shortmid"{marking}, Rightarrow, no head, from=4-2, to=4-3]
      \arrow["\otimes"', from=4-2, to=5-2]
      \arrow["\otimes", from=4-3, to=5-3]
      \arrow[""{name=2, anchor=center, inner sep=0}, "\shortmid"{marking}, Rightarrow, no head, from=5-2, to=5-3]
      \arrow["{\Delta_x}"', from=1-2, to=2-2]
      \arrow["{\Delta_x}", from=1-3, to=2-3]
      \arrow[""{name=3, anchor=center, inner sep=0}, "\shortmid"{marking}, Rightarrow, no head, from=1-2, to=1-3]
      \arrow[""{name=4, anchor=center, inner sep=0}, "\shortmid"{marking}, Rightarrow, no head, from=1-1, to=1-2]
      \arrow[""{name=5, anchor=center, inner sep=0}, "\shortmid"{marking}, Rightarrow, no head, from=5-1, to=5-2]
      \arrow[Rightarrow, no head, from=1-1, to=5-1]
      \arrow[""{name=6, anchor=center, inner sep=0}, "\shortmid"{marking}, Rightarrow, no head, from=3-3, to=3-4]
      \arrow["{\Delta_x}", from=1-4, to=2-4]
      \arrow["{1_x \Delta_x}", from=2-4, to=3-4]
      \arrow[""{name=7, anchor=center, inner sep=0}, "\shortmid"{marking}, Rightarrow, no head, from=1-3, to=1-4]
      \arrow["{1_x \otimes}", from=3-4, to=4-4]
      \arrow["\otimes", from=4-4, to=5-4]
      \arrow[""{name=8, anchor=center, inner sep=0}, "\shortmid"{marking}, Rightarrow, no head, from=5-3, to=5-4]
      \arrow["{\delta\; \mathrm{id}_{1_x}}"{description}, draw=none, from=0, to=1]
      \arrow["{\mathrm{id}_\otimes}"{description}, draw=none, from=1, to=2]
      \arrow["{\mathrm{id}_{\Delta_x}}"{description}, draw=none, from=3, to=0]
      \arrow["\delta"{description}, draw=none, from=4, to=5]
      \arrow["{\mathrm{id}}"{description}, draw=none, from=7, to=6]
      \arrow["\alpha"{description}, draw=none, from=6, to=8]
    \end{tikzcd}
    \quad=\quad
    \begin{tikzcd}[row sep=small]
      x & x & x \\
      & {x^2} & {x^2} \\
      && {x^3} \\
      & {x^2} & {x^2} \\
      x & x & x
      \arrow[""{name=0, anchor=center, inner sep=0}, "\shortmid"{marking}, Rightarrow, no head, from=2-2, to=2-3]
      \arrow["{1_x \Delta_x}", from=2-3, to=3-3]
      \arrow["{1_x \otimes}", from=3-3, to=4-3]
      \arrow[Rightarrow, no head, from=2-2, to=4-2]
      \arrow[""{name=1, anchor=center, inner sep=0}, "\shortmid"{marking}, Rightarrow, no head, from=4-2, to=4-3]
      \arrow["\otimes"', from=4-2, to=5-2]
      \arrow["\otimes", from=4-3, to=5-3]
      \arrow[""{name=2, anchor=center, inner sep=0}, "\shortmid"{marking}, Rightarrow, no head, from=5-2, to=5-3]
      \arrow["{\Delta_x}"', from=1-2, to=2-2]
      \arrow["{\Delta_x}", from=1-3, to=2-3]
      \arrow[""{name=3, anchor=center, inner sep=0}, "\shortmid"{marking}, Rightarrow, no head, from=1-2, to=1-3]
      \arrow[""{name=4, anchor=center, inner sep=0}, "\shortmid"{marking}, Rightarrow, no head, from=1-1, to=1-2]
      \arrow[""{name=5, anchor=center, inner sep=0}, "\shortmid"{marking}, Rightarrow, no head, from=5-1, to=5-2]
      \arrow[Rightarrow, no head, from=1-1, to=5-1]
      \arrow["{\mathrm{id}_{1_x}\; \delta}"{description}, draw=none, from=0, to=1]
      \arrow["{\mathrm{id}_\otimes}"{description}, draw=none, from=1, to=2]
      \arrow["{\mathrm{id}_{\Delta_x}}"{description}, draw=none, from=3, to=0]
      \arrow["\delta"{description}, draw=none, from=4, to=5]
    \end{tikzcd}
  \end{equation*}
  and commutativity law is
  \begin{equation*}
    \begin{tikzcd}[row sep=small]
      x & x & x \\
      & {x^2} & {x^2} \\
      && {x^2} \\
      x & x & x
      \arrow[""{name=0, anchor=center, inner sep=0}, "\shortmid"{marking}, Rightarrow, no head, from=1-1, to=1-2]
      \arrow["{\Delta_x}"', from=1-2, to=2-2]
      \arrow["\otimes"', from=2-2, to=4-2]
      \arrow[Rightarrow, no head, from=1-1, to=4-1]
      \arrow[""{name=1, anchor=center, inner sep=0}, "\shortmid"{marking}, Rightarrow, no head, from=4-1, to=4-2]
      \arrow["\otimes", from=3-3, to=4-3]
      \arrow["{\sigma_{x,x}}", from=2-3, to=3-3]
      \arrow[""{name=2, anchor=center, inner sep=0}, "\shortmid"{marking}, Rightarrow, no head, from=2-2, to=2-3]
      \arrow[""{name=3, anchor=center, inner sep=0}, "\shortmid"{marking}, Rightarrow, no head, from=4-2, to=4-3]
      \arrow[""{name=4, anchor=center, inner sep=0}, "\shortmid"{marking}, Rightarrow, no head, from=1-2, to=1-3]
      \arrow["{\Delta_x}", from=1-3, to=2-3]
      \arrow["\delta"{description}, draw=none, from=0, to=1]
      \arrow["\sigma"{description}, draw=none, from=2, to=3]
      \arrow["{\mathrm{id}_{\Delta_x}}"{description}, draw=none, from=4, to=2]
    \end{tikzcd}
    \quad=\quad
    \begin{tikzcd}[row sep=small]
      x & x \\
      & {x^2} \\
      \\
      x & x
      \arrow[""{name=0, anchor=center, inner sep=0}, "\shortmid"{marking}, Rightarrow, no head, from=1-1, to=1-2]
      \arrow["{\Delta_x}", from=1-2, to=2-2]
      \arrow["\otimes", from=2-2, to=4-2]
      \arrow[Rightarrow, no head, from=1-1, to=4-1]
      \arrow[""{name=1, anchor=center, inner sep=0}, "\shortmid"{marking}, Rightarrow, no head, from=4-1, to=4-2]
      \arrow["\delta"{description}, draw=none, from=0, to=1]
    \end{tikzcd}.
  \end{equation*}
  The coherence axioms assert that the comultiplication $\delta$ and counit
  $\varepsilon$ commute (in the strict direction) with the monoidal product and
  unit. For example, coherence between the counit and the monoidal product is
  the equation
  \begin{equation*}
    \begin{tikzcd}[row sep=small]
      {x^2} & {x^2} \\
      x & x \\
      & {\mathsf{1}} \\
      x & x
      \arrow[""{name=0, anchor=center, inner sep=0}, "\shortmid"{marking}, Rightarrow, no head, from=1-1, to=1-2]
      \arrow["\otimes"', from=1-1, to=2-1]
      \arrow["\otimes", from=1-2, to=2-2]
      \arrow[""{name=1, anchor=center, inner sep=0}, "\shortmid"{marking}, Rightarrow, no head, from=2-1, to=2-2]
      \arrow["{!}", from=2-2, to=3-2]
      \arrow["I", from=3-2, to=4-2]
      \arrow[Rightarrow, no head, from=2-1, to=4-1]
      \arrow[""{name=2, anchor=center, inner sep=0}, "\shortmid"{marking}, Rightarrow, no head, from=4-1, to=4-2]
      \arrow["\varepsilon"{description}, draw=none, from=1, to=2]
      \arrow["{\mathrm{id}_\otimes}"{description}, draw=none, from=0, to=1]
    \end{tikzcd}
    \quad=\quad
    \begin{tikzcd}[row sep=small]
      {x^2} & {x^2} & {x^2} \\
      & x & x \\
      {x^2} & {x^2} \\
      x & x & x
      \arrow["{! \cdot I}"', from=1-2, to=2-2]
      \arrow["{I \times 1_x}"', from=2-2, to=3-2]
      \arrow[""{name=0, anchor=center, inner sep=0}, "\shortmid"{marking}, Rightarrow, no head, from=1-1, to=1-2]
      \arrow[Rightarrow, no head, from=1-1, to=3-1]
      \arrow[""{name=1, anchor=center, inner sep=0}, "\shortmid"{marking}, Rightarrow, no head, from=3-1, to=3-2]
      \arrow["\otimes"', from=3-1, to=4-1]
      \arrow["\otimes", from=3-2, to=4-2]
      \arrow[""{name=2, anchor=center, inner sep=0}, "\shortmid"{marking}, Rightarrow, no head, from=4-1, to=4-2]
      \arrow["{! \cdot I}", from=1-3, to=2-3]
      \arrow[Rightarrow, no head, from=2-3, to=4-3]
      \arrow[""{name=3, anchor=center, inner sep=0}, "\shortmid"{marking}, Rightarrow, no head, from=4-2, to=4-3]
      \arrow[""{name=4, anchor=center, inner sep=0}, "\shortmid"{marking}, Rightarrow, no head, from=2-2, to=2-3]
      \arrow[""{name=5, anchor=center, inner sep=0}, "\shortmid"{marking}, Rightarrow, no head, from=1-2, to=1-3]
      \arrow["{\varepsilon^2}"{description}, draw=none, from=0, to=1]
      \arrow["{\mathrm{id}_\otimes}"{description}, draw=none, from=1, to=2]
      \arrow["\lambda"{description}, draw=none, from=4, to=3]
      \arrow["{\mathrm{id}}"{description}, draw=none, from=5, to=4]
    \end{tikzcd}
    \quad=\quad
    \begin{tikzcd}[row sep=small]
      {x^2} & {x^2} & {x^2} \\
      & x & x \\
      {x^2} & {x^2} \\
      x & x & x
      \arrow["{! \cdot I}"', from=1-2, to=2-2]
      \arrow["{1_x \times I}"', from=2-2, to=3-2]
      \arrow[""{name=0, anchor=center, inner sep=0}, "\shortmid"{marking}, Rightarrow, no head, from=1-1, to=1-2]
      \arrow[Rightarrow, no head, from=1-1, to=3-1]
      \arrow[""{name=1, anchor=center, inner sep=0}, "\shortmid"{marking}, Rightarrow, no head, from=3-1, to=3-2]
      \arrow["\otimes"', from=3-1, to=4-1]
      \arrow["\otimes", from=3-2, to=4-2]
      \arrow[""{name=2, anchor=center, inner sep=0}, "\shortmid"{marking}, Rightarrow, no head, from=4-1, to=4-2]
      \arrow["{! \cdot I}", from=1-3, to=2-3]
      \arrow[Rightarrow, no head, from=2-3, to=4-3]
      \arrow[""{name=3, anchor=center, inner sep=0}, "\shortmid"{marking}, Rightarrow, no head, from=4-2, to=4-3]
      \arrow[""{name=4, anchor=center, inner sep=0}, "\shortmid"{marking}, Rightarrow, no head, from=2-2, to=2-3]
      \arrow[""{name=5, anchor=center, inner sep=0}, "\shortmid"{marking}, Rightarrow, no head, from=1-2, to=1-3]
      \arrow["{\varepsilon^2}"{description}, draw=none, from=0, to=1]
      \arrow["{\mathrm{id}_\otimes}"{description}, draw=none, from=1, to=2]
      \arrow["\rho"{description}, draw=none, from=4, to=3]
      \arrow["{\mathrm{id}}"{description}, draw=none, from=5, to=4]
    \end{tikzcd}
  \end{equation*}
  and coherence between the comultiplication and the monoidal product is the
  equation
  \begin{equation*}
    \begin{tikzcd}[row sep=small]
      {x^2} & {x^2} \\
      x & x \\
      & {x^2} \\
      x & x
      \arrow[""{name=0, anchor=center, inner sep=0}, "\shortmid"{marking}, Rightarrow, no head, from=2-1, to=2-2]
      \arrow["{\Delta_x}", from=2-2, to=3-2]
      \arrow["\otimes", from=3-2, to=4-2]
      \arrow[Rightarrow, no head, from=2-1, to=4-1]
      \arrow[""{name=1, anchor=center, inner sep=0}, "\shortmid"{marking}, Rightarrow, no head, from=4-1, to=4-2]
      \arrow["\otimes"', from=1-1, to=2-1]
      \arrow["\otimes", from=1-2, to=2-2]
      \arrow[""{name=2, anchor=center, inner sep=0}, "\shortmid"{marking}, Rightarrow, no head, from=1-1, to=1-2]
      \arrow["\delta"{description}, draw=none, from=0, to=1]
      \arrow["{\mathrm{id}_\otimes}"{description}, draw=none, from=2, to=0]
    \end{tikzcd}
    \quad=\quad
    \begin{tikzcd}[row sep=small]
      {x^2} & {x^2} & {x^2} \\
      & {x^4} & x \\
      {x^2} & {x^2} & {x^2} \\
      x & x & x
      \arrow[""{name=0, anchor=center, inner sep=0}, "\shortmid"{marking}, Rightarrow, no head, from=1-1, to=1-2]
      \arrow["{\Delta_x^2}", from=1-2, to=2-2]
      \arrow["{\otimes^2}", from=2-2, to=3-2]
      \arrow[Rightarrow, no head, from=1-1, to=3-1]
      \arrow[""{name=1, anchor=center, inner sep=0}, "\shortmid"{marking}, Rightarrow, no head, from=3-1, to=3-2]
      \arrow["\otimes"', from=3-1, to=4-1]
      \arrow["\otimes", from=3-2, to=4-2]
      \arrow[""{name=2, anchor=center, inner sep=0}, "\shortmid"{marking}, Rightarrow, no head, from=4-1, to=4-2]
      \arrow[""{name=3, anchor=center, inner sep=0}, "\shortmid"{marking}, Rightarrow, no head, from=1-2, to=1-3]
      \arrow["\otimes", from=1-3, to=2-3]
      \arrow["{\Delta_x}", from=2-3, to=3-3]
      \arrow["\otimes", from=3-3, to=4-3]
      \arrow[""{name=4, anchor=center, inner sep=0}, "\shortmid"{marking}, Rightarrow, no head, from=4-2, to=4-3]
      \arrow["{\delta^2}"{description}, draw=none, from=0, to=1]
      \arrow["{\mathrm{id}_\otimes}"{description}, draw=none, from=1, to=2]
      \arrow["\cong"{description}, draw=none, from=3, to=4]
    \end{tikzcd},
  \end{equation*}
  where the elided cell on the right-hand side is formed from associators and a
  braiding.

  A model of the theory of internal commutative comonoids has been called, among
  other things, a symmetric monoidal category with a ``homomorphic supply of
  commutative comonoids'' \cite{fong2019}. By Fox's theorem \cite{fox1976}, this
  is equivalent to a cartesian monoidal category.
\end{theory}

Of course, both preceding theories can be dualized, giving theories of
\emph{cocartesian objects} and \emph{internal commutative monoids} whose models
are cocartesian monoidal categories. The next theory combines the theories
behind monoidal copresheaves and cartesian monoidal categories.

\begin{theory}[Algebraic theories and models]
  Consider the cartesian double theory that is the pushout of the theories of
  symmetric pseudomonoid actions (\cref{th:symmetric-pseudomonoid-action}) and
  internal commutative comonoids (\cref{th:internal-commutative-comonoid}) along
  their common copy of the theory of symmetric pseudomonoids
  (\cref{th:symmetric-pseudomonoid}). By construction, a model of this theory is
  a cartesian monoidal category $\cat{C}$ together with a lax symmetric monoidal
  functor $P\colon \cat{C}\to\Set$. But since $P$ is automatically oplax by the
  universal property of products, $P$ is actually strong, hence is a
  finite-product-preserving functor.

  Double categorical semantics thus reproduces \emph{algebraic type theories}
  along with all of their models \cite[Chapter 3]{crole1993}, \cite{adamek2010},
  which are known to include monoids, groups, and semilattices, among other
  algebraic structures. Here we invoke the dictionary of categorical logic: that
  every algebraic type theory has a classifying category with finite products;
  that every cartesian category has an internal algebraic type theory; and that
  these correspond in the sense that
  \begin{enumerate}[nosep]
    \item every cartesian category is the classifying category of its internal
    algebraic type theory \cite[Theorem 3.9.3]{crole1993}; and
    \item every algebraic type theory is the internal algebraic type theory of
    its classifying category \cite[Theorem 3.9.6]{crole1993}.
  \end{enumerate}
  Moreover, the classifying category of an algebraic type theory is universal in
  the sense that set-valued models of the theory are in one-to-one
  correspondence with finite-product-preserving set-valued functors on the
  classifying category \cite[Theorem 3.8.6]{crole1993}. In this way, we regard
  cartesian categories and cartesian set-valued functors as syntax-independent
  \emph{finite-product theories}, generalizing Lawvere theories, and see that
  they are subsumed by the present machinery of cartesian lax double functors.
\end{theory}

Formulating certain doctrines as double theories seems to require restricting
proarrows along their source or target, as in an equipment. It is tempting to
define such a double theory to be a small cartesian equipment, analogous to how
a cartesian double theory is just a small cartesian double category, but
allowing arbitrary restrictions causes of an explosion of proarrows in the
theory whose laxators in the models are difficult to control. Instead, we
introduce restriction cells via a sketch, a technique familiar from categorical
logic \cite[\S{D2}]{johnstone2002}. The following notion of restriction sketch
should be considered preliminary; our inclusion of it is motivated by several
compelling examples, but correctly using such sketches appears to be subtle. As
we discuss in the final section, \emph{virtual} equipments may offer an
alternative approach.

\begin{definition}[Restriction sketch]
  A \define{restriction sketch} is a cartesian double theory $\dbl{T}$ together
  with a distinguished subset of cells, denoted $\dbl{T}_\res$, closed under the
  operations in \cref{lem:closure-restrictions}. A \define{morphism} of
  restriction sketches from $(\dbl{T}, \dbl{T}_\res)$ to
  $(\dbl{T}', \dbl{T}'_\res)$ is a morphism $F: \dbl{T} \to \dbl{T}'$ of
  cartesian double theories that preserves the sketched restrictions in the
  sense that $F_1(\dbl{T}_\res) \subseteq \dbl{T}'_\res$.

  A \define{model} of a restriction sketch $(\dbl{T}, \dbl{T}_\res)$ in a
  cartesian equipment $\dbl{S}$ is a model $F: \dbl{T} \to \dbl{S}$ of the
  theory such that for every cell $\alpha$ belonging to $\dbl{T}_\res$, its
  image $F(\alpha)$ is a restriction cell in $\dbl{S}$.
\end{definition}

In particular, the sketched restriction cells $\dbl{T}_\res$ form a wide and
replete subcategory of $\dbl{T}_1$. Requiring that the sketched restrictions
respect the closure properties in \cref{lem:closure-restrictions} gives a
somewhat more invariant notion. More importantly, as the examples below show,
the closure properties are needed to control the laxators associated with
restricted proarrows.

As a first example, we show that restriction cells can be used to express the
universal property of finite products.

\begin{theory}[Categories with finite products] \label{th:finite-products}
  The \define{theory of finite products} $\Th{\cat{Fp}}$ is the restriction
  sketch generated by
  \begin{itemize}[noitemsep]
    \item an object $x$;
    \item arrows $\otimes: x^2 \to x$ and $I: 1 \to x$;
    \item proarrows $\otimes^*: x \proto x^2$ and $I^*: x \proto 1$, along with
      sketched restriction cells
      \begin{equation*}
        \begin{tikzcd}
          x & {x^2} \\
          x & {x }
          \arrow["\otimes", from=1-2, to=2-2]
          \arrow[""{name=0, anchor=center, inner sep=0}, "{\id_x}"', "\shortmid"{marking}, from=2-1, to=2-2]
          \arrow[""{name=1, anchor=center, inner sep=0}, "{\otimes^*}", "\shortmid"{marking}, from=1-1, to=1-2]
          \arrow[Rightarrow, no head, from=1-1, to=2-1]
          \arrow["{\mathrm{res}}"{description}, draw=none, from=1, to=0]
        \end{tikzcd}
        \qquad\text{and}\qquad
        \begin{tikzcd}
          x & 1 \\
          x & x
          \arrow[""{name=0, anchor=center, inner sep=0}, "{\id_x}"', "\shortmid"{marking}, from=2-1, to=2-2]
          \arrow["I", from=1-2, to=2-2]
          \arrow[""{name=1, anchor=center, inner sep=0}, "{I^*}", "\shortmid"{marking}, from=1-1, to=1-2]
          \arrow[Rightarrow, no head, from=1-1, to=2-1]
          \arrow["\res"{description}, draw=none, from=1, to=0]
        \end{tikzcd};
      \end{equation*}
    \item \define{projection} and \define{deletion} cells,
      \begin{equation*}
        \begin{tikzcd}
          {x^2} & {x^2} \\
          x & x
          \arrow[""{name=0, anchor=center, inner sep=0}, "{\mathrm{id}_x^2}", "\shortmid"{marking}, from=1-1, to=1-2]
          \arrow["\otimes"', from=1-1, to=2-1]
          \arrow["{\pi_{x,x}}", from=1-2, to=2-2]
          \arrow[""{name=1, anchor=center, inner sep=0}, "{\mathrm{id}_x}"', "\shortmid"{marking}, from=2-1, to=2-2]
          \arrow["\pi"{description}, draw=none, from=0, to=1]
        \end{tikzcd}
        \qquad
        \begin{tikzcd}
          {x^2} & {x^2} \\
          x & x
          \arrow[""{name=0, anchor=center, inner sep=0}, "{\mathrm{id}_x^2}", "\shortmid"{marking}, from=1-1, to=1-2]
          \arrow["\otimes"', from=1-1, to=2-1]
          \arrow["{\pi_{x,x}'}", from=1-2, to=2-2]
          \arrow[""{name=1, anchor=center, inner sep=0}, "{\mathrm{id}_x}"', "\shortmid"{marking}, from=2-1, to=2-2]
          \arrow["{\pi'}"{description}, draw=none, from=0, to=1]
        \end{tikzcd}
        \qquad
        \begin{tikzcd}[row sep=small]
          x & x \\
          & 1 \\
          x & x
          \arrow["{!}", from=1-2, to=2-2]
          \arrow["I", from=2-2, to=3-2]
          \arrow[""{name=0, anchor=center, inner sep=0}, "{\mathrm{id}_x}", "\shortmid"{marking}, from=1-1, to=1-2]
          \arrow[Rightarrow, no head, from=1-1, to=3-1]
          \arrow[""{name=1, anchor=center, inner sep=0}, "{\mathrm{id}_x}"', "\shortmid"{marking}, from=3-1, to=3-2]
          \arrow["\varepsilon"{description}, draw=none, from=0, to=1]
        \end{tikzcd}.
      \end{equation*}
  \end{itemize}
  These satisfy the following two axioms.
  \begin{itemize}
    \item Universal property of binary products: the cell below is a sketched
      restriction.
      \begin{equation*}
        \begin{tikzcd}
          x && {x^2} \\
          {x^2} && {x^2}
          \arrow[""{name=0, anchor=center, inner sep=0}, "{\otimes^*}", "\shortmid"{marking}, from=1-1, to=1-3]
          \arrow["{\Delta_x}"', from=1-1, to=2-1]
          \arrow[Rightarrow, no head, from=1-3, to=2-3]
          \arrow[""{name=1, anchor=center, inner sep=0}, "{\id_x^2}"', "\shortmid"{marking}, from=2-1, to=2-3]
          \arrow["{\langle \res \odot \pi, \res \odot \pi' \rangle}"{description}, draw=none, from=0, to=1]
        \end{tikzcd}
      \end{equation*}
    \item Universal property of terminal object: the cell below is a sketched
      restriction.
      \begin{equation*}
        \begin{tikzcd}[row sep=small]
          x & 1 \\
          x & x \\
          & 1 \\
          x & x
          \arrow["{!}", from=2-2, to=3-2]
          \arrow["I", from=3-2, to=4-2]
          \arrow[""{name=0, anchor=center, inner sep=0}, "{\mathrm{id}_x}"', "\shortmid"{marking}, from=2-1, to=2-2]
          \arrow[Rightarrow, no head, from=2-1, to=4-1]
          \arrow[""{name=1, anchor=center, inner sep=0}, "{\mathrm{id}_x}"', "\shortmid"{marking}, from=4-1, to=4-2]
          \arrow["I", from=1-2, to=2-2]
          \arrow[""{name=2, anchor=center, inner sep=0}, "{I^*}", "\shortmid"{marking}, from=1-1, to=1-2]
          \arrow[Rightarrow, no head, from=1-1, to=2-1]
          \arrow["\varepsilon"{description}, draw=none, from=0, to=1]
          \arrow["\res"{description}, draw=none, from=2, to=0]
        \end{tikzcd}
        \qedhere
      \end{equation*}
  \end{itemize}
\end{theory}

\begin{proposition}
  A model of the theory of finite products is precisely a cartesian category
  with chosen binary products and a chosen terminal object.
\end{proposition}
\begin{proof}
  A model $F: \Th{\cat{Fp}} \to \Span$ of the theory of finite products includes
  the data of a small category $\cat{C}$, functors
  $\otimes: \cat{C} \times \cat{C} \to \cat{C}$ and $I: 1 \to \cat{C}$, and
  natural families of projection maps $\pi_{c,c'}: c \otimes c' \to c$ and
  $\pi_{c,c'}: c \otimes c' \to c'$ and deletion maps $\varepsilon_c: c \to I$.
  These satisfiy the expected universal properties, as we now show.

  The cell $F(\res \odot \pi): F(\otimes^*) \to F \id_x$ maps morphisms in
  $\cat{C}$ of the form $h: a \to b \otimes c$ to the composite
  $h \cdot \pi_{b,c}: a \to b$, since we have
  \begin{equation*}
    \begin{tikzcd}
      Fx & {Fx^2} \\
      Fx & Fx
      \arrow[""{name=0, anchor=center, inner sep=0}, "{F\otimes^*}", "\shortmid"{marking}, from=1-1, to=1-2]
      \arrow["{F\pi_{x,x}}", from=1-2, to=2-2]
      \arrow[""{name=1, anchor=center, inner sep=0}, "{F\id_x}"', "\shortmid"{marking}, from=2-1, to=2-2]
      \arrow[Rightarrow, no head, from=1-1, to=2-1]
      \arrow["{F(\res \odot \pi)}"{description}, draw=none, from=0, to=1]
    \end{tikzcd}
    \quad=\quad
    \begin{tikzcd}
      Fx & {Fx^2} & {Fx^2} \\
      Fx & {Fx^2} & {Fx^2} \\
      Fx & Fx & Fx \\
      Fx && Fx
      \arrow[""{name=0, anchor=center, inner sep=0}, "{F\otimes^*}", "\shortmid"{marking}, from=2-1, to=2-2]
      \arrow[""{name=1, anchor=center, inner sep=0}, "{F \id_{x^2}}", "\shortmid"{marking}, from=2-2, to=2-3]
      \arrow[""{name=2, anchor=center, inner sep=0}, "{F \otimes^*}", "\shortmid"{marking}, from=1-1, to=1-2]
      \arrow["F\otimes"{description}, from=2-2, to=3-2]
      \arrow[""{name=3, anchor=center, inner sep=0}, "{F\id_x}"', "\shortmid"{marking}, from=3-1, to=3-2]
      \arrow[""{name=4, anchor=center, inner sep=0}, "{F\id_x}"', "\shortmid"{marking}, from=3-2, to=3-3]
      \arrow[""{name=5, anchor=center, inner sep=0}, "{F\id_x}"', "\shortmid"{marking}, from=4-1, to=4-3]
      \arrow[Rightarrow, no head, from=3-1, to=4-1]
      \arrow[Rightarrow, no head, from=3-3, to=4-3]
      \arrow["{F\pi_{x,x}}", from=2-3, to=3-3]
      \arrow[Rightarrow, no head, from=2-1, to=3-1]
      \arrow[""{name=6, anchor=center, inner sep=0}, "{\id_{Fx^2}}", "\shortmid"{marking}, from=1-2, to=1-3]
      \arrow[Rightarrow, no head, from=1-2, to=2-2]
      \arrow[Rightarrow, no head, from=1-3, to=2-3]
      \arrow[Rightarrow, no head, from=1-1, to=2-1]
      \arrow["{F_{x,x}}"{description}, draw=none, from=3-2, to=5]
      \arrow["F\pi"{description}, draw=none, from=1, to=4]
      \arrow["{F(\res)}"{description}, draw=none, from=0, to=3]
      \arrow["{F_{x^2}}"{description, pos=0.4}, draw=none, from=6, to=1]
      \arrow["1"{description, pos=0.4}, draw=none, from=2, to=0]
    \end{tikzcd}
  \end{equation*}
  by the naturality and unitality of the laxators. Similarly, the cell
  $F(\res \otimes \pi')$ maps morphisms of form $h: a \to b \otimes c$ to
  $h \cdot \pi_{b,c}': a \to c$. Thus, the first axiom asserts that any
  morphisms $f: a \to b$ and $g: a \to c$ with a common domain are equal to
  $h \cdot \pi_{b,c}$ and $h \cdot \pi_{b,c}'$ for a unique morphism
  $h: a \to b \otimes c$, which is precisely the universal property of a binary
  product.

  Moreover, the cell $F(\res \cdot \varepsilon): F(I^*) \to F \id_x$ maps
  morphisms in $\cat{C}$ of the form $f: a \to I$ to the deletion map
  $\varepsilon_x = f \cdot \varepsilon_I: x \to I$. It being a restriction cell
  along $F(I \circ !_x \circ I) = F(I)$ is then equivalent to the universal
  property of a terminal object.

  We have shown that a model of the theory gives a category with chosen binary
  and nullary products. Conversely, it is clear that a category $\cat{C}$ with
  chosen binary and nullary products determines functors
  $\otimes: \cat{C} \times \cat{C} \to \cat{C}$ and $I: 1 \to \cat{C}$ and
  natural transformations $\pi, \pi'$ and $\varepsilon$ that, in turn, determine
  a lax functor $F: \Th{\cat{Fp}} \to \Span$ at least on the generating data. It
  remains to check that the laxators of the lax functor $F$ do not require any
  extra data that is not already determined by the foregoing. This follows from
  \cref{lem:laxators-restrictions} and the closure properties of a restriction
  sketch (cf.\ \cref{lem:closure-restrictions}).

  A few examples will illustrate the general situation. First, the laxators
  $F_{x,\otimes^*}$ and $F_{\otimes^*,x^2}$ are uniquely determined by the
  laxator $F_{x,x}$ given the composition in $\cat{C}$. Moreover, since the
  cells
  \begin{equation*}
    \begin{tikzcd}
      x & {x^2} & {x^3} \\
      & {x^2} & {x^2} \\
      x & x & x
      \arrow[""{name=0, anchor=center, inner sep=0}, "{\otimes^*}", "\shortmid"{marking}, from=1-1, to=1-2]
      \arrow[""{name=1, anchor=center, inner sep=0}, "{\id_x \otimes^*}", "\shortmid"{marking}, from=1-2, to=1-3]
      \arrow[Rightarrow, no head, from=1-1, to=3-1]
      \arrow[Rightarrow, no head, from=1-2, to=2-2]
      \arrow["{1_x \otimes}", from=1-3, to=2-3]
      \arrow[""{name=2, anchor=center, inner sep=0}, "{\id_x^2}"', "\shortmid"{marking}, from=2-2, to=2-3]
      \arrow["\otimes"', from=2-2, to=3-2]
      \arrow[""{name=3, anchor=center, inner sep=0}, "{\id_x}"', "\shortmid"{marking}, from=3-1, to=3-2]
      \arrow[""{name=4, anchor=center, inner sep=0}, "{\id_x}"', "\shortmid"{marking}, from=3-2, to=3-3]
      \arrow["\otimes", from=2-3, to=3-3]
      \arrow["\res"{description}, draw=none, from=0, to=3]
      \arrow["{1_{\id_x}\; \res}"{description}, draw=none, from=1, to=2]
      \arrow["{\id_\otimes}"{description, pos=0.6}, draw=none, from=2, to=4]
    \end{tikzcd}
    \qquad\text{and}\qquad
    \begin{tikzcd}
      x & {x^2} & {x^3} \\
      & {x^2} & {x^2} \\
      x & x & x
      \arrow[""{name=0, anchor=center, inner sep=0}, "{\otimes^*}", "\shortmid"{marking}, from=1-1, to=1-2]
      \arrow[""{name=1, anchor=center, inner sep=0}, "{\otimes^* \id_x}", "\shortmid"{marking}, from=1-2, to=1-3]
      \arrow[Rightarrow, no head, from=1-1, to=3-1]
      \arrow[Rightarrow, no head, from=1-2, to=2-2]
      \arrow["{\otimes 1_x}", from=1-3, to=2-3]
      \arrow[""{name=2, anchor=center, inner sep=0}, "{\id_x^2}"', "\shortmid"{marking}, from=2-2, to=2-3]
      \arrow["\otimes"', from=2-2, to=3-2]
      \arrow[""{name=3, anchor=center, inner sep=0}, "{\id_x}"', "\shortmid"{marking}, from=3-1, to=3-2]
      \arrow[""{name=4, anchor=center, inner sep=0}, "{\id_x}"', "\shortmid"{marking}, from=3-2, to=3-3]
      \arrow["\otimes", from=2-3, to=3-3]
      \arrow["\res"{description}, draw=none, from=0, to=3]
      \arrow["{\res\; 1_{\id_x}}"{description}, draw=none, from=1, to=2]
      \arrow["{\id_\otimes}"{description, pos=0.6}, draw=none, from=2, to=4]
    \end{tikzcd}
  \end{equation*}
  are both sketched restrictions, corresponding to the conjoints
  $(\otimes \circ (1_x \times \otimes))^*$ and
  $(\otimes \circ (\otimes \times 1_x))^*$, the laxators
  $F_{\otimes^*, \id_x \otimes^*}$ and $F_{\otimes^*, \otimes^* \id_x}$ are also
  uniquely determined by $F_{x,x}$. Similar reasoning applies to any proarrows
  formed by products and external composites of $\otimes^*$, $I^*$, and $\id_x$.
\end{proof}

Using restriction sketches, we can treat symmetric multicategories, cartesian
multicategories, and other generalized multicategories
\cite[\S{2.6}]{shulman2016}. Unlike all of the others presented here, the
following double theory has a countably infinite, rather than a finite,
presentation.

\begin{theory}[Cartesian multicategories]
  Let $\cat{F}$ be the skeleton of $\FinSet$ spanned by the sets
  $[n] \coloneqq \{1,\dots,n\}$ for all $n \in \N$. The \define{theory of
    cartesian promonoids} is the restriction sketch augmenting the theory of
  promonoids (\cref{th:promonoid}) with, for every map $\sigma: [m] \to [n]$
  in $\cat{F}$,
  \begin{itemize}[noitemsep]
    \item a proarrow $x^\sigma_!: x^n \proto x^m$ along with a sketched
      restriction cell
      \begin{equation*}
        \begin{tikzcd}
          {x^n} & {x^m} \\
          {x^m} & {x^m}
          \arrow["{x^\sigma}"', from=1-1, to=2-1]
          \arrow[""{name=0, anchor=center, inner sep=0}, "{\id_x^m}"', "\shortmid"{marking}, from=2-1, to=2-2]
          \arrow[Rightarrow, no head, from=1-2, to=2-2]
          \arrow[""{name=1, anchor=center, inner sep=0}, "{x^\sigma_!}", "\shortmid"{marking}, from=1-1, to=1-2]
          \arrow["\res"{description}, draw=none, from=1, to=0]
        \end{tikzcd};
      \end{equation*}
    \item a \define{$\sigma$-action cell}
      \begin{equation*}
        \begin{tikzcd}
          {x^n} & x \\
          {x^n} & x
          \arrow[Rightarrow, no head, from=1-1, to=2-1]
          \arrow[""{name=0, anchor=center, inner sep=0}, "{p_m^\sigma}", "\shortmid"{marking}, from=1-1, to=1-2]
          \arrow[Rightarrow, no head, from=1-2, to=2-2]
          \arrow[""{name=1, anchor=center, inner sep=0}, "{p_n}"', "\shortmid"{marking}, from=2-1, to=2-2]
          \arrow["{\rho(\sigma)}"{description}, draw=none, from=0, to=1]
        \end{tikzcd}
        \quad\coloneqq\quad
        \begin{tikzcd}
          {x^n} & {x^m} & x \\
          {x^n} && x
          \arrow[Rightarrow, no head, from=1-1, to=2-1]
          \arrow["{x^\sigma_!}", "\shortmid"{marking}, from=1-1, to=1-2]
          \arrow["{p_m}", "\shortmid"{marking}, from=1-2, to=1-3]
          \arrow[Rightarrow, no head, from=1-3, to=2-3]
          \arrow[""{name=0, anchor=center, inner sep=0}, "{p_n}"', "\shortmid"{marking}, from=2-1, to=2-3]
          \arrow["{\rho(\sigma)}"{description, pos=0.4}, draw=none, from=1-2, to=0]
        \end{tikzcd},
      \end{equation*}
      where we abbreviate the composite proarrow $x^\sigma_! \odot p_m$ by
      $p_m^\sigma$.
  \end{itemize}
  These satisfy the following axioms.
  \begin{itemize}
    \item Functorality of action: for every pair of composable maps
      $[m] \xto{\sigma} [n] \xto{\tau} [q]$ in $\cat{F}$, we have
      $x^\tau_! \odot x^\sigma_! = x^{\sigma \cdot \tau}_!$ and
      \begin{equation*}
        \begin{tikzcd}
          {x^q} & {x^n} & {x^m} & x \\
          {x^q} & {x^n} && x \\
          {x^q} &&& x
          \arrow[Rightarrow, no head, from=1-2, to=2-2]
          \arrow[""{name=0, anchor=center, inner sep=0}, "{x^\sigma_!}", "\shortmid"{marking}, from=1-2, to=1-3]
          \arrow["{p_m}", "\shortmid"{marking}, from=1-3, to=1-4]
          \arrow[Rightarrow, no head, from=1-4, to=2-4]
          \arrow[""{name=1, anchor=center, inner sep=0}, "{p_n}"', "\shortmid"{marking}, from=2-2, to=2-4]
          \arrow[Rightarrow, no head, from=2-4, to=3-4]
          \arrow[""{name=2, anchor=center, inner sep=0}, "{x^\tau_!}", "\shortmid"{marking}, from=1-1, to=1-2]
          \arrow[Rightarrow, no head, from=1-1, to=2-1]
          \arrow[""{name=3, anchor=center, inner sep=0}, "{x^\tau_!}"', "\shortmid"{marking}, from=2-1, to=2-2]
          \arrow[""{name=4, anchor=center, inner sep=0}, "{p_q}"', "\shortmid"{marking}, from=3-1, to=3-4]
          \arrow[Rightarrow, no head, from=2-1, to=3-1]
          \arrow["{\rho(\sigma)}"{description}, draw=none, from=1-3, to=1]
          \arrow["1"{description}, draw=none, from=2, to=3]
          \arrow["{\rho(\tau)}"{description, pos=0.8}, draw=none, from=0, to=4]
        \end{tikzcd}
        \quad=\quad
        \begin{tikzcd}
          {x^q} & {x^m} & x \\
          {x^q} && x
          \arrow[Rightarrow, no head, from=1-1, to=2-1]
          \arrow["{x^{\sigma \cdot \tau}_!}", "\shortmid"{marking}, from=1-1, to=1-2]
          \arrow["{p_m}", "\shortmid"{marking}, from=1-2, to=1-3]
          \arrow[Rightarrow, no head, from=1-3, to=2-3]
          \arrow[""{name=0, anchor=center, inner sep=0}, "{p_q}"', "\shortmid"{marking}, from=2-1, to=2-3]
          \arrow["{\rho(\sigma \cdot \tau)}"{description}, draw=none, from=1-2, to=0]
        \end{tikzcd},
      \end{equation*}
      and for every $n \in \N$, we have $x^{1_{[n]}}_! = \id_x$ and
      $\rho(1_{[n]}) = 1_{p_n}$.
    \item Naturality of action (i): for all $k \in \N$ and maps
      $\sigma_i: [m_i] \to [n_i]$, $1 \leq i \leq k$, in $\cat{F}$, we have
      \begin{equation} \label{eq:multicategory-naturality-1}
        \begin{tikzcd}
          {x^n} && {x^k} & x \\
          {x^n} && {x^k} & x
          \arrow[""{name=0, anchor=center, inner sep=0}, "{p_{m_1}^{\sigma_1} \times \cdots \times p_{m_k}^{\sigma_k}}", "\shortmid"{marking}, from=1-1, to=1-3]
          \arrow[""{name=1, anchor=center, inner sep=0}, "{p_k}", "\shortmid"{marking}, from=1-3, to=1-4]
          \arrow[""{name=2, anchor=center, inner sep=0}, "{p_k}"', "\shortmid"{marking}, from=2-3, to=2-4]
          \arrow[Rightarrow, no head, from=1-3, to=2-3]
          \arrow[Rightarrow, no head, from=1-4, to=2-4]
          \arrow[Rightarrow, no head, from=1-1, to=2-1]
          \arrow[""{name=3, anchor=center, inner sep=0}, "{p_{n_1} \times \cdots \times p_{n_k}}"', "\shortmid"{marking}, from=2-1, to=2-3]
          \arrow["{1_{p_k}}"{description}, draw=none, from=1, to=2]
          \arrow["{\rho(\sigma_1) \times \cdots \times \rho(\sigma_k)}"{description}, draw=none, from=0, to=3]
        \end{tikzcd}
        \quad=\quad
        \begin{tikzcd}
          {x^n} && x \\
          {x^n} && x
          \arrow[""{name=0, anchor=center, inner sep=0}, "{p_n}"', "\shortmid"{marking}, from=2-1, to=2-3]
          \arrow[""{name=1, anchor=center, inner sep=0}, "{p_m^{\sigma_1 + \cdots + \sigma_k}}", "\shortmid"{marking}, from=1-1, to=1-3]
          \arrow[Rightarrow, no head, from=1-1, to=2-1]
          \arrow[Rightarrow, no head, from=1-3, to=2-3]
          \arrow["{\rho(\sigma_1 + \cdots + \sigma_k)}"{description}, draw=none, from=1, to=0]
        \end{tikzcd},
      \end{equation}
      where $m \coloneqq m_1 + \cdots + m_k$ and
      $n \coloneqq n_1 + \cdots + n_k$, which is well-typed since
      \begin{equation*}
        x^{\sigma_1}_! \times \cdots \times x^{\sigma_k}_!
          = (x^{\sigma_1} \times \cdots \times x^{\sigma_k})_!
          = x^{\sigma_1 + \cdots + \sigma_k}_!.
      \end{equation*}
    \item Naturality of action (ii): for every map $\sigma: [k] \to [\ell]$ in
      $\cat{F}$ and all numbers $n_i \in \N$, $1 \leq i \leq \ell$, we have
      \begin{equation} \label{eq:multicategory-naturality-2}
        \begin{tikzcd}
          {x^n} && {x^\ell} & x \\
          {x^n} && {x^\ell} & x
          \arrow[""{name=0, anchor=center, inner sep=0}, "{p_{n_1} \times \cdots \times p_{n_\ell}}", "\shortmid"{marking}, from=1-1, to=1-3]
          \arrow[""{name=1, anchor=center, inner sep=0}, "{p_{n_1} \times \cdots \times p_{n_\ell}}"', "\shortmid"{marking}, from=2-1, to=2-3]
          \arrow[Rightarrow, no head, from=1-1, to=2-1]
          \arrow[Rightarrow, no head, from=1-3, to=2-3]
          \arrow[""{name=2, anchor=center, inner sep=0}, "{p_k^\sigma}", "\shortmid"{marking}, from=1-3, to=1-4]
          \arrow[""{name=3, anchor=center, inner sep=0}, "{p_\ell}"', "\shortmid"{marking}, from=2-3, to=2-4]
          \arrow[Rightarrow, no head, from=1-4, to=2-4]
          \arrow["1"{description}, draw=none, from=0, to=1]
          \arrow["{\rho(\sigma)}"{description}, draw=none, from=2, to=3]
        \end{tikzcd}
        \quad=\quad
        \begin{tikzcd}
          {x^n} && x \\
          {x^n} && x
          \arrow[""{name=0, anchor=center, inner sep=0}, "{p_m^{\sigma \wr (n_1,\dots,n_\ell)}}", "\shortmid"{marking}, from=1-1, to=1-3]
          \arrow[""{name=1, anchor=center, inner sep=0}, "{p_n}"', "\shortmid"{marking}, from=2-1, to=2-3]
          \arrow[Rightarrow, no head, from=1-1, to=2-1]
          \arrow[Rightarrow, no head, from=1-3, to=2-3]
          \arrow["{\rho(\sigma \wr (n_1,\dots,n_\ell))}"{description}, draw=none, from=0, to=1]
        \end{tikzcd},
      \end{equation}
      where $m \coloneqq n_{\sigma 1} + \cdots + n_{\sigma k}$ and
      $n \coloneqq n_1 + \cdots + n_\ell$, and the map
      $\sigma \wr (n_1,\dots,n_\ell): [m] \to [n]$ applies $\sigma$ blockwise
      and is the identity within each block $n_i$, and we also have the side
      equations
      \begin{equation} \label{eq:multicategory-naturality-2-1}
        (p_{n_1} \times \cdots \times p_{n_\ell}) \odot x^\sigma_! =
          x^{\sigma \wr (n_1,\dots,n_\ell)}_! \odot (p_{n_{\sigma 1}} \times \cdots \times p_{n_{\sigma k}})
      \end{equation}
      and
      \begin{equation} \label{eq:multicategory-naturality-2-2}
        \begin{tikzcd}
          {x^n} && {x^\ell} & {x^k} \\
          {x^m} && {x^k} & {x^k}
          \arrow[""{name=0, anchor=center, inner sep=0}, "{p_{n_1} \times \cdots \times p_{n_\ell}}", "\shortmid"{marking}, from=1-1, to=1-3]
          \arrow[""{name=1, anchor=center, inner sep=0}, "{x^\sigma_!}", "\shortmid"{marking}, from=1-3, to=1-4]
          \arrow["{x^{\sigma \wr (n_1,\dots,n_\ell)}}"', from=1-1, to=2-1]
          \arrow["{x^\sigma}"', from=1-3, to=2-3]
          \arrow[""{name=2, anchor=center, inner sep=0}, "{p_{n_{\sigma 1}} \times \cdots \times p_{n_{\sigma k}}}"', "\shortmid"{marking}, from=2-1, to=2-3]
          \arrow[Rightarrow, no head, from=1-4, to=2-4]
          \arrow[""{name=3, anchor=center, inner sep=0}, "{\id_x^k}"', "\shortmid"{marking}, from=2-3, to=2-4]
          \arrow["{\sigma_{p_1,\dots,p_\ell}}"{description}, draw=none, from=0, to=2]
          \arrow["\res"{description}, draw=none, from=1, to=3]
        \end{tikzcd}
        =
        \begin{tikzcd}
          {x^n} & {x^m} && {x^k} \\
          {x^m} & {x^m} && {x^k}
          \arrow["{x^{\sigma \wr (n_1,\dots,n_\ell)}}"', from=1-1, to=2-1]
          \arrow[""{name=0, anchor=center, inner sep=0}, "{\id_x^m}"', "\shortmid"{marking}, from=2-1, to=2-2]
          \arrow[Rightarrow, no head, from=1-2, to=2-2]
          \arrow[""{name=1, anchor=center, inner sep=0}, "{x^{\sigma \wr (n_1,\dots,n_\ell)}_!}", "\shortmid"{marking}, from=1-1, to=1-2]
          \arrow[""{name=2, anchor=center, inner sep=0}, "{p_{n_{\sigma 1}} \times \cdots \times p_{n_{\sigma k}}}"', "\shortmid"{marking}, from=2-2, to=2-4]
          \arrow[""{name=3, anchor=center, inner sep=0}, "{p_{n_{\sigma 1}} \times \cdots \times p_{n_{\sigma k}}}", "\shortmid"{marking}, from=1-2, to=1-4]
          \arrow[Rightarrow, no head, from=1-4, to=2-4]
          \arrow["1"{description}, draw=none, from=3, to=2]
          \arrow["\res"{description}, draw=none, from=1, to=0]
        \end{tikzcd}.
        \qedhere
      \end{equation}
    \end{itemize}
\end{theory}

\begin{proposition}
  A model of the theory of cartesian promonoids is precisely a cartesian
  multicategory.
\end{proposition}
\begin{proof}
  A model $F$ of the theory in $\Mat$ includes the data of a multicategory
  $\cat{C}$ equipped with, for each map $\sigma: [m] \to [n]$ in $\cat{F}$, a
  family of operations
  \begin{equation*}
    \sigma^* \coloneqq F(\rho(\sigma)):
      \cat{C}(c_{\sigma(1)}, \dots, c_{\sigma(m)}; c) \to
      \cat{C}(c_1,\dots,c_n; c), \qquad
    c_1, \dots, c_n, c \in \cat{C}.
  \end{equation*}
  We have used that $F(p_m^\sigma) = F(x_!^\sigma \odot p_m)$ is the restriction
  of $F(p_m)$ along $x^\sigma: [n] \to [m]$, by \cref{lem:closure-restrictions}.
  Moreover, by the functorality axiom, we have
  $\tau^* \circ \sigma^* = (\tau \circ \sigma)^*$ for composable maps $\sigma$
  and $\tau$ in $\cat{F}$, and also that $(1_{[n]})^*$ is the identity for each
  $n \in \N$. To establish that $\cat{C}$ is cartesian, we must show that the
  two naturality axioms prove the remaining two axioms of a cartesian
  multicategory \cite[Definition 2.6.4]{shulman2016}, namely
  \begin{equation} \label{eq:multicategory-naturality-1-alt}
    g \circ (\sigma_1^* f_1, \dots, \sigma_k^* f_k) =
      (\sigma_1 + \dots + \sigma_k)^*(g \circ (f_1,\dots,f_k))
  \end{equation}
  for maps $\sigma_i: [m_i] \to [n_i]$, $1 \leq i \leq k$, in $\cat{F}$, and
  \begin{equation} \label{eq:multicategory-naturality-2-alt}
    (\sigma^* g) \circ (f_1,\dots,f_\ell) =
      (\sigma \wr (n_1,\dots,n_\ell))^*(g \circ (f_{\sigma(1)}, \dots, f_{\sigma(k)}))
  \end{equation}
  for each map $\sigma: [k] \to [\ell]$ in $\cat{F}$, where these equations
  range over all well-typed $m_i$-ary (resp.\ $n_i$-ary) multimorphisms $f_i$,
  for $i \in [k]$ (resp.\ $i \in [\ell]$), and $k$-ary multimorphisms $g$ in
  $\cat{C}$.

  By the naturality axiom \eqref{eq:multicategory-naturality-1} for the action,
  as well as the naturality of the laxators, we have
  \begin{equation*}
    \begin{tikzcd}
      {Fx^n} && {Fx^k} & Fx \\
      {Fx^n} && {Fx^k} & Fx \\
      {Fx^ n} &&& Fx
      \arrow[""{name=0, anchor=center, inner sep=0}, "{Fp_k}", "\shortmid"{marking}, from=1-3, to=1-4]
      \arrow[""{name=1, anchor=center, inner sep=0}, "{Fp_{m_1}^{\sigma_1} \times \cdots \times Fp_{m_k}^{\sigma_k}}", "\shortmid"{marking}, from=1-1, to=1-3]
      \arrow[""{name=2, anchor=center, inner sep=0}, "{Fp_k}"', "\shortmid"{marking}, from=2-3, to=2-4]
      \arrow[Rightarrow, no head, from=1-4, to=2-4]
      \arrow[Rightarrow, no head, from=1-3, to=2-3]
      \arrow[Rightarrow, no head, from=1-1, to=2-1]
      \arrow[""{name=3, anchor=center, inner sep=0}, "{Fp_{n_1} \times \cdots \times Fp_{n_k}}"', "\shortmid"{marking}, from=2-1, to=2-3]
      \arrow[Rightarrow, no head, from=2-4, to=3-4]
      \arrow[Rightarrow, no head, from=2-1, to=3-1]
      \arrow[""{name=4, anchor=center, inner sep=0}, "{Fp_n}"', "\shortmid"{marking}, from=3-1, to=3-4]
      \arrow["{\sigma_1^* \times \cdots \times \sigma_k^*}"{description}, draw=none, from=1, to=3]
      \arrow["1"{description}, draw=none, from=0, to=2]
      \arrow["{F_{(p_{n_1} \times \cdots \times p_{n_k}),p_k}}"{description}, draw=none, from=2-3, to=4]
    \end{tikzcd}
    \quad=\quad
    \begin{tikzcd}
      {Fx^n} && {Fx^k} & Fx \\
      {Fx^n} &&& Fx \\
      {Fx^ n} &&& Fx
      \arrow["{Fp_k}", "\shortmid"{marking}, from=1-3, to=1-4]
      \arrow[""{name=0, anchor=center, inner sep=0}, "{Fp_{m_1}^{\sigma_1} \times \cdots \times Fp_{m_k}^{\sigma_k}}", "\shortmid"{marking}, from=1-1, to=1-3]
      \arrow[Rightarrow, no head, from=1-4, to=2-4]
      \arrow[Rightarrow, no head, from=2-4, to=3-4]
      \arrow[""{name=1, anchor=center, inner sep=0}, "{Fp_n}"', "\shortmid"{marking}, from=3-1, to=3-4]
      \arrow[Rightarrow, no head, from=1-1, to=2-1]
      \arrow[Rightarrow, no head, from=2-1, to=3-1]
      \arrow[""{name=2, anchor=center, inner sep=0}, "{Fp_n^{\sigma_1 + \cdots + \sigma_k}}"', "\shortmid"{marking}, from=2-1, to=2-4]
      \arrow["{(\sigma_1 + \cdots + \sigma_k)^*}"{description, pos=0.6}, draw=none, from=2, to=1]
      \arrow["{F_{(p_{m_1}^{\sigma_1} \times \cdots \times p_{m_k}^{\sigma_k}),p_k}}"{description}, draw=none, from=0, to=2]
    \end{tikzcd}.
  \end{equation*}
  This is \cref{eq:multicategory-naturality-1-alt} provided that the laxator on
  the right-hand side has the expected behavior. But, by
  \cref{lem:closure-restrictions}, the cell
  \begin{equation*}
    \begin{tikzcd}
      {x^n} && {x^k} & x \\
      {x^m} && {x^k} & x
      \arrow["{x^{\sigma_1} \times \cdots \times x^{\sigma_k}}"', from=1-1, to=2-1]
      \arrow[""{name=0, anchor=center, inner sep=0}, "{p_{m_1}^{\sigma_1} \times \cdots \times p_{m_k}^{\sigma_k}}", "\shortmid"{marking}, from=1-1, to=1-3]
      \arrow[Rightarrow, no head, from=1-3, to=2-3]
      \arrow[""{name=1, anchor=center, inner sep=0}, "{p_k}", "\shortmid"{marking}, from=1-3, to=1-4]
      \arrow[Rightarrow, no head, from=1-4, to=2-4]
      \arrow[""{name=2, anchor=center, inner sep=0}, "{p_k}"', "\shortmid"{marking}, from=2-3, to=2-4]
      \arrow[""{name=3, anchor=center, inner sep=0}, "{p_{m_1} \times \cdots \times p_{m_k}}"', "\shortmid"{marking}, from=2-1, to=2-3]
      \arrow["1"{description}, draw=none, from=1, to=2]
      \arrow["{\prod_{i=1}^k \res}"{description}, draw=none, from=0, to=3]
    \end{tikzcd}
  \end{equation*}
  is a sketched restriction and hence by \cref{lem:laxators-restrictions} the
  laxator for $p_{m_1}^{\sigma_1} \times \cdots \times p_{m_k}^{\sigma_k}$ and
  $p_k$ is uniquely determined by the laxator for
  $p_{m_1} \times \cdots \times p_{m_k}$ and $p_k$, which is a composition
  operation in the multicategory.

  Similarly, the naturality axiom \eqref{eq:multicategory-naturality-2} for the
  action along with the naturality of the laxators implies that
  \begin{equation*}
    \begin{tikzcd}
      {Fx^n} && {Fx^\ell} & Fx \\
      {Fx^n} && {Fx^\ell} & Fx \\
      {Fx^n} &&& Fx
      \arrow[""{name=0, anchor=center, inner sep=0}, "{Fp_k^\sigma}", "\shortmid"{marking}, from=1-3, to=1-4]
      \arrow[""{name=1, anchor=center, inner sep=0}, "{Fp_\ell}"', "\shortmid"{marking}, from=2-3, to=2-4]
      \arrow[Rightarrow, no head, from=1-4, to=2-4]
      \arrow[Rightarrow, no head, from=1-3, to=2-3]
      \arrow[""{name=2, anchor=center, inner sep=0}, "{Fp_{n_1} \times \cdots \times Fp_{n_\ell}}", "\shortmid"{marking}, from=1-1, to=1-3]
      \arrow[Rightarrow, no head, from=1-1, to=2-1]
      \arrow[""{name=3, anchor=center, inner sep=0}, "{Fp_{n_1} \times \cdots \times F_{n_\ell}}"', "\shortmid"{marking}, from=2-1, to=2-3]
      \arrow[Rightarrow, no head, from=2-4, to=3-4]
      \arrow[Rightarrow, no head, from=2-1, to=3-1]
      \arrow[""{name=4, anchor=center, inner sep=0}, "{Fp_n}"', "\shortmid"{marking}, from=3-1, to=3-4]
      \arrow["1"{description}, draw=none, from=2, to=3]
      \arrow["{\sigma^*}"{description}, draw=none, from=0, to=1]
      \arrow["{F_{(p_{n_1} \times \cdots \times p_{n_\ell}),p_\ell}}"{description}, draw=none, from=2-3, to=4]
    \end{tikzcd}
    \quad=\quad
    \begin{tikzcd}
      {Fx^n} && {Fx^\ell} & Fx \\
      {Fx^n} &&& Fx \\
      {Fx^n} &&& Fx
      \arrow["{Fp_k^\sigma}", "\shortmid"{marking}, from=1-3, to=1-4]
      \arrow[Rightarrow, no head, from=1-4, to=2-4]
      \arrow["{Fp_{n_1} \times \cdots \times Fp_{n_\ell}}", "\shortmid"{marking}, from=1-1, to=1-3]
      \arrow[Rightarrow, no head, from=1-1, to=2-1]
      \arrow[Rightarrow, no head, from=2-4, to=3-4]
      \arrow[Rightarrow, no head, from=2-1, to=3-1]
      \arrow[""{name=0, anchor=center, inner sep=0}, "{Fp_n}"', "\shortmid"{marking}, from=3-1, to=3-4]
      \arrow[""{name=1, anchor=center, inner sep=0}, "{Fp_m^{\sigma \wr (n_1,\dots,n_\ell)}}"', "\shortmid"{marking}, from=2-1, to=2-4]
      \arrow["{F_{(p_{n_1} \times \cdots \times p_{n_\ell}),p_k^\sigma}}"{description, pos=0.4}, draw=none, from=1-3, to=1]
      \arrow["{(\sigma \wr (n_1,\dots,n_\ell))^*}"{description, pos=0.6}, draw=none, from=1, to=0]
    \end{tikzcd},
  \end{equation*}
  where we have used \cref{eq:multicategory-naturality-2-1} on the right-hand
  side. This is \cref{eq:multicategory-naturality-2-alt} provided that the
  laxator on the right-hand side has the desired effect. Now,
  \cref{eq:multicategory-naturality-2-2} implies that
  \begin{equation*}
    \begin{tikzcd}
      {x^n} && {x^\ell} & x \\
      {x^m} && {x^k} & x
      \arrow[""{name=0, anchor=center, inner sep=0}, "{p_{n_1} \times \cdots \times p_{n_\ell}}", "\shortmid"{marking}, from=1-1, to=1-3]
      \arrow[""{name=1, anchor=center, inner sep=0}, "{p_k^\sigma}", "\shortmid"{marking}, from=1-3, to=1-4]
      \arrow["{x^{\sigma \wr (n_1,\dots,n_\ell)}}"', from=1-1, to=2-1]
      \arrow["{x^\sigma}"', from=1-3, to=2-3]
      \arrow[""{name=2, anchor=center, inner sep=0}, "{p_{n_{\sigma 1}} \times \cdots \times p_{n_{\sigma k}}}"', "\shortmid"{marking}, from=2-1, to=2-3]
      \arrow[Rightarrow, no head, from=1-4, to=2-4]
      \arrow[""{name=3, anchor=center, inner sep=0}, "{p_k}"', "\shortmid"{marking}, from=2-3, to=2-4]
      \arrow["{\sigma_{p_1,\dots,p_\ell}}"{description}, draw=none, from=0, to=2]
      \arrow["\res"{description}, draw=none, from=1, to=3]
    \end{tikzcd}
    \quad=\quad
    \begin{tikzcd}
      {x^n} & {x^m} & x \\
      {x^m} & {x^m} & x
      \arrow["{x^{\sigma \wr (n_1,\dots,n_\ell)}}"', from=1-1, to=2-1]
      \arrow[""{name=0, anchor=center, inner sep=0}, "{\id_x^m}"', "\shortmid"{marking}, from=2-1, to=2-2]
      \arrow[Rightarrow, no head, from=1-2, to=2-2]
      \arrow[""{name=1, anchor=center, inner sep=0}, "{x^{\sigma \wr (n_1,\dots,n_\ell)}_!}", "\shortmid"{marking}, from=1-1, to=1-2]
      \arrow[""{name=2, anchor=center, inner sep=0}, "{p_m}", "\shortmid"{marking}, from=1-2, to=1-3]
      \arrow[Rightarrow, no head, from=1-3, to=2-3]
      \arrow[""{name=3, anchor=center, inner sep=0}, "{p_m}"', "\shortmid"{marking}, from=2-2, to=2-3]
      \arrow["\res"{description}, draw=none, from=1, to=0]
      \arrow["1"{description}, draw=none, from=2, to=3]
    \end{tikzcd},
  \end{equation*}
  where, by \cref{lem:closure-restrictions}, the right-hand side is a sketched
  restriction and hence so is the left-hand side. Thus, by
  \cref{lem:laxators-restrictions}, the laxator for
  $p_{n_1} \times \cdots \times p_{n_\ell}$ and $p_k^\sigma$ is uniquely
  determined by the laxator for
  $p_{n_{\sigma(1)}} \times \cdots \times p_{n_{\sigma(k)}}$ and $p_k$, which is
  again a composition operation in the multicategory.

  We have shown that a model $F$ of the theory of cartesian promonoids gives a
  cartesian multicategory. To establish the converse, that a cartesian
  multicategory uniquely determines a model, requires checking that the laxators
  of $F$ do not contain any further data that is not already uniquely determined
  by the cartesian multicategory operations. This follows from
  \cref{lem:laxators-restrictions,lem:closure-restrictions} as in the previous
  proposition.
\end{proof}

By confining the action by the category $\cat{F}$ to wide subcategories of
$\cat{F}$ that are closed under the cocartesian monoidal product and the wreath
operation (called ``faithful cartesian clubs'' in \cite[Definition
2.6.3]{shulman2016}), we can define restriction sketches whose models are other
kinds of generalized multicategories. Most importantly, taking the core of
$\cat{F}$, which is the category of finite sets and bijections, yields symmetric
multicategories.

\section{Lax transformations}
\label{sec:lax-transformations}

Having seen numerous examples of simple and cartesian double theories and their
models, we begin the task of constructing the virtual double category of models
of a double theory. In this section, we construct a mere category of models,
focusing on the natural transformations between lax functors that will play the
role of homomorphisms between models.

Natural transformations between double functors, be they strict, pseudo, or lax,
are standard \cite[Definition 3.5.4]{grandis2019}, and they give the correct
definition of strict maps between models of a double theory. However, it is also
typical in two-dimensional category theory to consider lax and oplax maps, which
requires notions of \emph{lax} and \emph{oplax} natural transformations. The
following definition of a lax transformation generalizes what Grandis calls a
\emph{pseudo transformation} of lax double functors \cite[Definition
3.8.1]{grandis2019} but to our knowledge does not appear in the literature. We
therefore present a detailed account.

\begin{definition}[Lax transformation] \label{def:lax-transformation}
  A \define{lax natural transformation} $\alpha: F \To G$ of lax double functors
  $F,G: \dbl{D} \to \dbl{E}$ consists of
  \begin{itemize}[noitemsep]
    \item for every object $x \in \dbl{D}$, its \define{component} at $x$, an
      arrow $\alpha_x: Fx \to Gx$ in $\dbl{E}$;
    \item for every proarrow $m: x \proto y$ in $\dbl{D}$, its
      \define{component} at $m$, a cell in $\dbl{E}$
    \begin{equation*}
      \begin{tikzcd}
        Fx & Fy \\
        Gx & Gy
        \arrow["{\alpha_x}"', from=1-1, to=2-1]
        \arrow[""{name=0, anchor=center, inner sep=0}, "Fm", "\shortmid"{marking}, from=1-1, to=1-2]
        \arrow["{\alpha_y}", from=1-2, to=2-2]
        \arrow[""{name=1, anchor=center, inner sep=0}, "{G_m}"', "\shortmid"{marking}, from=2-1, to=2-2]
        \arrow["{\alpha_m}"{description}, draw=none, from=0, to=1]
      \end{tikzcd};
    \end{equation*}
    \item for every arrow $f: x \to y$ in $\dbl{D}$, its \define{naturality
      comparison} at $f$, a cell in $\dbl{E}$
    \begin{equation*}
      \begin{tikzcd}
        Fx & Fx \\
        Gx & Fy \\
        Gy & Gy
        \arrow[""{name=0, anchor=center, inner sep=0}, "{\mathrm{id}_{Fx}}", "\shortmid"{marking}, from=1-1, to=1-2]
        \arrow["{\alpha_x}"', from=1-1, to=2-1]
        \arrow["Gf"', from=2-1, to=3-1]
        \arrow["Ff", from=1-2, to=2-2]
        \arrow["{\alpha_y}", from=2-2, to=3-2]
        \arrow[""{name=1, anchor=center, inner sep=0}, "{G(\mathrm{id}_y)}"', "\shortmid"{marking}, from=3-1, to=3-2]
        \arrow["{\alpha_f}"{description}, draw=none, from=0, to=1]
      \end{tikzcd};
    \end{equation*}
  \end{itemize}
  such that the following axioms are satisfied:
  \begin{itemize}
    \item Naturality with respect to cells: for every cell
      $\stdInlineCell{\gamma}$ in $\dbl{D}$,
    \begin{equation*}
      \begin{tikzcd}
        Fx & Fx & Fy \\
        Gw & Fw & Fz \\
        Gw & Gw & Gz \\
        Gw && Gz
        \arrow[""{name=0, anchor=center, inner sep=0}, "{\mathrm{id}_{Fx}}", "\shortmid"{marking}, from=1-1, to=1-2]
        \arrow["{\alpha_x}"', from=1-1, to=2-1]
        \arrow["Gf"', from=2-1, to=3-1]
        \arrow["Ff"', from=1-2, to=2-2]
        \arrow["{\alpha_w}"', from=2-2, to=3-2]
        \arrow[""{name=1, anchor=center, inner sep=0}, "{G \mathrm{id}_w}"', "\shortmid"{marking}, from=3-1, to=3-2]
        \arrow["Fg", from=1-3, to=2-3]
        \arrow[""{name=2, anchor=center, inner sep=0}, "Fm", "\shortmid"{marking}, from=1-2, to=1-3]
        \arrow[""{name=3, anchor=center, inner sep=0}, "Fn"', "\shortmid"{marking}, from=2-2, to=2-3]
        \arrow[""{name=4, anchor=center, inner sep=0}, "Gn"', "\shortmid"{marking}, from=3-2, to=3-3]
        \arrow["{\alpha_z}", from=2-3, to=3-3]
        \arrow[Rightarrow, no head, from=3-1, to=4-1]
        \arrow[Rightarrow, no head, from=3-3, to=4-3]
        \arrow[""{name=5, anchor=center, inner sep=0}, "Gn"', "\shortmid"{marking}, from=4-1, to=4-3]
        \arrow["{\alpha_f}"{description}, draw=none, from=0, to=1]
        \arrow["F\gamma"{description}, draw=none, from=2, to=3]
        \arrow["{\alpha_n}"{description}, draw=none, from=3, to=4]
        \arrow["{G_{w,n}}"{description}, draw=none, from=3-2, to=5]
      \end{tikzcd}
      \quad=\quad
      \begin{tikzcd}
        Fx & Fy & Fy \\
        Gx & Gy & Fz \\
        Gw & Gz & Gz \\
        Gw && Gz
        \arrow[""{name=0, anchor=center, inner sep=0}, "{\mathrm{id}_{Fy}}", "\shortmid"{marking}, from=1-2, to=1-3]
        \arrow["{\alpha_y}", from=1-2, to=2-2]
        \arrow["Gg", from=2-2, to=3-2]
        \arrow["Fg", from=1-3, to=2-3]
        \arrow["{\alpha_z}", from=2-3, to=3-3]
        \arrow[""{name=1, anchor=center, inner sep=0}, "{G \mathrm{id}_z}"', "\shortmid"{marking}, from=3-2, to=3-3]
        \arrow[""{name=2, anchor=center, inner sep=0}, "Fm", "\shortmid"{marking}, from=1-1, to=1-2]
        \arrow["{\alpha_x}"', from=1-1, to=2-1]
        \arrow[""{name=3, anchor=center, inner sep=0}, "Gm"', "\shortmid"{marking}, from=2-1, to=2-2]
        \arrow[""{name=4, anchor=center, inner sep=0}, "Gn"', "\shortmid"{marking}, from=3-1, to=3-2]
        \arrow["Gf"', from=2-1, to=3-1]
        \arrow[Rightarrow, no head, from=3-1, to=4-1]
        \arrow[Rightarrow, no head, from=3-3, to=4-3]
        \arrow[""{name=5, anchor=center, inner sep=0}, "Gn"', "\shortmid"{marking}, from=4-1, to=4-3]
        \arrow["{\alpha_g}"{description}, draw=none, from=0, to=1]
        \arrow["{\alpha_m}"{description}, draw=none, from=2, to=3]
        \arrow["G\gamma"{description}, draw=none, from=3, to=4]
        \arrow["{G_{n,z}}"{description}, draw=none, from=3-2, to=5]
      \end{tikzcd}
    \end{equation*}
    \item External functorality: for every consecutive pair of proarrows
    $x \xproto{m} y \xproto{n} z$ in $\dbl{D}$,
    \begin{equation*}
      \begin{tikzcd}
        Fx & Fy & Fz \\
        Fx && Fz \\
        Gx && Gz
        \arrow["{\alpha_x}"', from=2-1, to=3-1]
        \arrow["{\alpha_z}", from=2-3, to=3-3]
        \arrow[Rightarrow, no head, from=1-1, to=2-1]
        \arrow[Rightarrow, no head, from=1-3, to=2-3]
        \arrow[""{name=0, anchor=center, inner sep=0}, "{F(m \odot n)}"', "\shortmid"{marking}, from=2-1, to=2-3]
        \arrow["Fm", "\shortmid"{marking}, from=1-1, to=1-2]
        \arrow["Fn", "\shortmid"{marking}, from=1-2, to=1-3]
        \arrow[""{name=1, anchor=center, inner sep=0}, "{G(m \odot n)}"', "\shortmid"{marking}, from=3-1, to=3-3]
        \arrow["{F_{m,n}}"{description}, draw=none, from=1-2, to=0]
        \arrow["{\alpha_{m \odot n}}"{description}, draw=none, from=0, to=1]
      \end{tikzcd}
      \quad=\quad
      \begin{tikzcd}
        Fx & Fy & Fz \\
        Gx & Gy & Gz \\
        Gx && Gz
        \arrow[""{name=0, anchor=center, inner sep=0}, "Fm", "\shortmid"{marking}, from=1-1, to=1-2]
        \arrow["{\alpha_x}"', from=1-1, to=2-1]
        \arrow["{\alpha_y}", from=1-2, to=2-2]
        \arrow[""{name=1, anchor=center, inner sep=0}, "Fn", "\shortmid"{marking}, from=1-2, to=1-3]
        \arrow["{\alpha_z}", from=1-3, to=2-3]
        \arrow[""{name=2, anchor=center, inner sep=0}, "Gm"', "\shortmid"{marking}, from=2-1, to=2-2]
        \arrow[""{name=3, anchor=center, inner sep=0}, "Gn"', "\shortmid"{marking}, from=2-2, to=2-3]
        \arrow[Rightarrow, no head, from=2-1, to=3-1]
        \arrow[Rightarrow, no head, from=2-3, to=3-3]
        \arrow[""{name=4, anchor=center, inner sep=0}, "{G(m \odot n)}"', "\shortmid"{marking}, from=3-1, to=3-3]
        \arrow["{\alpha_m}"{description}, draw=none, from=0, to=2]
        \arrow["{\alpha_n}"{description}, draw=none, from=1, to=3]
        \arrow["{G_{m,n}}"{description}, draw=none, from=2-2, to=4]
      \end{tikzcd},
    \end{equation*}
    and for every object $x \in \dbl{D}$,
    \begin{equation*}
      \begin{tikzcd}
        Fx & Fx \\
        Fx & Fx \\
        Gx & Gx
        \arrow[""{name=0, anchor=center, inner sep=0}, "{\mathrm{id}_{Fx}}", "\shortmid"{marking}, from=1-1, to=1-2]
        \arrow[Rightarrow, no head, from=1-1, to=2-1]
        \arrow[""{name=1, anchor=center, inner sep=0}, "{F \mathrm{id}_x}"', "\shortmid"{marking}, from=2-1, to=2-2]
        \arrow[Rightarrow, no head, from=1-2, to=2-2]
        \arrow[""{name=2, anchor=center, inner sep=0}, "{G \mathrm{id}_x}"', "\shortmid"{marking}, from=3-1, to=3-2]
        \arrow["{\alpha_x}"', from=2-1, to=3-1]
        \arrow["{\alpha_x}", from=2-2, to=3-2]
        \arrow["{F_x}"{description}, draw=none, from=0, to=1]
        \arrow["{\alpha_{\mathrm{id}_x}}"{description}, draw=none, from=1, to=2]
      \end{tikzcd}
      \quad=\quad
      \begin{tikzcd}
        Fx & Fx \\
        Gx & Gx \\
        Gx & Gx
        \arrow[""{name=0, anchor=center, inner sep=0}, "{\mathrm{id}_{Fx}}", "\shortmid"{marking}, from=1-1, to=1-2]
        \arrow[""{name=1, anchor=center, inner sep=0}, "{\mathrm{id}_{Gx}}"', "\shortmid"{marking}, from=2-1, to=2-2]
        \arrow["{\alpha_x}"', from=1-1, to=2-1]
        \arrow["{\alpha_x}", from=1-2, to=2-2]
        \arrow[Rightarrow, no head, from=2-1, to=3-1]
        \arrow[Rightarrow, no head, from=2-2, to=3-2]
        \arrow[""{name=2, anchor=center, inner sep=0}, "{G \mathrm{id}_x}"', "\shortmid"{marking}, from=3-1, to=3-2]
        \arrow["{\mathrm{id}_{\alpha_x}}"{description}, draw=none, from=0, to=1]
        \arrow["{G_x}"{description}, draw=none, from=1, to=2]
      \end{tikzcd}
    \end{equation*}
    \item Functorality of naturality comparisons: for every consecutive pair of
      arrows $x \xto{f} y \xto{g} z$ in $\dbl{D}$,
    \begin{equation*}
      \begin{tikzcd}
        Fx & Fx & Fx \\
        Gx & Fy & Fy \\
        Gy & Gy & Fz \\
        Gz & Gz & Gz \\
        Gz && Gz
        \arrow[""{name=0, anchor=center, inner sep=0}, "{\mathrm{id}_{Fx}}", "\shortmid"{marking}, from=1-1, to=1-2]
        \arrow[""{name=1, anchor=center, inner sep=0}, "{\mathrm{id}_{Fx}}", from=1-2, to=1-3]
        \arrow["{\alpha_x}"', from=1-1, to=2-1]
        \arrow["Gf"', from=2-1, to=3-1]
        \arrow["Ff"', from=1-2, to=2-2]
        \arrow["{\alpha_y}"', from=2-2, to=3-2]
        \arrow[""{name=2, anchor=center, inner sep=0}, "{G \mathrm{id}_y}", "\shortmid"{marking}, from=3-1, to=3-2]
        \arrow[""{name=3, anchor=center, inner sep=0}, "{\mathrm{id}_{Fy}}"', "\shortmid"{marking}, from=2-2, to=2-3]
        \arrow["Ff", from=1-3, to=2-3]
        \arrow["Gg", from=3-2, to=4-2]
        \arrow["Gg"', from=3-1, to=4-1]
        \arrow[""{name=4, anchor=center, inner sep=0}, "{G \mathrm{id}_z}"', "\shortmid"{marking}, from=4-1, to=4-2]
        \arrow["Fg", from=2-3, to=3-3]
        \arrow["{\alpha_z}", from=3-3, to=4-3]
        \arrow[""{name=5, anchor=center, inner sep=0}, "{G \mathrm{id}_z}"', "\shortmid"{marking}, from=4-2, to=4-3]
        \arrow[Rightarrow, no head, from=4-1, to=5-1]
        \arrow[Rightarrow, no head, from=4-3, to=5-3]
        \arrow[""{name=6, anchor=center, inner sep=0}, "{G \mathrm{id}_z}"', "\shortmid"{marking}, from=5-1, to=5-3]
        \arrow["{\mathrm{id}_{Ff}}"{description}, draw=none, from=1, to=3]
        \arrow["{\alpha_f}"{description}, draw=none, from=0, to=2]
        \arrow["{G \mathrm{id}_g}"{description}, draw=none, from=2, to=4]
        \arrow["{\alpha_g}"{description}, draw=none, from=3, to=5]
        \arrow["{G_{z,z}}"{description}, draw=none, from=4-2, to=6]
      \end{tikzcd}
      \quad=\quad
      \begin{tikzcd}
        Fx & Fx \\
        Gx & Fz \\
        Gz & Gz
        \arrow[""{name=0, anchor=center, inner sep=0}, "{\mathrm{id}_{Fx}}", "\shortmid"{marking}, from=1-1, to=1-2]
        \arrow["{\alpha_x}"', from=1-1, to=2-1]
        \arrow["{G(f \cdot g)}"', from=2-1, to=3-1]
        \arrow["{F(f \cdot g)}", from=1-2, to=2-2]
        \arrow["{\alpha_z}", from=2-2, to=3-2]
        \arrow[""{name=1, anchor=center, inner sep=0}, "{G \mathrm{id}_z}"', "\shortmid"{marking}, from=3-1, to=3-2]
        \arrow["{\alpha_{f \cdot g}}"{description}, draw=none, from=0, to=1]
      \end{tikzcd}
    \end{equation*}
    and for every object $x \in \dbl{D}$,
    \begin{equation*}
      \begin{tikzcd}
        Fx & Fx \\
        Gx & Gx \\
        Gx & Gx
        \arrow[""{name=0, anchor=center, inner sep=0}, "{\mathrm{id}_{Fx}}", "\shortmid"{marking}, from=1-1, to=1-2]
        \arrow["{\alpha_x}"', from=1-1, to=2-1]
        \arrow["{\alpha_x}", from=1-2, to=2-2]
        \arrow[""{name=1, anchor=center, inner sep=0}, "{\mathrm{id}_{Gx}}"', "\shortmid"{marking}, from=2-1, to=2-2]
        \arrow[Rightarrow, no head, from=2-1, to=3-1]
        \arrow[Rightarrow, no head, from=2-2, to=3-2]
        \arrow[""{name=2, anchor=center, inner sep=0}, "{G \mathrm{id}_x}"', "\shortmid"{marking}, from=3-1, to=3-2]
        \arrow["{\mathrm{id}_{\alpha_x}}"{description}, draw=none, from=0, to=1]
        \arrow["{G_x}"{description}, draw=none, from=1, to=2]
      \end{tikzcd}
      \quad=\quad
      \begin{tikzcd}
        Fx & Fx \\
        Gx & Fx \\
        Gx & Gx
        \arrow[""{name=0, anchor=center, inner sep=0}, "{\mathrm{id}_{Fx}}", "\shortmid"{marking}, from=1-1, to=1-2]
        \arrow["{G(1_x)}"', from=2-1, to=3-1]
        \arrow["{\alpha_x}"', from=1-1, to=2-1]
        \arrow["{F(1_x)}", from=1-2, to=2-2]
        \arrow["{\alpha_x}", from=2-2, to=3-2]
        \arrow[""{name=1, anchor=center, inner sep=0}, "{G \mathrm{id}_x}"', "\shortmid"{marking}, from=3-1, to=3-2]
        \arrow["{\alpha_{1_x}}"{description}, draw=none, from=0, to=1]
      \end{tikzcd}.
    \end{equation*}
  \end{itemize}

  The lax transformation $\alpha: F \To G$ is \define{pseudo} if for every arrow
  $f: x \to y$ in $\dbl{D}$, there exists a cell in $\dbl{E}$ (necessarily
  unique)
  \begin{equation} \label{eq:pseudo-transformation-inverse}
    \begin{tikzcd}
      Fx & Fx \\
      Fy & Gx \\
      Gy & Gy
      \arrow[""{name=0, anchor=center, inner sep=0}, "{\id_{Fx}}", "\shortmid"{marking}, from=1-1, to=1-2]
      \arrow["{\alpha_x}", from=1-2, to=2-2]
      \arrow["Gf", from=2-2, to=3-2]
      \arrow["Ff"', from=1-1, to=2-1]
      \arrow["{\alpha_y}"', from=2-1, to=3-1]
      \arrow[""{name=1, anchor=center, inner sep=0}, "{G \id_y}"', "\shortmid"{marking}, from=3-1, to=3-2]
      \arrow["{\alpha_f^{-1}}"{description}, draw=none, from=0, to=1]
    \end{tikzcd}
  \end{equation}
  such that the equations
  \begin{equation*}
    (\alpha_f \odot \alpha_f^{-1}) \cdot G_{y,y} = \id_{\alpha_x \cdot Gf} \cdot G_y
    \qquad\text{and}\qquad
    (\alpha_f^{-1} \odot \alpha_f) \cdot G_{y,y} = \id_{Ff \cdot \alpha_y} \cdot G_y
  \end{equation*}
  hold. Finally, the lax transformation $\alpha: F \To G$ is \define{strict} if,
  for every arrow $f: x \to y$ in $\dbl{D}$, the naturaliy equation
  $\alpha_x \cdot Gf = Ff \cdot \alpha_y$ holds and the naturality comparison at
  $f$ factorizes as
  \begin{equation} \label{eq:strict-transformation-as-lax}
    \begin{tikzcd}
      Fx & Fx \\
      Gx & Fy \\
      Gy & Gy
      \arrow[""{name=0, anchor=center, inner sep=0}, "{\id_{Fx}}", "\shortmid"{marking}, from=1-1, to=1-2]
      \arrow["{\alpha_x}"', from=1-1, to=2-1]
      \arrow["Gf"', from=2-1, to=3-1]
      \arrow["Ff", from=1-2, to=2-2]
      \arrow["{\alpha_y}", from=2-2, to=3-2]
      \arrow[""{name=1, anchor=center, inner sep=0}, "{G \id_y}"', "\shortmid"{marking}, from=3-1, to=3-2]
      \arrow["{\alpha_f}"{description}, draw=none, from=0, to=1]
    \end{tikzcd}
    \quad=\quad
    \begin{tikzcd}[row sep=small]
      Fx & Fx \\
      Gx & Fy \\
      Gy & Gy \\
      Gy & Gy
      \arrow["{\alpha_x}"', from=1-1, to=2-1]
      \arrow["Gf"', from=2-1, to=3-1]
      \arrow[""{name=0, anchor=center, inner sep=0}, "{\id_{Fx}}", "\shortmid"{marking}, from=1-1, to=1-2]
      \arrow["Ff", from=1-2, to=2-2]
      \arrow["{\alpha_y}", from=2-2, to=3-2]
      \arrow[Rightarrow, no head, from=3-1, to=4-1]
      \arrow[Rightarrow, no head, from=3-2, to=4-2]
      \arrow[""{name=1, anchor=center, inner sep=0}, "{G \id_y}"', "\shortmid"{marking}, from=4-1, to=4-2]
      \arrow[""{name=2, anchor=center, inner sep=0}, "{\id_{Gy}}", "\shortmid"{marking}, from=3-1, to=3-2]
      \arrow["{G_y}"{description}, draw=none, from=2, to=1]
      \arrow["\id"{description}, draw=none, from=0, to=2]
    \end{tikzcd}.
  \end{equation}
\end{definition}

\define{Oplax} transformations between lax double functors are defined
similarly, except that the direction of the naturality comparisons is reversed.
Unless otherwise clarified, all results in this paper stated for lax
transformations apply equally to oplax transformations.

\begin{remark}[Special cases]
  \label{remark:special-cases-of-lax-transfs}
  A lax natural transformation between lax double functors that is strict in our 
  sense is equivalent to a natural transformation between lax functors in the 
  usual sense. However, a pseudo natural transformation in our sense is more 
  general than a pseudo transformation in Grandis' sense 
  \cite[Definition 3.8.1]{grandis2019}.
\end{remark}

\begin{remark}[Lax transformations between 2-functors]
  \label{remark:lax-trans-dbl-func-reduce-to-lax-trans-of-2-func}
  Given two \emph{unitary} lax double functors
  $F,G\colon \dbl{D} \rightrightarrows \dbl{E}$, any lax transformation
  $\alpha\colon F \To G$ induces a lax transformation
  $\VerTwoCat(\alpha) \colon \VerTwoCat(F) \To \VerTwoCat(F)$ between the
  vertical 2-functors
  $\VerTwoCat(F), \VerTwoCat(G) \colon \VerTwoCat(\dbl{D}) \rightrightarrows \VerTwoCat(\dbl{E})$
  that agrees with the standard notion of lax natural transformation in
  2-category theory \cite[\S{4.2}]{johnson2021}. Conversely, any lax
  transformation $\beta\colon G\To H$ of 2-functors
  $G,H\colon \bicat{A}\rightrightarrows\bicat{B}$ induces a lax transformation
  $\VerDbl(\beta)\colon \VerDbl(H)\To\VerDbl(K)$ of unitary lax double functors
  $\VerDbl(H),\VerDbl(K)\colon \VerDbl(\bicat{A})\rightrightarrows \VerDbl(\bicat{B})$
  between horizontally trivial double categories. These constructions also work
  if everywhere \emph{unitary} is replaced with \emph{normal}, provided that
  unitor isomorphisms are inserted where appropriate.
\end{remark}

Lax functors and lax transformations between two fixed double categories form a
category.

\begin{proposition}[Category of lax functors] \label{prop:lax-functor-category}
  For any double categories $\dbl{D}$ and $\dbl{E}$, lax double functors
  $\dbl{D} \to \dbl{E}$ and lax natural transformations between them form a
  category $\LaxFunOne_\ell(\dbl{D},\dbl{E})$. In this category, the composite
  of lax transformations $\alpha: F \To G$ and $\beta: G \To H$ has components
  \begin{equation*}
    (\alpha \cdot \beta)_x \coloneqq \alpha_x \cdot \beta_x, \qquad
    (\alpha \cdot \beta)_m \coloneqq \alpha_m \cdot \beta_m,
  \end{equation*}
  and naturality comparisons
  \begin{equation*}
    \begin{tikzcd}
      Fx & Fx \\
      Hx & Fy \\
      Hy & Hy
      \arrow["{(\alpha \cdot \beta)_x}"', from=1-1, to=2-1]
      \arrow["Hf"', from=2-1, to=3-1]
      \arrow[""{name=0, anchor=center, inner sep=0}, "{\mathrm{id}_{Fx}}", "\shortmid"{marking}, from=1-1, to=1-2]
      \arrow["Ff", from=1-2, to=2-2]
      \arrow["{(\alpha \cdot \beta)_y}", from=2-2, to=3-2]
      \arrow[""{name=1, anchor=center, inner sep=0}, "{H \mathrm{id}_y}"', "\shortmid"{marking}, from=3-1, to=3-2]
      \arrow["{(\alpha \cdot \beta)_f}"{description}, draw=none, from=0, to=1]
    \end{tikzcd}
    \quad\coloneqq\quad
    \begin{tikzcd}
      Fx & Fx & Fx \\
      Gx & Gx & Fy \\
      {H x} & Gy & Gy \\
      Hy & Hy & Hy \\
      Hy && Hy
      \arrow[""{name=0, anchor=center, inner sep=0}, "{\mathrm{id}_{Fx}}", "\shortmid"{marking}, from=1-1, to=1-2]
      \arrow[""{name=1, anchor=center, inner sep=0}, "{\mathrm{id}_{Fx}}", "\shortmid"{marking}, from=1-2, to=1-3]
      \arrow["{\alpha_x}", from=1-2, to=2-2]
      \arrow["{\alpha_x}"', from=1-1, to=2-1]
      \arrow[""{name=2, anchor=center, inner sep=0}, "{\mathrm{id}_{Gx}}"', "\shortmid"{marking}, from=2-1, to=2-2]
      \arrow["{\beta_x}"', from=2-1, to=3-1]
      \arrow["Hf"', from=3-1, to=4-1]
      \arrow["{\beta_y}"', from=3-2, to=4-2]
      \arrow["Gf"', from=2-2, to=3-2]
      \arrow[""{name=3, anchor=center, inner sep=0}, "{H \mathrm{id}_y}"', "\shortmid"{marking}, from=4-1, to=4-2]
      \arrow["Ff", from=1-3, to=2-3]
      \arrow["{\alpha_y}", from=2-3, to=3-3]
      \arrow[""{name=4, anchor=center, inner sep=0}, "{G \mathrm{id}_y}", "\shortmid"{marking}, from=3-2, to=3-3]
      \arrow[""{name=5, anchor=center, inner sep=0}, "{H \mathrm{id}_y}"', from=4-2, to=4-3]
      \arrow["{\beta_y}", from=3-3, to=4-3]
      \arrow[""{name=6, anchor=center, inner sep=0}, "{H \mathrm{id}_y}"', "\shortmid"{marking}, from=5-1, to=5-3]
      \arrow[Rightarrow, no head, from=4-1, to=5-1]
      \arrow[Rightarrow, no head, from=4-3, to=5-3]
      \arrow["{\mathrm{id}_{\alpha_x}}"{description}, draw=none, from=0, to=2]
      \arrow["{\beta_f}"{description}, draw=none, from=2, to=3]
      \arrow["{\beta_{\mathrm{id}_y}}"{description}, draw=none, from=4, to=5]
      \arrow["{\alpha_f}"{description}, draw=none, from=1, to=4]
      \arrow["{H_{y,y}}"{description}, draw=none, from=4-2, to=6]
    \end{tikzcd}.
  \end{equation*}
  The identity lax transformation $1_F: F \To F$ has components
  $(1_F)_x \coloneqq 1_{Fx}$ and $(1_F)_m \coloneqq 1_{Fm}$ and naturality
  comparisons
  \begin{equation*}
    \begin{tikzcd}
      Fx & Fx \\
      Fx & Fy \\
      Fy & Fy
      \arrow[Rightarrow, no head, from=1-1, to=2-1]
      \arrow["Ff"', from=2-1, to=3-1]
      \arrow[""{name=0, anchor=center, inner sep=0}, "{\mathrm{id}_{Fx}}", "\shortmid"{marking}, from=1-1, to=1-2]
      \arrow["Fy", from=1-2, to=2-2]
      \arrow[Rightarrow, no head, from=2-2, to=3-2]
      \arrow[""{name=1, anchor=center, inner sep=0}, "{F \mathrm{id}_y}"', "\shortmid"{marking}, from=3-1, to=3-2]
      \arrow["{(1_F)_f}"{description}, draw=none, from=0, to=1]
    \end{tikzcd}
    \quad\coloneqq\quad
    \begin{tikzcd}
      Fx & Fx \\
      Fy & Fy \\
      Fy & Fy
      \arrow[""{name=0, anchor=center, inner sep=0}, "{\mathrm{id}_{Fx}}", "\shortmid"{marking}, from=1-1, to=1-2]
      \arrow[""{name=1, anchor=center, inner sep=0}, "{\mathrm{id}_{Fy}}"', "\shortmid"{marking}, from=2-1, to=2-2]
      \arrow["Ff"', from=1-1, to=2-1]
      \arrow["Ff", from=1-2, to=2-2]
      \arrow[""{name=2, anchor=center, inner sep=0}, "{F \mathrm{id}_y}"', "\shortmid"{marking}, from=3-1, to=3-2]
      \arrow[Rightarrow, no head, from=2-1, to=3-1]
      \arrow[Rightarrow, no head, from=2-2, to=3-2]
      \arrow["{\mathrm{id}_{Ff}}"{description}, draw=none, from=0, to=1]
      \arrow["{F_y}"{description}, draw=none, from=1, to=2]
    \end{tikzcd}.
  \end{equation*}
\end{proposition}
\begin{proof}
  It is easy to show that the associativity and unitality laws hold, using the
  corresponding laws for the double category $\dbl{E}$ and the laxators of $F$.
  Proving that composites and identities obey the lax natural transformation
  axioms is a long series of calculations. We show what is perhaps the most
  involved one: that given two lax transformations $\alpha: F \To G$ and
  $\beta: G \To H$, the composite transformation $\alpha \beta: F \To H$ has
  naturality comparisons $(\alpha \beta)_f$ satisfying the naturality axiom.

  Fixing a cell $\stdInlineCell{\gamma}$ in $\dbl{D}$, we calculate
  \begin{align*}
    \begin{dblArray}{cc}
      \Block{2-1}{(\alpha \beta)_f} & F\gamma \\
      & (\alpha \beta)_n \\
      \Block{1-2}{H_{w,n}} &
    \end{dblArray} &=
    \begin{dblArray}{ccc}
      \id_{\alpha_x} & \Block{2-1}{\alpha_f} & F\gamma \\
      \Block{2-1}{\beta_f} & & \alpha_n \\
      & \beta_{\id_w} & \beta_n \\
      \Block{1-2}{H_{w,w}} & & 1 \\
      \Block{1-3}{H_{w,n}} & &
    \end{dblArray} =
    \begin{dblArray}{ccc}
      \id_{\alpha_x} & \Block{2-1}{\alpha_f} & F\gamma \\
      \Block{2-1}{\beta_f} & & \alpha_n \\
      & \beta_{\id_w} & \beta_n \\
      1 & \Block{1-2}{H_{w,n}} & \\
      \Block{1-3}{H_{w,n}} & &
    \end{dblArray} =
    \begin{dblArray}{ccc}
      \id_{\alpha_x} & \Block{2-1}{\alpha_f} & F\gamma \\
      \Block{2-1}{\beta_f} & & \alpha_n \\
      & \Block{1-2}{G_{w,n}} &  \\
      1 & \Block{1-2}{\beta_n} & \\
      \Block{1-3}{H_{w,n}} & &
    \end{dblArray} \\
    &= \begin{dblArray}{ccc}
      \id_{\alpha_x} & \alpha_m & \Block{2-1}{\alpha_g} \\
      \Block{2-1}{\beta_f} & G\gamma & \\
      & \Block{1-2}{G_{n,z}} &  \\
      1 & \Block{1-2}{\beta_n} & \\
      \Block{1-3}{H_{w,n}} & &
    \end{dblArray} =
    \begin{dblArray}{ccc}
      \id_{\alpha_x} & \alpha_m & \Block{2-1}{\alpha_g} \\
      \Block{2-1}{\beta_f} & G\gamma & \\
      & \beta_n & \beta_{\id_z} \\
      1 & \Block{1-2}{H_{n,z}} & \\
      \Block{1-3}{H_{w,n}} & &
    \end{dblArray} =
    \begin{dblArray}{ccc}
      \id_{\alpha_x} & \alpha_m & \Block{2-1}{\alpha_g} \\
      \Block{2-1}{\beta_f} & G\gamma & \\
      & \beta_n & \beta_{\id_z} \\
      \Block{1-2}{H_{w,n}} & & 1 \\
      \Block{1-3}{H_{n,z}} & &
    \end{dblArray} \\
    &= \begin{dblArray}{ccc}
      \alpha_m & \id_{\alpha_y} & \Block{2-1}{\alpha_g} \\
      \beta_m & \Block{2-1}{\beta_g} & \\
      H\gamma & & \beta_{\id_z} \\
      \Block{1-2}{H_{n,z}} & & 1 \\
      \Block{1-3}{H_{n,z}} & &
    \end{dblArray} =
    \begin{dblArray}{ccc}
      \alpha_m & \id_{\alpha_y} & \Block{2-1}{\alpha_g} \\
      \beta_m & \Block{2-1}{\beta_g} & \\
      H\gamma & & \beta_{\id_z} \\
      1 & \Block{1-2}{H_{z,z}} & \\
      \Block{1-3}{H_{n,z}} & &
    \end{dblArray} =
    \begin{dblArray}{cc}
      (\alpha \beta)_m & \Block{2-1}{(\alpha \beta)_g} \\
      H\gamma & \\
      \Block{1-2}{H_{n,z}} &
    \end{dblArray}\ . \qedhere
  \end{align*}
\end{proof}

The reduction to lax transformations between 2-functors noted in
\cref{remark:lax-trans-dbl-func-reduce-to-lax-trans-of-2-func} extends the
correspondence of \cref{cor:lax-double-fs-corresp-to-2-fs-upto-iso} to an
equivalence of categories. Setting the notation, for fixed 2-categories
$\bicat{A}$ and $\bicat{B}$, strict 2-functors $\bicat{A}\to\bicat{B}$ and lax
natural transformations between them form a category
$\TwoCatOne_\ell(\bicat{A},\bicat{B})$. The full subcategories of
$\LaxFunOne_\ell(\dbl{D},\dbl{E})$ spanned by \emph{normal} and \emph{unitary}
lax functors are denoted $\LaxFunOne_{\ell,n}(\dbl{D},\dbl{E})$ and
$\LaxFunOne_{\ell,u}(\dbl{D},\dbl{E})$, respectively.

\begin{corollary} \label{cor:characterize-unitary-models-on-veritical-dbl-cat-I}
  For any 2-category $\bicat{A}$ and double category $\dbl{E}$, there is an
  isomorphism of categories
  \begin{equation*}
    \LaxFunOne_{\ell, u}(\VerDbl(\bicat{A}),\dbl{E}) \cong
      \TwoCatOne_\ell(\bicat{A},\VerTwoCat(\dbl{E}))
  \end{equation*}
  and an equivalence of categories
  \begin{equation*}
    \LaxFunOne_{\ell, n}(\VerDbl(\bicat{A}),\dbl{E}) \simeq
      \TwoCatOne_\ell(\bicat{A},\VerTwoCat(\dbl{E})).
  \end{equation*}
\end{corollary}
  
Having defined the category of lax double functors and lax transformations, our
next goal is to extend the correspondence of
\cref{prop:unitalizationoflaxdoublefunctor} to an isomorphism of categories.
That is, when $\dbl{D}$ is a double category and $\dbl{E}$ is an equipment with
local coequalizers, the assignment \emph{postcomposition with $\epsilon$}
extends to a functor
\begin{equation*}
  \LaxFunOne_{\ell, u}(\dbl{D},\Module{\dbl{E}}) \to \LaxFunOne_\ell(\dbl{D},\dbl{E}),
  \qquad H\mapsto \epsilon H.
\end{equation*}
We shall see that this functor is in fact an isomorphism. We know already that
it is a bijection on objects. It is worth remarking that this is something of an
extraordinary situation, as the assignment on morphisms, namely, sending a lax
transformation $\alpha$ to its whiskering $\epsilon \ast\alpha$ does not always
result in a well-defined lax transformation \cite{shulman2009}. For this reason,
we shall recount this development in some detail.

\begin{lemma} \label{lem:whiskering-with-mod-counit}
  When $\dbl{D}$ is a double category and $\dbl{E}$ is an equipment with local
  coequalizers, there is a functor
  \begin{equation*}
    \epsilon \circ (-): \LaxFunOne_{\ell, u}(\dbl{D},\Module{\dbl{E}}) \to \LaxFunOne_\ell(\dbl{D},\dbl{E}),
    \qquad H\mapsto \epsilon H,\quad \alpha\mapsto \epsilon\ast\alpha,
  \end{equation*}
  extending the mapping on objects from
  \cref{prop:unitalizationoflaxdoublefunctor}.
\end{lemma}
\begin{proof}
  Suppose that $H,K\colon \dbl{D} \rightrightarrows \Module{\dbl{E}}$ are
  unitary lax functors and $\alpha\colon H \To K$ is a lax natural
  transformation. We need extra notation to describe the proposed whiskered
  transformation $\epsilon\ast \alpha$. For each object $x$ in $\dbl{D}$, the
  image $Hx$ has the structure of a category object, which we take to name the
  proarrow, so that $H(x)_0$ denotes the underlying object and
  $Hx\colon H(x)_0\proto H(x)_0$ is the associated proarrow. We also take $\mu$
  and $\upsilon$ as generic symbols for the structure cells. Likewise, given a
  morphism $f\colon x\to y$ of $\dbl{D}$, the image $Hf$ is a functor with
  object map $H(f)_0\colon H(x)_0\to H(x)_0$ and with morphism map
  $Hf\colon Hx \To Hx$, whose external source and target are $H(f)_0$.

  The extended action of the mapping $H\mapsto \epsilon H$ essentially takes the
  underlying structure with only a few modifications. For the action on lax
  transformations, note that a lax transformation $\alpha$ has components
  $\alpha_x$ and $\alpha_m$ and naturality comparisons $\alpha_f$. The first are
  functors between category objects and the latter two are both maps between
  profunctor objects. The proposed transformation $\epsilon\ast\alpha$ is
  defined to have components
  \begin{equation*}
    (\epsilon\ast \alpha)_x \coloneqq (\alpha_x)_0 \colon H(x)_0\to K(x)_0,
    \qquad\qquad
    (\epsilon\ast \alpha)_m \coloneqq \alpha_m,
  \end{equation*}
  and naturality comparisons $(\epsilon\ast \alpha)_f$ given by
  \begin{equation*}
    \begin{tikzcd}
      {H(x)_0} & {H(x)_0} \\
      {H(x)_0} & {H(x)_0} \\
      {K(y)_0} & {K(y)_0}
      \arrow[""{name=0, anchor=center, inner sep=0}, "Hx", "\shortmid"{marking}, from=2-1, to=2-2]
      \arrow[from=2-2, to=3-2]
      \arrow[from=2-1, to=3-1]
      \arrow[""{name=1, anchor=center, inner sep=0}, "{K\id_y}"', "\shortmid"{marking}, from=3-1, to=3-2]
      \arrow[Rightarrow, no head, from=1-1, to=2-1]
      \arrow[""{name=2, anchor=center, inner sep=0}, "{\id_{H(x)_0}}", "\shortmid"{marking}, from=1-1, to=1-2]
      \arrow[Rightarrow, no head, from=1-2, to=2-2]
      \arrow["{\upsilon_{Hx}}"{description, pos=0.4}, draw=none, from=2, to=0]
      \arrow["{\alpha_f}"{description}, draw=none, from=0, to=1]
    \end{tikzcd},
  \end{equation*}
  that is, precomposing $\alpha_f$ with the unit cell $\upsilon_{Hx}$ coming 
  with the category $Hx$. Note that this is well-typed due to the fact 
  that, since $H$ is unitary, the proarrow domain of $\alpha_f$ is the unit on $Hx$, namely, $Hx$ viewed as a bimodule over itself.

  Well-definition is essentially trivial owing to the fact that the assignment
  starts with a lax transformation and merely precomposes with a well-behaved
  cell. That the correspondence is functorial involves verifying composition and
  identities are preserved as they are defined in
  \cref{prop:lax-functor-category}. But this just involves some straightforward
  computations using the unit conditions for lax functors and lax
  transformations.
\end{proof}

Now, to extend the inverse correspondence of
\cref{prop:unitalizationoflaxdoublefunctor}, consider a lax transformation
$\alpha\colon F\To G$ between lax double functors
$F,G\colon \dbl{D}\rightrightarrows \dbl{E}$. We shall define components of a
lax transformation $\bar\alpha$ between \emph{unitary} lax double functors
$\Module{F}\eta \To \Module{G}\eta$. The data is as follows. First, for objects
and proarrows, take
\begin{enumerate}[noitemsep]
  \item $\bar \alpha_x \coloneqq (\alpha_x, \alpha_{\id_x})$
  \item $\bar \alpha_m \coloneqq \alpha_m$.
\end{enumerate}
That is, the proarrow component is the given one, while the component at an
object takes the underlying arrow component $\alpha_x\colon Fx\to Gx$ along with
the cell $\alpha_{\id_x}$ of the form
\begin{equation*}
  \begin{tikzcd}
    Fx & Fx \\
    Gx & Gx
    \arrow[""{name=0, anchor=center, inner sep=0}, "{F\id_x}", "\shortmid"{marking}, from=1-1, to=1-2]
    \arrow["{\alpha_x}", from=1-2, to=2-2]
    \arrow["{\alpha_x}"', from=1-1, to=2-1]
    \arrow[""{name=1, anchor=center, inner sep=0}, "{G\id_x}"', "\shortmid"{marking}, from=2-1, to=2-2]
    \arrow["{\alpha_{\id_x}}"{description}, draw=none, from=0, to=1]
  \end{tikzcd}.
\end{equation*}
This defines a functor between category objects. One can also check that each
cell $\bar \alpha_m$ is a morphism of bimodules over the images of trivial
categories.

The naturality comparison $\bar\alpha_f$ is more involved to construct. If
$f\colon x\to y$ is an arrow in $\dbl{D}$, define a proposed laxity cell
$\bar\alpha_f$ by either side of the equation
\begin{equation*}
  \begin{tikzcd}
    Fx & Fx & Fx \\
    Gx & Fy & Fy \\
    Gy & Gy & Gy \\
    Gy && Gy
    \arrow[""{name=0, anchor=center, inner sep=0}, "{\mathrm{id}_{Fx}}", "\shortmid"{marking}, from=1-1, to=1-2]
    \arrow["{\alpha_x}"', from=1-1, to=2-1]
    \arrow["Gf"', from=2-1, to=3-1]
    \arrow["Ff"', from=1-2, to=2-2]
    \arrow["{\alpha_y}"', from=2-2, to=3-2]
    \arrow[""{name=1, anchor=center, inner sep=0}, "{G \id_y}"', "\shortmid"{marking}, from=3-1, to=3-2]
    \arrow["Ff", from=1-3, to=2-3]
    \arrow[""{name=2, anchor=center, inner sep=0}, "{F\id_x}", "\shortmid"{marking}, from=1-2, to=1-3]
    \arrow[""{name=3, anchor=center, inner sep=0}, "{F\id_y}"', "\shortmid"{marking}, from=2-2, to=2-3]
    \arrow[""{name=4, anchor=center, inner sep=0}, "{G\id_y}"', "\shortmid"{marking}, from=3-2, to=3-3]
    \arrow["{\alpha_y}", from=2-3, to=3-3]
    \arrow[Rightarrow, no head, from=3-1, to=4-1]
    \arrow[Rightarrow, no head, from=3-3, to=4-3]
    \arrow[""{name=5, anchor=center, inner sep=0}, "{G\id_y}"', "\shortmid"{marking}, from=4-1, to=4-3]
    \arrow["{\alpha_f}"{description}, draw=none, from=0, to=1]
    \arrow["{F\id_f}"{description}, draw=none, from=2, to=3]
    \arrow["{G_{y,y}}"{description}, draw=none, from=3-2, to=5]
    \arrow["{\alpha_{\id_y}}"{description}, draw=none, from=3, to=4]
  \end{tikzcd}
  \quad=\quad
  \begin{tikzcd}
    Fx & Fx & Fx \\
    Gx & Gx & Fy \\
    Gy & Gy & Gy \\
    Gy && Gy
    \arrow[""{name=0, anchor=center, inner sep=0}, "{F\id_x}", "\shortmid"{marking}, from=1-1, to=1-2]
    \arrow[""{name=1, anchor=center, inner sep=0}, "{\id_{Fx}}", "\shortmid"{marking}, from=1-2, to=1-3]
    \arrow["{\alpha_x}"', from=1-1, to=2-1]
    \arrow[""{name=2, anchor=center, inner sep=0}, "\shortmid"{marking}, from=2-1, to=2-2]
    \arrow["Gf"', from=2-1, to=3-1]
    \arrow[""{name=3, anchor=center, inner sep=0}, "{G\id_y}"', "\shortmid"{marking}, from=3-1, to=3-2]
    \arrow["Gf", from=2-2, to=3-2]
    \arrow["{\alpha_x}", from=1-2, to=2-2]
    \arrow[""{name=4, anchor=center, inner sep=0}, "{G\id_y}"', "\shortmid"{marking}, from=3-2, to=3-3]
    \arrow[Rightarrow, no head, from=3-1, to=4-1]
    \arrow[""{name=5, anchor=center, inner sep=0}, "{G\id_y}"', "\shortmid"{marking}, from=4-1, to=4-3]
    \arrow[Rightarrow, no head, from=3-3, to=4-3]
    \arrow["Ff", from=1-3, to=2-3]
    \arrow["{\alpha_y}", from=2-3, to=3-3]
    \arrow["{G_{y,y}}"{description}, draw=none, from=3-2, to=5]
    \arrow["{\alpha_f}"{description}, draw=none, from=1, to=4]
    \arrow["{G\id_f}"{description}, draw=none, from=2, to=3]
    \arrow["{\alpha_{\id_x}}"{description}, draw=none, from=0, to=2]
  \end{tikzcd}
\end{equation*}
given by instantiating the cell naturality axiom of
\cref{def:lax-transformation} at the external identity cell $\id_f$. Our claim
is that with this definition, $\bar \alpha$ is indeed a lax natural
transformation in the sense of \cref{def:lax-transformation}. To see this, it is
helpful to have a preliminary result.

\begin{lemma} \label{lemma:lax-transf-data-bijection}
  Let $\alpha: F \To G$ be a lax transformation between lax double functors
  $F, G: \dbl{D} \rightrightarrows \dbl{E}$. For any arrow $f\colon x\to y$ in
  $\dbl{D}$, there is an equality
  \begin{equation*}
    \begin{tikzcd}
      Fx & Fx \\
      Fx & Fx \\
      Gy & Gy
      \arrow[""{name=0, anchor=center, inner sep=0}, "\shortmid"{marking}, from=2-1, to=2-2]
      \arrow[from=2-2, to=3-2]
      \arrow[from=2-1, to=3-1]
      \arrow[""{name=1, anchor=center, inner sep=0}, "{G\id_y}"', "\shortmid"{marking}, from=3-1, to=3-2]
      \arrow[Rightarrow, no head, from=1-1, to=2-1]
      \arrow[""{name=2, anchor=center, inner sep=0}, "{\id_x}", "\shortmid"{marking}, from=1-1, to=1-2]
      \arrow[Rightarrow, no head, from=1-2, to=2-2]
      \arrow["{\bar\alpha_f}"{description}, draw=none, from=0, to=1]
      \arrow["{F_x}"{description}, draw=none, from=2, to=0]
    \end{tikzcd}
    \quad=\quad
    \begin{tikzcd}
      Fx & Fx \\
      Gx & Fy \\
      Gy & Gy
      \arrow["{\alpha_X}"', from=1-1, to=2-1]
      \arrow["Gf"', from=2-1, to=3-1]
      \arrow[""{name=0, anchor=center, inner sep=0}, "{\id_x}", "\shortmid"{marking}, from=1-1, to=1-2]
      \arrow["Ff", from=1-2, to=2-2]
      \arrow["{\alpha_y}", from=2-2, to=3-2]
      \arrow[""{name=1, anchor=center, inner sep=0}, "{G\id_y}"', "\shortmid"{marking}, from=3-1, to=3-2]
      \arrow["{\alpha_f}"{description}, draw=none, from=0, to=1]
    \end{tikzcd}
  \end{equation*}
  when $\bar\alpha_f$ is defined as above. Likewise, if $\beta\colon H\To K$ is
  a lax transformation of module-valued, unitary lax functors
  $H,K\colon \dbl{D}\rightrightarrows \Module{\dbl{E}}$, then the comparison
  $\beta_f$ satisfies
  \begin{equation*}
    \begin{tikzcd}
      {H(x)_0} & {H(x)_0} & {H(x)_0} \\
      {H(x)_0} & {H(x)_0} & {H(x)_0} \\
      {K(x)_0} & {H(y)_0} & {H(y)_0} \\
      {K(y)_0} & {K(y)_0} & {K(y)_0} \\
      {K(y)_0} && {K(y)_0}
      \arrow[""{name=0, anchor=center, inner sep=0}, "Hx", "\shortmid"{marking}, from=2-1, to=2-2]
      \arrow["{\alpha_x}"', from=2-1, to=3-1]
      \arrow["Kf"', from=3-1, to=4-1]
      \arrow["Hf"', from=2-2, to=3-2]
      \arrow["{\alpha_y}"', from=3-2, to=4-2]
      \arrow[""{name=1, anchor=center, inner sep=0}, "Ky"', "\shortmid"{marking}, from=4-1, to=4-2]
      \arrow["Hf", from=2-3, to=3-3]
      \arrow[""{name=2, anchor=center, inner sep=0}, "Hx", "\shortmid"{marking}, from=2-2, to=2-3]
      \arrow[""{name=3, anchor=center, inner sep=0}, "Hx"', "\shortmid"{marking}, from=3-2, to=3-3]
      \arrow[""{name=4, anchor=center, inner sep=0}, "Ky"', "\shortmid"{marking}, from=4-2, to=4-3]
      \arrow["{\alpha_y}", from=3-3, to=4-3]
      \arrow[""{name=5, anchor=center, inner sep=0}, Rightarrow, no head, from=4-1, to=5-1]
      \arrow[""{name=6, anchor=center, inner sep=0}, Rightarrow, no head, from=4-3, to=5-3]
      \arrow["Ky"', "\shortmid"{marking}, from=5-1, to=5-3]
      \arrow[Rightarrow, no head, from=1-1, to=2-1]
      \arrow[""{name=7, anchor=center, inner sep=0}, "{\id_{H(x)_0}}", "\shortmid"{marking}, from=1-1, to=1-2]
      \arrow[Rightarrow, no head, from=1-2, to=2-2]
      \arrow[""{name=8, anchor=center, inner sep=0}, "Hx", "\shortmid"{marking}, from=1-2, to=1-3]
      \arrow[Rightarrow, no head, from=1-3, to=2-3]
      \arrow["{K_{y,y}}"{description}, draw=none, from=5, to=6]
      \arrow["1"{description}, draw=none, from=8, to=2]
      \arrow["{\beta_f}"{description}, draw=none, from=0, to=1]
      \arrow["{\upsilon_{Hx}}"{description}, draw=none, from=7, to=0]
      \arrow["{\beta_{\id_y}}"{description}, draw=none, from=3, to=4]
      \arrow["{H\id_f}"{description}, draw=none, from=2, to=3]
    \end{tikzcd}
    \quad=\quad
    \begin{tikzcd}
        {H(x)_0} & {H(x)_0} \\
        {K(x)_0} & {H(y)_0} \\
        {K(y)_0} & {K(y)_0}
        \arrow[""{name=0, anchor=center, inner sep=0}, "Hx", "\shortmid"{marking}, from=1-1, to=1-2]
        \arrow["{\alpha_x}"', from=1-1, to=2-1]
        \arrow["Kf"', from=2-1, to=3-1]
        \arrow["Hf", from=1-2, to=2-2]
        \arrow["{\alpha_y}", from=2-2, to=3-2]
        \arrow[""{name=1, anchor=center, inner sep=0}, "Ky"', "\shortmid"{marking}, from=3-1, to=3-2]
        \arrow["{\beta_f}"{description}, draw=none, from=0, to=1]
    \end{tikzcd}.
  \end{equation*}
\end{lemma}
\begin{proof}
  This proof follows closely that of the bijection in
  \cref{prop:natural-transformation-in-dbl}. On the one hand, calculate that
  \begin{equation*}
    \begin{dblArray}{cc}
      \Block{1-2}{F_x} & \\
      \Block{1-2}{\cong} & \\
      \alpha_{\id_x} & \Block{2-1}{\alpha_f} \\
      G\id_f &  \\
      \Block{1-2}{G_{y,y}} &
    \end{dblArray} =
    \begin{dblArray}{cc}
      \Block{1-2}{\cong} & \\
      F_x & 1 \\
      \alpha_{\id_x} & \Block{2-1}{\alpha_f} \\
      G\id_f &  \\
      \Block{1-2}{G_{y,y}} &
    \end{dblArray} =
    \begin{dblArray}{cc}
      \Block{1-2}{\cong} & \\
      \Block{2-1}{\alpha_{1_x}} & 1 \\
      & \Block{2-1}{\alpha_f} \\
      G\id_f &  \\
      \Block{1-2}{G_{y,y}} &
    \end{dblArray} =
    \begin{dblArray}{c}
      \alpha_f
    \end{dblArray}
  \end{equation*}
  using first unitality in \cref{def:lax-functor}, then both of the unit
  conditions in \cref{def:lax-transformation}, and finally the functoriality of
  naturality comparisons in \cref{def:lax-transformation} with one of the two
  morphisms an identity. On the other hand, calculate
  \begin{equation*}
    \begin{dblArray}{cc}
      \Block{1-2}{\cong} & \\
      \upsilon & 1 \\
      \Block{2-1}{\beta_f} & H\id_x  \\
      & \beta_{\id_x}  \\
      \Block{1-2}{K_{x,x}} &
    \end{dblArray} =
    \begin{dblArray}{cc}
      \Block{1-2}{\cong} & \\
      \upsilon & 1 \\
      \Block{1-2}{\mu} & \\
      \Block{1-2}{\beta_f} &
    \end{dblArray} =
    \begin{dblArray}{c}
      \beta_f
    \end{dblArray}
  \end{equation*}
  using the fact that $\beta_f$ is a transformation in $\dbl{E}$ and therefore
  satisfies the equivariance axiom of \cref{def:profunctor-object} applied to the 
  category structure $\upsilon$, $\mu$ of $Hx$.
\end{proof}

\begin{proposition}[Unitalization of lax functors, one-dimensional]
  \label{prop:unitalizationisoof1categories}
  If $\dbl{D}$ is a double category and $\dbl{E}$ is an equipment with local
  coequalizers, then the functor $H\mapsto \epsilon H$ has the functor
  $F\mapsto \Module{F}\eta$, $\alpha\mapsto\bar\alpha$ as its inverse, inducing
  an isomorphism of categories
  \begin{equation*}
    \LaxFunOne_{\ell,u}(\dbl{D},\Module{\dbl{E}}) \xrightarrow{\cong}
    \LaxFunOne_\ell(\dbl{D},\dbl{E})
  \end{equation*}
  that extends the bijection on objects from
  \cref{prop:unitalizationoflaxdoublefunctor}.
\end{proposition}
\begin{proof}
  We first show that $\bar\alpha$ as defined above is a well-defined lax natural
  transformation. We verify the cell naturality condition of
  \cref{def:lax-transformation}, omitting the other verifications. What needs to
  be shown for this condition is that for any cell
  \begin{equation*}
    \begin{tikzcd}
      x & z \\
      y & w
      \arrow[""{name=0, anchor=center, inner sep=0}, "m", "\shortmid"{marking}, from=1-1, to=1-2]
      \arrow["g", from=1-2, to=2-2]
      \arrow["f"', from=1-1, to=2-1]
      \arrow[""{name=1, anchor=center, inner sep=0}, "n"', "\shortmid"{marking}, from=2-1, to=2-2]
      \arrow["\theta"{description}, draw=none, from=0, to=1]
    \end{tikzcd}
  \end{equation*}
  we have an equality of composites
  \begin{equation*}
    \begin{tikzcd}
      Fx && Fz \\
      Gy && Gw \\
      Gy && Gw
      \arrow[""{name=0, anchor=center, inner sep=0}, "{F\id_x\otimes_{Fx}Fm}", "\shortmid"{marking}, from=1-1, to=1-3]
      \arrow[from=1-1, to=2-1]
      \arrow[""{name=1, anchor=center, inner sep=0}, "{G\id_y\otimes_{Gy}Gn}"', "\shortmid"{marking}, from=2-1, to=2-3]
      \arrow[from=1-3, to=2-3]
      \arrow[Rightarrow, no head, from=2-1, to=3-1]
      \arrow[""{name=2, anchor=center, inner sep=0}, "Gn"', "\shortmid"{marking}, from=3-1, to=3-3]
      \arrow[Rightarrow, no head, from=2-3, to=3-3]
      \arrow["{G_{y,n}^\otimes}"{description, pos=0.6}, draw=none, from=1, to=2]
      \arrow["{\bar\alpha_f\otimes\alpha_nF\theta}"{description}, draw=none, from=0, to=1]
    \end{tikzcd}
    \quad=\quad
    \begin{tikzcd}
      Fx && Fz \\
      Gy && Gw \\
      Gy && Gw
      \arrow[from=1-1, to=2-1]
      \arrow[Rightarrow, no head, from=2-1, to=3-1]
      \arrow[""{name=0, anchor=center, inner sep=0}, "Gn"', "\shortmid"{marking}, from=3-1, to=3-3]
      \arrow[""{name=1, anchor=center, inner sep=0}, "{Gn\otimes_{Gy}G\id_w}"', "\shortmid"{marking}, from=2-1, to=2-3]
      \arrow[Rightarrow, no head, from=2-3, to=3-3]
      \arrow[from=1-3, to=2-3]
      \arrow[""{name=2, anchor=center, inner sep=0}, "{Fm\otimes_{Fz}F\id_z}", "\shortmid"{marking}, from=1-1, to=1-3]
      \arrow["{G_{n,w}^\otimes}"{description, pos=0.6}, draw=none, from=1, to=0]
      \arrow["{G\theta\alpha_m\otimes \bar\alpha_g}"{description}, draw=none, from=2, to=1]
    \end{tikzcd}
  \end{equation*}
  in $\Module{\dbl{E}}$. We will show that this holds by precomposing with
  certain universal morphisms, using the previous lemma, and appealing to
  uniqueness. We have on the one hand that
  \begin{equation*}
    \begin{dblArray}{c}
      \lambda \\
      \cong \\
      \bar\alpha_f\otimes\alpha_nF\theta \\
      G_{y,n}^\otimes
    \end{dblArray} =
    \begin{dblArray}{c}
      F_x \odot 1\\
      \mathrm{coeq} \\
      \bar\alpha_f\otimes\alpha_nF\theta\\
      G_{y,n}^\otimes
    \end{dblArray} =
    \begin{dblArray}{cc}
      F_x & 1\\
      \Block{2-1}{\bar\alpha_f} & F\theta \\
      & \alpha_n \\
      \Block{1-2}{G_{y,n}}
    \end{dblArray} =
    \begin{dblArray}{cc}
      \Block{2-1}{\alpha_f} & F\theta \\
      & \alpha_n \\
      \Block{1-2}{G_{y,n}}
    \end{dblArray}
  \end{equation*}
  using \cref{lemma:lax-transf-data-bijection} for the last equality. Here the
  isomorphism is the canonical right unit isomorphism in $\Module{\dbl{E}}$. So,
  the first equality follows by its construction. The middle equality is
  \cref{equation:auxiliaryequationforexternalcompositeofmodulations}. Of course,
  on the other hand, we can analogously compute that
  \begin{equation*}
    \begin{dblArray}{c}
      \rho \\
      \cong \\
      G\theta\alpha_m\otimes\bar\alpha_g \\
      G_{n,w}^\otimes
    \end{dblArray} =
    \begin{dblArray}{c}
      1\odot F_z\\
      \mathrm{coeq} \\
      G\theta\alpha_m\otimes\bar\alpha_g\\
      G_{n,w}^\otimes
    \end{dblArray} =
    \begin{dblArray}{cc}
      1 & F_x \\
      \alpha_m & \Block{2-1}{\bar\alpha_f} \\
      G\theta & \\
      \Block{1-2}{G_{n,w}}
    \end{dblArray} =
    \begin{dblArray}{cc}
      \alpha_m & \Block{2-1}{\alpha_g}\\
      G\theta & \\
      \Block{1-2}{G_{n,w}}
    \end{dblArray}.
  \end{equation*}
  But the right-most side of each of the last two displays are equal by the cell
  naturality condition assumed for $\alpha$. So, the desired equation holds
  modulo the precomposed isomorphisms $\lambda$ and $\rho$. But these are easily
  cancelled, so the desired equality holds.

  Now, \cref{lemma:lax-transf-data-bijection} proves that the mapping is a
  bijection on morphism since it suffices to check that the component cells of
  lax transformations on each side are in bijection via the defined
  correspondence. We need to see then that $\alpha\mapsto\bar\alpha$ is
  functorial using the definition of composition in
  \cref{prop:lax-functor-category}. For this, fix composable lax transformations
  $\alpha\colon F\To G$ and $\beta\colon G\To H$. On the one hand, the component
  $(\overline{\beta\alpha})_f$ is computed as
  \begin{equation*}
    \begin{dblArray}{ccc}
      \id_{\alpha_x} & \Block{2-1}{\alpha_f} & F\id_f \\
      \Block{2-1}{\beta_f} & & \Block{2-1}{(\beta\alpha)_{\id_y}} \\
      & \beta_{\id_y} & \\
      \Block{1-2}{H_{y,y}} & & 1 \\
      \Block{1-3}{H_{y,y}} & &
    \end{dblArray}=
    \begin{dblArray}{ccc}
      \id_{\alpha_x} & \Block{2-1}{\alpha_f} & F\id_f \\
      \Block{2-1}{\beta_f} & & \alpha_{\id_y}  \\
      & \beta_{\id_y} & \beta_{\id_y} \\
      1 & \Block{1-2}{H_{y,y}} & \\
      \Block{1-3}{H_{y,y}} & &
    \end{dblArray}=
    \begin{dblArray}{ccc}
      \id_{\alpha_x} & \Block{2-1}{\alpha_f} & F\id_f \\
      \Block{2-1}{\beta_f} & & \alpha_{\id_y}  \\
      & \Block{1-2}{H_{y,y}} & \\
      1 & \Block{1-2}{\beta_{\id_y}}  & \\
      \Block{1-3}{H_{y,y}} & &
    \end{dblArray}=
    (\beta\bar\alpha)_f
  \end{equation*}
  just using laxator associativity and external functoriality of the proarrow
  components of $\beta$. On the other hand, $(\bar\beta\bar\alpha)_f$ requires
  composing in $\Module{\dbl{E}}$. To see how this works we need to be more
  careful with inserting the usually suppressed canonical unit isomorphisms. We
  have
  \begin{equation*}
    \begin{dblArray}{c}
      \lambda \\
      \cong \\
      \bar\beta_f\otimes\bar\alpha_f \\
      H^\otimes_{y,y}
    \end{dblArray}=
    \begin{dblArray}{c}
      H_x \odot 1 \\
      \mathrm{coeq} \\
      \bar\beta_f \otimes \bar\alpha_f \\
      H^\otimes_{y,y}
    \end{dblArray}=
    \begin{dblArray}{cc}
      H_x & 1 \\
      \bar\beta_f & \bar\alpha_f \\
      \Block{1-2}{H_{y,y}} & \\
    \end{dblArray}=
    \begin{dblArray}{cc}
      \beta_f & \bar\alpha_f \\
      \Block{1-2}{H_{y,y}} & \\
    \end{dblArray}
  \end{equation*}
  using \cref{lemma:lax-transf-data-bijection} for $\beta$ in the last step.
  Preservation of units is easier to check and is omitted.
\end{proof}

\begin{corollary} \label{cor:1-categories-of-models-on-2-category}
  For any 2-category $\bicat{A}$ and any equipment $\dbl{E}$ with local
  coequalizers, there is an equivalence of categories
  \begin{equation*}
    \LaxFunOne_\ell(\VerDbl(\bicat{A}),\dbl{E}) \simeq
      \TwoCatOne_\ell(\bicat{A},\Cat(\dbl{E})).
  \end{equation*}
  In particular, there is an equivalence of categories
  \begin{equation*}
    \LaxFunOne_\ell(\VerDbl(\bicat{A}),\Span) \simeq
      \TwoCatOne_\ell(\bicat{A},\Cat).
  \end{equation*}
\end{corollary}
\begin{proof}
  Compute that 
  \begin{align}
      \LaxFunOne_\ell(\VerDbl(\bicat{A}),\dbl{E})
      &\cong \LaxFunOne_{\ell,u}(\VerDbl(\bicat{A}),\Module{\dbl{E}})\qquad
        &\text{(\cref{prop:unitalizationisoof1categories})}\notag\\
      &\cong \TwoCatOne_\ell(\bicat{A},\VerTwoCat(\Module{\dbl{E}}))\qquad
        &\text{(\cref{cor:characterize-unitary-models-on-veritical-dbl-cat-I})}\notag\\
      &\cong \TwoCatOne_\ell(\bicat{A},\Cat(\dbl{E})),\notag
  \end{align}
  where the last isomorphism holds since the 2-category underlying
  $\Module{\dbl{E}}$ is just the 2-category of categories, functors and
  transformations in $\dbl{E}$. The last statement follows by choosing
  $\dbl{E}=\Span$.
\end{proof}

We now turn to lax transformations between \emph{cartesian} double functors. As
is well known, ordinary natural transformations automatically preserve products,
which has the consequence that no extra conditions are needed in forming the
2-category of cartesian categories, cartesian functors, and natural
transformations. For lax natural transformations, we need an extra axiom
asserting that products are preserved strictly.

\begin{definition}[Cartesian lax transformation]
  \label{def:cartesian-lax-transformation}
  Let $\dbl{D}$ and $\dbl{E}$ be precartesian double categories and let
  $F,G: \dbl{D} \to \dbl{E}$ be (not necessarily cartesian) lax double functors.
  A lax natural transformation $\alpha: F \To G$ is \define{cartesian} if it is
  \emph{strictly} natural with respect to projections, meaning that the squares
  \begin{equation*}
    \begin{tikzcd}
      {F(x \times x')} & Fx \\
      {G(x \times x')} & Gx
      \arrow["{F\pi_{x,x'}}", from=1-1, to=1-2]
      \arrow["{\alpha_{x \times x'}}"', from=1-1, to=2-1]
      \arrow["{\alpha_x}", from=1-2, to=2-2]
      \arrow["{G \pi_{x,x'}}"', from=2-1, to=2-2]
    \end{tikzcd}
    \qquad\text{and}\qquad
    \begin{tikzcd}
      {F(x \times x')} & {Fx'} \\
      {G(x \times x')} & {Gx'}
      \arrow["{F\pi_{x,x'}'}", from=1-1, to=1-2]
      \arrow["{\alpha_{x \times x'}}"', from=1-1, to=2-1]
      \arrow["{\alpha_{x'}}", from=1-2, to=2-2]
      \arrow["{G \pi_{x,x'}'}"', from=2-1, to=2-2]
    \end{tikzcd}
  \end{equation*}
  commute for all objects $x$ and $x'$ in $\dbl{D}$ and that the corresponding
  naturality comparisons have the form of
  \cref{eq:strict-transformation-as-lax}.
\end{definition}

In other contexts, the phrase ``cartesian natural transformation'' is often
taken to mean a natural transformation whose naturality squares are pullbacks.
We shall not employ that usage.

Cartesian lax transformations are clearly closed under the composition defined
in \cref{prop:lax-functor-category}. Let $\CartLaxFunOne_\ell(\dbl{D}, \dbl{E})$
denote the category of cartesian lax double functors and cartesian lax
transformations between two fixed precartesian double categories $\dbl{D}$ and
$\dbl{E}$. Similarly, let $\CartLaxFunOne_{\ell,n}(\dbl{D},\dbl{E})$ and
$\CartLaxFunOne_{\ell,u}(\dbl{D},\dbl{E})$ be the full subcategories spanned by
normal and unitary cartesian lax double functors.

The naturality comparisons of a cartesian lax natural transformation preserve
finite products, as we now show in the case of binary products.

\begin{lemma}[Naturality comparisons for products]
  Let $\dbl{D}$ and $\dbl{E}$ be precartesian double categories, let
  $F,G: \dbl{D} \to \dbl{E}$ be lax double functors, and let $\alpha: F \To G$
  be a cartesian lax natural transformation. Then for any arrows $f: x \to y$
  and $f': x' \to y'$ in $\dbl{D}$,
  \begin{equation*}
    \begin{tikzcd}[column sep=large]
      {F(x \times x')} & {F(x \times x')} \\
      {G(x \times x')} & {F(y \times y')} \\
      {G(y \times y')} & {G(y \times y')} \\
      {G(y \times y')} & {G(y \times y')} \\
      {Gy \times Gy'} & {Gy \times Gy'}
      \arrow["{\alpha_{x \times x'}}"', from=1-1, to=2-1]
      \arrow["{G(f \times f')}"', from=2-1, to=3-1]
      \arrow[""{name=0, anchor=center, inner sep=0}, "{\mathrm{id}_{F(x \times x')}}", "\shortmid"{marking}, from=1-1, to=1-2]
      \arrow["{F(f \times f')}", from=1-2, to=2-2]
      \arrow["{\alpha_{y \times y'}}", from=2-2, to=3-2]
      \arrow[""{name=1, anchor=center, inner sep=0}, "{G \mathrm{id}_{y \times y'}}", "\shortmid"{marking}, from=3-1, to=3-2]
      \arrow[""{name=2, anchor=center, inner sep=0}, "{G(\mathrm{id}_y \times \mathrm{id}_y)}"', "\shortmid"{marking}, from=4-1, to=4-2]
      \arrow[Rightarrow, no head, from=3-1, to=4-1]
      \arrow[Rightarrow, no head, from=3-2, to=4-2]
      \arrow["{\Psi_{y,y'}}"', from=4-1, to=5-1]
      \arrow["{\Psi_{y,y'}}", from=4-2, to=5-2]
      \arrow[""{name=3, anchor=center, inner sep=0}, "{G \mathrm{id}_y \times G \mathrm{id}_{y'}}"', "\shortmid"{marking}, from=5-1, to=5-2]
      \arrow["{\alpha_{f \times f'}}"{description}, draw=none, from=0, to=1]
      \arrow["{G \times_{(y,y')}}"{description}, draw=none, from=1, to=2]
      \arrow["{\Psi_{\mathrm{id}_y, \mathrm{id}_{y'}}}"{description, pos=0.6}, draw=none, from=2, to=3]
    \end{tikzcd}
    \quad=\quad
    \begin{tikzcd}[column sep=large]
      {F(x \times x')} & {F(x \times x')} \\
      {Fx \times Fx'} & {Fx \times Fx'} \\
      {Fx \times Fx'} & {Fx \times Fx'} \\
      {Gx \times Gx'} & {Fy \times Fy'} \\
      {Gy \times Gy'} & {Gy \times Gy'}
      \arrow[""{name=0, anchor=center, inner sep=0}, "{\mathrm{id}_{F(x \times x')}}", "\shortmid"{marking}, from=1-1, to=1-2]
      \arrow["{\Phi_{x,x'}}"', from=1-1, to=2-1]
      \arrow["{\Phi_{x,x'}}", from=1-2, to=2-2]
      \arrow[""{name=1, anchor=center, inner sep=0}, "{\mathrm{id}_{Fx \times Fx'}}"', "\shortmid"{marking}, from=2-1, to=2-2]
      \arrow[Rightarrow, no head, from=2-2, to=3-2]
      \arrow[Rightarrow, no head, from=2-1, to=3-1]
      \arrow[""{name=2, anchor=center, inner sep=0}, "{\mathrm{id}_{Fx} \times \mathrm{id}_{Fx'}}"', "\shortmid"{marking}, from=3-1, to=3-2]
      \arrow["{\alpha_x \times \alpha_{x'}}"', from=3-1, to=4-1]
      \arrow["{Ff \times Ff'}", from=3-2, to=4-2]
      \arrow[""{name=3, anchor=center, inner sep=0}, "{G \mathrm{id}_y \times G \mathrm{id}_{y'}}"', "\shortmid"{marking}, from=5-1, to=5-2]
      \arrow["{Gf \times Gf'}"', from=4-1, to=5-1]
      \arrow["{\alpha_y \times \alpha_{y'}}", from=4-2, to=5-2]
      \arrow["{\mathrm{id}_{\Phi_{x,x'}}}"{description}, draw=none, from=0, to=1]
      \arrow["{\times_{(Fx, Fx')}}"{description}, draw=none, from=1, to=2]
      \arrow["{\alpha_f \times \alpha_{f'}}"{description}, draw=none, from=2, to=3]
    \end{tikzcd}.
  \end{equation*}
  In particular, when the double category $\dbl{D}$ and lax functor $G$ are both
  cartesian, the comparison cell $\alpha_{f \times f'}$ is completely determined
  by the product of the comparisons $\alpha_f$ and $\alpha_{f'}$.
\end{lemma}
\begin{proof}
  Since $(f \times f') \cdot \pi_{y,y'} = \pi_{x,x'} \cdot f$, we have trivially
  that $\alpha_{(f \times f') \cdot \pi_{y,y'}} = \alpha_{\pi_{x,x'} \cdot f}$.
  This equation between cells expands to
  \begin{equation*}
    \begin{tikzcd}[row sep=scriptsize, column sep=large]
      {F(x \times x')} & {F(x \times x')} \\
      {G(x \times x')} & {F(y \times y')} \\
      {G(y \times y')} & {G(y \times y')} \\
      Gy & Gy
      \arrow["{\alpha_{x \times x'}}"', from=1-1, to=2-1]
      \arrow["{G(f \times f')}"', from=2-1, to=3-1]
      \arrow[""{name=0, anchor=center, inner sep=0}, "{\id_{F(x \times x')}}", "\shortmid"{marking}, from=1-1, to=1-2]
      \arrow["{F(f \times f')}", from=1-2, to=2-2]
      \arrow["{\alpha_{y \times y'}}", from=2-2, to=3-2]
      \arrow[""{name=1, anchor=center, inner sep=0}, "{G\id_{y \times y'}}", "\shortmid"{marking}, from=3-1, to=3-2]
      \arrow[""{name=2, anchor=center, inner sep=0}, "{G\id_y}"', "\shortmid"{marking}, from=4-1, to=4-2]
      \arrow["{G\pi_{y,y'}}"', from=3-1, to=4-1]
      \arrow["{G\pi_{y,y'}}", from=3-2, to=4-2]
      \arrow["{G \id_{\pi_{y,y'}}}"{description}, draw=none, from=1, to=2]
      \arrow["{\alpha_{f \times f'}}"{description}, draw=none, from=0, to=1]
    \end{tikzcd}
    \quad=\quad
    \begin{tikzcd}[row sep=scriptsize, column sep=large]
      {F(x \times x')} & {F(x \times x')} \\
      Fx & Fx \\
      Gx & Fy \\
      Gy & Gy
      \arrow[""{name=0, anchor=center, inner sep=0}, "{\id_{F(x \times x')}}", "\shortmid"{marking}, from=1-1, to=1-2]
      \arrow[""{name=1, anchor=center, inner sep=0}, "{\id_{Fx}}"', "\shortmid"{marking}, from=2-1, to=2-2]
      \arrow["{\alpha_x}"', from=2-1, to=3-1]
      \arrow["Ff", from=2-2, to=3-2]
      \arrow[""{name=2, anchor=center, inner sep=0}, "{G\id_y}"', "\shortmid"{marking}, from=4-1, to=4-2]
      \arrow["Gf"', from=3-1, to=4-1]
      \arrow["{\alpha_y}", from=3-2, to=4-2]
      \arrow["{F\pi_{x,x'}}"', from=1-1, to=2-1]
      \arrow["{F\pi_{x,x'}}", from=1-2, to=2-2]
      \arrow["{\id_{F\pi_{x,x'}}}"{description}, draw=none, from=0, to=1]
      \arrow["{\alpha_f}"{description}, draw=none, from=1, to=2]
    \end{tikzcd}
  \end{equation*}
  using the functorality of the naturality comparisons, the assumption that the
  lax transformation $\alpha$ is strictly natural with respect to projections,
  and the naturality of unitors. Starting from the equation
  $(f \times f') \cdot \pi_{y,y'}' = \pi_{x,x'}' \cdot f'$ yields a similar
  equation relating $\alpha_{f \times f'}$ and $\alpha_{f'}$. The pairing of
  these two equations is
  \begin{equation*}
    \begin{tikzcd}[row sep=scriptsize, column sep=large]
      {F(x \times x')} & {F(x \times x')} \\
      {G(x \times x')} & {F(y \times y')} \\
      {G(y \times y')} & {G(y \times y')} \\
      {G(y \times y')} & {G(y \times y')}
      \arrow["{\alpha_{x \times x'}}"', from=1-1, to=2-1]
      \arrow["{G(f \times f')}"', from=2-1, to=3-1]
      \arrow[""{name=0, anchor=center, inner sep=0}, "{\mathrm{id}_{F(x \times x')}}", "\shortmid"{marking}, from=1-1, to=1-2]
      \arrow["{F(f \times f')}", from=1-2, to=2-2]
      \arrow["{\alpha_{y \times y'}}", from=2-2, to=3-2]
      \arrow[""{name=1, anchor=center, inner sep=0}, "{G \mathrm{id}_{y \times y'}}", "\shortmid"{marking}, from=3-1, to=3-2]
      \arrow[""{name=2, anchor=center, inner sep=0}, "{G(\mathrm{id}_y \times \mathrm{id}_y)}"', "\shortmid"{marking}, from=4-1, to=4-2]
      \arrow["{\Psi_{y,y'}}"', from=3-1, to=4-1]
      \arrow["{\Psi_{y,y'}}", from=3-2, to=4-2]
      \arrow["{\alpha_{f \times f'}}"{description}, draw=none, from=0, to=1]
      \arrow["{\langle G\id_{\pi_{y,y'}}, G\id_{\pi_{y,y'}'}\rangle}"{description}, draw=none, from=1, to=2]
    \end{tikzcd}
    \quad=\quad
    \begin{tikzcd}[row sep=scriptsize, column sep=large]
      {Fx \times Fx'} & {Fx \times Fx'} \\
      {Fx \times Fx'} & {Fx \times Fx'} \\
      {Gx \times Gx'} & {Fy \times Fy'} \\
      {Gy \times Gy'} & {Gy \times Gy'}
      \arrow[""{name=0, anchor=center, inner sep=0}, "{\mathrm{id}_{Fx \times Fx'}}", "\shortmid"{marking}, from=1-1, to=1-2]
      \arrow["{\Phi_{x,x'}}", from=1-2, to=2-2]
      \arrow["{\Phi_{x,x'}}"', from=1-1, to=2-1]
      \arrow[""{name=1, anchor=center, inner sep=0}, "{\mathrm{id}_{Fx} \times \mathrm{id}_{Fx'}}"', "\shortmid"{marking}, from=2-1, to=2-2]
      \arrow["{\alpha_x \times \alpha_{x'}}"', from=2-1, to=3-1]
      \arrow["{Ff \times Ff'}", from=2-2, to=3-2]
      \arrow[""{name=2, anchor=center, inner sep=0}, "{G \mathrm{id}_y \times G \mathrm{id}_{y'}}"', "\shortmid"{marking}, from=4-1, to=4-2]
      \arrow["{Gf \times Gf'}"', from=3-1, to=4-1]
      \arrow["{\alpha_y \times \alpha_{y'}}", from=3-2, to=4-2]
      \arrow["{\langle \id_{F\pi_{x,x'}}, \id_{F\pi_{x,x'}'} \rangle}"{description}, draw=none, from=0, to=1]
      \arrow["{\alpha_f \times \alpha_{f'}}"{description}, draw=none, from=1, to=2]
    \end{tikzcd}.
  \end{equation*}
  In view of
  \cref{eq:functor-product-id-comparison-1,eq:functor-product-id-comparison-2},
  the lemma is proved.
\end{proof}

By a similar argument, cartesian lax transformations can be shown to respect
diagonals, even though they are only defined to be strictly natural with respect
to projections.

The equivalences of \cref{cor:1-categories-of-models-on-2-category} specialize
to the cartesian setting, which will be helpful in the analysis of several
examples below. To state the result, for any \emph{cartesian} 2-categories
$\bicat{A}$ and $\bicat{B}$, let $\CartTwoCatOne_l(\bicat{A},\bicat{B})$ denote
the category of cartesian strict 2-functors and cartesian lax transformations.

\begin{corollary} \label{cor:1-categories-of-cartesian-models-on-2-category}
  For any cartesian 2-category $\bicat{A}$ and cartesian equipment $\dbl{E}$
  with local coequalizers, there is an isomorphism of categories
  \begin{equation*}
    \CartLaxFunOne_\ell(\VerDbl(\bicat{A}),\dbl{E}) \cong
      \CartTwoCatOne_\ell(\bicat{A},\Cat(\dbl{E})).
  \end{equation*}
  In particular, there is an isomorphism of categories
  \begin{equation*}
    \CartLaxFunOne_\ell(\VerDbl(\bicat{A}),\Span) \cong
      \CartTwoCatOne_\ell(\bicat{A},\Cat).
  \end{equation*}
\end{corollary}

To conclude this section, we verify, for a number of different double theories,
that cartesian lax natural transformations give the correct notion of lax
morphism between models.

\begin{example}[Lax monoidal functors]
  \label{example:lax-transfs-are-lax-monoidal-funcs}
  Let $\dbl{T}$ be the cartesian double theory of monoids (\cref{th:monoid}). A
  model, that is, a cartesian, unitary lax double functor
  $M\colon \dbl{T}\to\Prof$, amounts to a strict monoidal category
  $\cat{M} \coloneqq Mx$. Since the theory $\dbl{T}$ has the form
  $\VerDbl(\bicat{T})$ for a 2-category $\bicat{T}$ and is thus horizontally
  trivial, models of the theory amount to cartesian 2-functors
  $M\colon \bicat{T} \to \Cat$ and this characterization affords a more
  straightforward confirmation that such models are precisely monoidal
  categories.

  Cartesian lax transformations of models are then precisely monoidal functors.
  Let $\phi\colon M\To N$ denote such a transformation. It consists of a
  component $\phi_x\colon \cat{M}\to\cat{N}$, which is an arrow of $\Cat$, hence
  is a functor. Laxity implies that there are comparison cells
  \begin{equation*}
    \begin{tikzcd}
      {\cat{M}^2} & {\cat{M}} \\
      {\cat{N}^2} & {\cat{N}}
      \arrow["{\otimes_{\cat{M}}}", from=1-1, to=1-2]
      \arrow[""{name=0, anchor=center, inner sep=0}, "{\phi_x}", from=1-2, to=2-2]
      \arrow[""{name=1, anchor=center, inner sep=0}, "{\phi_x^2}"', from=1-1, to=2-1]
      \arrow["{\otimes_{\cat{N}}}"', from=2-1, to=2-2]
      \arrow["{\phi_\otimes}", shorten <=8pt, shorten >=8pt, Rightarrow, from=1, to=0]
    \end{tikzcd}
    \qquad\text{and}\qquad
    \begin{tikzcd}
      1 & {\cat{M}} \\
      1 & {\cat{N}}
      \arrow[""{name=0, anchor=center, inner sep=0}, Rightarrow, no head, from=1-1, to=2-1]
      \arrow["{I_{\cat{N}}}"', from=2-1, to=2-2]
      \arrow["{I_{\cat{M}}}", from=1-1, to=1-2]
      \arrow[""{name=1, anchor=center, inner sep=0}, "{\phi_x}", from=1-2, to=2-2]
      \arrow["{\phi_I}", shorten <=7pt, shorten >=7pt, Rightarrow, from=0, to=1]
    \end{tikzcd}
  \end{equation*}
  corresponding to the arrows $\otimes \colon x^2\to x$ and $I\colon 1\to x$ in
  the theory. These satisfy the usual associativity and unitality conditions of
  a lax monoidal functor, as in \cite[\S{11.2}]{maclane1998}. Associativity
  follows by instantiating the arrow functoriality equation of
  \cref{def:lax-transformation} for each side of the associativity law
  $\otimes (1\times \otimes) = \otimes(\otimes\times 1)$. Likewise for the
  unitality axioms. Conversely, every lax monoidal functor defines a lax natural
  transformation of the unitary lax functors $\bicat{T}\to\Cat$ defined by
  the strict monoidal categories.

  These considerations also apply to the (non-strict) monoidal categories that
  are the models of \cref{th:pseudomonoid}. That is, monoidal categories are the
  same as cartesian lax functors $\dbl{T}\to\Span$ or cartesian unitary lax
  functors $\bicat{T}\to\Cat$, where $\bicat{T}$ is the 2-category underlying
  the theory. Again, cartesian lax transformations of such theories are
  precisely lax monoidal functors as in
  \cref{example:lax-transfs-are-lax-monoidal-funcs} above. The difference in the
  proofs is merely that in this case the non-trivial associators and unitors of
  the theory must be accounted for. However, the associators and unitors in the
  monoidal category are images of those from the theory; so the extra data can
  be seen to obey the required associativity and unitality laws using arrow
  functoriality in conjunction with cell naturality in
  \cref{def:lax-transformation}.
\end{example}

\begin{example}[Cartesian functors]
  In particular, \cref{example:lax-transfs-are-lax-monoidal-funcs} implies that
  in the case of cartesian monoidal categories (\cref{th:cart-mon-cat-I} and
  \cref{th:internal-commutative-comonoid}), pseudo natural transformations
  between models (i.e., cartesian categories) are precisely cartesian (i.e.,
  product-preserving) functors.
\end{example}

\begin{example}[Monad functors] \label{ex:lax-transfs-are-monad-functors}
  Now consider the theory of monads (\cref{th:monad}). Lax transformations of
  models are precisely the \emph{monad functors} described by Street
  \cite[\S{1}]{street1972}. If $\phi$ is such a transformation, then the
  naturality comparison at $t\colon x\to x$ provides the required cell
  \begin{equation*}
    \begin{tikzcd}
      {\cat{C}} & {\cat{C}} \\
      {\cat{D}} & {\cat{D}}
      \arrow["S", from=1-1, to=1-2]
      \arrow[""{name=0, anchor=center, inner sep=0}, "U"', from=1-1, to=2-1]
      \arrow["T"', from=2-1, to=2-2]
      \arrow[""{name=1, anchor=center, inner sep=0}, "U", from=1-2, to=2-2]
      \arrow["\phi", shorten <=7pt, shorten >=7pt, Rightarrow, from=0, to=1]
    \end{tikzcd}
  \end{equation*}
  which satisfies the unit condition
  \begin{equation*}
    \begin{tikzcd}
      {\cat{C}} & {\cat{D}} \\
      {\cat{C}} & {\cat{D}}
      \arrow[""{name=0, anchor=center, inner sep=0}, "U", from=1-1, to=1-2]
      \arrow[""{name=1, anchor=center, inner sep=0}, "T"', from=1-2, to=2-2]
      \arrow["S"', from=1-1, to=2-1]
      \arrow[""{name=2, anchor=center, inner sep=0}, "U"', from=2-1, to=2-2]
      \arrow[""{name=3, anchor=center, inner sep=0}, "1", curve={height=-18pt}, from=1-2, to=2-2]
      \arrow["\phi"', shorten <=4pt, shorten >=4pt, Rightarrow, from=0, to=2]
      \arrow["\eta"'{pos=0.6}, shift left, shorten <=4pt, shorten >=4pt, Rightarrow, from=3, to=1]
    \end{tikzcd}
    \quad=\quad
    \begin{tikzcd}
      {\cat{C}} & {\cat{C}} & {\cat{D}}
      \arrow[""{name=0, anchor=center, inner sep=0}, "1", curve={height=-12pt}, from=1-1, to=1-2]
      \arrow[""{name=1, anchor=center, inner sep=0}, "S"', curve={height=12pt}, from=1-1, to=1-2]
      \arrow["U", from=1-2, to=1-3]
      \arrow["\eta", shorten <=3pt, shorten >=3pt, Rightarrow, from=0, to=1]
    \end{tikzcd}
  \end{equation*}
  for a monad functor by instantiating the cell naturality condition for a lax
  transformation at the generating unit cell $\eta$. It also satisfies the
  multiplication condition
  \begin{equation*}
    \begin{tikzcd}
      {\cat{C}} & {\cat{C}} \\
      {\cat{D}} & {\cat{D}}
      \arrow["S", from=1-1, to=1-2]
      \arrow[""{name=0, anchor=center, inner sep=0}, "U", from=1-2, to=2-2]
      \arrow[""{name=1, anchor=center, inner sep=0}, "U"', from=1-1, to=2-1]
      \arrow[""{name=2, anchor=center, inner sep=0}, "TT"', curve={height=18pt}, from=2-1, to=2-2]
      \arrow[""{name=3, anchor=center, inner sep=0}, "T", from=2-1, to=2-2]
      \arrow["\phi", shorten <=7pt, shorten >=7pt, Rightarrow, from=1, to=0]
      \arrow["\mu"', shorten <=2pt, shorten >=2pt, Rightarrow, from=2, to=3]
    \end{tikzcd}
    \quad=\quad
    \begin{tikzcd}
      {\cat{C}} & {\cat{C}} & {\cat{C}} \\
      {\cat{D}} & {\cat{D}} & {\cat{D}}
      \arrow["S", from=1-1, to=1-2]
      \arrow[""{name=0, anchor=center, inner sep=0}, "U"', from=1-1, to=2-1]
      \arrow["T", from=2-1, to=2-2]
      \arrow["S", from=1-2, to=1-3]
      \arrow["T"', from=2-2, to=2-3]
      \arrow[""{name=1, anchor=center, inner sep=0}, "U", from=1-3, to=2-3]
      \arrow[""{name=2, anchor=center, inner sep=0}, "U", from=1-2, to=2-2]
      \arrow[""{name=3, anchor=center, inner sep=0}, "S", curve={height=-24pt}, from=1-1, to=1-3]
      \arrow["\phi", shorten <=11pt, shorten >=7pt, Rightarrow, from=2, to=1]
      \arrow["\phi", shorten <=7pt, shorten >=7pt, Rightarrow, from=0, to=2]
      \arrow["\mu"{pos=0.3}, shorten >=2pt, Rightarrow, from=1-2, to=3]
    \end{tikzcd}
  \end{equation*}
  by instantiating both cell naturality and arrow functoriality at the
  generating cell $\mu$.
\end{example}

Lest it appear that all double theories of interest are merely 2-categorical in
nature, we examine the model homomorphisms for the theory of promonoids
(\cref{th:promonoid}). The axiomatization of this theory is purely
bicategorical, i.e., has only trivial arrows, but the double-categorical
structure is needed to get the correct morphisms, namely, multifunctors.

\begin{example}[Multifunctors]
  \label{example:lax-transfs-are-multifunctors}
  Letting $\dbl{T}$ be the theory of promonoids (\cref{th:promonoid}),
  span-valued models are precisely multicategories. Cartesian lax
  transformations are then precisely multifunctors in the usual sense
  \cite[Definition 2.1.9]{leinster2004}. Given such a transformation $\phi$
  between models $\cat{C}$ and $\cat{D}$, for each proarrow
  $p_n\colon x^n\proto x$ in $\dbl{T}$, there is a cell amounting to a map of
  spans of the form
  \begin{equation*}
    \begin{tikzcd}
      {\cat{C}_0^n} & {\cat{C}(...;.)} & {\cat{C}_0} \\
      {\cat{D}_0^n} & {\cat{D}(...;.)} & {\cat{D}_0}
      \arrow["{\phi_x}", from=1-3, to=2-3]
      \arrow["{\phi_x^n}"', from=1-1, to=2-1]
      \arrow[from=1-2, to=1-1]
      \arrow[from=1-2, to=1-3]
      \arrow[from=2-2, to=2-1]
      \arrow[from=2-2, to=2-3]
      \arrow["{\phi_{p_n}}", from=1-2, to=2-2]
    \end{tikzcd}.
  \end{equation*}
  Thus, for any objects $a_1,\dots, a_n, a$ in $\cat{C}$, we have a function
  \begin{equation*}
    \phi_{p_n}\colon \cat{C}(a_1,\dots, a_n;a)\to \cat{D}(\phi_x(a_1),\dots, \phi_x(a_n);\phi_x(a)))
  \end{equation*}
  satisfying the required functoriality conditions by the functoriality of
  naturality comparisons axiom in \cref{def:lax-transformation}. Conversely,
  given a multifunctor from $\cat{C}$ to $\cat{D}$, a cartesian lax
  transformation can be defined between the models.

  The correspondence of span-valued models with unitary profunctor-valued lax
  functors (\cref{prop:unitalizationisoof1categories}) gives another way to look
  at such transformations. In this case, a cartesian lax transformation $\phi$
  amounts to a family of morphisms of profunctors
  \begin{equation*}
    \begin{tikzcd}
      {(\cat{C}^n)^{\op}\times \cat{C}} && \Set \\
      {(\cat{D}^n)^{\op}\times \cat{D}} && \Set
      \arrow[""{name=0, anchor=center, inner sep=0}, "{\cat{C}(...;.)}", from=1-1, to=1-3]
      \arrow["{(\phi_x)^n\times \phi_x}"', from=1-1, to=2-1]
      \arrow[""{name=1, anchor=center, inner sep=0}, "{\cat{D}(...;.)}"', from=2-1, to=2-3]
      \arrow[Rightarrow, no head, from=1-3, to=2-3]
      \arrow["{\phi_{p_n}}"', shorten <=4pt, shorten >=4pt, Rightarrow, from=0, to=1]
    \end{tikzcd}
  \end{equation*}
  for each $n \geq 0$, where $\cat{C}$ and $\cat{D}$ are the underlying
  categories of unary morphisms associated with the multicategory source and
  target of the multifunctor $\phi_x$. Each component $\phi_{p_n}$ of the
  transformation $\phi$ thus yields the required function of multihomsets as in
  the display above. Again, the required functoriality conditions are recovered
  by the functoriality conditions of \cref{def:lax-transformation}.
\end{example}

\begin{example}[Lax maps of monoidal copresheaves]
  The theory of pseudomonoid actions (\cref{th:pseudomonoid-action}) affords
  another example that nontrivially uses the double-categorical structure. In
  fact, genuine double-categorical structure is needed to obtain both the models
  and the correct model homomorphisms. Span-valued models of the theory are
  monoidal copresheaves on monoidal categories. Cartesian lax transformations
  between models are pairs consisting of a lax monoidal functor
  $\phi_x: \cat{M} \to \cat{N}$ and a monoidal natural transformation
  \begin{equation*}
    \begin{tikzcd}
      {\cat{M}} & \Set \\
      {\cat{N}} & \Set
      \arrow[""{name=0, anchor=center, inner sep=0}, "P", from=1-1, to=1-2]
      \arrow["{\phi_x}"', from=1-1, to=2-1]
      \arrow[""{name=1, anchor=center, inner sep=0}, "Q"', from=2-1, to=2-2]
      \arrow[Rightarrow, no head, from=1-2, to=2-2]
      \arrow["{\phi_p}"', shorten <=4pt, shorten >=4pt, Rightarrow, from=0, to=1]
    \end{tikzcd}
  \end{equation*}
  of lax monoidal functors $P$ and $Q \circ \phi_x$. On the one hand, given a
  cartesian lax transformation $\phi$ of models $P\colon\cat{M}\to\Set$ and
  $Q\colon\cat{N}\to\Set$, we have again ordinary natural transformations
  \begin{equation*}
    \begin{tikzcd}
      {\cat{M}^2} & {\cat{M}} \\
      {\cat{N}^2} & {\cat{N}}
      \arrow["{\otimes_{\cat{M}}}", from=1-1, to=1-2]
      \arrow[""{name=0, anchor=center, inner sep=0}, "{\phi_x}", from=1-2, to=2-2]
      \arrow[""{name=1, anchor=center, inner sep=0}, "{\phi_x^2}"', from=1-1, to=2-1]
      \arrow["{\otimes_{\cat{N}}}"', from=2-1, to=2-2]
      \arrow["{\phi_\otimes}", shorten <=8pt, shorten >=8pt, Rightarrow, from=1, to=0]
    \end{tikzcd}
    \qquad\text{and}\qquad
    \begin{tikzcd}
      1 & {\cat{M}} \\
      1 & {\cat{N}}
      \arrow[""{name=0, anchor=center, inner sep=0}, Rightarrow, no head, from=1-1, to=2-1]
      \arrow["{I_{\cat{N}}}"', from=2-1, to=2-2]
      \arrow["{I_{\cat{M}}}", from=1-1, to=1-2]
      \arrow[""{name=1, anchor=center, inner sep=0}, "{\phi_x}", from=1-2, to=2-2]
      \arrow["{\phi_I}", shorten <=7pt, shorten >=7pt, Rightarrow, from=0, to=1]
    \end{tikzcd}
  \end{equation*}
  as in \cref{example:lax-transfs-are-lax-monoidal-funcs}, corresponding to the
  arrows $\otimes \colon x^2\to x$ and $I\colon 1\to x$ in the theory and making
  $\phi_x$ into a lax monoidal functor. Additionally, corresponding to the
  generating proarrow $p\colon 1\proto x$, there is the further natural
  transformation $\phi_p$ displayed above satisfying several conditions. In
  particular, cell naturality in \cref{def:lax-transformation} instantiated at
  the generating cell $\mu$ amounts to the equation
  \begin{equation*}
    \begin{tikzcd}
      1 & {\cat{M}^2} \\
      1 & {\cat{M}} \\
      1 & {\cat{N}}
      \arrow[""{name=0, anchor=center, inner sep=0}, "{P^2}", "\shortmid"{marking}, from=1-1, to=1-2]
      \arrow["{\otimes_{\cat{M}}}", from=1-2, to=2-2]
      \arrow[Rightarrow, no head, from=1-1, to=2-1]
      \arrow[""{name=1, anchor=center, inner sep=0}, "P"', "\shortmid"{marking}, from=2-1, to=2-2]
      \arrow["{\phi_x}", from=2-2, to=3-2]
      \arrow[Rightarrow, no head, from=2-1, to=3-1]
      \arrow[""{name=2, anchor=center, inner sep=0}, "Q"', "\shortmid"{marking}, from=3-1, to=3-2]
      \arrow["\mu"{description}, draw=none, from=0, to=1]
      \arrow["{\phi_p}"{description}, draw=none, from=1, to=2]
    \end{tikzcd}
    \quad=\quad
    \begin{tikzcd}
      1 & {\cat{M}^2} & {\cat{M}^2} \\
      1 & {\cat{N}^2} & {\cat{M}} \\
      1 & {\cat{N}} & {\cat{N}}
      \arrow[""{name=0, anchor=center, inner sep=0}, "{P^2}", "\shortmid"{marking}, from=1-1, to=1-2]
      \arrow[""{name=1, anchor=center, inner sep=0}, "\shortmid"{marking}, Rightarrow, from=1-2, to=1-3]
      \arrow["{\otimes_{\cat{M}}}", from=1-3, to=2-3]
      \arrow[Rightarrow, no head, from=1-1, to=2-1]
      \arrow[""{name=2, anchor=center, inner sep=0}, "{Q^2}"', "\shortmid"{marking}, from=2-1, to=2-2]
      \arrow["{\phi_x^2}", from=1-2, to=2-2]
      \arrow["{\cat{N}^2}", from=2-2, to=3-2]
      \arrow[Rightarrow, no head, from=2-1, to=3-1]
      \arrow[""{name=3, anchor=center, inner sep=0}, "Q"', "\shortmid"{marking}, from=3-1, to=3-2]
      \arrow[""{name=4, anchor=center, inner sep=0}, "\shortmid"{marking}, Rightarrow, from=3-2, to=3-3]
      \arrow["{\phi_x}", from=2-3, to=3-3]
      \arrow["\nu"{description}, draw=none, from=2, to=3]
      \arrow["{\phi_p^2}"{description}, draw=none, from=0, to=2]
      \arrow["{\phi_\otimes}"{description}, draw=none, from=1, to=4]
    \end{tikzcd},
  \end{equation*}
  where we have written $\mu$ and $\nu$ for the monoidal product comparisons
  coming with the lax monoidal functors $P$ and $Q$. Unpacked in $\Cat$ at a
  component $(a,b)\in\cat{M}^2$, this equation is just the statement that the
  diagram
  \begin{equation*}
    \begin{tikzcd}
      {Pa\times Pb} &&&& {P(a\otimes_{\cat{M}}b)} \\
      {Q(\phi(a))\times Q(\phi(b))} && {Q(\phi(a)\otimes_{\cat{N}}\phi(b))} && {Q(\phi(a\otimes_{\cat{M}}b))}
      \arrow["{\mu_{a,b}}", from=1-1, to=1-5]
      \arrow["{(\phi_p)_a\times(\phi_p)_b}"', from=1-1, to=2-1]
      \arrow["{\nu_{\phi(a),\phi(b)}}"', from=2-1, to=2-3]
      \arrow["{(\phi_p)_{a\otimes_{\cat{M}}b}}", from=1-5, to=2-5]
      \arrow["{Q((\phi_\otimes)_{a,b})}"', from=2-3, to=2-5]
    \end{tikzcd}
  \end{equation*}
  commutes. This equation is, in turn, the first condition for $\phi_p$ to be a
  monoidal natural transformation \cite[XI.2.(5)]{maclane1998}. The other
  condition follows similarly. Conversely, any monoidal transformation
  determines a cartesian lax transformation when defined suitably on the
  generating structure of the theory.
\end{example}

\section{2-categories of models}
\label{sec:model-2-categories}

Our objective is now to upgrade the category of lax double functors to a
2-category and obtain a 2-category of models of a double theory. For the
purposes of double-categorical logic, the correct notion of morphism between lax
transformations is not a modification but rather a \emph{modulation}, as
introduced by Paré \cite{pare2011}. In general, a modulation is a square-shaped
cell bounded by two transformations and two modules of lax double functors, but
we have not yet said what a module is. For now, we state a simplified version of
the definition where both modules involved are identities, which suffices to
construct a 2-category.

\begin{definition}[Modulation, special case] \label{def:modulation-special}
  A \define{modulation} $\mu: \alpha \Tto \beta$ between lax natural
  transformations $\alpha, \beta: F \To G$ of lax double functors
  $F,G: \dbl{D} \to \dbl{E}$ consists of, for every object $x$ of $\dbl{D}$, a
  cell in $\dbl{E}$
  \begin{equation*}
    \begin{tikzcd}
      Fx & Fx \\
      Gx & Gx
      \arrow[""{name=0, anchor=center, inner sep=0}, "{\mathrm{id}_{Fx}}", "\shortmid"{marking}, from=1-1, to=1-2]
      \arrow["{\alpha_x}"', from=1-1, to=2-1]
      \arrow["{\beta_x}", from=1-2, to=2-2]
      \arrow[""{name=1, anchor=center, inner sep=0}, "{G \mathrm{id}_x}"', "\shortmid"{marking}, from=2-1, to=2-2]
      \arrow["{\mu_x}"{description}, draw=none, from=0, to=1]
    \end{tikzcd},
  \end{equation*}
  the \define{component} of $\mu$ at $x$, satisfying the following two axioms.
  \begin{itemize}
    \item Equivariance: for every proarrow $m: x \proto y$ in $\dbl{D}$,
      \begin{equation} \label{eq:modulation-special-equivariance}
        \begin{aligned}
          \begin{tikzcd}[ampersand replacement=\&]
            Fx \& Fx \& Fy \\
            Gx \& Gx \& Gy \\
            Gx \&\& Gy
            \arrow[""{name=0, anchor=center, inner sep=0}, "{\mathrm{id}_{Fx}}", "\shortmid"{marking}, from=1-1, to=1-2]
            \arrow["{\alpha_x}"', from=1-1, to=2-1]
            \arrow["{\beta_x}", from=1-2, to=2-2]
            \arrow[""{name=1, anchor=center, inner sep=0}, "{G \mathrm{id}_x}"', "\shortmid"{marking}, from=2-1, to=2-2]
            \arrow[""{name=2, anchor=center, inner sep=0}, "Fm", "\shortmid"{marking}, from=1-2, to=1-3]
            \arrow["{\beta_y}", from=1-3, to=2-3]
            \arrow[""{name=3, anchor=center, inner sep=0}, "Gm"', "\shortmid"{marking}, from=2-2, to=2-3]
            \arrow[""{name=4, anchor=center, inner sep=0}, "Gm"', "\shortmid"{marking}, from=3-1, to=3-3]
            \arrow[Rightarrow, no head, from=2-1, to=3-1]
            \arrow[Rightarrow, no head, from=2-3, to=3-3]
            \arrow["{\mu_x}"{description}, draw=none, from=0, to=1]
            \arrow["{G_{x,m}}"{description}, draw=none, from=2-2, to=4]
            \arrow["{\beta_m}"{description}, draw=none, from=2, to=3]
          \end{tikzcd}
          &\eqqcolon
          \begin{tikzcd}[ampersand replacement=\&]
            Fx \& Fy \\
            Gx \& Gy
            \arrow[""{name=0, anchor=center, inner sep=0}, "Fm", "\shortmid"{marking}, from=1-1, to=1-2]
            \arrow["{\alpha_x}"', from=1-1, to=2-1]
            \arrow["{\beta_y}", from=1-2, to=2-2]
            \arrow[""{name=1, anchor=center, inner sep=0}, "Gm"', "\shortmid"{marking}, from=2-1, to=2-2]
            \arrow["{\mu_m}"{description}, draw=none, from=0, to=1]
          \end{tikzcd}
          \\
          &=
          \begin{tikzcd}[ampersand replacement=\&]
            Fx \& Fy \& Fy \\
            Gx \& Gy \& Gy \\
            Gx \&\& Gy
            \arrow[""{name=0, anchor=center, inner sep=0}, "Fm", "\shortmid"{marking}, from=1-1, to=1-2]
            \arrow["{\alpha_x}"', from=1-1, to=2-1]
            \arrow["{\alpha_y}"', from=1-2, to=2-2]
            \arrow[""{name=1, anchor=center, inner sep=0}, "Gm"', "\shortmid"{marking}, from=2-1, to=2-2]
            \arrow[""{name=2, anchor=center, inner sep=0}, "{G \mathrm{id}_y}"', "\shortmid"{marking}, from=2-2, to=2-3]
            \arrow[""{name=3, anchor=center, inner sep=0}, "Gm"', "\shortmid"{marking}, from=3-1, to=3-3]
            \arrow["{\beta_y}", from=1-3, to=2-3]
            \arrow[""{name=4, anchor=center, inner sep=0}, "{\mathrm{id}_{Fy}}", "\shortmid"{marking}, from=1-2, to=1-3]
            \arrow[Rightarrow, no head, from=2-1, to=3-1]
            \arrow[Rightarrow, no head, from=2-3, to=3-3]
            \arrow["{G_{m,y}}"{description}, draw=none, from=2-2, to=3]
            \arrow["{\alpha_m}"{description}, draw=none, from=0, to=1]
            \arrow["{\mu_y}"{description}, draw=none, from=4, to=2]
          \end{tikzcd}.
        \end{aligned}
      \end{equation}
      The \define{component} of $\mu$ at $m$ is defined by either side of this
      equation.
    \item Naturality: for every cell $\stdInlineCell{\gamma}$ in $\dbl{D}$,
      \begin{equation} \label{eq:modulation-naturality}
        \begin{tikzcd}
          Fx & Fx & Fy \\
          Gx & Fw & Fz \\
          Gw & Gw & Gz \\
          Gw && Gz
          \arrow["Ff"', from=1-2, to=2-2]
          \arrow["Fg", from=1-3, to=2-3]
          \arrow[""{name=0, anchor=center, inner sep=0}, "Fm", "\shortmid"{marking}, from=1-2, to=1-3]
          \arrow[""{name=1, anchor=center, inner sep=0}, "Fn"', "\shortmid"{marking}, from=2-2, to=2-3]
          \arrow[""{name=2, anchor=center, inner sep=0}, "{\mathrm{id}_{Fx}}", "\shortmid"{marking}, from=1-1, to=1-2]
          \arrow["{\alpha_x}"', from=1-1, to=2-1]
          \arrow["Gf"', from=2-1, to=3-1]
          \arrow["{\alpha_w}"', from=2-2, to=3-2]
          \arrow[""{name=3, anchor=center, inner sep=0}, "{G \mathrm{id}_w}"', "\shortmid"{marking}, from=3-1, to=3-2]
          \arrow["{\beta_z}", from=2-3, to=3-3]
          \arrow[""{name=4, anchor=center, inner sep=0}, "Gn"', "\shortmid"{marking}, from=3-2, to=3-3]
          \arrow[""{name=5, anchor=center, inner sep=0}, "Gn"', "\shortmid"{marking}, from=4-1, to=4-3]
          \arrow[Rightarrow, no head, from=3-1, to=4-1]
          \arrow[Rightarrow, no head, from=3-3, to=4-3]
          \arrow["F\gamma"{description}, draw=none, from=0, to=1]
          \arrow["{\mu_n}"{description}, draw=none, from=1, to=4]
          \arrow["{G_{w,n}}"{description}, draw=none, from=3-2, to=5]
          \arrow["{\alpha_f}"{description}, draw=none, from=2, to=3]
        \end{tikzcd}
        \quad=\quad
        \begin{tikzcd}
          Fx & Fy & Fy \\
          Gx & Gy & Gz \\
          Gw & Gz & Gz \\
          Gw && Gz
          \arrow[""{name=0, anchor=center, inner sep=0}, "Fm", "\shortmid"{marking}, from=1-1, to=1-2]
          \arrow["{\alpha_x}"', from=1-1, to=2-1]
          \arrow["{\beta_y}", from=1-2, to=2-2]
          \arrow[""{name=1, anchor=center, inner sep=0}, "Gm"', "\shortmid"{marking}, from=2-1, to=2-2]
          \arrow["Gf"', from=2-1, to=3-1]
          \arrow["Gg", from=2-2, to=3-2]
          \arrow[""{name=2, anchor=center, inner sep=0}, "Gn"', "\shortmid"{marking}, from=3-1, to=3-2]
          \arrow[""{name=3, anchor=center, inner sep=0}, "{\mathrm{id}_{Fy}}", "\shortmid"{marking}, from=1-2, to=1-3]
          \arrow["Gg", from=1-3, to=2-3]
          \arrow["{\beta_z}", from=2-3, to=3-3]
          \arrow[""{name=4, anchor=center, inner sep=0}, "{G \mathrm{id}_z}"', "\shortmid"{marking}, from=3-2, to=3-3]
          \arrow[Rightarrow, no head, from=3-1, to=4-1]
          \arrow[Rightarrow, no head, from=3-3, to=4-3]
          \arrow[""{name=5, anchor=center, inner sep=0}, "Gn"', "\shortmid"{marking}, from=4-1, to=4-3]
          \arrow["{\mu_m}"{description}, draw=none, from=0, to=1]
          \arrow["G\gamma"{description}, draw=none, from=1, to=2]
          \arrow["{G_{n,z}}"{description}, draw=none, from=3-2, to=5]
          \arrow["{\beta_g}"{description}, draw=none, from=3, to=4]
        \end{tikzcd}.
      \end{equation}
  \end{itemize}
\end{definition}

\begin{remark}[Special cases]
  \label{rmk:modulations-reduce-to-modications}
  As with lax transformations between 2-functors of vertical 2-categories in
  \cref{remark:lax-trans-dbl-func-reduce-to-lax-trans-of-2-func}, any modulation
  $\mu\colon \alpha \Tto \beta$ of lax natural transformations restricts to an
  ordinary \emph{modification} between lax transformations of 2-functors
  $\VerTwoCat(\mu)\colon \VerTwoCat(\alpha) \Tto \VerTwoCat(\beta)$, as defined
  in \cite[\S 4.2]{johnson2021}. In particular, the first condition of
  \cref{def:modulation-special} is rendered trivial by such restriction whereas
  the second condition is precisely the familiar \emph{cylinder} condition for a
  modification.
\end{remark}

\begin{example}[Monad functor transformations]
  \label{ex:modulations-are-monad-functor-transfs}
  Returning to \cref{ex:lax-transfs-are-monad-functors}, we examine the 2-cells
  between the \emph{monad functors} that were seen to be precisely the lax
  transformations between models of the theory of monads (\cref{th:monad}). In
  this case, modulations, which amount to modifications since the theory is
  2-categorical, are precisely the \emph{monad functor transformations}
  considered by Street \cite[\S{1}]{street1972}. Such a modification
  $\sigma\colon (U,\phi) \Rrightarrow (V,\psi)$ comes with a cell
  $\sigma_x\colon U\Rightarrow V$ corresponding to the single generating object
  $x$. By the cell naturality condition in \cref{def:modulation-special}
  instantiated at the identity cell on the generating arrow $t\colon x\to x$, we
  have the equation
  \begin{equation*}
    \begin{tikzcd}
      {\cat{C}} & {\cat{D}} \\
      {\cat{C}} & {\cat{D}}
      \arrow["T", from=1-2, to=2-2]
      \arrow[""{name=0, anchor=center, inner sep=0}, "U", from=2-1, to=2-2]
      \arrow[""{name=1, anchor=center, inner sep=0}, "V"', curve={height=18pt}, from=2-1, to=2-2]
      \arrow["S"', from=1-1, to=2-1]
      \arrow[""{name=2, anchor=center, inner sep=0}, "U", from=1-1, to=1-2]
      \arrow["\sigma"{pos=0.4}, shift right, shorten <=2pt, shorten >=2pt, Rightarrow, from=0, to=1]
      \arrow["\phi"{pos=0.4}, shorten <=4pt, shorten >=9pt, Rightarrow, from=2, to=0]
    \end{tikzcd}
    \quad=\quad
    \begin{tikzcd}
      {\cat{C}} & {\cat{D}} \\
      {\cat{C}} & {\cat{D}}
      \arrow["T", from=1-2, to=2-2]
      \arrow[""{name=0, anchor=center, inner sep=0}, "V"', from=2-1, to=2-2]
      \arrow["S"', from=1-1, to=2-1]
      \arrow[""{name=1, anchor=center, inner sep=0}, "V"', from=1-1, to=1-2]
      \arrow[""{name=2, anchor=center, inner sep=0}, "U", curve={height=-18pt}, from=1-1, to=1-2]
      \arrow["\psi"{pos=0.6}, shorten <=9pt, shorten >=4pt, Rightarrow, from=1, to=0]
      \arrow["\sigma", shift right, shorten <=2pt, shorten >=2pt, Rightarrow, from=2, to=1]
    \end{tikzcd},
  \end{equation*}
  which is precisely the required compatibility condition for a monad functor
  transformation as described in the reference.
\end{example}

Having defined a modulation, the category of lax double functors
(\cref{prop:lax-functor-category}) upgrades to a 2-category:

\begin{theorem}[2-category of lax functors] \label{thm:lax-functor-2-category}
  For any double categories $\dbl{D}$ and $\dbl{E}$, there is a 2-category
  $\LaxFun_\ell(\dbl{D},\dbl{E})$ whose objects are lax double functors
  $\dbl{D} \to \dbl{E}$, morphisms are lax natural transformations, and
  2-morphisms are modulations.

  In each hom-category $\LaxFun_\ell(\dbl{D},\dbl{E})(F,G)$, the composite of
  modulations $\mu: \alpha \Tto \beta$ and $\nu: \beta \Tto \gamma$ is the
  modulation $\mu \cdot \nu: \alpha \Tto \gamma$ with components
  \begin{equation*} \label{eq:homwise-composition-modulations}
    \begin{tikzcd}
      Fx & Fx \\
      Gx & Gx
      \arrow[""{name=0, anchor=center, inner sep=0}, "{G \mathrm{id}_x}"', "\shortmid"{marking}, from=2-1, to=2-2]
      \arrow[""{name=1, anchor=center, inner sep=0}, "{\mathrm{id}_{Fx}}", "\shortmid"{marking}, from=1-1, to=1-2]
      \arrow["{\alpha_x}"', from=1-1, to=2-1]
      \arrow["{\gamma_x}", from=1-2, to=2-2]
      \arrow["{(\mu \cdot \nu)_x}"{description}, draw=none, from=1, to=0]
    \end{tikzcd}
    \quad\coloneqq\quad
    \begin{tikzcd}
      Fx & Fx & Fx \\
      Gx & Gx & Gx \\
      Gx && Gx
      \arrow[""{name=0, anchor=center, inner sep=0}, "{G \mathrm{id}_x}"', "\shortmid"{marking}, from=3-1, to=3-3]
      \arrow[""{name=1, anchor=center, inner sep=0}, "{G \mathrm{id}_x}"', "\shortmid"{marking}, from=2-1, to=2-2]
      \arrow[""{name=2, anchor=center, inner sep=0}, "{G \mathrm{id}_x}"', "\shortmid"{marking}, from=2-2, to=2-3]
      \arrow[Rightarrow, no head, from=2-1, to=3-1]
      \arrow[Rightarrow, no head, from=2-3, to=3-3]
      \arrow[""{name=3, anchor=center, inner sep=0}, "{\mathrm{id}_{Fx}}", "\shortmid"{marking}, from=1-1, to=1-2]
      \arrow[""{name=4, anchor=center, inner sep=0}, "{\mathrm{id}_{Fx}}", "\shortmid"{marking}, from=1-2, to=1-3]
      \arrow["{\alpha_x}"', from=1-1, to=2-1]
      \arrow["{\beta_x}"', from=1-2, to=2-2]
      \arrow["{\gamma_x}"', from=1-3, to=2-3]
      \arrow["{G_{x,x}}"{description}, draw=none, from=2-2, to=0]
      \arrow["{\mu_x}"{description}, draw=none, from=3, to=1]
      \arrow["{\nu_x}"{description}, draw=none, from=4, to=2]
    \end{tikzcd},
  \end{equation*}
  and the identity modulation $1_\alpha: \alpha \To \alpha$ has components
  $(1_\alpha)_x \coloneqq \id_{\alpha_x} \cdot G_x$.

  The composite of modulations $\mu: \alpha \Tto \beta: F \To G$ and
  $\nu: \gamma \Tto \delta: G \To H$ is the modulation
  $(\mu * \nu): (\alpha \cdot \gamma) \Tto (\beta \cdot \delta): F \To H$ with
  components
  \begin{equation*}
    \begin{tikzcd}
      Fx & Fx \\
      Hx & Hx
      \arrow[""{name=0, anchor=center, inner sep=0}, "{\mathrm{id}_{Fx}}", "\shortmid"{marking}, from=1-1, to=1-2]
      \arrow["{(\alpha \cdot \gamma)_x}"', from=1-1, to=2-1]
      \arrow["{(\beta \cdot \delta)_x}", from=1-2, to=2-2]
      \arrow[""{name=1, anchor=center, inner sep=0}, "{H \mathrm{id}_x}"', "\shortmid"{marking}, from=2-1, to=2-2]
      \arrow["{(\mu * \nu)_x}"{description}, draw=none, from=0, to=1]
    \end{tikzcd}
    \quad\coloneqq\quad
    \begin{tikzcd}
      Fx & Fx \\
      Gx & Gx \\
      Hx & Hx
      \arrow[""{name=0, anchor=center, inner sep=0}, "{\mathrm{id}_{Fx}}", "\shortmid"{marking}, from=1-1, to=1-2]
      \arrow["{\alpha_x}"', from=1-1, to=2-1]
      \arrow["{\beta_x}", from=1-2, to=2-2]
      \arrow[""{name=1, anchor=center, inner sep=0}, "{G \mathrm{id}_x}"', "\shortmid"{marking}, from=2-1, to=2-2]
      \arrow[""{name=2, anchor=center, inner sep=0}, "{H \mathrm{id}_x}"', "\shortmid"{marking}, from=3-1, to=3-2]
      \arrow["{\delta_x}", from=2-2, to=3-2]
      \arrow["{\gamma_x}"', from=2-1, to=3-1]
      \arrow["{\nu_{\mathrm{id}_x}}"{description}, draw=none, from=1, to=2]
      \arrow["{\mu_x}"{description}, draw=none, from=0, to=1]
    \end{tikzcd}
  \end{equation*}
  using the convention of \cref{eq:modulation-special-equivariance}.
\end{theorem}

\begin{proof}
  The associativity, unitality, and interchange laws for modulations follow from
  the corresponding laws for double categories and lax double functors. As for
  the proof that composite and identity modulations obey the axioms, we show
  only the longest in a series of calculations, namely that given modulations
  $\mu: \alpha \Tto \beta: F \To G$ and $\nu: \gamma \Tto \delta: G \To H$, the
  composite modulation
  \begin{equation*}
    (\mu * \nu): (\alpha \cdot \gamma) \Tto (\beta \cdot \delta):
      F \To H: \dbl{D} \to \dbl{E}
  \end{equation*}
  satisfies the naturality axiom.

  Fixing a cell $\stdInlineCell{\phi}$ in $\dbl{D}$, we begin by calculating
  \begin{align*}
    \begin{dblArray}{cc}
      (\mu * \nu)_m & \Block{2-1}{(\beta\delta)_g} \\
      H \phi & \\
      \Block{1-2}{H_{n,z}} &
    \end{dblArray} &=
    \begin{dblArray}{ccc}
      (\mu * \nu)_x & (\beta\delta)_m & \Block{3-1}{(\beta\delta)_g} \\
      \Block{1-2}{H_{x,m}} & & \\
      \Block{1-2}{H \phi} & & \\
      \Block{1-3}{H_{n,z}} & &
    \end{dblArray} =
    \begin{dblArray}{cccc}
      \mu_x & \beta_m & \id_{\beta_y} & \Block{2-1}{\beta_g} \\
      \nu_{\id_x} & \delta_m & \Block{2-1}{\delta_g} & \\
      H \id_f & H \phi & & \delta_{\id_z} \\
      \Block{1-4}{H_{w,n,z,z}} & & &
    \end{dblArray} \\ &=
    \begin{dblArray}{cccc}
      \mu_x & \id_{\beta_x} & \beta_m & \Block{2-1}{\beta_g} \\
      \nu_{\id_x} & \Block{2-1}{\delta_f} & G\phi & \\
      H \id_f & & \delta_n & \delta_{\id_z} \\
      \Block{1-4}{H_{w,w,n,z}} & & &
    \end{dblArray} =
    \begin{dblArray}{cccc}
      \id_{\alpha_x} & \mu_x & \beta_m & \Block{2-1}{\beta_g} \\
      \Block{2-1}{\gamma_f} & G \id_f & G\phi & \\
      & \nu_{\id_w} & \delta_n & \delta_{\id_z} \\
      \Block{1-4}{H_{w,w,n,z}} & & &
    \end{dblArray}\ .
  \end{align*}
  Focusing on the middle segment, we have
  \begin{align*}
    \begin{dblArray}{cc}
      \mu_x & \beta_m \\
      G \id_f & G \phi \\
      \nu_{\id_w} & \delta_n \\
      \Block{1-2}{H_{w,n}}
    \end{dblArray} &=
    \begin{dblArray}{ccc}
      \Block{1-2}{\mu_x} & & \beta_m \\
      \Block{1-2}{G \id_f} & & G\phi \\
      \nu_w & \delta_{\id_w} & \delta_n \\
      \Block{1-3}{H_{w,w,n}}
    \end{dblArray} =
    \begin{dblArray}{cc}
      \mu_x & \beta_m \\
      G \id_f & G\phi \\
      \Block{1-2}{G_{w,n}} \\
      \nu_w & \delta_n \\
      \Block{1-2}{H_{w,n}}
    \end{dblArray} =
    \begin{dblArray}{cc}
      \mu_x & \beta_m \\
      \Block{1-2}{G_{x,m}} \\
      \Block{1-2}{G\phi} \\
      \nu_w & \delta_n \\
      \Block{1-2}{H_{w,n}}
    \end{dblArray} \\&=
    \begin{dblArray}{cc}
      \alpha_m & \mu_y \\
      \Block{1-2}{G_{m,y}} \\
      \Block{1-2}{G\phi} \\
      \gamma_n & \nu_z \\
      \Block{1-2}{H_{n,z}}
    \end{dblArray} =
    \begin{dblArray}{cc}
      \alpha_m & \mu_y \\
      G\phi & G\id_g \\
      \Block{1-2}{G_{n,z}} \\
      \gamma_n & \nu_z \\
      \Block{1-2}{H_{n,z}}
    \end{dblArray} =
    \begin{dblArray}{ccc}
      \alpha_m & \Block{1-2}{\mu_y} & \\
      G\phi & \Block{1-2}{G\id_g} & \\
      \gamma_n & \gamma_{\id_z} & \nu_z \\
      \Block{1-3}{H_{n,z,z}}
    \end{dblArray}
  \end{align*}\ .
  Finally, substituting this equation into the previous one and using
  \cref{eq:modulation-id-cells}, we obtain
  \begin{align*}
    \begin{dblArray}{cc}
      (\mu * \nu)_m & \Block{2-1}{(\beta\delta)_g} \\
      H \phi & \\
      \Block{1-2}{H_{n,z}} &
    \end{dblArray} &=
    \begin{dblArray}{ccccc}
      \id_{\alpha_x} & \alpha_m & F_y \cdot \mu_{\id_y} & \Block{2-2}{\beta_g} \\
      \Block{2-1}{\gamma_f} & G \phi & G \id_g & \\
      & \gamma_n & \gamma_{\id_z} & \gamma_{\id_z} & \nu_z \\
      \Block{1-5}{H_{w,n,z,z,z}} & & & &
    \end{dblArray} =
    \begin{dblArray}{ccccc}
      \id_{\alpha_x} & \alpha_m & \Block{2-1}{\alpha_g} & \Block{1-2}{F_y \cdot F \id_g} & \\
      \Block{2-1}{\gamma_f} & G \phi & & \Block{1-2}{\mu_{\id_z}} & \\
      & \gamma_n & \gamma_{\id_z} & \gamma_{\id_z} & \nu_z \\
      \Block{1-5}{H_{w,n,z,z,z}} & & & &
    \end{dblArray} \\ &=
    \begin{dblArray}{cccc}
      \id_{\alpha_x} & \alpha_m & \Block{2-1}{\alpha_g} & \id_{Fg} \\
      \Block{2-1}{\gamma_f} & G \phi & & \mu_z \\
      & \gamma_n & \gamma_{\id_z} & \nu_{\id_z} \\
      \Block{1-4}{H_{w,n,z,z}} & & &
    \end{dblArray} =
    \begin{dblArray}{cccc}
      \id_{\alpha_x} & \Block{2-1}{\alpha_f} & F\phi & \id_{Fg} \\
      \Block{2-1}{\gamma_f} & & \alpha_n & \mu_z \\
      & \gamma_{\id_w} & \gamma_n & \nu_{\id_z} \\
      \Block{1-4}{H_{w,w,n,z}} & & &
    \end{dblArray} \\ &=
    \begin{dblArray}{ccc}
      \Block{2-1}{(\alpha\gamma)_f} & \Block{1-2}{F\phi} & \\
      & (\alpha\gamma)_n & (\mu * \nu)_z \\
      \Block{1-3}{H_{w,n,z}} & &
    \end{dblArray} =
    \begin{dblArray}{cc}
      \Block{2-1}{(\alpha\gamma)_f} & F\phi \\
      & (\mu*\nu)_n \\
      \Block{1-2}{H_{w,n}}
    \end{dblArray}\ . \qedhere
  \end{align*}
\end{proof}

Now that we have a 2-category of lax functors, our goal is to extend the
unitalization result of \cref{prop:unitalizationisoof1categories} from an
isomorphism of categories to an isomorphism of 2-categories. We first extend the
functor given by post-composition with the counit
$\epsilon: \Module{\dbl{E}} \to \dbl{E}$ to a 2-functor.

\begin{lemma}
  When $\dbl{D}$ is a double category and $\dbl{E}$ is an equipment with local
  coequalizers, there is a 2-functor
  \begin{equation*}
    \epsilon \circ (-): \LaxFun_{\ell,u}(\dbl{D},\Module{\dbl{E}}) \to 
    \LaxFun_\ell(\dbl{D},\dbl{E})
  \end{equation*}
  extending the functor from \cref{lem:whiskering-with-mod-counit}.
\end{lemma}
\begin{proof}
  Given a modulation $\mu\colon \alpha\Rrightarrow\beta$ of lax transformations
  $\alpha,\beta\colon H\Rightarrow K$ between unitary lax functors
  $H,K\colon \dbl{D} \rightrightarrows \Module{\dbl{E}}$, define a new
  modulation
  \begin{equation*}
    \epsilon\ast\mu \colon \epsilon\ast\alpha \Rrightarrow \epsilon\ast\beta
  \end{equation*}
  between the whiskered lax transformations having the component
  $(\epsilon\ast\mu)_x$ defined as the composite
  \begin{equation*}
    \begin{tikzcd}
      {H(x)_0} & {H(x)_0} \\
      {H(x)_0} & {H(x)_0} \\
      {K(x)_0} & {K(x)_0}
      \arrow[""{name=0, anchor=center, inner sep=0}, "{\id_{H(x)_0}}", "\shortmid"{marking}, from=1-1, to=1-2]
      \arrow[Rightarrow, no head, from=1-2, to=2-2]
      \arrow[Rightarrow, no head, from=1-1, to=2-1]
      \arrow[""{name=1, anchor=center, inner sep=0}, "Hx"', "\shortmid"{marking}, from=2-1, to=2-2]
      \arrow["{\alpha_x}"', from=2-1, to=3-1]
      \arrow[""{name=2, anchor=center, inner sep=0}, "Kx"', "\shortmid"{marking}, from=3-1, to=3-2]
      \arrow["{\beta_x}", from=2-2, to=3-2]
      \arrow["{\mu_x}"{description}, draw=none, from=1, to=2]
      \arrow["{\upsilon_{Hx}}"{description}, draw=none, from=0, to=1]
    \end{tikzcd}.
  \end{equation*}
  With this definition, it is immediate that $\epsilon\ast\mu$ is a modulation
  as in \cref{def:modulation-special} just using the corresponding properties
  for the given modulation $\mu$.
\end{proof}

We now construct the 2-categorical inverse to this 2-functor. To that end, let
$\mu\colon\alpha\Rrightarrow\beta$ be a modulation of lax transformations
$\alpha\colon F\Rightarrow G$ and $\beta\colon F\Rightarrow G$, where
$F,G\colon \dbl{D} \rightrightarrows \dbl{E}$ are lax double functors. The
component $\bar\mu_x$ of a proposed modulation
$\bar\mu\colon \bar\alpha \Rrightarrow\bar\beta$ is given as $\mu_{\id_x}$ that
is, as either composite in the equation
\begin{equation*}
  \begin{tikzcd}
    Fx & Fx & Fx \\
    Gx & Gx & Gx \\
    Gx && Gx
    \arrow[""{name=0, anchor=center, inner sep=0}, "{\id_{Fx}}", "\shortmid"{marking}, from=1-1, to=1-2]
    \arrow[from=1-2, to=2-2]
    \arrow[from=1-1, to=2-1]
    \arrow[""{name=1, anchor=center, inner sep=0}, "{G\id_x}"', "\shortmid"{marking}, from=2-1, to=2-2]
    \arrow[""{name=2, anchor=center, inner sep=0}, "{F\id_x}", "\shortmid"{marking}, from=1-2, to=1-3]
    \arrow[from=1-3, to=2-3]
    \arrow[""{name=3, anchor=center, inner sep=0}, "{G\id_x}"', "\shortmid"{marking}, from=2-2, to=2-3]
    \arrow[from=2-1, to=3-1]
    \arrow[""{name=4, anchor=center, inner sep=0}, "{G\id_x}"', "\shortmid"{marking}, from=3-1, to=3-3]
    \arrow[from=2-3, to=3-3]
    \arrow["{G_{x,x}}"{description}, draw=none, from=2-2, to=4]
    \arrow["{\mu_x}"{description}, draw=none, from=0, to=1]
    \arrow["{\beta_{\id_x}}"{description}, draw=none, from=2, to=3]
  \end{tikzcd}
  \quad=\quad
  \begin{tikzcd}
    Fx & Fx & Fx \\
    Gx & Gx & Gx \\
    Gx && Gx
    \arrow[from=1-1, to=2-1]
    \arrow[""{name=0, anchor=center, inner sep=0}, "{G\id_x}"', "\shortmid"{marking}, from=2-1, to=2-2]
    \arrow[""{name=1, anchor=center, inner sep=0}, "{F\id_x}", "\shortmid"{marking}, from=1-1, to=1-2]
    \arrow[""{name=2, anchor=center, inner sep=0}, "{\id_{Fx}}", "\shortmid"{marking}, from=1-2, to=1-3]
    \arrow[from=1-3, to=2-3]
    \arrow[from=1-2, to=2-2]
    \arrow[""{name=3, anchor=center, inner sep=0}, "{G\id_x}"', "\shortmid"{marking}, from=2-2, to=2-3]
    \arrow[from=2-1, to=3-1]
    \arrow[""{name=4, anchor=center, inner sep=0}, "{G\id_x}"', "\shortmid"{marking}, from=3-1, to=3-3]
    \arrow[from=2-3, to=3-3]
    \arrow["{G_{x,x}}"{description}, draw=none, from=2-2, to=4]
    \arrow["{\mu_x}"{description}, draw=none, from=2, to=3]
    \arrow["{\alpha_{\id_x}}"{description}, draw=none, from=1, to=0]
  \end{tikzcd}
\end{equation*}
which holds by the equivariance condition in \cref{def:modulation-special}.

A lemma will help with subsequent computations. It will show, roughly, that the
component cells of modulations of transformations of $\dbl{E}$-valued lax
functors are in bijective correspondence with component cells of modulations of
transformations of $\Module{\dbl{E}}$-valued unitary lax functors. The
correspondence is simply precomposing with an appropriate unitor.

\begin{lemma} \label{lemma:bijection-modulation-data}
  Let $\dbl{D}$ and $\dbl{E}$ be double categories. For any modulation
  $\mu\colon \alpha\Rrightarrow \beta$ of lax transformations
  $\alpha,\beta\colon F\Rightarrow G$ between lax functors
  $F,G\colon \dbl{D}\rightrightarrows \dbl{E}$, the defined cell $\bar\mu_x$
  above satisfies
  \begin{equation*}
    \begin{tikzcd}
      Fx & Fx \\
      Fx & Fx \\
      Gx & Gx
      \arrow[Rightarrow, no head, from=1-1, to=2-1]
      \arrow[""{name=0, anchor=center, inner sep=0}, "\shortmid"{marking}, from=2-1, to=2-2]
      \arrow[""{name=1, anchor=center, inner sep=0}, "{\id_{Fx}}", "\shortmid"{marking}, from=1-1, to=1-2]
      \arrow[Rightarrow, no head, from=1-2, to=2-2]
      \arrow["{\alpha_x}"', from=2-1, to=3-1]
      \arrow[""{name=2, anchor=center, inner sep=0}, "{G\id_x}"', "\shortmid"{marking}, from=3-1, to=3-2]
      \arrow["{\beta_x}", from=2-2, to=3-2]
      \arrow["{\bar\mu_x}"{description}, draw=none, from=0, to=2]
      \arrow["{F_x}"{description}, draw=none, from=1, to=0]
    \end{tikzcd}
    \quad=\quad
    \begin{tikzcd}
      Fx & Fx \\
      Gx & Gx
      \arrow[""{name=0, anchor=center, inner sep=0}, "{\id_{Fx}}", "\shortmid"{marking}, from=1-1, to=1-2]
      \arrow["{\alpha_x}"', from=1-1, to=2-1]
      \arrow[""{name=1, anchor=center, inner sep=0}, "{G\id_x}"', "\shortmid"{marking}, from=2-1, to=2-2]
      \arrow["{\beta_x}", from=1-2, to=2-2]
      \arrow["{\mu_x}"{description}, draw=none, from=0, to=1]
    \end{tikzcd}.
  \end{equation*}
  Likewise, if $\dbl{E}$ is an equipment with local coequalizers, then for any
  modulation $\nu\colon \gamma\Rrightarrow \delta$ of lax transformations
  $\gamma,\delta\colon H\Rightarrow K$ between unitary lax functors
  $H,K\colon \dbl{D} \rightrightarrows \Module{\dbl{E}}$, the component cell
  $\nu_x$ satisfies
  \begin{equation*}
    \begin{tikzcd}
      {H(x)_0} & {H(x)_0} & {H(x)_0} \\
      {H(x)_0} & {H(x)_0} & {H(x)_0} \\
      {K(x)_0} & {K(x)_0} & {K(x)_0} \\
      {K(x)_0} && {K(x)_0}
      \arrow[""{name=0, anchor=center, inner sep=0}, "\shortmid"{marking}, from=2-1, to=2-2]
      \arrow["{\delta_x}"', from=2-2, to=3-2]
      \arrow["{\gamma_x}"', from=2-1, to=3-1]
      \arrow[""{name=1, anchor=center, inner sep=0}, "Kx"', "\shortmid"{marking}, from=3-1, to=3-2]
      \arrow[Rightarrow, no head, from=1-1, to=2-1]
      \arrow[""{name=2, anchor=center, inner sep=0}, "{\id_{H(x)_0}}", "\shortmid"{marking}, from=1-1, to=1-2]
      \arrow[Rightarrow, no head, from=1-2, to=2-2]
      \arrow[""{name=3, anchor=center, inner sep=0}, "Hx", "\shortmid"{marking}, from=1-2, to=1-3]
      \arrow[Rightarrow, no head, from=1-3, to=2-3]
      \arrow[""{name=4, anchor=center, inner sep=0}, "\shortmid"{marking}, from=2-2, to=2-3]
      \arrow["{\delta_x}", from=2-3, to=3-3]
      \arrow[""{name=5, anchor=center, inner sep=0}, "Kx"', "\shortmid"{marking}, from=3-2, to=3-3]
      \arrow[""{name=6, anchor=center, inner sep=0}, Rightarrow, no head, from=3-1, to=4-1]
      \arrow["Kx"', "\shortmid"{marking}, from=4-1, to=4-3]
      \arrow[""{name=7, anchor=center, inner sep=0}, Rightarrow, no head, from=3-3, to=4-3]
      \arrow["{K_{x,x}}"{description}, draw=none, from=6, to=7]
      \arrow["{\nu_x}"{description}, draw=none, from=0, to=1]
      \arrow["{\delta_{id_x}}"{description}, draw=none, from=4, to=5]
      \arrow["{\upsilon_{Hx}}"{description}, draw=none, from=2, to=0]
      \arrow["1"{description}, draw=none, from=3, to=4]
    \end{tikzcd}
    \quad=\quad
    \begin{tikzcd}
      {H(x)_0} & {H(x)_0} \\
      {K(x)_0} & {K(x)_0}
      \arrow[""{name=0, anchor=center, inner sep=0}, "Hx", "\shortmid"{marking}, from=1-1, to=1-2]
      \arrow["{\delta_x}", from=1-2, to=2-2]
      \arrow["{\gamma_x}"', from=1-1, to=2-1]
      \arrow[""{name=1, anchor=center, inner sep=0}, "Kx"', "\shortmid"{marking}, from=2-1, to=2-2]
      \arrow["{\nu_x}"{description}, draw=none, from=0, to=1]
    \end{tikzcd}.
  \end{equation*}
\end{lemma}
\begin{proof}
  This proof is straightforward in that, as in the proof of 
  \cref{lemma:lax-transf-data-bijection}, it follows the pattern of that of
  the bijection in \cref{prop:natural-transformation-in-dbl}. On the one hand,
  we have
  \begin{equation*}
    \begin{dblArray}{cc}
      \Block{1-2}{F_x} & \\
      \Block{1-2}{\cong} & \\
      \mu_x & \beta_{\id_x}  \\
      \Block{1-2}{G_{x,x}} &
    \end{dblArray} =
    \begin{dblArray}{cc}
      \Block{1-2}{\cong} & \\
      1 & F_x \\
      \mu_x & \beta_{\id_x}  \\
      \Block{1-2}{G_{x,x}} &
    \end{dblArray} =
    \begin{dblArray}{cc}
      \mu_x & \id_{\alpha_x}  \\
      1 & G_x \\
      \Block{1-2}{G_{x,x}} &
    \end{dblArray} =
    \begin{dblArray}{c}
      \mu_x
    \end{dblArray}
  \end{equation*}
  using the first unit condition of \cref{def:lax-transformation} and then
  unitor coherence for a lax functor as in \cref{def:lax-functor}. Note that the
  unlabeled isomorphisms are the canonical unit isomorphisms required for the
  compositions to make sense. On the other hand, we have
  \begin{equation*}
    \begin{dblArray}{cc}
      \Block{1-2}{\cong} & \\
      \upsilon & 1 \\
      \nu_x & \delta_{\id_x}  \\
      \Block{1-2}{K_{x,x}} &
    \end{dblArray} =
    \begin{dblArray}{cc}
      \Block{1-2}{\cong} & \\
      \upsilon & 1 \\
      \Block{1-2}{\mu} & \\
      \Block{1-2}{\nu_x} &
    \end{dblArray} =
    \begin{dblArray}{c}
      \nu_x
    \end{dblArray}
  \end{equation*}
  by the unit laws in $\dbl{E}$ and the fact that $\nu_x$ is a modulation in
  $\dbl{E}$ and thus satisfies the equivariance axiom in
  \cref{def:profunctor-object}. Again $\upsilon$ and $\mu$ denote the structure 
  cells coming with $Hx$.
\end{proof}

\begin{remark}[Parameterizing modulations]
  In the definition of a modulation, the cells $\mu_m$, parameterized by proarrows
  $m: x \proto y$, are not part of the data of the modulation but are derived
  from the cells $\mu_x$ or $\mu_y$ parameterized by objects. \cref{lemma:bijection-modulation-data} shows that, conversely, $\mu_x$ can be recovered from $\mu_{\id_x}$ via the equation:
  \begin{equation} \label{eq:modulation-id-cells}
    \begin{tikzcd}
      Fx & Fx \\
      Gx & Gx
      \arrow[""{name=0, anchor=center, inner sep=0}, "{\mathrm{id}_{Fx}}", "\shortmid"{marking}, from=1-1, to=1-2]
      \arrow["{\alpha_x}"', from=1-1, to=2-1]
      \arrow["{\beta_x}", from=1-2, to=2-2]
      \arrow[""{name=1, anchor=center, inner sep=0}, "{G \mathrm{id}_x}"', "\shortmid"{marking}, from=2-1, to=2-2]
      \arrow["{\mu_x}"{description}, draw=none, from=0, to=1]
    \end{tikzcd}
    \quad=\quad
    \begin{tikzcd}
      Fx & Fx \\
      Fx & Fx \\
      Gx & Gx
      \arrow[""{name=0, anchor=center, inner sep=0}, "{\mathrm{id}_{Fx}}", "\shortmid"{marking}, from=1-1, to=1-2]
      \arrow[""{name=1, anchor=center, inner sep=0}, "{F \mathrm{id}_x}"', "\shortmid"{marking}, from=2-1, to=2-2]
      \arrow[Rightarrow, no head, from=1-1, to=2-1]
      \arrow[Rightarrow, no head, from=1-2, to=2-2]
      \arrow["{\alpha_x}"', from=2-1, to=3-1]
      \arrow["{\beta_x}", from=2-2, to=3-2]
      \arrow[""{name=2, anchor=center, inner sep=0}, "{G \mathrm{id}_x}"', "\shortmid"{marking}, from=3-1, to=3-2]
      \arrow["{F_x}"{description}, draw=none, from=0, to=1]
      \arrow["{\mu_{\mathrm{id}_x}}"{description}, draw=none, from=1, to=2]
    \end{tikzcd}.
  \end{equation}
  This relationship is the multiobject version of the two descriptions of a
  natural transformation in a double category
  (\cref{prop:natural-transformation-in-dbl}). It will be useful later when
  comparing with the general definition of modulation
  (\cref{lem:modulation-definitions}).
\end{remark}

\begin{theorem}[Unitalization of lax functors, two-dimensional]
  \label{th:unitalizationFINALFORM}
  If $\dbl{D}$ is a double category and $\dbl{E}$ is an equipment with local
  coequalizers, then the assignment $\mu\mapsto \bar \mu$ results in a 2-functor
  \begin{equation*}
    \LaxFun_\ell(\dbl{D},\dbl{E}) \to \LaxFun_{\ell,u}(\dbl{D}, \Module{\dbl{E}}), \qquad
    F\mapsto \Module{F}\eta,\ \alpha\mapsto\bar\alpha,\ \mu\mapsto\bar\mu,
  \end{equation*}
  extending the functor from \cref{prop:unitalizationisoof1categories}.
  Furthermore, this 2-functor is an isomorphism of 2-categories, whose inverse
  is the above 2-functor
  \begin{equation*}
    \epsilon \circ (-): \LaxFun_{\ell,u}(\dbl{D},\Module{\dbl{E}}) \xrightarrow{\cong}
      \LaxFun_\ell(\dbl{D},\dbl{E}).
  \end{equation*}
\end{theorem}
\begin{proof}
  We first need to verify that with components $\bar\mu_x$ defined above,
  $\bar \mu$ is a well-defined modulation $\bar\alpha\Rrightarrow \bar\beta$.
  There are two conditions of substance to check. We will verify the
  equivariance condition. To this end, fix a proarrow $m\colon x\proto y$ in
  $\dbl{D}$. We compute on the one hand that
  \begin{equation*}
    \begin{dblArray}{c}
      \lambda \\
      \cong \\
      \bar\mu_x \otimes \beta_m  \\
      G^\otimes_{x,m}
    \end{dblArray} =
    \begin{dblArray}{cc}
      F_x & 1 \\
      \bar\mu_x & \beta_m  \\
      \Block{1-2}{G_{x,m}} &
    \end{dblArray} =
    \begin{dblArray}{ccc}
      1 & F_x & 1 \\
      \mu_x & \beta_{\id_x} & \beta_m  \\
      \Block{1-2}{G_{x,x}} && 1 \\
      \Block{1-3}{G_{x,m}} &&
    \end{dblArray} =
    \begin{dblArray}{ccc}
      \mu_x & G_x & \beta_m  \\
      \Block{1-2}{G_{x,x}} && 1 \\
      \Block{1-3}{G_{x,m}} &&
    \end{dblArray} =
    \begin{dblArray}{cc}
      \mu_x & \beta_m  \\
      \Block{1-2}{G_{x,m}} &
    \end{dblArray}
  \end{equation*}
  using first the now standard trick that $\lambda$ with the canonical
  comparison introduces the coequalizer defining modules composition and
  translates these back to composites in $\dbl{E}$. The next equality uses the
  definition of $\bar\mu_x$. The penultimate one uses the first unit condition
  for the lax transformation $\beta$. The last equality is the lax functor
  unitor condition. Note that we have suppressed a few coherence isomorphisms
  for readability. On the other hand, by an analogous sequence of computations,
  we have that
  \begin{equation*}
    \begin{dblArray}{c}
      \rho \\
      \cong \\
      \alpha_m \otimes \bar\mu_y  \\
      G^\otimes_{m,y}
    \end{dblArray} =
    \begin{dblArray}{cc}
      \alpha_m & \mu_y  \\
      \Block{1-2}{G_{m,y}} &
    \end{dblArray}
  \end{equation*}
  but the right-most composites of each of the last two displays are equal by
  the equivariance condition assumed for $\mu$. Thus, cancelling the canonical
  comparisons and unitors, we conclude that
  \begin{equation*}
    \begin{dblArray}{c}
      \bar\mu_x \otimes \beta_m  \\
      G^\otimes_{x,m}
    \end{dblArray} =
    \begin{dblArray}{c}
      \alpha_m \otimes \bar\mu_y  \\
      G^\otimes_{m,y}
    \end{dblArray}
  \end{equation*}
  which is precisely the required equivariance condition. Naturality follows by
  a similar argument that we omit.

  As a result of \cref{prop:unitalizationisoof1categories}, it remains only to
  check that the assignment on modulations is 2-functorial and bijective. It
  suffices to check that each is the case for the modulation data, namely, the
  component cells $\mu_x$. That the correspondence is a bijection was proved in
  \cref{lemma:bijection-modulation-data}. So, we need to see that the two types
  of cell composition are preserved by the map $\mu\mapsto \bar \mu$. We shall
  show that homwise composition in \cref{eq:homwise-composition-modulations} is
  preserved and omit the other since it is easier. On the one hand,
  $(\overline{\nu\mu})_x$ is given as the composite
  \begin{equation*}
    \begin{tikzcd}
      Fx && Fx & Fx \\
      Fx & Fx & Fx & Fx \\
      Gx & Gx & Gx & Gx \\
      Gx && Gx & Gx \\
      Gx &&& Gx
      \arrow[""{name=0, anchor=center, inner sep=0}, "{F\id_x}", "\shortmid"{marking}, from=1-3, to=1-4]
      \arrow["{\id_{Fx}}", "\shortmid"{marking}, from=1-1, to=1-3]
      \arrow[""{name=1, anchor=center, inner sep=0}, Rightarrow, no head, from=1-1, to=2-1]
      \arrow[""{name=2, anchor=center, inner sep=0}, "{\id_{Fx}}", "\shortmid"{marking}, from=2-1, to=2-2]
      \arrow[""{name=3, anchor=center, inner sep=0}, "{\id_{Fx}}", "\shortmid"{marking}, from=2-2, to=2-3]
      \arrow[""{name=4, anchor=center, inner sep=0}, Rightarrow, no head, from=1-3, to=2-3]
      \arrow[Rightarrow, no head, from=1-4, to=2-4]
      \arrow[""{name=5, anchor=center, inner sep=0}, "{F\id_x}"', "\shortmid"{marking}, from=2-3, to=2-4]
      \arrow[from=2-1, to=3-1]
      \arrow[""{name=6, anchor=center, inner sep=0}, "{G\id_x}"', "\shortmid"{marking}, from=3-1, to=3-2]
      \arrow[""{name=7, anchor=center, inner sep=0}, "{G\id_x}"', "\shortmid"{marking}, from=3-2, to=3-3]
      \arrow[from=2-3, to=3-3]
      \arrow[""{name=8, anchor=center, inner sep=0}, "{G\id_x}"', "\shortmid"{marking}, from=3-3, to=3-4]
      \arrow[from=2-4, to=3-4]
      \arrow[""{name=9, anchor=center, inner sep=0}, Rightarrow, no head, from=3-1, to=4-1]
      \arrow["{G\id_x}"', "\shortmid"{marking}, from=4-1, to=4-3]
      \arrow[""{name=10, anchor=center, inner sep=0}, Rightarrow, no head, from=3-3, to=4-3]
      \arrow[""{name=11, anchor=center, inner sep=0}, "{G\id_x}"', "\shortmid"{marking}, from=4-3, to=4-4]
      \arrow[Rightarrow, no head, from=3-4, to=4-4]
      \arrow[""{name=12, anchor=center, inner sep=0}, Rightarrow, no head, from=4-1, to=5-1]
      \arrow["{G\id_x}"', "\shortmid"{marking}, from=5-1, to=5-4]
      \arrow[""{name=13, anchor=center, inner sep=0}, Rightarrow, no head, from=4-4, to=5-4]
      \arrow[from=2-2, to=3-2]
      \arrow["{G_{x,x}}"{description}, draw=none, from=12, to=13]
      \arrow["1"{description}, draw=none, from=8, to=11]
      \arrow["{\gamma_{\id_x}}"{description}, draw=none, from=5, to=8]
      \arrow["{\mu_x}"{description}, draw=none, from=2, to=6]
      \arrow["{\nu_x}"{description}, draw=none, from=3, to=7]
      \arrow["1"{description}, draw=none, from=0, to=5]
      \arrow["{G_{x,x}}"{description}, draw=none, from=9, to=10]
      \arrow["\cong"{description}, draw=none, from=1, to=4]
    \end{tikzcd}.
  \end{equation*}
  By an application of associativity of the laxators, this is the same as
  \begin{equation*}
    \begin{tikzcd}
      Fx & Fx & Fx \\
      Gx & Gx & Gx \\
      Gx && Gx
      \arrow[""{name=0, anchor=center, inner sep=0}, "{\id_{Fx}}", "\shortmid"{marking}, from=1-1, to=1-2]
      \arrow[""{name=1, anchor=center, inner sep=0}, "{F\id_x}", "\shortmid"{marking}, from=1-2, to=1-3]
      \arrow[from=1-1, to=2-1]
      \arrow[""{name=2, anchor=center, inner sep=0}, "{G\id_x}"', "\shortmid"{marking}, from=2-1, to=2-2]
      \arrow[""{name=3, anchor=center, inner sep=0}, "{G\id_x}"', "\shortmid"{marking}, from=2-2, to=2-3]
      \arrow[from=1-3, to=2-3]
      \arrow[""{name=4, anchor=center, inner sep=0}, Rightarrow, no head, from=2-1, to=3-1]
      \arrow["{G\id_x}"', "\shortmid"{marking}, from=3-1, to=3-3]
      \arrow[""{name=5, anchor=center, inner sep=0}, Rightarrow, no head, from=2-3, to=3-3]
      \arrow[from=1-2, to=2-2]
      \arrow["{\mu_x}"{description}, draw=none, from=0, to=2]
      \arrow["{\bar\nu_x}"{description}, draw=none, from=1, to=3]
      \arrow["{G_{x,x}}"{description}, draw=none, from=4, to=5]
    \end{tikzcd}.
  \end{equation*}
  On the other hand, the composite $(\bar\nu \bar\mu)_x$ is taken in modules. A
  calculation shows it is the same as the composite immediately above:
  \begin{equation*}
    \begin{dblArray}{c}
      \lambda \\
      \cong \\
      \bar \mu_x \otimes \bar\nu_x  \\
      G^\otimes_{x,x}
    \end{dblArray} =
    \begin{dblArray}{c}
      F_x\odot 1 \\
      \mathrm{coeq} \\
      \bar \mu_x \otimes \bar\nu_x  \\
      G^\otimes_{x,x}
    \end{dblArray} =
    \begin{dblArray}{cc}
      F_x & 1 \\
      \bar \mu_x & \bar\nu_x  \\
      \Block{1-2}{G_{x,x}} &
    \end{dblArray} =
    \begin{dblArray}{cc}
      \mu_x & \bar\nu_x  \\
      \Block{1-2}{G_{x,x}} &
    \end{dblArray}
  \end{equation*}
  using \cref{lemma:bijection-modulation-data} in the last step. The other
  verification is, as mentioned above, easier. Setting up the proposed equation
  $(\overline{\mu\ast\nu})_x = (\bar\mu\ast\bar\nu)_x$ one would just observe
  that a single application of proarrow functoriality eliminates the only
  ostensible extra cells on the right side of the equation.
\end{proof}

\begin{corollary} \label{cor:unitalize-span-valued-2-cat-iso-FINAL}
  For any double category $\dbl{D}$, the assignment $H\mapsto \Ob H$ induces an
  isomorphism of 2-categories
  \begin{equation*}
    \LaxFun_{\ell,u}(\dbl{D},\Prof) \xrightarrow{\cong}
      \LaxFun_\ell(\dbl{D},\Span) .
  \end{equation*}
\end{corollary}
\begin{proof}
  Instantiate \cref{th:unitalizationFINALFORM} for $\Prof = \Module{\Span}$.
\end{proof}

Boosting the result of \cref{cor:1-categories-of-models-on-2-category} to 2-categories, we have the following.

\begin{corollary} \label{cor:2-categories-of-models-non-cartesian}
  For any 2-category $\bicat{A}$, there is an isomorphism of 2-categories
  \begin{equation*}
    \LaxFun_\ell(\VerDbl(\bicat{A}),\dbl{E})
    \xrightarrow{\cong}
    \TwoCat_\ell(\bicat{A},\Cat(\dbl{E})).
  \end{equation*}
  In particular, there is an isomorphism of 2-categories
  \begin{equation*}
    \LaxFun_\ell(\VerDbl(\bicat{A}),\Span)
    \xrightarrow{\cong}
    \TwoCat_\ell(\bicat{A},\Cat).
  \end{equation*}
\end{corollary}
\begin{proof}
  The argument is formally the same computation as in the proof of
  \cref{cor:1-categories-of-models-on-2-category} except with 2-categories in
  place of 1-categories. The identifications of 1-categories involved there have
  been shown to work for the corresponding 2-categories except the equivalence
  of \cref{cor:1-categories-of-models-on-2-category}. But this follows too since
  the special case of modulations (\cref{def:modulation-special}) amounts to
  ordinary modifications when the parameterizing double category is horizontally
  trivial (\cref{rmk:modulations-reduce-to-modications}).
\end{proof}

\cref{cor:2-categories-of-models-non-cartesian} shows that when double theories
are both simple and purely 2-categorical, we can identify the 2-categories of
models with familiar constructions from the theory of 2-categories. That is,
span-valued models are profunctor-valued 2-functors; transformations of models
are 2-natural transformations; and finally modulations are ordinary
modifications. Accordingly, several of the simple double theories given in
\cref{sec:simple-double-theories} have easily characterized 2-categories of
models.

\begin{example}[Categories]
  As a basic check, taking the unit theory (\cref{th:unit-theory}), we recover
  an equivalence
  \begin{equation*}
    \LaxFun_\ell(\dbl{1},\Span)
    \xrightarrow{\simeq}
    \TwoCat_\ell(\bicat{1},\Cat)
    \xrightarrow{\simeq}
    \Cat
  \end{equation*}
  with the usual 2-category of categories, functors, and natural transformations.
  If $\dbl{T}$ is the walking arrow theory (\cref{th:walking-arrow}), we have 
  \begin{equation*}
    \LaxFun_\ell(\dbl{T},\Span)
    \xrightarrow{\simeq}
    \TwoCat_\ell(\bicat{T},\Cat)
    \xrightarrow{\simeq}
    \Cat^{\mathbf 2}
  \end{equation*}
  recovering the usual arrow 2-category of $\Cat$.
\end{example}

\begin{example}[Monads]
  Returning to the monad functor transformations of
  \cref{ex:modulations-are-monad-functor-transfs}, the corollary now implies
  that the 2-category of span-valued models of the theory of monads
  (\cref{th:monad}) is precisely the 2-category of monads, monad functors, and
  monad functor transformations in $\Cat$ as axiomatized by Street
  \cite{street1972}.
\end{example}

Of course, not all double theories of interest are simple. So, we need a further
analysis of the 2-categorical structure of cartesian models, model
homomorphisms, and their transformations. To this end, we observe that the
components of a modulation between cartesian lax natural transformations
automatically preserve products, as the following lemma shows.

\begin{lemma}[Modulation components for products] \label{lemma:modulations-preserve-products}
  Let $\dbl{D}$ and $\dbl{E}$ be precartesian double categories, let
  $F,G: \dbl{D} \to \dbl{E}$ be lax double functors, let
  $\alpha, \beta: F \To G$ be cartesian lax natural transformations, and let
  $\mu: \alpha \Tto \beta$ be a modulation. Then for any proarrows
  $m: x \proto y$ and $m': x' \proto y'$ in $\dbl{D}$,
  \begin{equation*}
    \begin{tikzcd}[column sep=large]
      {F(x \times x')} & {F(y \times y')} \\
      {G(x \times x')} & {G(y \times y')} \\
      {Gx \times Gx'} & {Gy \times Gy'}
      \arrow[""{name=0, anchor=center, inner sep=0}, "{F(m \times m')}", "\shortmid"{marking}, from=1-1, to=1-2]
      \arrow["{\alpha_{x \times x'}}"', from=1-1, to=2-1]
      \arrow["{\beta_{y \times y'}}", from=1-2, to=2-2]
      \arrow[""{name=1, anchor=center, inner sep=0}, "{G(m \times m')}"', "\shortmid"{marking}, from=2-1, to=2-2]
      \arrow["{\Psi_{x,x'}}"', from=2-1, to=3-1]
      \arrow[""{name=2, anchor=center, inner sep=0}, "{Gm \times Gm'}"', "\shortmid"{marking}, from=3-1, to=3-2]
      \arrow["{\Phi_{y,y'}}", from=2-2, to=3-2]
      \arrow["{\mu_{m \times m'}}"{description}, draw=none, from=0, to=1]
      \arrow["{\Psi_{m,m'}}"{description}, draw=none, from=1, to=2]
    \end{tikzcd}
    \quad=\quad
    \begin{tikzcd}[column sep=large]
      {F(x \times x')} & {F(y \times y')} \\
      {Fx \times Fx'} & {Fy \times Fy'} \\
      {Gx \times Gx'} & {Gy \times Gy'}
      \arrow[""{name=0, anchor=center, inner sep=0}, "{F(m \times m')}", "\shortmid"{marking}, from=1-1, to=1-2]
      \arrow["{\Phi_{x,x'}}"', from=1-1, to=2-1]
      \arrow["{\Phi_{y,y'}}", from=1-2, to=2-2]
      \arrow[""{name=1, anchor=center, inner sep=0}, "{Fm \times Fm'}"', "\shortmid"{marking}, from=2-1, to=2-2]
      \arrow["{\alpha_x \times \alpha_{x'}}"', from=2-1, to=3-1]
      \arrow["{\beta_y \times \beta_{y'}}", from=2-2, to=3-2]
      \arrow[""{name=2, anchor=center, inner sep=0}, "{Gm \times Gm'}"', "\shortmid"{marking}, from=3-1, to=3-2]
      \arrow["{\Phi_{m,m'}}"{description}, draw=none, from=0, to=1]
      \arrow["{\mu_m \times \mu_{m'}}"{description}, draw=none, from=1, to=2]
    \end{tikzcd}.
  \end{equation*}
  Also, for any objects $x$ and $x'$ in $\dbl{D}$,
  \begin{equation*}
    \begin{tikzcd}[column sep=large]
      {F(x \times x')} & {F(x \times x')} \\
      {G(x \times x')} & {G(x \times x')} \\
      {G(x \times x')} & {G(x \times x')} \\
      {Gx \times Gx'} & {Gx \times Gx'}
      \arrow[""{name=0, anchor=center, inner sep=0}, "{\id_{F(x \times x')}}", "\shortmid"{marking}, from=1-1, to=1-2]
      \arrow["{\alpha_{x \times x'}}"', from=1-1, to=2-1]
      \arrow["{\beta_{x \times x'}}", from=1-2, to=2-2]
      \arrow[""{name=1, anchor=center, inner sep=0}, "{G \id_{x \times x'}}", "\shortmid"{marking}, from=2-1, to=2-2]
      \arrow[Rightarrow, no head, from=2-1, to=3-1]
      \arrow[Rightarrow, no head, from=2-2, to=3-2]
      \arrow[""{name=2, anchor=center, inner sep=0}, "{G(\id_x \times \id_{x'})}", "\shortmid"{marking}, from=3-1, to=3-2]
      \arrow["{\Psi_{x,x'}}"', from=3-1, to=4-1]
      \arrow["{\Psi_{x,x'}}", from=3-2, to=4-2]
      \arrow[""{name=3, anchor=center, inner sep=0}, "{G\id_x \times G\id_{x'}}"', "\shortmid"{marking}, from=4-1, to=4-2]
      \arrow["{\mu_{x \times x'}}"{description, pos=0.4}, draw=none, from=0, to=1]
      \arrow["{G \times_{(x,x')}}"{description, pos=0.4}, draw=none, from=1, to=2]
      \arrow["{\Psi_{\id_x, \id_{x'}}}"{description}, draw=none, from=2, to=3]
    \end{tikzcd}
    \quad=\quad
    \begin{tikzcd}[column sep=large]
      {F(x \times x')} & {F(x \times x')} \\
      {Fx \times Fx'} & {Fx \times Fx'} \\
      {Fx \times Fx'} & {Fx \times Fx'} \\
      {Gx \times Gx'} & {Gx \times Gx'}
      \arrow[""{name=0, anchor=center, inner sep=0}, "{\mathrm{id}_{F(x \times x')}}", "\shortmid"{marking}, from=1-1, to=1-2]
      \arrow["{\Phi_{x,x'}}"', from=1-1, to=2-1]
      \arrow["{\Phi_{x,x'}}", from=1-2, to=2-2]
      \arrow[""{name=1, anchor=center, inner sep=0}, "{\mathrm{id}_{Fx \times Fx'}}", "\shortmid"{marking}, from=2-1, to=2-2]
      \arrow[Rightarrow, no head, from=2-2, to=3-2]
      \arrow[Rightarrow, no head, from=2-1, to=3-1]
      \arrow[""{name=2, anchor=center, inner sep=0}, "{\mathrm{id}_{Fx} \times \mathrm{id}_{Fx'}}", "\shortmid"{marking}, from=3-1, to=3-2]
      \arrow["{\alpha_x \times \alpha_{x'}}"', from=3-1, to=4-1]
      \arrow["{\beta_x \times \beta_{x'}}", from=3-2, to=4-2]
      \arrow[""{name=3, anchor=center, inner sep=0}, "{G\id_x \times G\id_{x'}}"', "\shortmid"{marking}, from=4-1, to=4-2]
      \arrow["{\mathrm{id}_{\Phi_{x,x'}}}"{description, pos=0.4}, draw=none, from=0, to=1]
      \arrow["{\mu_x \times \mu_{x'}}"{description}, draw=none, from=2, to=3]
      \arrow["{\times_{(Fx, Fx')}}"{description, pos=0.4}, draw=none, from=1, to=2]
    \end{tikzcd}.
  \end{equation*}
  In particular, when the double category $\dbl{D}$ and lax functor $G$ are both
  cartesian, the cell $\mu_{m \times m'}$ is completely determined by the
  product of the cells $\mu_m$ and $\mu_{m'}$, and the cell $\mu_{x \times x'}$
  is determined by the product of the cells $\mu_x$ and $\mu_{x'}$.
\end{lemma}
\begin{proof}
  Applying the naturality axiom \eqref{eq:modulation-naturality} of a modulation
  to the cell $\pi_{m,m'}: m \times m' \to m$ in $\dbl{D}$ and using that
  cartesian lax transformations strictly preserve products, we have
  \begin{equation*}
    \begin{tikzcd}[row sep=scriptsize, column sep=large]
      {F(x \times x')} & {F(y \times y')} \\
      {G(x \times x')} & {G(y \times y')} \\
      Gx & Gy
      \arrow[""{name=0, anchor=center, inner sep=0}, "{F(m \times m')}", "\shortmid"{marking}, from=1-1, to=1-2]
      \arrow["{\alpha_{x \times x'}}"', from=1-1, to=2-1]
      \arrow["{\beta_{y \times y'}}", from=1-2, to=2-2]
      \arrow[""{name=1, anchor=center, inner sep=0}, "{G(m \times m')}"', "\shortmid"{marking}, from=2-1, to=2-2]
      \arrow["{G \pi_{x,x'}}"', from=2-1, to=3-1]
      \arrow["{G\pi_{y,y'}}", from=2-2, to=3-2]
      \arrow[""{name=2, anchor=center, inner sep=0}, "Gm"', "\shortmid"{marking}, from=3-1, to=3-2]
      \arrow["{\mu_{m \times m'}}"{description}, draw=none, from=0, to=1]
      \arrow["{G\pi_{m,m'}}"{description}, draw=none, from=1, to=2]
    \end{tikzcd}
    \quad=\quad
    \begin{tikzcd}[row sep=scriptsize, column sep=large]
      {F(x \times x')} & {F(y \times y')} \\
      Fx & Fy \\
      Gx & Gy
      \arrow[""{name=0, anchor=center, inner sep=0}, "{F(m \times m')}", "\shortmid"{marking}, from=1-1, to=1-2]
      \arrow[""{name=1, anchor=center, inner sep=0}, "Fm"', "\shortmid"{marking}, from=2-1, to=2-2]
      \arrow["{F\pi_{x,x'}}"', from=1-1, to=2-1]
      \arrow["{F\pi_{y,y'}}", from=1-2, to=2-2]
      \arrow["{\alpha_x}"', from=2-1, to=3-1]
      \arrow["{\beta_y}", from=2-2, to=3-2]
      \arrow[""{name=2, anchor=center, inner sep=0}, "Gm"', "\shortmid"{marking}, from=3-1, to=3-2]
      \arrow["{F\pi_{m,m'}}"{description}, draw=none, from=0, to=1]
      \arrow["{\mu_m}"{description}, draw=none, from=1, to=2]
    \end{tikzcd}.
  \end{equation*}
  Similarly, applying the naturality axiom to the cell
  $\pi_{m,m'}: m \times m' \to m'$ yields an equation relating the cells
  $\mu_{m \times m'}$ and $\mu_{m'}$. The pairing of these two equations is the
  first equation in the lemma statement.

  Next, for any objects $x$ and $x'$ in $\dbl{D}$, we have
  \begin{equation*}
    \begin{tikzcd}[row sep=scriptsize, column sep=large]
      {F(x \times x')} & {F(x \times x')} \\
      {G(x \times x')} & {G(x \times x')} \\
      Gx & Gx
      \arrow[""{name=0, anchor=center, inner sep=0}, "{F \id_{x \times x'}}", "\shortmid"{marking}, from=1-1, to=1-2]
      \arrow["{\alpha_{x \times x'}}"', from=1-1, to=2-1]
      \arrow["{\beta_{x \times x'}}", from=1-2, to=2-2]
      \arrow[""{name=1, anchor=center, inner sep=0}, "{G \id_{x \times x'}}", "\shortmid"{marking}, from=2-1, to=2-2]
      \arrow[""{name=2, anchor=center, inner sep=0}, "{G \id_x}"', "\shortmid"{marking}, from=3-1, to=3-2]
      \arrow["{G\pi_{x,x'}}"', from=2-1, to=3-1]
      \arrow["{G\pi_{x,x'}}", from=2-2, to=3-2]
      \arrow["{\mu_{\id_{x \times x'}}}"{description, pos=0.4}, draw=none, from=0, to=1]
      \arrow["{G\id_{\pi_{x,x'}}}"{description}, draw=none, from=1, to=2]
    \end{tikzcd}
    \quad=\quad
    \begin{tikzcd}[row sep=scriptsize, column sep=large]
      {F(x \times x')} & {F(x \times x')} \\
      Fx & Fx \\
      Gx & Gx
      \arrow[""{name=0, anchor=center, inner sep=0}, "{F\id_x}", "\shortmid"{marking}, from=2-1, to=2-2]
      \arrow[""{name=1, anchor=center, inner sep=0}, "{G\id_x}"', "\shortmid"{marking}, from=3-1, to=3-2]
      \arrow["{\alpha_x}"', from=2-1, to=3-1]
      \arrow["{\beta_x}", from=2-2, to=3-2]
      \arrow["{F\pi_{x,x'}}"', from=1-1, to=2-1]
      \arrow["{F\pi_{x,x'}}", from=1-2, to=2-2]
      \arrow[""{name=2, anchor=center, inner sep=0}, "{F\id_{x \times x'}}", "\shortmid"{marking}, from=1-1, to=1-2]
      \arrow["{F\id_{\pi_{x,x'}}}"{description, pos=0.4}, draw=none, from=2, to=0]
      \arrow["{\mu_{\id_x}}"{description}, draw=none, from=0, to=1]
    \end{tikzcd}
  \end{equation*}
  by the naturality axiom \eqref{eq:modulation-naturality} at the cell
  $\id_{\pi_{x,x'}}: \id_{x \times x'} \to \id_x$ in $\dbl{D}$. Precomposing
  with the unitor $F_{x \times x'}$ and then using the naturality of unitors on
  the right-hand side, we obtain
  \begin{equation*}
    \begin{tikzcd}[row sep=scriptsize, column sep=large]
      {F(x \times x')} & {F(x \times x')} \\
      {G(x \times x')} & {G(x \times x')} \\
      Gx & Gx
      \arrow[""{name=0, anchor=center, inner sep=0}, "{\id_{F(x \times x')}}", "\shortmid"{marking}, from=1-1, to=1-2]
      \arrow["{\alpha_{x \times x'}}"', from=1-1, to=2-1]
      \arrow["{\beta_{x \times x'}}", from=1-2, to=2-2]
      \arrow[""{name=1, anchor=center, inner sep=0}, "{G \id_{x \times x'}}", "\shortmid"{marking}, from=2-1, to=2-2]
      \arrow[""{name=2, anchor=center, inner sep=0}, "{G \id_x}"', "\shortmid"{marking}, from=3-1, to=3-2]
      \arrow["{G\pi_{x,x'}}"', from=2-1, to=3-1]
      \arrow["{G\pi_{x,x'}}", from=2-2, to=3-2]
      \arrow["{\mu_{x \times x'}}"{description, pos=0.4}, draw=none, from=0, to=1]
      \arrow["{G\id_{\pi_{x,x'}}}"{description}, draw=none, from=1, to=2]
    \end{tikzcd}
    \quad=\quad
    \begin{tikzcd}[row sep=scriptsize, column sep=large]
      {F(x \times x')} & {F(x \times x')} \\
      Fx & Fx \\
      Gx & Gx
      \arrow[""{name=0, anchor=center, inner sep=0}, "{\id_{Fx}}", "\shortmid"{marking}, from=2-1, to=2-2]
      \arrow[""{name=1, anchor=center, inner sep=0}, "{G\id_x}"', "\shortmid"{marking}, from=3-1, to=3-2]
      \arrow["{\alpha_x}"', from=2-1, to=3-1]
      \arrow["{\beta_x}", from=2-2, to=3-2]
      \arrow["{F\pi_{x,x'}}"', from=1-1, to=2-1]
      \arrow["{F\pi_{x,x'}}", from=1-2, to=2-2]
      \arrow[""{name=2, anchor=center, inner sep=0}, "{\id_{F(x \times x')}}", "\shortmid"{marking}, from=1-1, to=1-2]
      \arrow["{\mu_x}"{description}, draw=none, from=0, to=1]
      \arrow["{\id_{F\pi_{x,x'}}}"{description, pos=0.4}, draw=none, from=2, to=0]
    \end{tikzcd}.
  \end{equation*}
  Starting from the cell $\id_{\pi_{x,x'}'}: \id_{x \times x'} \to \id_{x'}$ in
  $\dbl{D}$ yields a similar equation relating $\mu_{x \times x'}$ and
  $\mu_{x'}$. The pairing these two equations is the second equation in the
  lemma statement, after applying
  \cref{eq:functor-product-id-comparison-1,eq:functor-product-id-comparison-2}.
\end{proof}

Modulations recover a number of well-known examples of 2-cells between morphisms
of structures that have already been seen to be lax transformations of models of
cartesian double theories.

\begin{example}[Monoidal transformations] \label{ex:modulations-are-monoidal-transformations}
  For the various axiomatiziations of monoidal categories (namely,
  \cref{th:monoid}, \cref{th:pseudomonoid}, and
  \cref{th:internal-commutative-comonoid}), a modulation between cartesian lax
  transformations amounts to a monoidal transformation between lax monoidal
  functors, and conversely. This follows readily by the lemma and
  \cref{rmk:modulations-reduce-to-modications}, owing to the fact that all these
  theories are vertically trivial double categories. In more detail, a
  modulation $\mu\colon \phi\Rrightarrow\psi$ is effectively an ordinary
  modification of lax transformations of the associated 2-functors. The single
  object $x$ yields a component $\mu_x\colon \phi_x\Rightarrow\psi_x$ which is a
  cell in $\Cat$, hence an ordinary natural transformation. Now, in particular
  the modification condition applied to the arrow $\otimes\colon x^2\to x$
  yields
  \begin{equation*}
    \begin{tikzcd}
      {\cat{M}^2} && {\cat{M}} \\
      {\cat{N}^2} && {\cat{N}}
      \arrow["{\otimes_{\cat{M}}}", from=1-1, to=1-3]
      \arrow[""{name=0, anchor=center, inner sep=0}, "{\psi_x}", from=1-3, to=2-3]
      \arrow[""{name=1, anchor=center, inner sep=0}, "{\phi^2_x}"', curve={height=12pt}, from=1-1, to=2-1]
      \arrow["{\otimes_{\cat{N}}}"', from=2-1, to=2-3]
      \arrow[""{name=2, anchor=center, inner sep=0}, "{\psi^2_x}", curve={height=-12pt}, from=1-1, to=2-1]
      \arrow["{\psi_\otimes}", shorten <=18pt, shorten >=12pt, Rightarrow, from=2, to=0]
      \arrow["{\mu_{x^2}}", shorten <=5pt, shorten >=5pt, Rightarrow, from=1, to=2]
    \end{tikzcd}
    \quad=\quad
    \begin{tikzcd}
      {\cat{M}^2} && {\cat{M}} \\
      {\cat{N}^2} && {\cat{N}}
      \arrow[""{name=0, anchor=center, inner sep=0}, "{\phi^2_x}"', from=1-1, to=2-1]
      \arrow["{\otimes_{\cat{M}}}", from=1-1, to=1-3]
      \arrow["{\otimes_{\cat{N}}}"', from=2-1, to=2-3]
      \arrow[""{name=1, anchor=center, inner sep=0}, "{\psi_x}", curve={height=-12pt}, from=1-3, to=2-3]
      \arrow[""{name=2, anchor=center, inner sep=0}, "{\phi_x}"', curve={height=12pt}, from=1-3, to=2-3]
      \arrow["{\phi_\otimes}", shorten <=12pt, shorten >=18pt, Rightarrow, from=0, to=2]
      \arrow["{\mu_x}", shorten <=5pt, shorten >=5pt, Rightarrow, from=2, to=1]
    \end{tikzcd}
  \end{equation*}
  which at each component is precisely the equation stating that the components
  of $\mu$ interact correctly with the monoidal product as in the definition of
  a monoidal transformation \cite[\S{XI.2}]{maclane1998}. Applied to the unit
  morphism, the unit condition follows as well. Conversely, each monoidal
  natural transformation determines a modulation.
\end{example}

\begin{example}[Multinaturality] \label{ex:multinaturality}
  \cref{example:lax-transfs-are-multifunctors} showed that cartesian lax
  transformations of models of the theory of promonoids are precisely
  multifunctors between multicategories. Modulations of such transformations are
  precisely natural transformations of multifunctors satisfying the
  \emph{multinaturality} condition described by Hermida \cite[Definition
  6.6]{hermida2000} in the definition of the 2-category of multicategories.
  Given such a modulation $\mu\colon \phi\Rrightarrow \psi$, viewing the models
  as span-valued, we have a cell
  \begin{equation*}
    \begin{tikzcd}
      {\cat{C}_0^n} & {\cat{C}_0^n} & {\cat{C}_0^n} \\
      {\cat{D}_0^n} & {\cat{D}_1^n} & {\cat{D}_0^n}
      \arrow[from=1-2, to=1-1]
      \arrow["{\mu_n}", from=1-2, to=2-2]
      \arrow["{\phi_{x^n}}"', from=1-1, to=2-1]
      \arrow[from=2-2, to=2-1]
      \arrow[from=1-2, to=1-3]
      \arrow["{\psi_{x^n}}", from=1-3, to=2-3]
      \arrow[from=2-2, to=2-3]
    \end{tikzcd}
  \end{equation*}
  which has the effect of associating to every $n$-tuple of objects of
  $\cat{C}$, say $a_1,\dots a_n$, an $n$-tuple of unary morphisms
  \begin{equation*}
    \langle \mu_1,\dots,\mu_n\rangle\colon \langle \phi(a_1),\dots,\phi(a_n)\rangle
    \to \langle \psi(a_1),\dots,\psi(a_n)\rangle.
  \end{equation*}
  The equivariance axiom in \cref{def:modulation-special} applied at the
  proarrow $p_n\colon x^n\proto x$ then ensures the multinaturality condition,
  namely, that the square
  \begin{equation*}
    \begin{tikzcd}
      {\langle \phi(a_1),\dots,\phi(a_n)\rangle} && {\phi(a)} \\
      {\langle \psi(a_1),\dots,\psi(a_n)\rangle} && {\psi(a)}
      \arrow["{\phi(m)}", from=1-1, to=1-3]
      \arrow["{\mu_a}", from=1-3, to=2-3]
      \arrow["{\langle \mu_1,\dots,\mu_n\rangle}"', from=1-1, to=2-1]
      \arrow["{\psi(m)}"', from=2-1, to=2-3]
    \end{tikzcd}
  \end{equation*}
  commutes for any multimorphism $m\colon a_1,\dots, a_n \to a$ of $\cat{C}$.
\end{example}

The fact, proved in \cref{lemma:modulations-preserve-products}, that components
of modulations between cartesian lax transformations preserve products motivates
the following definition and theorem.

\begin{theorem}[2-category of cartesian lax functors]
  \label{thm:cartesian-lax-functor-2-category}
  For any cartesian double categories $\dbl{D}$ and $\dbl{E}$, there is a
  2-category $\CartLaxFun_\ell(\dbl{D}, \dbl{E})$ whose objects are cartesian
  lax double functors $\dbl{D} \to \dbl{E}$, morphisms are cartesian lax natural
  transformations, and 2-morphisms are modulations.
\end{theorem}
\begin{proof}
  In view of \cref{thm:lax-functor-2-category}, we just need to check that
  cartesian lax natural transformations are closed under composition. This is
  clear by pasting the strict naturality squares in
  \cref{def:cartesian-lax-transformation}.
\end{proof}

We immediately obtain the sought after 2-category of models of a cartesian
double theory.

\begin{corollary}[2-category of models] \label{cor:models-2-category}
  Let $\dbl{T}$ be a cartesian double theory and let $\dbl{S}$ be a cartesian
  double category. A \define{lax} (resp.\ \define{pseudo}, resp.\
  \define{strict}) \define{map} between models $M$ and $M'$ of $\dbl{T}$ in
  $\dbl{S}$ is a cartesian lax (resp.\ pseudo, resp. strict) natural
  transformation $M \To M'$. A \define{transformation} between maps of models is
  a modulation.

  There is a 2-category whose objects are models of the theory $\dbl{T}$ in
  $\dbl{S}$, morphisms are (lax, pseudo, or strict) maps between models, and
  2-morphisms are transformations between maps.
\end{corollary}

Similarly, a restriction sketch has a 2-category of models in any cartesian
equipment, which is a full sub-2-category of the 2-category of models of its
underlying cartesian double theory. Although it seems impossible to prove these
definitions ``correct'' in all cases, the following example offers favorable
evidence. A further corollary helps with the calculations.

\begin{corollary}
  If $\dbl{D}$ is a double category and $\dbl{E}$ is a cartesian equipment, then
  the bijection of \cref{cor:unitalizationofCartesianLaxFunctors} extends to an
  isomophism of 2-categories
    \begin{equation*}
      \CartLaxFun_\ell(\dbl{D}, \dbl{E}) \xrightarrow{\cong} 
      \CartLaxFun_{\ell,u}(\dbl{D}, \Module{\dbl{E}})
    \end{equation*}
  specializing that of \cref{th:unitalizationFINALFORM}. In particular, there is
  an isomorphism of 2-categories
    \begin{equation*}
      \CartLaxFun_\ell(\dbl{D},\Span) \xrightarrow{\cong}
      \CartLaxFun_{\ell,u}(\dbl{D},\Prof)
    \end{equation*}
  specializing \cref{cor:unitalize-span-valued-2-cat-iso-FINAL}.
\end{corollary}
\begin{proof}
  There is nothing to prove at the level of modulations by
  \cref{lemma:modulations-preserve-products}. But cartesian lax transformations
  on one side of the purported isomorphism correspond to cartesian lax
  transformations on the other side by the construction of products in
  $\Module{\dbl{E}}$. The final statement follows again by taking
  $\dbl{E} = \Span$.
\end{proof}

\begin{example}[2-category of multicategories]
  Returning to \cref{ex:multinaturality} about natural transformations between
  multifunctors, we have that the 2-category of models of the theory of
  promonoids (\cref{th:promonoid}) is precisely the 2-category of
  multicategories considered by Hermida \cite{hermida2000}.
\end{example}

Extending \cref{cor:1-categories-of-cartesian-models-on-2-category},
we have the following.

\begin{corollary} \label{cor:2-categories-of-models-cartesian}
  For any cartesian 2-category $\bicat{A}$ and cartesian equipment $\dbl{E}$
  with local coequalizers, there is an isomorphism of 2-categories
  \begin{equation*}
    \CartLaxFun_\ell(\VerDbl(\bicat{A}),\dbl{E})
    \xrightarrow{\cong}
    \CartTwoCat_\ell(\bicat{A},\Cat(\dbl{E})).
  \end{equation*}
  In particular, there is an isomorphism of 2-categories
  \begin{equation*}
    \CartLaxFun_\ell(\VerDbl(\bicat{A}),\Span)
    \xrightarrow{\cong}
    \CartTwoCat_\ell(\bicat{A},\Cat).
  \end{equation*}
\end{corollary}

A paradigmatic example of a cartesian double theory is that of pseudomonoids
(\cref{th:pseudomonoid}), which is 2-categorical in nature, having no
non-identity proarrows. Thus, the corollary confirms the following example,
concluding the thread began with
\cref{example:lax-transfs-are-lax-monoidal-funcs} and
\cref{ex:modulations-are-monoidal-transformations}.

\begin{example}[2-category of monoidal categories]
  The 2-category of models of the theory of pseudomonoids
  (\cref{th:pseudomonoid}) with lax (resp.\ pseudo, resp.\ strict) maps is
  precisely the 2-category of monoidal categories, lax (resp.\ pseudo, resp.\
  strict) monoidal functors, and monoidal natural transformations.
\end{example}

\section{Modules and modulations}
\label{sec:modules}

Modules between lax double functors, and modulations between a square of
transformations and modules, were introduced by Paré \cite{pare2011,pare2013},
generalizing the corresponding definitions for bicategories \cite{cockett2003}.
In this section, we review the definitions of modules and modulations, define
the notion of a cartesian module between cartesian lax functors, and consider
several examples of cartesian modules between models of cartesian double
theories. Denoting by $\dbl{I}$ the strict double category freely generated by a
single proarrow $0 \proto 1$, a module can be succinctly defined as follows
\cite[Remark 3.5]{pare2011}.

\begin{definition}[Module] \label{def:module}
  A \define{module} $M: F \proTo G$ between two lax double functors
  $F,G: \dbl{D} \to \dbl{E}$ is a lax double functor
  $M: \dbl{D} \times \dbl{I} \to \dbl{E}$ such that $M(-,0) = F$ and
  $M(-,1) = G$.
\end{definition}

When this definition is fully unpacked as in \cite[Definition 3.2]{pare2011}, a
module $M: F \proTo G$ is seen to consist of
\begin{itemize}
  \item for every proarrow $m: x \proto y$ in $\dbl{D}$, a proarrow
  $M(m): Fx \proto Gy$ in $\dbl{E}$;
  \item for every cell $\stdInlineCell{\alpha}$ in $\dbl{D}$, a cell in
    $\dbl{E}$
  \begin{equation*}
    \begin{tikzcd}
      Fx & Gy \\
      Fw & Gz
      \arrow["Ff"', from=1-1, to=2-1]
      \arrow[""{name=0, anchor=center, inner sep=0}, "{M(m)}", "\shortmid"{marking}, from=1-1, to=1-2]
      \arrow[""{name=1, anchor=center, inner sep=0}, "{M(n)}"', "\shortmid"{marking}, from=2-1, to=2-2]
      \arrow["Gg", from=1-2, to=2-2]
      \arrow["{M(\alpha)}"{description}, draw=none, from=0, to=1]
    \end{tikzcd};
  \end{equation*}
  \item for every consecutive pair of proarrows $x \xproto{m} y \xproto{n} z$ in
  $\dbl{D}$, globular cells in $\dbl{E}$
  \begin{equation*}
    \begin{tikzcd}
      Fx & Fy & Gz \\
      Fx && Gz
      \arrow["{F(m)}", "\shortmid"{marking}, from=1-1, to=1-2]
      \arrow["{M(n)}", "\shortmid"{marking}, from=1-2, to=1-3]
      \arrow[Rightarrow, no head, from=1-1, to=2-1]
      \arrow[Rightarrow, no head, from=1-3, to=2-3]
      \arrow[""{name=0, anchor=center, inner sep=0}, "{M(m \odot n)}"', "\shortmid"{marking}, from=2-1, to=2-3]
      \arrow["{M^\ell_{m,n}}"{description}, draw=none, from=1-2, to=0]
    \end{tikzcd}
    \qquad\text{and}\qquad
    \begin{tikzcd}
      Fx & Gy & Gz \\
      Fx && Gz
      \arrow["{M(m)}", "\shortmid"{marking}, from=1-1, to=1-2]
      \arrow["{G(n)}", "\shortmid"{marking}, from=1-2, to=1-3]
      \arrow[Rightarrow, no head, from=1-1, to=2-1]
      \arrow[Rightarrow, no head, from=1-3, to=2-3]
      \arrow[""{name=0, anchor=center, inner sep=0}, "{M(m \odot n)}"', "\shortmid"{marking}, from=2-1, to=2-3]
      \arrow["{M^r_{m,n}}"{description}, draw=none, from=1-2, to=0]
    \end{tikzcd},
  \end{equation*}
  the \define{left} and \define{right actions} of $M$.
\end{itemize}
The following axioms must be satisfied.
\begin{itemize}
  \item Functorality on cells: for any cells
    $\inlineCell{x}{x'}{y}{y'}{m}{n}{f}{f'}{\alpha}$ and
    $\inlineCell{y}{y'}{z}{z'}{n}{p}{g}{g'}{\beta}$ in $\dbl{D}$,
    \begin{equation*}
      \begin{tikzcd}
        Fx & {Gx'} \\
        Fy & {Gy'} \\
        Fz & {Gz'}
        \arrow["Ff"', from=1-1, to=2-1]
        \arrow[""{name=0, anchor=center, inner sep=0}, "{M(m)}", "\shortmid"{marking}, from=1-1, to=1-2]
        \arrow[""{name=1, anchor=center, inner sep=0}, "{M(n)}"', "\shortmid"{marking}, from=2-1, to=2-2]
        \arrow["{Gf'}", from=1-2, to=2-2]
        \arrow["Fg"', from=2-1, to=3-1]
        \arrow["{Gg'}", from=2-2, to=3-2]
        \arrow[""{name=2, anchor=center, inner sep=0}, "{M(p)}"', "\shortmid"{marking}, from=3-1, to=3-2]
        \arrow["{M(\alpha)}"{description}, draw=none, from=0, to=1]
        \arrow["{M(\beta)}"{description}, draw=none, from=1, to=2]
      \end{tikzcd}
      \quad=\quad
      \begin{tikzcd}
        Fx & {Gx'} \\
        Fy & {Gy'} \\
        Fz & {Gz'}
        \arrow["Ff"', from=1-1, to=2-1]
        \arrow[""{name=0, anchor=center, inner sep=0}, "{M(m)}", "\shortmid"{marking}, from=1-1, to=1-2]
        \arrow["{Gf'}", from=1-2, to=2-2]
        \arrow["Fg"', from=2-1, to=3-1]
        \arrow["{Gg'}", from=2-2, to=3-2]
        \arrow[""{name=1, anchor=center, inner sep=0}, "{M(p)}"', "\shortmid"{marking}, from=3-1, to=3-2]
        \arrow["{M(\alpha \cdot \beta)}"{description}, draw=none, from=0, to=1]
      \end{tikzcd},
    \end{equation*}
    and $M(1_m) = 1_{M(m)}$ for any proarrow $m: x \proto y$ in $\dbl{D}$.
  \item Naturality of actions: for any cells
    $\inlineCell{x}{y}{x'}{y'}{m}{m'}{f}{g}{\alpha}$ and
    $\inlineCell{y}{z}{y'}{z'}{n}{n'}{g}{h}{\beta}$ in $\dbl{D}$,
    \begin{equation*}
      \begin{tikzcd}
        Fx & Fy & Gz \\
        {Fx'} & {Fy'} & {Gz'} \\
        {Fx'} && {Gz'}
        \arrow[""{name=0, anchor=center, inner sep=0}, "{F(m')}"', "\shortmid"{marking}, from=2-1, to=2-2]
        \arrow[""{name=1, anchor=center, inner sep=0}, "{M(n')}"', "\shortmid"{marking}, from=2-2, to=2-3]
        \arrow[Rightarrow, no head, from=2-1, to=3-1]
        \arrow[Rightarrow, no head, from=2-3, to=3-3]
        \arrow[""{name=2, anchor=center, inner sep=0}, "{M(m' \odot n')}"', "\shortmid"{marking}, from=3-1, to=3-3]
        \arrow["Fg"{description}, from=1-2, to=2-2]
        \arrow["Ff"', from=1-1, to=2-1]
        \arrow[""{name=3, anchor=center, inner sep=0}, "{F(m)}", "\shortmid"{marking}, from=1-1, to=1-2]
        \arrow[""{name=4, anchor=center, inner sep=0}, "{M(n)}", "\shortmid"{marking}, from=1-2, to=1-3]
        \arrow["Gh", from=1-3, to=2-3]
        \arrow["{M^\ell_{m',n'}}"{description}, draw=none, from=2-2, to=2]
        \arrow["{F(\alpha)}"{description}, draw=none, from=3, to=0]
        \arrow["{M(\beta)}"{description}, draw=none, from=4, to=1]
      \end{tikzcd}
      \quad=\quad
      \begin{tikzcd}
        Fx & Fy & Gz \\
        Fx && Gz \\
        {Fx'} && {Gz'}
        \arrow["{F(m)}", "\shortmid"{marking}, from=1-1, to=1-2]
        \arrow["{M(n)}", "\shortmid"{marking}, from=1-2, to=1-3]
        \arrow[Rightarrow, no head, from=1-1, to=2-1]
        \arrow[Rightarrow, no head, from=1-3, to=2-3]
        \arrow[""{name=0, anchor=center, inner sep=0}, "{M(m \odot n)}"', "\shortmid"{marking}, from=2-1, to=2-3]
        \arrow["Ff"', from=2-1, to=3-1]
        \arrow["Gh", from=2-3, to=3-3]
        \arrow[""{name=1, anchor=center, inner sep=0}, "{M(m' \odot n')}"', "\shortmid"{marking}, from=3-1, to=3-3]
        \arrow["{M^\ell_{m,n}}"{description}, draw=none, from=1-2, to=0]
        \arrow["{M(\alpha \odot \beta)}"{description, pos=0.6}, draw=none, from=0, to=1]
      \end{tikzcd}
    \end{equation*}
    and
    \begin{equation*}
      \begin{tikzcd}
        Fx & Gy & Gz \\
        {Fx'} & {Gy'} & {Gz'} \\
        {Fx'} && Gz
        \arrow[""{name=0, anchor=center, inner sep=0}, "{M(m')}"', "\shortmid"{marking}, from=2-1, to=2-2]
        \arrow[""{name=1, anchor=center, inner sep=0}, "{G(n')}"', "\shortmid"{marking}, from=2-2, to=2-3]
        \arrow[Rightarrow, no head, from=2-1, to=3-1]
        \arrow[Rightarrow, no head, from=2-3, to=3-3]
        \arrow[""{name=2, anchor=center, inner sep=0}, "{M(m' \odot n')}"', "\shortmid"{marking}, from=3-1, to=3-3]
        \arrow["Ff"', from=1-1, to=2-1]
        \arrow["Gg"{description}, from=1-2, to=2-2]
        \arrow["Gh", from=1-3, to=2-3]
        \arrow[""{name=3, anchor=center, inner sep=0}, "{G(n)}", "\shortmid"{marking}, from=1-2, to=1-3]
        \arrow[""{name=4, anchor=center, inner sep=0}, "{M(m)}", "\shortmid"{marking}, from=1-1, to=1-2]
        \arrow["{M^r_{m',n'}}"{description}, draw=none, from=2-2, to=2]
        \arrow["{G(\beta)}"{description}, draw=none, from=3, to=1]
        \arrow["{M(\alpha)}"{description}, draw=none, from=4, to=0]
      \end{tikzcd}
      \quad=\quad
      \begin{tikzcd}
        Fx & Gy & Gz \\
        Fx && Gz \\
        {Fx'} && {Gz'}
        \arrow["{M(m)}", "\shortmid"{marking}, from=1-1, to=1-2]
        \arrow["{G(n)}", "\shortmid"{marking}, from=1-2, to=1-3]
        \arrow[Rightarrow, no head, from=1-1, to=2-1]
        \arrow[Rightarrow, no head, from=1-3, to=2-3]
        \arrow[""{name=0, anchor=center, inner sep=0}, "{M(m \odot n)}"', "\shortmid"{marking}, from=2-1, to=2-3]
        \arrow["Ff"', from=2-1, to=3-1]
        \arrow["Gh", from=2-3, to=3-3]
        \arrow[""{name=1, anchor=center, inner sep=0}, "{M(m' \odot n')}"', "\shortmid"{marking}, from=3-1, to=3-3]
        \arrow["{M^r_{m,n}}"{description}, draw=none, from=1-2, to=0]
        \arrow["{M(\alpha \odot \beta)}"{description, pos=0.6}, draw=none, from=0, to=1]
      \end{tikzcd}.
    \end{equation*}
  \item Associativity and unitality of left actions: for any triple of proarrows
    $w \xproto{m} x \xproto{n} y \xproto{p} z$,
    \begin{equation*}
      \begin{tikzcd}
        Fw & Fx & Fy & Gz \\
        Fw && Fy & Gz \\
        Fw &&& Gz
        \arrow["Fm", "\shortmid"{marking}, from=1-1, to=1-2]
        \arrow[""{name=0, anchor=center, inner sep=0}, "Fn", "\shortmid"{marking}, from=1-2, to=1-3]
        \arrow[""{name=1, anchor=center, inner sep=0}, "Mp", "\shortmid"{marking}, from=1-3, to=1-4]
        \arrow[""{name=2, anchor=center, inner sep=0}, "Mp"', "\shortmid"{marking}, from=2-3, to=2-4]
        \arrow[Rightarrow, no head, from=1-3, to=2-3]
        \arrow[Rightarrow, no head, from=1-4, to=2-4]
        \arrow[""{name=3, anchor=center, inner sep=0}, "{F(m \odot n)}"', "\shortmid"{marking}, from=2-1, to=2-3]
        \arrow[Rightarrow, no head, from=1-1, to=2-1]
        \arrow[""{name=4, anchor=center, inner sep=0}, "{M(m \odot n \odot p)}"', "\shortmid"{marking}, from=3-1, to=3-4]
        \arrow[Rightarrow, no head, from=2-4, to=3-4]
        \arrow[Rightarrow, no head, from=2-1, to=3-1]
        \arrow["{1_{Mp}}"{description}, draw=none, from=1, to=2]
        \arrow["{F_{m,n}}"{description}, draw=none, from=1-2, to=3]
        \arrow["{M^\ell_{m \odot n, p}}"{description, pos=0.8}, draw=none, from=0, to=4]
      \end{tikzcd}
      \quad=\quad
      \begin{tikzcd}
        Fw & Fx & Fy & Gz \\
        Fw & Fx && Gz \\
        Fw &&& Gz
        \arrow[""{name=0, anchor=center, inner sep=0}, "Fm", "\shortmid"{marking}, from=1-1, to=1-2]
        \arrow[""{name=1, anchor=center, inner sep=0}, "Fn", "\shortmid"{marking}, from=1-2, to=1-3]
        \arrow["Mp", "\shortmid"{marking}, from=1-3, to=1-4]
        \arrow[""{name=2, anchor=center, inner sep=0}, "{M(n \odot p)}"', "\shortmid"{marking}, from=2-2, to=2-4]
        \arrow[Rightarrow, no head, from=1-1, to=2-1]
        \arrow[""{name=3, anchor=center, inner sep=0}, "{M(m \odot n \odot p)}"', "\shortmid"{marking}, from=3-1, to=3-4]
        \arrow[Rightarrow, no head, from=2-4, to=3-4]
        \arrow[Rightarrow, no head, from=2-1, to=3-1]
        \arrow[""{name=4, anchor=center, inner sep=0}, "Fm"', "\shortmid"{marking}, from=2-1, to=2-2]
        \arrow[Rightarrow, no head, from=1-2, to=2-2]
        \arrow[Rightarrow, no head, from=1-4, to=2-4]
        \arrow["{M^\ell_{m, n \odot p}}"{description, pos=0.8}, draw=none, from=1, to=3]
        \arrow["{1_{Fm}}"{description}, draw=none, from=0, to=4]
        \arrow["{M^\ell_{n,p}}"{description}, draw=none, from=1-3, to=2]
      \end{tikzcd},
    \end{equation*}
    and for every proarrow $m: x \proto y$,
    \begin{equation*}
      \begin{tikzcd}
        Fx & Fx & Gy \\
        Fx & Fx & Gy \\
        Fx && Gy
        \arrow[""{name=0, anchor=center, inner sep=0}, "{F \id_x}"', "\shortmid"{marking}, from=2-1, to=2-2]
        \arrow[""{name=1, anchor=center, inner sep=0}, "Mm"', "\shortmid"{marking}, from=2-2, to=2-3]
        \arrow[Rightarrow, no head, from=2-1, to=3-1]
        \arrow[Rightarrow, no head, from=2-3, to=3-3]
        \arrow[""{name=2, anchor=center, inner sep=0}, "Mm"', "\shortmid"{marking}, from=3-1, to=3-3]
        \arrow[""{name=3, anchor=center, inner sep=0}, "{\id_{Fx}}", "\shortmid"{marking}, from=1-1, to=1-2]
        \arrow[Rightarrow, no head, from=1-3, to=2-3]
        \arrow[Rightarrow, no head, from=1-1, to=2-1]
        \arrow[Rightarrow, no head, from=1-2, to=2-2]
        \arrow[""{name=4, anchor=center, inner sep=0}, "Mm", "\shortmid"{marking}, from=1-2, to=1-3]
        \arrow["{M^\ell_{x,m}}"{description}, draw=none, from=2-2, to=2]
        \arrow["{F_x}"{description}, draw=none, from=3, to=0]
        \arrow["{1_{Mm}}"{description}, draw=none, from=4, to=1]
      \end{tikzcd}
      \quad=\quad
      \begin{tikzcd}
        Fx & Gy \\
        Fx & Gy
        \arrow[""{name=0, anchor=center, inner sep=0}, "Mm"', "\shortmid"{marking}, from=2-1, to=2-2]
        \arrow[Rightarrow, no head, from=1-2, to=2-2]
        \arrow[Rightarrow, no head, from=1-1, to=2-1]
        \arrow[""{name=1, anchor=center, inner sep=0}, "Mm", "\shortmid"{marking}, from=1-1, to=1-2]
        \arrow["{1_{Mm}}"{description}, draw=none, from=1, to=0]
      \end{tikzcd}.
    \end{equation*}
  \item Associativity and unitality of right actions, dual to the previous
    axiom.
  \item Compatibility of left and right actions: for any triple of proarrows
    $w \xproto{m} x \xproto{n} y \xproto{p} z$,
    \begin{equation*}
      \begin{tikzcd}
        Fw & Fx & Gy & Gz \\
        Fw && Gy & Gz \\
        Fw &&& Gz
        \arrow[""{name=0, anchor=center, inner sep=0}, "Mn", "\shortmid"{marking}, from=1-2, to=1-3]
        \arrow["Fm", "\shortmid"{marking}, from=1-1, to=1-2]
        \arrow[""{name=1, anchor=center, inner sep=0}, "Gp", "\shortmid"{marking}, from=1-3, to=1-4]
        \arrow[""{name=2, anchor=center, inner sep=0}, "{M(m \odot n \odot p)}"', "\shortmid"{marking}, from=3-1, to=3-4]
        \arrow[Rightarrow, no head, from=1-1, to=2-1]
        \arrow[Rightarrow, no head, from=1-3, to=2-3]
        \arrow[Rightarrow, no head, from=2-1, to=3-1]
        \arrow[Rightarrow, no head, from=1-4, to=2-4]
        \arrow[Rightarrow, no head, from=2-4, to=3-4]
        \arrow[""{name=3, anchor=center, inner sep=0}, "{M(m \odot n)}"', "\shortmid"{marking}, from=2-1, to=2-3]
        \arrow[""{name=4, anchor=center, inner sep=0}, "Gp"', "\shortmid"{marking}, from=2-3, to=2-4]
        \arrow["{M^\ell_{m,n}}"{description}, draw=none, from=1-2, to=3]
        \arrow["{1_{Gp}}"{description}, draw=none, from=1, to=4]
        \arrow["{M^r_{m \odot n, p}}"{description, pos=0.8}, draw=none, from=0, to=2]
      \end{tikzcd}
      \quad=\quad
      \begin{tikzcd}
        Fw & Fx & Gy & Gz \\
        Fw & Fx && Gz \\
        Fw &&& Gz
        \arrow[""{name=0, anchor=center, inner sep=0}, "Mn", "\shortmid"{marking}, from=1-2, to=1-3]
        \arrow[""{name=1, anchor=center, inner sep=0}, "Fm", "\shortmid"{marking}, from=1-1, to=1-2]
        \arrow["Gp", "\shortmid"{marking}, from=1-3, to=1-4]
        \arrow[""{name=2, anchor=center, inner sep=0}, "{M(m \odot n \odot p)}"', "\shortmid"{marking}, from=3-1, to=3-4]
        \arrow[Rightarrow, no head, from=1-1, to=2-1]
        \arrow[Rightarrow, no head, from=2-1, to=3-1]
        \arrow[Rightarrow, no head, from=1-4, to=2-4]
        \arrow[Rightarrow, no head, from=2-4, to=3-4]
        \arrow[Rightarrow, no head, from=1-2, to=2-2]
        \arrow[""{name=3, anchor=center, inner sep=0}, "Fm"', "\shortmid"{marking}, from=2-1, to=2-2]
        \arrow[""{name=4, anchor=center, inner sep=0}, "{M(n \odot p)}"', "\shortmid"{marking}, from=2-2, to=2-4]
        \arrow["{1_{Fm}}"{description}, draw=none, from=1, to=3]
        \arrow["{M^r_{n,p}}"{description}, draw=none, from=1-3, to=4]
        \arrow["{M^\ell_{m, n \odot p}}"{description, pos=0.8}, draw=none, from=0, to=2]
      \end{tikzcd}.
    \end{equation*}
\end{itemize}

\begin{example}[Profunctors]
  A module between lax double functors $\dbl{1} \to \Span$, which we recall are
  the same thing as categories, is a profunctor between the categories. Indeed,
  modules were originally conceived as a kind of ``multi-object profunctor''
  \cite{pare2011}.
\end{example}

A module contains slightly less data than it might seem at first glance, as its
assignment on any cell out of an identity proarrow is already determined by its
assignment on external identity cells, along with the left or right actions.

\begin{lemma} \label{lem:module-cells}
  Let $\dbl{D}$ be a strict double category and let
  $M: F \proTo G: \dbl{D} \to \dbl{E}$ be a module between lax double functors.
  For any cell of the form $\inlineCell{x}{x}{y}{z}{\id_x}{n}{f}{g}{\alpha}$ in
  $\dbl{D}$,
  \begin{equation*}
    \begin{tikzcd}
      Fx & Fx & Gx \\
      Fx & Fx & Gx \\
      Fy & Fz & Gz \\
      Fy && Gz
      \arrow[""{name=0, anchor=center, inner sep=0}, "{M \id_x}", "\shortmid"{marking}, from=2-2, to=2-3]
      \arrow[""{name=1, anchor=center, inner sep=0}, "{F\id_x}", "\shortmid"{marking}, from=2-1, to=2-2]
      \arrow[""{name=2, anchor=center, inner sep=0}, "{M\id_z}"', "\shortmid"{marking}, from=3-2, to=3-3]
      \arrow[""{name=3, anchor=center, inner sep=0}, "Fn"', "\shortmid"{marking}, from=3-1, to=3-2]
      \arrow["Ff"', from=2-1, to=3-1]
      \arrow["Fg"{description}, from=2-2, to=3-2]
      \arrow["Gg", from=2-3, to=3-3]
      \arrow[""{name=4, anchor=center, inner sep=0}, "Mn"', "\shortmid"{marking}, from=4-1, to=4-3]
      \arrow[Rightarrow, no head, from=3-1, to=4-1]
      \arrow[Rightarrow, no head, from=3-3, to=4-3]
      \arrow[Rightarrow, no head, from=1-1, to=2-1]
      \arrow[Rightarrow, no head, from=1-2, to=2-2]
      \arrow[Rightarrow, no head, from=1-3, to=2-3]
      \arrow[""{name=5, anchor=center, inner sep=0}, "{ M\id_x}", "\shortmid"{marking}, from=1-2, to=1-3]
      \arrow[""{name=6, anchor=center, inner sep=0}, "{\id_{Fx}}", "\shortmid"{marking}, from=1-1, to=1-2]
      \arrow["{F \alpha}"{description}, draw=none, from=1, to=3]
      \arrow["{M\id_g}"{description}, draw=none, from=0, to=2]
      \arrow["{M^\ell_{n,z}}"{description}, draw=none, from=3-2, to=4]
      \arrow["{F_x}"{description}, draw=none, from=6, to=1]
      \arrow["1"{description}, draw=none, from=5, to=0]
    \end{tikzcd}
    \quad=\quad
    \begin{tikzcd}
      Fx & Gx \\
      Fy & Gz
      \arrow[""{name=0, anchor=center, inner sep=0}, "{M\id_x}", "\shortmid"{marking}, from=1-1, to=1-2]
      \arrow["Ff"', from=1-1, to=2-1]
      \arrow["Gg", from=1-2, to=2-2]
      \arrow[""{name=1, anchor=center, inner sep=0}, "Mn"', "\shortmid"{marking}, from=2-1, to=2-2]
      \arrow["M\alpha"{description}, draw=none, from=0, to=1]
    \end{tikzcd}
    \quad=\quad
    \begin{tikzcd}
      Fx & Gx & Gx \\
      Fx & Gx & Gx \\
      Fy & Gy & Gz \\
      Fy && Gz
      \arrow[""{name=0, anchor=center, inner sep=0}, "{M\id_x}", "\shortmid"{marking}, from=2-1, to=2-2]
      \arrow[""{name=1, anchor=center, inner sep=0}, "{G\id_x}", "\shortmid"{marking}, from=2-2, to=2-3]
      \arrow[""{name=2, anchor=center, inner sep=0}, "{M\id_y}"', "\shortmid"{marking}, from=3-1, to=3-2]
      \arrow["Ff"', from=2-1, to=3-1]
      \arrow["Gf"{description}, from=2-2, to=3-2]
      \arrow["Gg", from=2-3, to=3-3]
      \arrow[""{name=3, anchor=center, inner sep=0}, "Gn"', "\shortmid"{marking}, from=3-2, to=3-3]
      \arrow[""{name=4, anchor=center, inner sep=0}, "Mn"', "\shortmid"{marking}, from=4-1, to=4-3]
      \arrow[Rightarrow, no head, from=3-1, to=4-1]
      \arrow[Rightarrow, no head, from=3-3, to=4-3]
      \arrow[Rightarrow, no head, from=1-1, to=2-1]
      \arrow[Rightarrow, no head, from=1-2, to=2-2]
      \arrow[Rightarrow, no head, from=1-3, to=2-3]
      \arrow[""{name=5, anchor=center, inner sep=0}, "{M\id_x}", "\shortmid"{marking}, from=1-1, to=1-2]
      \arrow[""{name=6, anchor=center, inner sep=0}, "{\id_{Gx}}", "\shortmid"{marking}, from=1-2, to=1-3]
      \arrow["G\alpha"{description}, draw=none, from=1, to=3]
      \arrow["{M\id_f}"{description}, draw=none, from=0, to=2]
      \arrow["{M^r_{y,n}}"{description}, draw=none, from=3-2, to=4]
      \arrow["1"{description}, draw=none, from=5, to=0]
      \arrow["{G_x}"{description}, draw=none, from=6, to=1]
    \end{tikzcd}.
  \end{equation*}
\end{lemma}
\begin{proof}
  Apply the naturality axiom of a modulation to the equation
  $\id_f \odot \alpha = \alpha = \alpha \odot \id_g$, then precompose with the
  unitors $F_x$ and $G_x$ and use the unitality of the left and right actions.
\end{proof}

Since modules are themselves lax functors, a module that preserves finite
products should be analogous to a cartesian lax functor
(\cref{def:cartesian-lax-functor}). A cartesian module
$M: F \proTo G: \dbl{D} \to \dbl{E}$ cannot be directly defined by requiring the
lax functor $M: \dbl{D} \times \dbl{I} \to \dbl{E}$ to be cartesian, since the
double category $\dbl{I} \coloneqq \{0 \xproto{i} 1\}$ does not have all finite
products. However, $\dbl{I}$ does have a few products: for every positive number
$n$, we have $0^n = 0$ and $1^n = 1$ in the discrete category $\dbl{I}_0$ and
similarly $\id_0^n = \id_0$, $\id_1^n = \id_1$, and $i^n = i$ in the discrete
category $\dbl{I}_1$, and these products are clearly preserved by the source and
target functors. Thus, we can define a module $M$ to be cartesian if the lax
functor $M: \dbl{D} \times \dbl{I} \to \dbl{E}$ preserves finite products
whenever they exist in $\dbl{D} \times \dbl{I}$. This leads to the following
definition.

\begin{definition}[Cartesian module] \label{def:cartesian-module}
  Let $F,G: \dbl{D} \to \dbl{E}$ be cartesian double functors between
  precartesian double categories $\dbl{D}$ and $\dbl{E}$. A module
  $M: F \proTo G$ is \define{cartesian} if for every pair of proarrows
  $m: x \proto y$ and $m': x' \proto y'$ in $\dbl{D}$, the canonical comparison
  cell
  \begin{equation*}
    \begin{tikzcd}[column sep=large]
      {F(x \times x')} & {G(y \times y')} \\
      {Fx \times Fx'} & {Gy \times Gy'}
      \arrow["{\Phi_{x,x'}}"', from=1-1, to=2-1]
      \arrow["{\Psi_{y,y'}}", from=1-2, to=2-2]
      \arrow[""{name=0, anchor=center, inner sep=0}, "{M(m \times m')}", "\shortmid"{marking}, from=1-1, to=1-2]
      \arrow[""{name=1, anchor=center, inner sep=0}, "{Mm \times Mm'}"', "\shortmid"{marking}, from=2-1, to=2-2]
      \arrow["{\mu_{m,m'}}"{description}, draw=none, from=0, to=1]
    \end{tikzcd}
    \quad\coloneqq\quad
    \begin{tikzcd}[column sep=large]
      {F(x \times x')} & {G(y \times y')} \\
      {Fx \times Fx'} & {Gy \times Gy'}
      \arrow["{\langle F\pi_{x,x'}, F\pi_{x,x'}' \rangle}"', from=1-1, to=2-1]
      \arrow["{\langle G\pi_{y,y'}, G\pi_{y,y'}' \rangle}", from=1-2, to=2-2]
      \arrow[""{name=0, anchor=center, inner sep=0}, "{M(m \times m')}", "\shortmid"{marking}, from=1-1, to=1-2]
      \arrow[""{name=1, anchor=center, inner sep=0}, "{Mm \times Mm'}"', "\shortmid"{marking}, from=2-1, to=2-2]
      \arrow["{\langle M \pi_{m,m'}, M \pi_{m,m'}' \rangle}"{description}, draw=none, from=0, to=1]
    \end{tikzcd}
  \end{equation*}
  is an isomorphism in $\dbl{E}_1$, and so is the unique cell
  $M(I_1) \xto{!} I_1$.
\end{definition}

The left and right actions of a cartesian module preserve finite products in a
sense analogous to \cref{lem:cartesian-laxators-unitors}, by an analogous
argument. We omit a detailed statement.

\begin{example}[Monoidal profunctors]
  A cartesian module between models of the theory of pseudomonoids
  (\cref{th:pseudomonoid}) is a \define{monoidal profunctor} between the
  corresponding monoidal categories. Monoidal profunctors between monoidal
  categories have been previously considered by Koudenburg \cite[Example
  3.23]{koudenburg2015} and by Spivak et al \cite[\S{3.2}]{spivak2017}.

  Unpacking the definition, a monoidal profunctor between monoidal categories
  $(\cat{C}, \otimes_{\cat{C}}, I_{\cat{C}})$ and
  $(\cat{D}, \otimes_{\cat{D}}, I_{\cat{D}})$ consists of a profunctor
  $M: \cat{C} \proto \cat{D}$ and maps
  \begin{equation*}
    \begin{tikzcd}
      {\cat{C}_0 \times \cat{C}_0} & {\cat{D}_0 \times \cat{D}_0} \\
      {\cat{C}_0} & {\cat{D}_0}
      \arrow[""{name=0, anchor=center, inner sep=0}, "{M \times M}", "\shortmid"{marking}, from=1-1, to=1-2]
      \arrow["{\otimes_{\cat{C}}}"', from=1-1, to=2-1]
      \arrow["{\otimes_{\cat{D}}}", from=1-2, to=2-2]
      \arrow[""{name=1, anchor=center, inner sep=0}, "M"', "\shortmid"{marking}, from=2-1, to=2-2]
      \arrow["{\otimes_M}"{description}, draw=none, from=0, to=1]
    \end{tikzcd}
    \qquad\text{and}\qquad
    \begin{tikzcd}
      1 & 1 \\
      {\cat{C}_0} & {\cat{D}_0}
      \arrow[""{name=0, anchor=center, inner sep=0}, "{\id_1}", "\shortmid"{marking}, from=1-1, to=1-2]
      \arrow["{I_{\cat{C}}}"', from=1-1, to=2-1]
      \arrow["{I_{\cat{D}}}", from=1-2, to=2-2]
      \arrow[""{name=1, anchor=center, inner sep=0}, "M"', "\shortmid"{marking}, from=2-1, to=2-2]
      \arrow["{I_M}"{description}, draw=none, from=0, to=1]
    \end{tikzcd}
  \end{equation*}
  that, by the naturality of the left and right actions, satisfy
  \begin{align*}
    (f_1 \otimes_{\cat{C}} f_2) \cdot (m_1 \otimes_M m_2)
      &= (f_1 \cdot m_1) \otimes_M (f_2 \cdot m_2) \\
    (m_1 \otimes_M m_2) \cdot (g_1 \otimes_{\cat{D}} g_2)
      &= (m_1 \cdot g_1) \otimes_M (m_2 \cdot g_2)
  \end{align*}
  for all heteromorphisms $m_i \in M(x_i, y_y)$ and morphisms $f_i: w_i \to x_i$
  in $\cat{C}$ and $g_i: y_i \to z_i$ in $\cat{D}$. In addition, by
  \cref{lem:module-cells}, the associators and unitors of the monoidal
  categories must be preserved in the sense that
  \begin{equation*}
    \alpha_{x_1,x_2,x_3}^{\cat{C}} \cdot (m_1 \otimes_M (m_2 \otimes_M m_3))
      = ((m_1 \otimes_M m_2) \otimes_M m_3) \cdot \alpha_{y_1,y_2,y_3}^{\cat{D}}
  \end{equation*}
  for all heteromorphisms $m_i \in M(x_i, y_i)$, and
  \begin{equation*}
    \lambda_x^{\cat{C}} \cdot m = (I_M \otimes_M m) \cdot \lambda_y^{\cat{D}}
    \qquad\text{and}\qquad
    \rho_x^{\cat{C}} \cdot m = (m \otimes_M I_M) \cdot \rho_y^{\cat{D}}
  \end{equation*}
  for all heteromorphisms $m \in M(x,y)$.
\end{example}

\begin{example}[Multiprofunctors]
  A cartesian module between models of the theory of promonoids
  (\cref{th:promonoid}) is a \define{multiprofunctor} between the corresponding
  multicategories. Thus, a multiprofunctor $P: \cat{C} \proto \cat{D}$ between
  multicategories $\cat{C}$ and $\cat{D}$ consists of, for every arity
  $n \in \N$, sets of ``$n$-ary heteromorphisms'' between $\cat{C}$ and
  $\cat{D}$
  \begin{equation*}
    P(x_1, \dots, x_n; y) \in \Set, \qquad
    x_1, \dots, x_n \in \cat{C},\ y \in \cat{D},
  \end{equation*}
  equipped with a left action by the multicategory $\cat{C}$
  \begin{equation*}
    \prod_{i=1}^k \cat{C}(w_{i,1}, \dots, w_{i,n_i}; x_i) \times
      P(x_1, \dots, x_k; y) \to P(w_{1,1}, \dots, w_{k,n_k}; y)
  \end{equation*}
  and a right action by the multicategory $\cat{D}$
  \begin{equation*}
    \prod_{i=1}^k P(x_{i,1}, \dots, x_{i,n_i}; y_i) \times
      \cat{D}(y_1, \dots, y_k; z) \to P(x_{1,1}, \dots, x_{k,n_k}; z),
  \end{equation*}
  subject to laws of associativity, unitality, and compatibility. A
  multiprofunctor is precisely what Leinster calls a module between
  multicategories \cite[Definition 2.3.6]{leinster2004}.
\end{example}

The existing definition of a modulation \cite[Definition 3.3]{pare2011} requires
a slight generalization to include \emph{lax} natural transformations.

\begin{definition}[Modulation] \label{def:modulation} Given lax double functors
  $F,G,H,K: \dbl{D} \to \dbl{E}$, lax natural transformations $\alpha: F \To H$
  and $\beta: G \To K$, and modules $M: F \proTo G$ and $N: H \proTo K$, a
  \define{modulation}
  \begin{equation*}
    \begin{tikzcd}
      F & G \\
      H & K
      \arrow[""{name=0, anchor=center, inner sep=0}, "M", "\shortmid"{marking}, from=1-1, to=1-2]
      \arrow[""{name=1, anchor=center, inner sep=0}, "N"', "\shortmid"{marking}, from=2-1, to=2-2]
      \arrow["\alpha"', from=1-1, to=2-1]
      \arrow["\beta", from=1-2, to=2-2]
      \arrow["\mu"{description}, draw=none, from=0, to=1]
    \end{tikzcd}
  \end{equation*}
  consists of, for every proarrow $m: x \proto y$ of $\dbl{D}$, a cell in
  $\dbl{E}$
  \begin{equation*}
    \begin{tikzcd}
      Fx & Gy \\
      Hx & Ky
      \arrow[""{name=0, anchor=center, inner sep=0}, "Mm", "\shortmid"{marking}, from=1-1, to=1-2]
      \arrow[""{name=1, anchor=center, inner sep=0}, "Nm"', "\shortmid"{marking}, from=2-1, to=2-2]
      \arrow["{\alpha_x}"', from=1-1, to=2-1]
      \arrow["{\beta_y}", from=1-2, to=2-2]
      \arrow["{\mu_m}"{description}, draw=none, from=0, to=1]
    \end{tikzcd},
  \end{equation*}
  the \define{component} of $\mu$ at $m$, satisfying the following two axioms.
  \begin{itemize}
    \item Equivariance: for every consecutive pair of proarrows
      $x \xproto{m} y \xproto{n} z$ in $\dbl{D}$,
      \begin{equation*}
        \begin{tikzcd}
          Fx & Gy & Gz \\
          Hx & Ky & Kz \\
          Hx && Kz
          \arrow[""{name=0, anchor=center, inner sep=0}, "Mm", "\shortmid"{marking}, from=1-1, to=1-2]
          \arrow[""{name=1, anchor=center, inner sep=0}, "Gn", "\shortmid"{marking}, from=1-2, to=1-3]
          \arrow["{\beta_y}", from=1-2, to=2-2]
          \arrow["{\beta_z}", from=1-3, to=2-3]
          \arrow[""{name=2, anchor=center, inner sep=0}, "Kn"', "\shortmid"{marking}, from=2-2, to=2-3]
          \arrow[""{name=3, anchor=center, inner sep=0}, "Nm"', "\shortmid"{marking}, from=2-1, to=2-2]
          \arrow["{\alpha_x}"', from=1-1, to=2-1]
          \arrow[""{name=4, anchor=center, inner sep=0}, "{N(m \odot n)}"', "\shortmid"{marking}, from=3-1, to=3-3]
          \arrow[Rightarrow, no head, from=2-3, to=3-3]
          \arrow[Rightarrow, no head, from=2-1, to=3-1]
          \arrow["{\mu_m}"{description}, draw=none, from=0, to=3]
          \arrow["{N^r_{m,n}}"{description}, draw=none, from=2-2, to=4]
          \arrow["{\beta_n}"{description}, draw=none, from=1, to=2]
        \end{tikzcd}
        \quad=\quad
        \begin{tikzcd}
          Fx & Gy & Gz \\
          Fx && Gz \\
          Hx && Kz
          \arrow["Mm", "\shortmid"{marking}, from=1-1, to=1-2]
          \arrow["Gn", "\shortmid"{marking}, from=1-2, to=1-3]
          \arrow[""{name=0, anchor=center, inner sep=0}, "{M(m \odot n)}"', "\shortmid"{marking}, from=2-1, to=2-3]
          \arrow[Rightarrow, no head, from=1-1, to=2-1]
          \arrow[Rightarrow, no head, from=1-3, to=2-3]
          \arrow[""{name=1, anchor=center, inner sep=0}, "{N(m \odot n)}"', "\shortmid"{marking}, from=3-1, to=3-3]
          \arrow["{\alpha_x}"', from=2-1, to=3-1]
          \arrow["{\beta_z}", from=2-3, to=3-3]
          \arrow["{M^r_{m,n}}"{description}, draw=none, from=1-2, to=0]
          \arrow["{\mu_{m \odot n}}"{description}, draw=none, from=0, to=1]
        \end{tikzcd}
      \end{equation*}
      and
      \begin{equation*}
        \begin{tikzcd}
          Fx & Fy & Gz \\
          Hx & Hy & Kz \\
          Hx && Kz
          \arrow[""{name=0, anchor=center, inner sep=0}, "Fm", "\shortmid"{marking}, from=1-1, to=1-2]
          \arrow[""{name=1, anchor=center, inner sep=0}, "Mn", "\shortmid"{marking}, from=1-2, to=1-3]
          \arrow["{\alpha_y}"', from=1-2, to=2-2]
          \arrow["{\beta_z}", from=1-3, to=2-3]
          \arrow[""{name=2, anchor=center, inner sep=0}, "Nn"', "\shortmid"{marking}, from=2-2, to=2-3]
          \arrow[""{name=3, anchor=center, inner sep=0}, "Hm"', "\shortmid"{marking}, from=2-1, to=2-2]
          \arrow["{\alpha_x}"', from=1-1, to=2-1]
          \arrow[""{name=4, anchor=center, inner sep=0}, "{N(m \odot n)}"', "\shortmid"{marking}, from=3-1, to=3-3]
          \arrow[Rightarrow, no head, from=2-3, to=3-3]
          \arrow[Rightarrow, no head, from=2-1, to=3-1]
          \arrow["{\alpha_m}"{description}, draw=none, from=0, to=3]
          \arrow["{N^\ell_{m,n}}"{description}, draw=none, from=2-2, to=4]
          \arrow["{\mu_n}"{description}, draw=none, from=1, to=2]
        \end{tikzcd}
        \quad=\quad
        \begin{tikzcd}
          Fx & Fy & Gz \\
          Fx && Gz \\
          Hx && Kz
          \arrow["Fm", "\shortmid"{marking}, from=1-1, to=1-2]
          \arrow["Mn", "\shortmid"{marking}, from=1-2, to=1-3]
          \arrow[""{name=0, anchor=center, inner sep=0}, "{M(m \odot n)}"', "\shortmid"{marking}, from=2-1, to=2-3]
          \arrow[Rightarrow, no head, from=1-1, to=2-1]
          \arrow[Rightarrow, no head, from=1-3, to=2-3]
          \arrow[""{name=1, anchor=center, inner sep=0}, "{N(m \odot n)}"', "\shortmid"{marking}, from=3-1, to=3-3]
          \arrow["{\alpha_x}"', from=2-1, to=3-1]
          \arrow["{\beta_z}", from=2-3, to=3-3]
          \arrow["{M^\ell_{m,n}}"{description}, draw=none, from=1-2, to=0]
          \arrow["{\mu_{m \odot n}}"{description}, draw=none, from=0, to=1]
        \end{tikzcd}.
      \end{equation*}
    \item Naturality: for every cell $\stdInlineCell{\gamma}$ in $\dbl{D}$,
      \begin{equation*}
        \begin{tikzcd}
          Fx & Fx & Gy \\
          Hx & Fw & Gz \\
          Hw & Hw & Kz \\
          Hw && Hz
          \arrow["Ff"', from=1-2, to=2-2]
          \arrow["Gg", from=1-3, to=2-3]
          \arrow[""{name=0, anchor=center, inner sep=0}, "Mm", "\shortmid"{marking}, from=1-2, to=1-3]
          \arrow[""{name=1, anchor=center, inner sep=0}, "Mn"', "\shortmid"{marking}, from=2-2, to=2-3]
          \arrow[""{name=2, anchor=center, inner sep=0}, "{\mathrm{id}_{Fx}}", "\shortmid"{marking}, from=1-1, to=1-2]
          \arrow["{\alpha_x}"', from=1-1, to=2-1]
          \arrow["Hf"', from=2-1, to=3-1]
          \arrow["{\alpha_w}"', from=2-2, to=3-2]
          \arrow[""{name=3, anchor=center, inner sep=0}, "{H \mathrm{id}_w}"', "\shortmid"{marking}, from=3-1, to=3-2]
          \arrow["{\beta_z}", from=2-3, to=3-3]
          \arrow[""{name=4, anchor=center, inner sep=0}, "Nn"', "\shortmid"{marking}, from=3-2, to=3-3]
          \arrow[""{name=5, anchor=center, inner sep=0}, "Nn"', "\shortmid"{marking}, from=4-1, to=4-3]
          \arrow[Rightarrow, no head, from=3-1, to=4-1]
          \arrow[Rightarrow, no head, from=3-3, to=4-3]
          \arrow["M\gamma"{description}, draw=none, from=0, to=1]
          \arrow["{\mu_n}"{description}, draw=none, from=1, to=4]
          \arrow["{N^\ell_{w,n}}"{description}, draw=none, from=3-2, to=5]
          \arrow["{\alpha_f}"{description}, draw=none, from=2, to=3]
        \end{tikzcd}
        \quad=\quad
        \begin{tikzcd}
          Fx & Gy & Gy \\
          Hx & Ky & Gz \\
          Hw & Kz & Kz \\
          Hw && Kz
          \arrow[""{name=0, anchor=center, inner sep=0}, "Mm", "\shortmid"{marking}, from=1-1, to=1-2]
          \arrow["{\alpha_x}"', from=1-1, to=2-1]
          \arrow["{\beta_y}", from=1-2, to=2-2]
          \arrow[""{name=1, anchor=center, inner sep=0}, "Nm"', "\shortmid"{marking}, from=2-1, to=2-2]
          \arrow["Hf"', from=2-1, to=3-1]
          \arrow["Kg", from=2-2, to=3-2]
          \arrow[""{name=2, anchor=center, inner sep=0}, "Nn"', "\shortmid"{marking}, from=3-1, to=3-2]
          \arrow[""{name=3, anchor=center, inner sep=0}, "{\mathrm{id}_{Gy}}", "\shortmid"{marking}, from=1-2, to=1-3]
          \arrow["Gg", from=1-3, to=2-3]
          \arrow["{\beta_z}", from=2-3, to=3-3]
          \arrow[""{name=4, anchor=center, inner sep=0}, "{K \mathrm{id}_z}"', "\shortmid"{marking}, from=3-2, to=3-3]
          \arrow[Rightarrow, no head, from=3-1, to=4-1]
          \arrow[Rightarrow, no head, from=3-3, to=4-3]
          \arrow[""{name=5, anchor=center, inner sep=0}, "Nn"', "\shortmid"{marking}, from=4-1, to=4-3]
          \arrow["{\mu_m}"{description}, draw=none, from=0, to=1]
          \arrow["N\gamma"{description}, draw=none, from=1, to=2]
          \arrow["{N^r_{n,z}}"{description}, draw=none, from=3-2, to=5]
          \arrow["{\beta_g}"{description}, draw=none, from=3, to=4]
        \end{tikzcd}.
      \end{equation*}
  \end{itemize}
\end{definition}

This definition of a modulation reduces to our previous one when the modules
involved are identities. Given any lax double functor $F$, there is a canonical
\define{identity module} $\id_F: F \proTo F$ where $\id_F(m) = F(m)$,
$\id_F(\alpha) = F(\alpha)$, and the left and right actions of $\id_F$ are both
defined by the laxators of $F$ \cite[Definition 5.1.1]{pare2013}.

\begin{proposition} \label{lem:modulation-definitions}
  Let $F,G: \dbl{D} \to \dbl{E}$ be lax double functors and let
  $\alpha, \beta: F \To G$ be lax natural transformations. A modulation
  $\mu: F \Tto G$ in the sense of \cref{def:modulation-special} is equivalent to
  a modulation
  \begin{equation*}
    \begin{tikzcd}
      F & F \\
      G & G
      \arrow[""{name=0, anchor=center, inner sep=0}, "{\mathrm{id}_F}", "\shortmid"{marking}, from=1-1, to=1-2]
      \arrow["\alpha"', from=1-1, to=2-1]
      \arrow["\beta", from=1-2, to=2-2]
      \arrow[""{name=1, anchor=center, inner sep=0}, "{\mathrm{id}_G}"', "\shortmid"{marking}, from=2-1, to=2-2]
      \arrow["\mu"{description}, draw=none, from=0, to=1]
    \end{tikzcd}
  \end{equation*}
  in the sense of \cref{def:modulation}.
\end{proposition}
\begin{proof}
  Note that the two kinds of modulations in question have components
  parameterized by objects and proarrows of $\dbl{D}$, respectively. Given a
  modulation $\mu = (\mu_x)_{x \in \dbl{D}_0}$ in the first sense, define a
  modulation in the second sense by \cref{eq:modulation-special-equivariance}.
  Conversely, given a modulation $\mu = (\mu_m)_{m \in \dbl{D}_1}$ in the second
  sense, define a modulation in the first sense by
  \cref{eq:modulation-id-cells}. It is straightforward to show that these
  operations are mutually inverse and that the two axioms called ``naturality''
  and ``equivariance'' in each case imply each other.

  Alternatively, this equivalence is recovered as a special case of
  \cref{thm:lax-functor-vdbl-category} below.
\end{proof}

\section{Double categories of models}
\label{sec:model-double-categories}

A double-categorical framework for doctrines might be expected to produce not
just a 2-category but a double category of models, having modules as proarrows
and modulations as cells. This prospect faces the obstacle that composing
modules between lax functors is difficult and subtle, with composites known to
exist only under sufficient conditions \cite{pare2013}. In this paper, except in
a certain important special case, we will sidestep these issues and settle for a
\emph{virtual} double category of models. Virtual double categories are related
to double categories in the same way that multicategories are to monoidal
categories. The phrase ``virtual double category'' was introduced by Cruttwell
and Shulman \cite[\S 2]{cruttwell2010}. Virtual double categories have been
studied earlier by Leinster \cite{leinster2004} under the name
``fc-multicategories'' and earlier still by Burroni \cite{burroni1971}.

To define a virtual double category of lax functors, we need the concept of a
multimodulation \cite{pare2011}, which like the concept of a modulation
(\cref{def:modulation}) we must generalize slightly to account for \emph{lax}
natural transformations.

\begin{definition}[Multimodulation] \label{def:multimodulation}
  Let $F_0, F_1, \dots, F_k, G, H: \dbl{D} \to \dbl{E}$ be lax double functors
  and let $\alpha: F_0 \To G$ and $\beta: F_k \To H$ be lax natural
  transformations. A \define{multimodulation} $\mu$ from a composable sequence
  of modules $F_0 \xproTo{M_1} F_1 \xproTo{M_2} \cdots \xproTo{M_k} F_k$ to a
  module $N: G \proTo H$, having source $\alpha$ and target $\beta$, is depicted
  as
  \begin{equation*}
    \begin{tikzcd}
      {F_0} & {F_1} & \cdots & {F_k} \\
      G &&& H
      \arrow["{M_1}", "\shortmid"{marking}, from=1-1, to=1-2]
      \arrow["\alpha"', from=1-1, to=2-1]
      \arrow["\beta", from=1-4, to=2-4]
      \arrow[""{name=0, anchor=center, inner sep=0}, "N"', "\shortmid"{marking}, from=2-1, to=2-4]
      \arrow[""{name=1, anchor=center, inner sep=0}, "{M_2}", "\shortmid"{marking}, from=1-2, to=1-3]
      \arrow["{M_k}", "\shortmid"{marking}, from=1-3, to=1-4]
      \arrow["\mu"{description}, draw=none, from=1, to=0]
    \end{tikzcd}
  \end{equation*}
  and consists of, for every sequence of proarrows
  $x_0 \xproto{m_1} x_1 \xproto{m_2} \cdots \xproto{m_k} x_k$ in $\dbl{D}$, a
  cell in $\dbl{E}$
  \begin{equation*}
    \begin{tikzcd}
      {F_0 x_0} & {F_1 x_1} & \cdots & {F_k x_k} \\
      {G x_0} &&& {H x_k}
      \arrow["{\alpha_{x_0}}"', from=1-1, to=2-1]
      \arrow["{M_1 m_1}", "\shortmid"{marking}, from=1-1, to=1-2]
      \arrow[""{name=0, anchor=center, inner sep=0}, "{M_2 m_2}", "\shortmid"{marking}, from=1-2, to=1-3]
      \arrow["{\beta_{x_k}}", from=1-4, to=2-4]
      \arrow[""{name=1, anchor=center, inner sep=0}, "{N(m_1 \odot \cdots \odot m_k)}"', "\shortmid"{marking}, from=2-1, to=2-4]
      \arrow["{M_k m_k}", "\shortmid"{marking}, from=1-3, to=1-4]
      \arrow["{\mu_{m_1,\dots,m_k}}"{description}, draw=none, from=0, to=1]
    \end{tikzcd}.
  \end{equation*}
  Note that an empty sequence of modules ($k=0$) is allowed: a nullary
  multimodulation
  \begin{equation*}
    \begin{tikzcd}
      & F \\
      G && H
      \arrow[""{name=0, anchor=center, inner sep=0}, "N"', "\shortmid"{marking}, from=2-1, to=2-3]
      \arrow["\alpha"', from=1-2, to=2-1]
      \arrow["\beta", from=1-2, to=2-3]
      \arrow["\mu"{description}, draw=none, from=1-2, to=0]
    \end{tikzcd}
    \qquad\leftrightsquigarrow\qquad
    \begin{tikzcd}
      Fx & Fx \\
      Gx & Hx
      \arrow[""{name=0, anchor=center, inner sep=0}, "{\id_{Fx}}", "\shortmid"{marking}, from=1-1, to=1-2]
      \arrow["{\alpha_x}"', from=1-1, to=2-1]
      \arrow[""{name=1, anchor=center, inner sep=0}, "{N \id_x}"', "\shortmid"{marking}, from=2-1, to=2-2]
      \arrow["{\beta_x}", from=1-2, to=2-2]
      \arrow["{\mu_x}"{description}, draw=none, from=0, to=1]
    \end{tikzcd},
    \quad x \in \dbl{D},
  \end{equation*}
  as on the left comprises a family of cells in $\dbl{E}$ parameterized by
  objects of $\dbl{D}$ as on the right.

  The following axioms must be satisfied.
  \begin{itemize}
    \item Left equivariance: if $k \geq 1$, then
      \begin{equation*}
        \begin{tikzcd}[row sep=scriptsize]
          {F_0 w} & {F_0 x_0} & {F_1 x_1} & \cdots & {F_k x_k} \\
          {G w} & {G x_0} &&& {H x_k} \\
          Gw &&&& {H_{x_k}}
          \arrow[""{name=0, anchor=center, inner sep=0}, "{F_0 m_0}", "\shortmid"{marking}, from=1-1, to=1-2]
          \arrow[""{name=1, anchor=center, inner sep=0}, "{ G m_0}"', "\shortmid"{marking}, from=2-1, to=2-2]
          \arrow["{\alpha_w}"', from=1-1, to=2-1]
          \arrow["{\alpha_{x_0}}", from=1-2, to=2-2]
          \arrow["{M_1 m_1}", "\shortmid"{marking}, from=1-2, to=1-3]
          \arrow[""{name=2, anchor=center, inner sep=0}, "{M_2 m_2}", "\shortmid"{marking}, from=1-3, to=1-4]
          \arrow["{M_k m_k}", "\shortmid"{marking}, from=1-4, to=1-5]
          \arrow["{\beta_{x_k}}", from=1-5, to=2-5]
          \arrow[""{name=3, anchor=center, inner sep=0}, "{N(m_1 \odot \cdots \odot m_k)}"', "\shortmid"{marking}, from=2-2, to=2-5]
          \arrow[Rightarrow, no head, from=2-1, to=3-1]
          \arrow[""{name=4, anchor=center, inner sep=0}, "{N(m_0 \odot m_1 \odot \cdots \odot m_k)}"', "\shortmid"{marking}, from=3-1, to=3-5]
          \arrow[Rightarrow, no head, from=2-5, to=3-5]
          \arrow["{\alpha_{m_0}}"{description}, draw=none, from=0, to=1]
          \arrow["{\mu_{m_1,\dots,m_k}}"{description}, draw=none, from=2, to=3]
          \arrow["{N^\ell_{m_0, m_1 \odot \cdots \odot m_k}}"{description, pos=0.8}, draw=none, from=1-3, to=4]
        \end{tikzcd}
      \end{equation*}
      is equal to
      \begin{equation*}
        \begin{tikzcd}[row sep=scriptsize]
          {F_0 w} & {F_0 x_0} & {F_1 x_1} & \cdots & {F_k x_k} \\
          {F_0 w} && {F_1 x_1} & \cdots & {F_k x_k} \\
          Gw &&&& {H_{x_k}}
          \arrow["{F_0 m_0}", "\shortmid"{marking}, from=1-1, to=1-2]
          \arrow["{M_1 m_1}", "\shortmid"{marking}, from=1-2, to=1-3]
          \arrow[""{name=0, anchor=center, inner sep=0}, "{M_2 m_2}", "\shortmid"{marking}, from=1-3, to=1-4]
          \arrow[""{name=1, anchor=center, inner sep=0}, "{M_k m_k}", "\shortmid"{marking}, from=1-4, to=1-5]
          \arrow[""{name=2, anchor=center, inner sep=0}, "{N(m_0 \odot m_1 \odot \cdots \odot m_k)}"', "\shortmid"{marking}, from=3-1, to=3-5]
          \arrow[Rightarrow, no head, from=1-1, to=2-1]
          \arrow["{\alpha_w}"', from=2-1, to=3-1]
          \arrow[Rightarrow, no head, from=1-3, to=2-3]
          \arrow[""{name=3, anchor=center, inner sep=0}, "{M_1(m_0 \odot m_1)}"', "\shortmid"{marking}, from=2-1, to=2-3]
          \arrow[Rightarrow, no head, from=1-5, to=2-5]
          \arrow[""{name=4, anchor=center, inner sep=0}, "{M_2 m_2}"', "\shortmid"{marking}, from=2-3, to=2-4]
          \arrow[""{name=5, anchor=center, inner sep=0}, "{M_k m_k}"', "\shortmid"{marking}, from=2-4, to=2-5]
          \arrow["{\beta_{x_k}}", from=2-5, to=3-5]
          \arrow["{(M_1)^\ell_{m_0, m_1}}"{description}, draw=none, from=1-2, to=3]
          \arrow["1"{description}, draw=none, from=0, to=4]
          \arrow["1"{description}, draw=none, from=1, to=5]
          \arrow["{\mu_{m_0 \odot m_1, m_2, \dots, m_k}}"{description}, draw=none, from=2-3, to=2]
        \end{tikzcd}.
      \end{equation*}
    \item Right equivariance, dual to the previous axiom.
    \item Inner equivariance: if $k \geq 2$, then for every $1 \leq i < k$,
      \begin{equation*}
        \begin{tikzcd}[sep=scriptsize]
          {F_0 x_0} & \cdots & {F_{i-1} x_{i-1}} & {F_i x_i} & {F_i x_i'} && {F_{i+1} x_{i+1}} & \cdots & {F_k x_k} \\
          {F_0 x_0} & \cdots & {F_{i-1} x_{i-1}} && {F_i x_i'} && {F_{i+1} x_{i+1}} & \cdots & {F_k x_k} \\
          {G x_0} &&&&&&&& {H x_k}
          \arrow["{M_i m_i}", "\shortmid"{marking}, from=1-3, to=1-4]
          \arrow["{F_i m_i'}", "\shortmid"{marking}, from=1-4, to=1-5]
          \arrow[""{name=0, anchor=center, inner sep=0}, "{M_{i+1} m_{i+1}}", "\shortmid"{marking}, from=1-5, to=1-7]
          \arrow[""{name=1, anchor=center, inner sep=0}, "{M_i(m_i \odot m_i')}"', "\shortmid"{marking}, from=2-3, to=2-5]
          \arrow[Rightarrow, no head, from=1-3, to=2-3]
          \arrow[Rightarrow, no head, from=1-1, to=2-1]
          \arrow[Rightarrow, no head, from=1-5, to=2-5]
          \arrow[Rightarrow, no head, from=1-7, to=2-7]
          \arrow[Rightarrow, no head, from=1-9, to=2-9]
          \arrow[""{name=2, anchor=center, inner sep=0}, "{M_{i+1} m_{i+1}}"', "\shortmid"{marking}, from=2-5, to=2-7]
          \arrow[""{name=3, anchor=center, inner sep=0}, "{N(m_1 \odot \cdots\odot m_i \odot m_i' \odot m_{i+1} \odot \cdots \odot m_k)}"', "\shortmid"{marking}, from=3-1, to=3-9]
          \arrow["{\beta_{x_k}}", from=2-9, to=3-9]
          \arrow["{\alpha_{x_0}}"', from=2-1, to=3-1]
          \arrow[""{name=4, anchor=center, inner sep=0}, "{M_1 m_1}", "\shortmid"{marking}, from=1-1, to=1-2]
          \arrow[""{name=5, anchor=center, inner sep=0}, "\shortmid"{marking}, from=1-2, to=1-3]
          \arrow[""{name=6, anchor=center, inner sep=0}, "{M_1 m_1}"', "\shortmid"{marking}, from=2-1, to=2-2]
          \arrow[""{name=7, anchor=center, inner sep=0}, "\shortmid"{marking}, from=2-2, to=2-3]
          \arrow[""{name=8, anchor=center, inner sep=0}, "\shortmid"{marking}, from=1-7, to=1-8]
          \arrow[""{name=9, anchor=center, inner sep=0}, "\shortmid"{marking}, from=2-7, to=2-8]
          \arrow[""{name=10, anchor=center, inner sep=0}, "{M_k m_k}", "\shortmid"{marking}, from=1-8, to=1-9]
          \arrow[""{name=11, anchor=center, inner sep=0}, "{M_k m_k}"', "\shortmid"{marking}, from=2-8, to=2-9]
          \arrow["{(M_i)^r_{m_i,m_i'}}"{description}, draw=none, from=1-4, to=1]
          \arrow["1"{description}, draw=none, from=0, to=2]
          \arrow["1"{description}, draw=none, from=4, to=6]
          \arrow["1"{description}, draw=none, from=5, to=7]
          \arrow["1"{description}, draw=none, from=10, to=11]
          \arrow["1"{description}, draw=none, from=8, to=9]
          \arrow["{\mu_{m_1, \dots, m_{i-1}, m_i \odot m_i', m_{i+1}, \dots, m_k}}"{description, pos=0.8}, draw=none, from=1-5, to=3]
        \end{tikzcd}
      \end{equation*}
      is equal to
      \begin{equation*}
        \begin{tikzcd}[sep=scriptsize]
          {F_0 x_0} & \cdots & {F_{i-1} x_{i-1}} & {F_i x_i} & {F_i x_i'} && {F_{i+1} x_{i+1}} & \cdots & {F_k x_k} \\
          {F_0 x_0} & \cdots & {F_{i-1} x_{i-1}} & {F_i x_i} &&& {F_{i+1} x_{i+1}} & \cdots & {F_k x_k} \\
          {G x_0} &&&&&&&& {H x_k}
          \arrow[""{name=0, anchor=center, inner sep=0}, "{M_i m_i}", "\shortmid"{marking}, from=1-3, to=1-4]
          \arrow[""{name=1, anchor=center, inner sep=0}, "{F_i m_i'}", "\shortmid"{marking}, from=1-4, to=1-5]
          \arrow["{M_{i+1} m_{i+1}}", "\shortmid"{marking}, from=1-5, to=1-7]
          \arrow[Rightarrow, no head, from=1-3, to=2-3]
          \arrow[Rightarrow, no head, from=1-1, to=2-1]
          \arrow[Rightarrow, no head, from=1-7, to=2-7]
          \arrow[Rightarrow, no head, from=1-9, to=2-9]
          \arrow[""{name=2, anchor=center, inner sep=0}, "{N(m_1 \odot \cdots \odot m_i \odot m_i' \odot m_{i+1} \odot \cdots \odot m_k)}"', "\shortmid"{marking}, from=3-1, to=3-9]
          \arrow["{\beta_{x_k}}", from=2-9, to=3-9]
          \arrow["{\alpha_{x_0}}"', from=2-1, to=3-1]
          \arrow[""{name=3, anchor=center, inner sep=0}, "{M_i m_i}"', "\shortmid"{marking}, from=2-3, to=2-4]
          \arrow[Rightarrow, no head, from=1-4, to=2-4]
          \arrow[""{name=4, anchor=center, inner sep=0}, "{M_{i+1}(m_i' \odot m_{i+1})}"', "\shortmid"{marking}, from=2-4, to=2-7]
          \arrow[""{name=5, anchor=center, inner sep=0}, "\shortmid"{marking}, from=1-7, to=1-8]
          \arrow[""{name=6, anchor=center, inner sep=0}, "\shortmid"{marking}, from=2-7, to=2-8]
          \arrow[""{name=7, anchor=center, inner sep=0}, "{M_k m_k}", "\shortmid"{marking}, from=1-8, to=1-9]
          \arrow[""{name=8, anchor=center, inner sep=0}, "{M_k m_k}"', "\shortmid"{marking}, from=2-8, to=2-9]
          \arrow[""{name=9, anchor=center, inner sep=0}, "{M_1 m_1}", "\shortmid"{marking}, from=1-1, to=1-2]
          \arrow[""{name=10, anchor=center, inner sep=0}, "{M_1 m_1}"', "\shortmid"{marking}, from=2-1, to=2-2]
          \arrow[""{name=11, anchor=center, inner sep=0}, "\shortmid"{marking}, from=1-2, to=1-3]
          \arrow[""{name=12, anchor=center, inner sep=0}, "\shortmid"{marking}, from=2-2, to=2-3]
          \arrow["{\mu_{m_1, \dots, m_i, m_i' \odot m_{i+1}, m_{i+2}, \dots, m_k}}"{description, pos=0.8}, draw=none, from=1, to=2]
          \arrow["1"{description}, draw=none, from=0, to=3]
          \arrow["{(M_{i+1})^\ell_{m_i', m_{i+1}}}"{description}, draw=none, from=1-5, to=4]
          \arrow["1"{description}, draw=none, from=7, to=8]
          \arrow["1"{description}, draw=none, from=5, to=6]
          \arrow["1"{description}, draw=none, from=9, to=10]
          \arrow["1"{description}, draw=none, from=11, to=12]
        \end{tikzcd}.
      \end{equation*}
    \item Nullary equivariance: if $k = 0$, then for every proarrow
      $m: x \proto y$ in $\dbl{D}$,
      \begin{equation*}
        \begin{tikzcd}[row sep=scriptsize]
          Fx & Fx & Fy \\
          Gx & Hx & Hy \\
          Gx && Hy
          \arrow[""{name=0, anchor=center, inner sep=0}, "{\id_{Fx}}", "\shortmid"{marking}, from=1-1, to=1-2]
          \arrow["{\alpha_x}"', from=1-1, to=2-1]
          \arrow[""{name=1, anchor=center, inner sep=0}, "{N \id_x}"', "\shortmid"{marking}, from=2-1, to=2-2]
          \arrow["{\beta_x}"{description}, from=1-2, to=2-2]
          \arrow[""{name=2, anchor=center, inner sep=0}, "Fm", "\shortmid"{marking}, from=1-2, to=1-3]
          \arrow["{\beta_y}", from=1-3, to=2-3]
          \arrow[""{name=3, anchor=center, inner sep=0}, "Hm"', "\shortmid"{marking}, from=2-2, to=2-3]
          \arrow[Rightarrow, no head, from=2-3, to=3-3]
          \arrow[Rightarrow, no head, from=2-1, to=3-1]
          \arrow[""{name=4, anchor=center, inner sep=0}, "Nm"', "\shortmid"{marking}, from=3-1, to=3-3]
          \arrow["{\mu_x}"{description}, draw=none, from=0, to=1]
          \arrow["{N_{x,m}^r}"{description}, draw=none, from=2-2, to=4]
          \arrow["{\beta_m}"{description}, draw=none, from=2, to=3]
        \end{tikzcd}
        \quad=\quad
        \begin{tikzcd}[row sep=scriptsize]
          Fx & Fy & Fy \\
          Gx & Gy & Hy \\
          Gx && Hy
          \arrow[""{name=0, anchor=center, inner sep=0}, "{\id_{Fy}}", "\shortmid"{marking}, from=1-2, to=1-3]
          \arrow["{\alpha_y}"{description}, from=1-2, to=2-2]
          \arrow[""{name=1, anchor=center, inner sep=0}, "{N \id_y}"', "\shortmid"{marking}, from=2-2, to=2-3]
          \arrow["{\beta_y}", from=1-3, to=2-3]
          \arrow[""{name=2, anchor=center, inner sep=0}, "Fm", "\shortmid"{marking}, from=1-1, to=1-2]
          \arrow["{\alpha_x}"', from=1-1, to=2-1]
          \arrow[""{name=3, anchor=center, inner sep=0}, "Gm"', "\shortmid"{marking}, from=2-1, to=2-2]
          \arrow[""{name=4, anchor=center, inner sep=0}, "Nm"', "\shortmid"{marking}, from=3-1, to=3-3]
          \arrow[Rightarrow, no head, from=2-1, to=3-1]
          \arrow[Rightarrow, no head, from=2-3, to=3-3]
          \arrow["{\mu_y}"{description}, draw=none, from=0, to=1]
          \arrow["{N_{m,y}^\ell}"{description}, draw=none, from=2-2, to=4]
          \arrow["{\alpha_m}"{description}, draw=none, from=2, to=3]
        \end{tikzcd}.
      \end{equation*}
    \item Naturality: for every sequence of cells
      $\inlineCell{x_{i-1}}{x_i}{y_{i-1}}{y_i}{m_i}{n_i}{f_{i-1}}{f_i}{\gamma_i}$,
      $i=1,\dots,k$, in $\dbl{D}$,
      \begin{equation*}
        \begin{tikzcd}
          {F_0 x_0} & {F_0 x_0} & \cdots & {F_k x_k} \\
          {G x_0} & {F_0 y_0} & \cdots & {F_k y_k} \\
          {G y_0} & {G y_0} && {H y_k} \\
          {G y_0} &&& {H y_k}
          \arrow["{\alpha_{y_0}}"', from=2-2, to=3-2]
          \arrow[""{name=0, anchor=center, inner sep=0}, "{M_1 n_1 }"', "\shortmid"{marking}, from=2-2, to=2-3]
          \arrow["{\beta_{y_k}}", from=2-4, to=3-4]
          \arrow[""{name=1, anchor=center, inner sep=0}, "{N(n_1 \odot \cdots \odot n_k)}"', "\shortmid"{marking}, from=3-2, to=3-4]
          \arrow[""{name=2, anchor=center, inner sep=0}, "{M_k n_k}"', "\shortmid"{marking}, from=2-3, to=2-4]
          \arrow["{F_0 f_0}"', from=1-2, to=2-2]
          \arrow[""{name=3, anchor=center, inner sep=0}, "{M_1 m_1}", "\shortmid"{marking}, from=1-2, to=1-3]
          \arrow["{F_k f_k}", from=1-4, to=2-4]
          \arrow[""{name=4, anchor=center, inner sep=0}, "{M_k m_k}", "\shortmid"{marking}, from=1-3, to=1-4]
          \arrow[""{name=5, anchor=center, inner sep=0}, "{\id_{F_0 x_0}}", "\shortmid"{marking}, from=1-1, to=1-2]
          \arrow[""{name=6, anchor=center, inner sep=0}, "{G \id_{y_0}}"', "\shortmid"{marking}, from=3-1, to=3-2]
          \arrow["{\alpha_{x_0}}"', from=1-1, to=2-1]
          \arrow["{G f_0}"', from=2-1, to=3-1]
          \arrow[Rightarrow, no head, from=3-1, to=4-1]
          \arrow[Rightarrow, no head, from=3-4, to=4-4]
          \arrow[""{name=7, anchor=center, inner sep=0}, "{N(n_1 \odot \cdots \odot n_k)}"', "\shortmid"{marking}, from=4-1, to=4-4]
          \arrow["{\mu_{n_1,\dots,n_k}}"{description}, draw=none, from=2-3, to=1]
          \arrow["{M_1 \gamma_1}"{description}, draw=none, from=3, to=0]
          \arrow["{M_k \gamma_k}"{description}, draw=none, from=4, to=2]
          \arrow["{\alpha_{f_0}}"{description}, draw=none, from=5, to=6]
          \arrow["{N_{y_0, n_1 \odot \cdots \odot n_k}^\ell}"{description}, draw=none, from=3-2, to=7]
        \end{tikzcd}
      \end{equation*}
      is equal to
      \begin{equation*}
        \begin{tikzcd}
          {F_0 x_0} & \cdots & {F_k x_k} & {F_k x_k} \\
          {G x_0} && {H x_k} & {F_k y_k} \\
          {G y_0} && {H y_k} & {H y_k} \\
          {G y_0} &&& {H y_k}
          \arrow["{\alpha_{x_0}}"', from=1-1, to=2-1]
          \arrow["{M_1 m_1 }", "\shortmid"{marking}, from=1-1, to=1-2]
          \arrow["{\beta_{x_k}}", from=1-3, to=2-3]
          \arrow[""{name=0, anchor=center, inner sep=0}, "{N(m_1 \odot \cdots \odot m_k)}"', "\shortmid"{marking}, from=2-1, to=2-3]
          \arrow["{M_k m_k}", "\shortmid"{marking}, from=1-2, to=1-3]
          \arrow["{G f_0}"', from=2-1, to=3-1]
          \arrow["{H f_k}", from=2-3, to=3-3]
          \arrow[""{name=1, anchor=center, inner sep=0}, "{N(n_1 \odot \cdots \odot n_k)}"', "\shortmid"{marking}, from=3-1, to=3-3]
          \arrow["{F_k f_k}", from=1-4, to=2-4]
          \arrow["{\beta_{y_k}}", from=2-4, to=3-4]
          \arrow[""{name=2, anchor=center, inner sep=0}, "{\id_{F_k x_k}}", "\shortmid"{marking}, from=1-3, to=1-4]
          \arrow[""{name=3, anchor=center, inner sep=0}, "{H \id_{y_k}}"', "\shortmid"{marking}, from=3-3, to=3-4]
          \arrow[""{name=4, anchor=center, inner sep=0}, "{N(n_1 \odot \cdots \odot n_k)}"', "\shortmid"{marking}, from=4-1, to=4-4]
          \arrow[Rightarrow, no head, from=3-1, to=4-1]
          \arrow[Rightarrow, no head, from=3-4, to=4-4]
          \arrow["{N(\gamma_1 \odot \cdots \odot \gamma_k)}"{description, pos=0.6}, draw=none, from=0, to=1]
          \arrow["{\mu_{m_1,\dots,m_k}}"{description}, draw=none, from=1-2, to=0]
          \arrow["{\beta_{f_k}}"{description}, draw=none, from=2, to=3]
          \arrow["{N_{n_1 \odot \cdots \odot n_k, y_k}^r}"{description}, draw=none, from=3-3, to=4]
        \end{tikzcd}.
      \end{equation*}
  \end{itemize}
\end{definition}

A unary multimodulation is precisely a modulation in the sense of
\cref{def:modulation}, whereas a nullary multimodulation reduces to a modulation
in the sense of \cref{def:modulation-special} when the target module $N$ is an
identity module. Thus, the new definition subsumes and generalizes both of the
previous ones.

Every double category has an underlying 2-category but this is not necessarily
true of a virtual double category, since it may lack the identity proarrows
needed to express the vertical composition of 2-cells in a 2-category. For the
virtual double category of models of this section to properly generalize the
2-category of models constructed in \cref{sec:lax-transformations}, it must
possess units in the following sense \cite[Definition 5.1]{cruttwell2010}.

\begin{definition}[Unit in virtual double category]
  \label{def:unit-virt-dbl-cat}
  A (strongly representable) \define{unit} or \define{identity} for an object
  $x$ in a virtual double category $\dbl{D}$ is a proarrow $\id_x: x \proto x$
  along with a globular nullary cell $\eta_x: ()_x \to \id_x$ satisfying the
  universal property that any multicell
  \begin{equation*}
    \begin{tikzcd}[row sep=scriptsize]
      {x_0} & \dots & {x_i} & \dots & {x_k} \\
      {x_0} &&&& {x_k}
      \arrow["{m_1}", "\shortmid"{marking}, from=1-1, to=1-2]
      \arrow["{m_i}", "\shortmid"{marking}, from=1-2, to=1-3]
      \arrow["{m_{i+1}}", "\shortmid"{marking}, from=1-3, to=1-4]
      \arrow["{m_k}", "\shortmid"{marking}, from=1-4, to=1-5]
      \arrow["f"', from=1-1, to=2-1]
      \arrow["g", from=1-5, to=2-5]
      \arrow[""{name=0, anchor=center, inner sep=0}, "n"', "\shortmid"{marking}, from=2-1, to=2-5]
      \arrow["\alpha"{description}, draw=none, from=1-3, to=0]
    \end{tikzcd}
  \end{equation*}
  in $\dbl{D}$ such that $x_i = x$ for any choice of $0 \leq i \leq k$
  factorizes uniquely as
  \begin{equation*}
    \begin{tikzcd}[row sep=scriptsize]
      {x_0} & \cdots &&& \cdots & {x_k} \\
      {x_0} & \dots & {x_i} & {x_i} & \dots & {x_k} \\
      {x_0} &&&&& {x_k}
      \arrow[""{name=0, anchor=center, inner sep=0}, "{m_1}"', "\shortmid"{marking}, from=2-1, to=2-2]
      \arrow[""{name=1, anchor=center, inner sep=0}, "{m_k}"', "\shortmid"{marking}, from=2-5, to=2-6]
      \arrow["f"', from=2-1, to=3-1]
      \arrow["g", from=2-6, to=3-6]
      \arrow[""{name=2, anchor=center, inner sep=0}, "n"', "\shortmid"{marking}, from=3-1, to=3-6]
      \arrow[""{name=3, anchor=center, inner sep=0}, "{\id_x}"', "\shortmid"{marking}, from=2-3, to=2-4]
      \arrow[""{name=4, anchor=center, inner sep=0}, "{m_i}"', "\shortmid"{marking}, from=2-2, to=2-3]
      \arrow[""{name=5, anchor=center, inner sep=0}, "{m_{i+1}}"', "\shortmid"{marking}, from=2-4, to=2-5]
      \arrow[Rightarrow, no head, from=1-1, to=2-1]
      \arrow[Rightarrow, no head, from=1-6, to=2-6]
      \arrow[""{name=6, anchor=center, inner sep=0}, "{m_k}", "\shortmid"{marking}, from=1-5, to=1-6]
      \arrow[""{name=7, anchor=center, inner sep=0}, "{m_1}", "\shortmid"{marking}, from=1-1, to=1-2]
      \arrow[""{name=8, anchor=center, inner sep=0}, "x"{description}, draw=none, from=1-2, to=1-5]
      \arrow["{1_{m_k}}"{description}, draw=none, from=6, to=1]
      \arrow["{\exists !}"{description, pos=0.6}, draw=none, from=3, to=2]
      \arrow["{1_{m_1}}"{description}, draw=none, from=7, to=0]
      \arrow[shorten <=4pt, Rightarrow, no head, from=8, to=2-3]
      \arrow[shorten <=4pt, Rightarrow, no head, from=8, to=2-4]
      \arrow["{\eta_x}"{description, pos=0.6}, draw=none, from=8, to=3]
      \arrow[""{name=9, anchor=center, inner sep=0}, "{m_i}", "\shortmid"{marking}, shorten >=8pt, from=1-2, to=8]
      \arrow[""{name=10, anchor=center, inner sep=0}, "{m_{i+1}}", "\shortmid"{marking}, shorten <=8pt, from=8, to=1-5]
      \arrow["{1_{m_i}}"{description}, draw=none, from=9, to=4]
      \arrow["{1_{m_{i+1}}}"{description}, draw=none, from=10, to=5]
    \end{tikzcd}.
  \end{equation*}
  A \define{unital virtual double category} is a virtual double category
  equipped with a choice of unit for each object.
\end{definition}

Lax functors, lax transformations, modules, and multimodulations assemble into a
virtual double categories with strongly representable units. This result has
been stated without proof for strict transformations in \cite[Theorem
4.3]{pare2011} and \cite[Theorem 1.2.5]{pare2013}.

\begin{theorem}[Virtual double category of lax functors]
  \label{thm:lax-functor-vdbl-category}
  For any double categories $\dbl{D}$ and $\dbl{E}$, there is a unital virtual
  double category $\vLaxFun_\ell(\dbl{D},\dbl{E})$ whose objects are lax double
  functors $\dbl{D} \to \dbl{E}$, arrows are lax natural transformations,
  proarrows are modules, and multicells are multimodulations.

  Multimodulations are composed by composing, in the same shape, their component
  cells in $\dbl{E}$. For example, a binary composite of binary multimodulations
  \begin{equation*}
    \begin{tikzcd}[row sep=scriptsize]
      {F_0} & {F_1} & {F_2} & {F_3} & {F_4} \\
      {G_0} && {G_1} && {G_2} \\
      {H_0} &&&& {H_1}
      \arrow["\alpha"', from=1-1, to=2-1]
      \arrow["\gamma"', from=2-1, to=3-1]
      \arrow[""{name=0, anchor=center, inner sep=0}, "P"', "\shortmid"{marking}, from=3-1, to=3-5]
      \arrow[""{name=1, anchor=center, inner sep=0}, "{N_1}"', "\shortmid"{marking}, from=2-1, to=2-3]
      \arrow[""{name=2, anchor=center, inner sep=0}, "{N_2}"', "\shortmid"{marking}, from=2-3, to=2-5]
      \arrow["{M_1}", "\shortmid"{marking}, from=1-1, to=1-2]
      \arrow["{M_2}", "\shortmid"{marking}, from=1-2, to=1-3]
      \arrow["{M_3}", "\shortmid"{marking}, from=1-3, to=1-4]
      \arrow["{M_4}", "\shortmid"{marking}, from=1-4, to=1-5]
      \arrow["\epsilon"', from=1-3, to=2-3]
      \arrow["\beta", from=1-5, to=2-5]
      \arrow["\delta", from=2-5, to=3-5]
      \arrow["\rho"{description}, draw=none, from=2-3, to=0]
      \arrow["\mu"{description}, draw=none, from=1-2, to=1]
      \arrow["\nu"{description}, draw=none, from=1-4, to=2]
    \end{tikzcd}
  \end{equation*}
  has components, for proarrows
  $x_0 \xproto{m_1} x_1 \xproto{m_2} x_2 \xproto{m_3} x_3 \xproto{m_4} x_4$ in
  $\dbl{D}$,
  \begin{equation*}
    \begin{tikzcd}[row sep=scriptsize]
      {F_0 x_0} & {F_1 x_1} & {F_2 x_2} & {F_3 x_3} & {F_4 x_4} \\
      {G_0 x_0} && {G_1 x_2} && {G_2 x_4} \\
      {H_0 x_0} &&&& {H_1 x_4}
      \arrow["{\alpha_{x_0}}"', from=1-1, to=2-1]
      \arrow["{\gamma_{x_0}}"', from=2-1, to=3-1]
      \arrow[""{name=0, anchor=center, inner sep=0}, "{P(m_1 \odot m_2 \odot m_3 \odot m_4)}"', "\shortmid"{marking}, from=3-1, to=3-5]
      \arrow[""{name=1, anchor=center, inner sep=0}, "{N_1 (m_1 \odot m_2)}"', "\shortmid"{marking}, from=2-1, to=2-3]
      \arrow[""{name=2, anchor=center, inner sep=0}, "{N_2 (m_3 \odot m_4)}"', "\shortmid"{marking}, from=2-3, to=2-5]
      \arrow["{M_1 m_1}", "\shortmid"{marking}, from=1-1, to=1-2]
      \arrow["{M_2 m_2}", "\shortmid"{marking}, from=1-2, to=1-3]
      \arrow["{M_3 m_3}", "\shortmid"{marking}, from=1-3, to=1-4]
      \arrow["{M_4 m_3}", "\shortmid"{marking}, from=1-4, to=1-5]
      \arrow["{\epsilon_{x_2}}"', from=1-3, to=2-3]
      \arrow["{\beta_{x_4}}", from=1-5, to=2-5]
      \arrow["{\delta_{x_4}}", from=2-5, to=3-5]
      \arrow["{\rho_{m_1 \odot m_2, m_3 \odot m_4}}"{description}, draw=none, from=2-3, to=0]
      \arrow["{\mu_{m_1,m_2}}"{description}, draw=none, from=1-2, to=1]
      \arrow["{\nu_{m_3,m_4}}"{description}, draw=none, from=1-4, to=2]
    \end{tikzcd}.
  \end{equation*}
  An identity modulation as on the left
  \begin{equation*}
    \begin{tikzcd}[row sep=scriptsize]
      F & G \\
      F & G
      \arrow[""{name=0, anchor=center, inner sep=0}, "M", "\shortmid"{marking}, from=1-1, to=1-2]
      \arrow[""{name=1, anchor=center, inner sep=0}, "M"', "\shortmid"{marking}, from=2-1, to=2-2]
      \arrow[Rightarrow, no head, from=1-1, to=2-1]
      \arrow[Rightarrow, no head, from=1-2, to=2-2]
      \arrow["{1_M}"{description}, draw=none, from=0, to=1]
    \end{tikzcd}
    \qquad\leftrightsquigarrow\qquad
    \begin{tikzcd}[row sep=scriptsize]
      Fx & Gy \\
      Fx & Gy
      \arrow[""{name=0, anchor=center, inner sep=0}, "Mm", "\shortmid"{marking}, from=1-1, to=1-2]
      \arrow[Rightarrow, no head, from=1-1, to=2-1]
      \arrow[""{name=1, anchor=center, inner sep=0}, "Mm"', "\shortmid"{marking}, from=2-1, to=2-2]
      \arrow[Rightarrow, no head, from=1-2, to=2-2]
      \arrow["{1_{Mm}}"{description}, draw=none, from=0, to=1]
    \end{tikzcd},
    \quad x \xproto{m} y \text{ in } \dbl{D}.
  \end{equation*}
  has components given by identities in $\dbl{E}_1$ as on the right.

  The units in the virtual double category $\vLaxFun_\ell(\dbl{D},\dbl{E})$ are
  the identity modules $\id_F$, defined by $\id_F(m) = F(m)$,
  $\id_F(\alpha) = F(\alpha)$, and the laxators of $F$, together with the cells
  $\eta_F: ()_F \to \id_F$, $(\eta_F)_x = F_x$, defined by the unitors of $F$.

  Similarly, there are unital virtual double categories
  $\vLaxFun_c(\dbl{D},\dbl{E})$ and $\vLaxFun_\pseudo(\dbl{D},\dbl{E})$ whose
  arrows are colax/oplax and pseudo natural transformations, respectively.
\end{theorem}
\begin{proof}
  We have already shown that lax functors and lax transformations form a
  category (\cref{prop:lax-functor-category}). The associativity and unitality
  laws of the virtual double category $\vLaxFun_\ell(\dbl{D},\dbl{E})$ follow
  directly from those of the double category $\dbl{E}$. We need to show that
  composite and identity multimodulations satisfy the equivariance and
  naturality axioms. We omit the proof of equivariance, which is the same for
  strict and lax transformations, and we prove that composite multimodulations
  are natural. To keep the notation manageable, we do this for the binary
  composite of binary multimodulations in the theorem statement.

  Given composable cells
  $\inlineCell{x_{i-1}}{x_i}{y_{i-1}}{y_i}{m_i}{n_i}{f_{i-1}}{f_i}{\phi_i}$,
  $i=1,\dots,4$, in $\dbl{D}$ and abbreviating
  $n \coloneqq n_1 \odot n_2 \odot n_3 \odot n_4$, we have by expanding the
  definitions
  \begin{equation*}
    \begin{dblArray}{ccccc}
      \Block{2-1}{(\alpha \gamma)_{f_0}} & M_1 \phi_1 & M_2 \phi_2 & M_3 \phi_3 & M_4 \phi_4 \\
      & \Block{1-4}{((\mu,\nu) \cdot \rho)_{n_1,n_2,n_3,n_4}} &&& \\
      \Block{1-5}{P^\ell_{y_0,n}} &&&&
    \end{dblArray}
    \quad=\quad
    \begin{dblArray}{cccccc}
      \id_{\alpha_{x_0}} & \Block{2-1}{\alpha_{f_0}} & M_1 \phi_1 & M_2 \phi_2 & M_3 \phi_3 & M_4 \phi_4 \\
      \Block{2-1}{\gamma_{f_0}} && \Block{1-2}{\mu_{n_1,n_2}} && \Block{1-2}{\nu_{n_3,n_4}} & \\
      & \gamma_{\id_{y_0}} & \Block{1-4}{\rho_{n_1 \odot n_2, n_3 \odot n_4}} &&& \\
      \Block{1-2}{(H_0)_{y_0,y_0}} && \Block{1-4}{1_{Pn}} &&& \\
      \Block{1-6}{P^\ell_{y_0,n}} &&&&&
    \end{dblArray}.
  \end{equation*}
  By the associativity of the module $P$, then the left equivariance of the
  multimodulation $\rho$, this is equal to the left-hand side
  \begin{equation*}
    \begin{dblArray}{cccccc}
      \id_{\alpha_{x_0}} & \Block{2-1}{\alpha_{f_0}} & M_1 \phi_1 & M_2 \phi_2 & M_3 \phi_3 & M_4 \phi_4 \\
      \Block{2-1}{\gamma_{f_0}} && \Block{1-2}{\mu_{n_1,n_2}} && \Block{1-2}{\nu_{n_3,n_4}} & \\
      & \Block{1-3}{(N_1)^\ell_{y_0, n_1 \odot n_2}} &&& \Block{1-2}{1_{N_2(n_3 \odot n_4)}} & \\
      1_{H_0 \id_{y_0}} & \Block{1-5}{\rho_{n_1 \odot n_2, n_3 \odot n_4}} &&&& \\
      \Block{1-6}{P^\ell_{y_0,n}} &&&&&
    \end{dblArray}
    \quad=\quad
    \begin{dblArray}{ccccc}
      \id_{\alpha_{x_0}} & \mu_{m_1,m_2} & \Block{2-1}{\epsilon_{f_2}} & M_3 \phi_3 & M_4 \phi_4 \\
      \Block{2-1}{\gamma_{f_0}} & N_1(\phi_1 \odot \phi_2) && \Block{1-2}{\nu_{n_3,n_4}} & \\
      & \Block{1-2}{(N_1)^r_{n_1 \odot n_2, y_2}} && \Block{1-2}{1_{N_2(n_3 \odot n_4)}} & \\
      1_{H_0 \id_{y_0}} & \Block{1-4}{\rho_{n_1 \odot n_2, n_3 \odot n_4}} &&& \\
      \Block{1-5}{P^\ell_{y_0,n}} &&&&
    \end{dblArray}
  \end{equation*}
  thence equal to the right-hand side by the naturality of the multimodulation
  $\mu$. Applying the inner equivariance of $\rho$ and then the naturality of
  $\nu$, we have
  \begin{equation*}
    \begin{dblArray}{ccccc}
      \id_{\alpha_{x_0}} & \mu_{m_1,m_2} & \Block{2-1}{\epsilon_{f_2}} & M_3 \phi_3 & M_4 \phi_4 \\
      \Block{2-1}{\gamma_{f_0}} & N_1(\phi_1 \odot \phi_2) && \Block{1-2}{\nu_{n_3,n_4}} & \\
      & 1_{N_1(n_1 \odot n_2)} & \Block{1-3}{N^\ell_{y_2, n_3 \odot n_4}} && \\
      1_{H_0 \id_{y_0}} & \Block{1-4}{\rho_{n_1 \odot n_2, n_3 \odot n_4}} &&& \\
      \Block{1-5}{P^\ell_{y_0,n}} &&&&
    \end{dblArray}
    \quad=\quad
    \begin{dblArray}{cccc}
      \id_{\alpha_{x_0}} & \mu_{m_1,m_2} & \nu_{m_3,m_4} & \Block{2-1}{\beta_{f_4}} \\
      \Block{2-1}{\gamma_{f_0}} & N_1(\phi_1 \odot \phi_2) & N_2(\phi_3 \odot \phi_4) & \\
      & 1_{N_1(n_1 \odot n_2)} & \Block{1-2}{N^r_{n_3 \odot n_4, y_4}} & \\
      1_{H_0 \id_{y_0}} & \Block{1-3}{\rho_{n_1 \odot n_2, n_3 \odot n_4}} && \\
      \Block{1-4}{P^\ell_{y_0,n}} &&&
    \end{dblArray}.
  \end{equation*}
  By the right equivariance of $\rho$ and the associativity of $P$, this is
  equal to
  \begin{equation*}
    \begin{dblArray}{cccc}
      \id_{\alpha_{x_0}} & \mu_{m_1,m_2} & \nu_{m_3,m_4} & \Block{2-1}{\beta_{f_4}} \\
      \Block{2-1}{\gamma_{f_0}} & N_1(\phi_1 \odot \phi_2) & N_2(\phi_3 \odot \phi_4) & \\
      & \Block{1-2}{\rho_{n_1 \odot n_2, n_3 \odot n_4}} && \beta_{\id_{y_4}} \\
      1_{H_0 \id_{y_0}} & \Block{1-3}{P^r_{n,y_4}} && \\
      \Block{1-4}{P^\ell_{y_0,n}} &&&
    \end{dblArray}
    \quad=\quad
    \begin{dblArray}{cccc}
      \id_{\alpha_{x_0}} & \mu_{m_1,m_2} & \nu_{m_3,m_4} & \Block{2-1}{\beta_{f_4}} \\
      \Block{2-1}{\gamma_{f_0}} & N_1(\phi_1 \odot \phi_2) & N_2(\phi_3 \odot \phi_4) & \\
      & \Block{1-2}{\rho_{n_1 \odot n_2, n_3 \odot n_4}} && \beta_{\id_{y_4}} \\
      \Block{1-3}{P^\ell_{y_0,n}} &&& 1_{G_2 \id_{y_4}} \\
      \Block{1-4}{P^r_{n,y_4}} &&&
    \end{dblArray}.
  \end{equation*}
  Finally, applying the associativity of $P$ yet again and collecting the
  definitions, we obtain the desired equality
  \begin{equation*}
    \begin{dblArray}{cccc}
      \mu_{m_1,m_2} & \nu_{m_3,m_4} & \id_{\delta_{x_4}} & \Block{2-1}{\beta_{f_4}} \\
      \Block{1-2}{\rho_{m_1 \odot m_2, m_3 \odot m_4}} && \Block{2-1}{\delta_{f_4}} & \\
      \Block{1-2}{P(\phi_1 \odot \cdots \odot \phi_4)} &&& \beta_{\id_{y_4}} \\
      \Block{1-3}{P^r_{n,y_4}} &&& 1_{G_2 \id_{y_4}} \\
      \Block{1-4}{P^r_{n,y_4}} &&&
    \end{dblArray}
    \quad=\quad
    \begin{dblArray}{cc}
      ((\mu,\nu) \cdot \rho)_{m_1,m_2,m_3,m_4} & \Block{2-1}{(\beta \delta)_{f_4}} \\
      P(\phi_1 \odot \phi_2 \odot \phi_3 \odot \phi_4) & \\
      \Block{1-2}{P^r_{n,y_4}} &
    \end{dblArray}.
  \end{equation*}

  It remains to prove that identity modules are units in the virtual double
  category $\vLaxFun_\ell(\dbl{D},\dbl{E})$. Given a multimodulation
  $\mu: (M_1,\dots,M_k) \to N$ and an index $0 \leq i \leq k$, the unique
  factorization
  \begin{equation*}
    \bar\mu: (M_1,\dots,M_i,\id_{F_i},M_{i+1},\dots,M_k) \to N
  \end{equation*}
  of $\mu$ through $\eta_{F_i}$ is given by either side of the corresponding
  equivariance axiom with the extra module being $\id_{F_i}$. Specifically, when
  $k = 0$, we use nullary equivariance; when $k \geq 1$ and $i = 0$, left
  equivariance; when $k \geq 1$ and $i = k$, right equivariance, and otherwise we
  use inner equivariance. The validity and uniqueness of these factorizations
  follow from the unitality axioms for natural transformations and modules.
\end{proof}

Paré has found sufficient conditions for the virtual double category of lax
functors, strict transformations, and modules to be a genuine double category,
such as the so-called condition \emph{AFP} \cite[Definition 3.1.1]{pare2013}.
This is a factorization condition for multicells with a binary proarrow domain
such that if the domain double category $\dbl{D}$ satisfies it and the codomain
$\dbl{E}$ is locally cocomplete (that is, each hom category has all colimits and
they are preserved by external composition in each argument), then the virtual
double category $\vLaxFun(\dbl{D},\dbl{E})$ is in fact a double category
\cite[Theorems 4.0.1 \& 5.1.10]{pare2013}. As examples, double categories
generated by either a 2-category or a bicategory satisfy the condition AFP
\cite[Corollaries 3.1.5 \& 3.1.7]{pare2013}.

Another approach is that of Cruttwell and Shulman \cite[Theorem
5.2]{cruttwell2010}, showing that a virtual double category is representable as
a double category whenever any sequence of proarrows satisfies the universal
property of a \emph{composite}, that is, has a certain proarrow and opcartesian
multicell associated with it \cite[Definition 5.1]{cruttwell2010}. Ultimately,
these references specify essentially the same universal property. Either
approach could conceivably be adapted to the case of \emph{lax} transformations,
but at least in the case when the domain double category
$\dbl{D} = \VerDbl(\bicat{A})$ is horizontally trivial, a direct proof is most
straightforward. We work with unitary lax functors to simplify the computations.

\begin{proposition}[Sufficient conditions for representability]
  \label{prop:lax-func-form-double-category}
  If $\bicat{A}$ is a 2-category and $\dbl{E}$ is a double category, then
  $\vLaxFun_{\ell,u}(\VerDbl(\bicat{A}),\dbl{E})$ is representable as a double
  category with modules as proarrows and modulations as cells. Likewise, if
  $\bicat{D}$ is a bicategory, then $\vLaxFun_{\ell}(\HorDbl(\bicat{D}),\Span)$
  is a double category.
\end{proposition}
\begin{proof}
  The last statement concerning bicategories follows from the known results
  mentioned above since in this case the double category $\HorDbl(\bicat{D})$
  has no non-identity arrows and so any lax transformation is automatically
  strict.

  For the first statement, suppose $\bicat{A}$ is a 2-category and consider two
  modules $M\colon F\proto G$ and $N\colon G\proto H$ between unitary lax
  functors $\dbl{D} \coloneqq \VerDbl(\bicat{A}) \to \dbl{E}$. The composite
  $M\odot N$ first needs two items of data specified, namely, its action on
  proarrows and its action on cells. But proarrows in $\dbl{D}$ are only
  identities. So, take as a definition
  \begin{equation*}
    (M\odot N)(\id_x) \coloneqq M\id_x\odot N\id_x
  \end{equation*}
  using the composition in $\dbl{E}$. Cells of $\dbl{D}$ are then of the form
  $\alpha\colon f \To g$ with proarrow identity domain $\id_x$ and codomain
  $\id_y$. Take as the definition of $(M\odot N)(\alpha)$ the composite cell
  \begin{equation*}
    (M\odot N)(\alpha)
    \quad\coloneqq\quad
    \begin{tikzcd}
      Fx &&& Hx \\
      Fx & Gx & Gx & Hx \\
      Fy & Gy & Gy & Hy \\
      Fy &&& Hy
      \arrow[""{name=0, anchor=center, inner sep=0}, "{M\id_x}", "\shortmid"{marking}, from=2-1, to=2-2]
      \arrow[""{name=1, anchor=center, inner sep=0}, "{\id_{Gx}}", "\shortmid"{marking}, from=2-2, to=2-3]
      \arrow[""{name=2, anchor=center, inner sep=0}, "{N\id_x}", "\shortmid"{marking}, from=2-3, to=2-4]
      \arrow["Hg", from=2-4, to=3-4]
      \arrow["Ff"', from=2-1, to=3-1]
      \arrow[""{name=3, anchor=center, inner sep=0}, "{M\id_y}"', "\shortmid"{marking}, from=3-1, to=3-2]
      \arrow[from=2-2, to=3-2]
      \arrow[""{name=4, anchor=center, inner sep=0}, "{\id_{Gy}}"', "\shortmid"{marking}, from=3-2, to=3-3]
      \arrow[from=2-3, to=3-3]
      \arrow[""{name=5, anchor=center, inner sep=0}, "{N\id_y}"', "\shortmid"{marking}, from=3-3, to=3-4]
      \arrow[Rightarrow, no head, from=3-1, to=4-1]
      \arrow[""{name=6, anchor=center, inner sep=0}, "{M(\id_y)\odot N(\id_y)}"', "\shortmid"{marking}, from=4-1, to=4-4]
      \arrow[Rightarrow, no head, from=3-4, to=4-4]
      \arrow[Rightarrow, no head, from=1-1, to=2-1]
      \arrow[""{name=7, anchor=center, inner sep=0}, "{M(\id_x)\odot N(\id_x)}", "\shortmid"{marking}, from=1-1, to=1-4]
      \arrow[Rightarrow, no head, from=1-4, to=2-4]
      \arrow["\cong"{description}, draw=none, from=4, to=6]
      \arrow["\cong"{description}, draw=none, from=7, to=1]
      \arrow["{N\id_g}"{description}, draw=none, from=2, to=5]
      \arrow["{M\id_f}"{description}, draw=none, from=0, to=3]
      \arrow["G\alpha"{description}, draw=none, from=1, to=4]
    \end{tikzcd}.
  \end{equation*}
  Since the lax functors involved are unitary, define the associated left and
  right actions on the proposed composite to be the globular cells given by the
  left and right unitors in $\dbl{E}$. It is straightforward to verify that with
  these definitions $M\odot N$ is a module as in \cref{def:module}. We have also
  to see that external compositions of modulations is well-defined. For this,
  suppose that
  \begin{equation*}
    \begin{tikzcd}
      F & G & H \\
      J & K & L
      \arrow[""{name=0, anchor=center, inner sep=0}, "M", "\shortmid"{marking}, from=1-1, to=1-2]
      \arrow[""{name=1, anchor=center, inner sep=0}, "N", "\shortmid"{marking}, from=1-2, to=1-3]
      \arrow["\gamma", from=1-3, to=2-3]
      \arrow["\beta", from=1-2, to=2-2]
      \arrow["\alpha"', from=1-1, to=2-1]
      \arrow[""{name=2, anchor=center, inner sep=0}, "P"', "\shortmid"{marking}, from=2-1, to=2-2]
      \arrow[""{name=3, anchor=center, inner sep=0}, "Q"', "\shortmid"{marking}, from=2-2, to=2-3]
      \arrow["\mu"{description}, draw=none, from=0, to=2]
      \arrow["\nu"{description}, draw=none, from=1, to=3]
    \end{tikzcd}
  \end{equation*}
  are modulations, where $\alpha$, $\beta$, and $\gamma$ are lax
  transformations. The required data for the purported composite $\mu\odot\nu$
  is a cell $(\mu\odot\nu)_{\id_x}$ associated to the proarrow
  $\id_x\colon x\proto x$, which we define as
  \begin{equation*}
    (\mu\odot\nu)_{\id_x}
    \qquad\coloneqq\qquad
    \begin{tikzcd}
      Fx &&& Hx \\
      Fx & Gx & Gx & Hx \\
      Jx & Kx & Kx & Lx \\
      Jx &&& Lx
      \arrow[""{name=0, anchor=center, inner sep=0}, "{M\id_x}", "\shortmid"{marking}, from=2-1, to=2-2]
      \arrow[""{name=1, anchor=center, inner sep=0}, "{\id_{Gx}}", "\shortmid"{marking}, from=2-2, to=2-3]
      \arrow[""{name=2, anchor=center, inner sep=0}, "{N\id_x}", "\shortmid"{marking}, from=2-3, to=2-4]
      \arrow["{\gamma_x}", from=2-4, to=3-4]
      \arrow[from=2-3, to=3-3]
      \arrow["{\alpha_x}"', from=2-1, to=3-1]
      \arrow[from=2-2, to=3-2]
      \arrow[""{name=3, anchor=center, inner sep=0}, "{\id_{Kx}}"', "\shortmid"{marking}, from=3-2, to=3-3]
      \arrow[""{name=4, anchor=center, inner sep=0}, "{Q\id_x}"', "\shortmid"{marking}, from=3-3, to=3-4]
      \arrow[Rightarrow, no head, from=3-1, to=4-1]
      \arrow[""{name=5, anchor=center, inner sep=0}, "{P\id_x\odot Q\id_x}"', "\shortmid"{marking}, from=4-1, to=4-4]
      \arrow[Rightarrow, no head, from=3-4, to=4-4]
      \arrow[Rightarrow, no head, from=1-1, to=2-1]
      \arrow[""{name=6, anchor=center, inner sep=0}, "{M\id_x\odot N\id_x}", "\shortmid"{marking}, from=1-1, to=1-4]
      \arrow[Rightarrow, no head, from=1-4, to=2-4]
      \arrow[""{name=7, anchor=center, inner sep=0}, "{P\id_x}"', "\shortmid"{marking}, from=3-1, to=3-2]
      \arrow["\cong"{description}, draw=none, from=6, to=1]
      \arrow["\cong"{description}, draw=none, from=3, to=5]
      \arrow["{\beta_{\id_x}}"{description}, draw=none, from=1, to=3]
      \arrow["{\nu_{\id_x}}"{description}, draw=none, from=2, to=4]
      \arrow["{\mu_{\id_x}}"{description}, draw=none, from=0, to=7]
    \end{tikzcd}.
  \end{equation*}
  The equivariance condition in \cref{def:modulation} follows immediately since
  the actions are given by unitors in $\dbl{E}$. We now check the naturality
  condition. Given a cell $\alpha\colon f \To g$ as in the considerations above,
  we can calculate that
  \begin{align*}
    \begin{dblArray}{cccc}
      \Block{2-1}{\alpha_f} & M\id_f & G\alpha & N\id_g \\
      & \mu_{\id_y} & \beta_{\id_y} & \nu_{\id_y}
    \end{dblArray}
    \quad&=\quad
    \begin{dblArray}{cccc}
      \mu_{\id_x} &\Block{2-1}{\beta_f} & G\alpha &N\id_g  \\
      P\id_f& & \beta_{\id_y} & \nu_{\id_g}
    \end{dblArray}\\
    \quad&=\quad
    \begin{dblArray}{cccc}
      \mu_{\id_x} &\beta_{\id_x} & \Block{2-1}{\beta_g} &N\id_g  \\
      P\id_f & K\alpha & & \nu_{\id_g}
    \end{dblArray}\\
    \quad&=\quad
    \begin{dblArray}{cccc}
      \mu_{\id_x} &\beta_{\id_x} & \nu_{\id_x} & \Block{2-1}{\gamma_g}  \\
      P\id_f & K\alpha & Q\id_g &
    \end{dblArray}
  \end{align*}
  using the naturality condition for the modulations $\mu$ and $\gamma$ in the
  first and last equalities, and the cell naturality condition for the lax
  transformation $\beta$ in the middle equality. Notice that this is where
  laxness manifests in the proof. It remains to see that the external
  composition of modules is associative up to (coherent) isomorphism. This is
  now relatively straightforward using pointwise definitions of the required
  associativity isomorphisms. Note that since these are by construction
  \emph{globular}, the required well-definition and naturality and coherence
  conditions follow almost immediately from the corresponding properties in
  $\dbl{E}$. The unitors and their coherence laws follow similarly.
\end{proof}

Owing to the fact that most of the simple double theories in
\cref{sec:simple-double-theories} are either 2- or bi-categorical, we have the
following result. Note that here we implicitly invoke
\cref{cor:unitalize-span-valued-2-cat-iso-FINAL} and \cref{cor:2-categories-of-models-non-cartesian} to provide the proof.

\begin{corollary} \label{cor:double-categories-of-models-examples}
  For each of the following simple theories, span-valued models are the objects
  of a double category of lax transformations, modules and modulations:
  \begin{enumerate}
    \item categories (\cref{th:unit-theory}), functors (\cref{th:walking-arrow}),
    transformations (\cref{th:special-cell});
    \item adjunctions (\cref{th:adjunctions}), monads (\cref{th:monad}),
    Frobenius monads (\cref{th:frob-monad}).
  \end{enumerate}
\end{corollary}

Multimodulations, like modulations, require no extra conditions in the presence
of finite products. We omit a precise statement and proceed directly to:

\begin{corollary}[Virtual double category of cartesian lax functors]
  \label{cor:cartesian-lax-functor-vdbl-category}
  For any cartesian double categories $\dbl{D}$ and $\dbl{E}$, there is a unital
  virtual double category $\vCartLaxFun_\ell(\dbl{D},\dbl{E})$ whose objects are
  cartesian lax functors $\dbl{D} \to \dbl{E}$, arrows are cartesian lax natural
  transformations, proarrows are cartesian modules, and multicells are
  multimodulations.

  Similarly, there are unital virtual double categories
  $\vCartLaxFun_c(\dbl{D},\dbl{E})$ and $\vCartLaxFun_\pseudo(\dbl{D},\dbl{E})$
  whose arrows are cartesian oplax transformations and cartesian pseudo
  transformations, respectively.
\end{corollary}
\begin{proof}
  This corollary follows directly from the previous
  \cref{thm:lax-functor-vdbl-category} since we have already observed that
  cartesian lax transformations are closed under composition, and identity
  modules $\id_F$ are plainly cartesian whenever $F$ is.
\end{proof}

\begin{corollary}
  For any cartesian 2-category $\dbl{A}$ and cartesian double category
  $\dbl{E}$, the virtual double category
  $\vCartLaxFun_{\ell,u}(\VerDbl(\bicat{A}),\dbl{E})$ is representable as a
  double category.
\end{corollary}
\begin{proof}
  It needs to be seen that cartesian modules are closed under external
  composition as defined in the proof of
  \cref{prop:lax-func-form-double-category}. Let $M\colon F\proto G$ and
  $N\colon G\proto H$ be cartesian modules between cartesian unitary lax functors
  $F,G,H\colon \VerDbl(\bicat{A})\to \dbl{E}$. We need to check the one
  condition of \cref{def:cartesian-module}, namely, that the cell
  \begin{equation*}
    \begin{tikzcd}
      {F(x\times x')} &&&&& {H(x\times x')} \\
      {Fx\times Fx'} &&&&& {Hx\times Hx'}
      \arrow[""{name=0, anchor=center, inner sep=0}, "{M\id_{x\times x'}\odot N\id_{x\times x'}}", "\shortmid"{marking}, from=1-1, to=1-6]
      \arrow[from=1-1, to=2-1]
      \arrow[""{name=1, anchor=center, inner sep=0}, "{M\id_x\odot N\id_x\times M\id_{x'}\odot N\id_{x'}}"', "\shortmid"{marking}, from=2-1, to=2-6]
      \arrow[from=1-6, to=2-6]
      \arrow["{\langle M\id_{\pi_{x,x'}}\odot N\id_{\pi_{x,x'}},M\id_{\pi'_{x,x'}}\odot N\id_{\pi'_{x,x'}}\rangle }"{description}, draw=none, from=0, to=1]
    \end{tikzcd}
  \end{equation*}
  is invertible. Note that this cell is well-typed and well-defined; for as the
  double category $\VerDbl(\bicat{A})$ is effectively a 2-category, each product
  $\id_x\times \id_{x'}$ is identically $\id_{x\times x'}$; accordingly, the
  projections are $\id_{\pi_{x,x'}}$ and $\id_{\pi'_{x,x'}}$. Now, to prove that
  the cell is invertible, note that it is in fact equal to the composite
  \begin{equation*}
    \begin{tikzcd}
      {F(x\times x')} &&& {G(x\times x')} &&& {H(x\times x')} \\
      {Fx\times Fx'} &&& {Gx\times Gx'} &&& {Hx\times Hx'} \\
      {Fx\times Fx'} &&&&&& {Hx\times Hx'}
      \arrow[from=1-1, to=2-1]
      \arrow[from=2-1, to=3-1]
      \arrow[""{name=0, anchor=center, inner sep=0}, "{M\id_{x\times x'}}", "\shortmid"{marking}, from=1-1, to=1-4]
      \arrow[from=1-4, to=2-4]
      \arrow[""{name=1, anchor=center, inner sep=0}, "{M\id_x\times M\id_{x'}}"', "\shortmid"{marking}, from=2-1, to=2-4]
      \arrow[""{name=2, anchor=center, inner sep=0}, "{N\id_{x\times x'}}", "\shortmid"{marking}, from=1-4, to=1-7]
      \arrow[from=1-7, to=2-7]
      \arrow[from=2-7, to=3-7]
      \arrow[""{name=3, anchor=center, inner sep=0}, "{N\id_x\times N\id_{x'}}"', "\shortmid"{marking}, from=2-4, to=2-7]
      \arrow[""{name=4, anchor=center, inner sep=0}, "{(M\id_x\odot N\id_x)\times (M\id_{x'}\odot N\id_{x'})}"', "\shortmid"{marking}, from=3-1, to=3-7]
      \arrow["{\langle\pi_{M\id_x,M\id_{x'}}\odot \pi_{N\id_x,N\id_{x'}},\pi'_{M\id_x,M\id_{x'}}\odot \pi'_{N\id_x,N\id_{x'}}\rangle}"{description}, draw=none, from=2-4, to=4]
      \arrow["{\langle M\id_{\pi_{x,x'}},M\id_{\pi'_{x,x'}}\rangle}"{description}, draw=none, from=0, to=1]
      \arrow["{\langle N\id_{\pi_{x,x'}},N\id_{\pi'_{x,x'}}\rangle}"{description}, draw=none, from=2, to=3]
    \end{tikzcd}
  \end{equation*}
  in which each constituent cell is invertible. The top two cells are invertible
  since $M$ and $N$ are cartesian modules; the bottom cell is an interchange
  cell in $\dbl{E}$ from \cref{eq:dbl-product-compose-comparison} which must be
  invertible since $\dbl{E}$ is cartesian. That the cells in the two displays
  are actually equal can be checked easily using the projections from the
  proarrow codomain of the above cell and appealing to the uniqueness property
  of morphisms into a product.
\end{proof}

The 2-category of models of a cartesian double theory
(\cref{cor:models-2-category}) now extends to a unital virtual double category.

\begin{corollary}[Virtual double category of models]
  Let $\dbl{T}$ be a cartesian double theory and let $\dbl{S}$ be a cartesian
  double category. A \define{bimodule} between models $M$ and $M'$ of $\dbl{T}$
  in $\dbl{S}$ is a cartesian module $M \proTo M'$. A
  \define{multitransformation} from a composable sequence of bimodules of models
  to another bimodule is a multimodulation.

  There is a unital virtual double category of models of the theory $\dbl{T}$ in
  $\dbl{S}$, (lax, oplax, pseudo, or strict) maps between models, bimodules
  between models, and multitransformations between those.

  Moreover, when the theory $\dbl{T}$ is generated by a 2-category or a
  bicategory (i.e., has only trivial arrows or trivial proarrows) and the
  semantics $\dbl{S}$ is a cartesian equipment with local coequalizers, the
  virtual double category of models is representable as a double category.
\end{corollary}

Similarly, a restriction sketch has a unital virtual double category of models
in any cartesian equipment. One might even expect that a double theory
interpreted in an equipment should have a virtual \emph{equipment} of models.
That turns to be true when the maps between models are taken to be pseudo or
strict. When the maps are lax or oplax, only certain restrictions exist.

\begin{theorem}[Virtual equipment of lax functors] \label{thm:lax-functor-veqp}
  Let $\dbl{D}$ be a double category and let $\dbl{E}$ be an equipment. Then:
  \begin{enumerate}[(i),noitemsep]
    \item $\vLaxFun_\ell(\dbl{D},\dbl{E})$ has restrictions along lax
      transformations $\alpha$ and $\beta$ whenever $\alpha$ is pseudo.
    \item $\vLaxFun_c(\dbl{D},\dbl{E})$ has restrictions along oplax
      transformations $\alpha$ and $\beta$ whenever $\beta$ is pseudo.
    \item $\vLaxFun_\pseudo(\dbl{D},\dbl{E})$ has restrictions along any pseudo
      transformations $\alpha$ and $\beta$.
  \end{enumerate}
  Consequently, the virtual double category $\vLaxFun_\pseudo(\dbl{D},\dbl{E})$
  of lax double functors $\dbl{D} \to \dbl{E}$, pseudo natural transformations,
  modules, and multimodulations is a virtual equipment.
\end{theorem}
\begin{proof}
  We have shown that the virtual double categories in question exist and have
  units in \cref{thm:lax-functor-vdbl-category}, so we need only show that they
  have the claimed restrictions, in the sense of \cite[\mbox{Definition
    7.1}]{cruttwell2010}. We will construct the restrictions ``pointwise'' using
  the restrictions in the target equipment $\dbl{E}$. We prove the first
  statement about lax transformations; the others are proved similarly.

  Given a niche as on the left
  \begin{equation*}
    \begin{tikzcd}
      F & G \\
      H & K
      \arrow["\alpha"', from=1-1, to=2-1]
      \arrow["\beta", from=1-2, to=2-2]
      \arrow["N"', "\shortmid"{marking}, from=2-1, to=2-2]
    \end{tikzcd}
    \qquad\leadsto\qquad
    \begin{tikzcd}
      F & G \\
      H & K
      \arrow["\alpha"', from=1-1, to=2-1]
      \arrow["\beta", from=1-2, to=2-2]
      \arrow[""{name=0, anchor=center, inner sep=0}, "N"', "\shortmid"{marking}, from=2-1, to=2-2]
      \arrow[""{name=1, anchor=center, inner sep=0}, "{N(\alpha,\beta)}", "\shortmid"{marking}, from=1-1, to=1-2]
      \arrow["\res"{description}, draw=none, from=1, to=0]
    \end{tikzcd}
  \end{equation*}
  where $\alpha$ is a pseudo transformation, $\beta$ is a lax transformation, and $N$ is
  a module between lax functors, we must construct a module $N(\alpha,\beta)$ along with
  a modulation as on the right. Define the module $N(\alpha,\beta)$ at a proarrow
  $m: x \proto y$ in $\dbl{D}$ to be the restricted proarrow
  $N(\alpha,\beta)(m) \coloneqq (Nm)(\alpha_x,\beta_y)$ in $\dbl{E}$. Next, for each cell
  $\stdInlineCell{\sigma}$ in $\dbl{D}$, define the cell
  \begin{equation*}
    \begin{tikzcd}[column sep=large]
      Fx & Gy \\
      Fw & Gz
      \arrow["Ff"', from=1-1, to=2-1]
      \arrow[""{name=0, anchor=center, inner sep=0}, "{N(\alpha,\beta)(m)}", "\shortmid"{marking}, from=1-1, to=1-2]
      \arrow["Gg", from=1-2, to=2-2]
      \arrow[""{name=1, anchor=center, inner sep=0}, "{N(\alpha,\beta)(n)}"', "\shortmid"{marking}, from=2-1, to=2-2]
      \arrow["{N(\alpha, \beta)(\sigma)}"{description}, draw=none, from=0, to=1]
    \end{tikzcd}
  \end{equation*}
  by applying the universal property of the restriction $N(\alpha,\beta)(n)$ in $\dbl{E}$:
  \begin{equation*}
    \begin{tikzcd}
      Fx & Fx & Gy & Gy \\
      Fw & Hx & Ky & Gz \\
      Hw & Hw & Kz & Kz \\
      Hw &&& Kz
      \arrow["{\alpha_x}"', from=1-2, to=2-2]
      \arrow["{\beta_y}", from=1-3, to=2-3]
      \arrow[""{name=0, anchor=center, inner sep=0}, "Nm", "\shortmid"{marking}, from=2-2, to=2-3]
      \arrow[""{name=1, anchor=center, inner sep=0}, "{N(\alpha,\beta)(m)}", "\shortmid"{marking}, from=1-2, to=1-3]
      \arrow["Hf"', from=2-2, to=3-2]
      \arrow["Kg", from=2-3, to=3-3]
      \arrow[""{name=2, anchor=center, inner sep=0}, "Nn"', "\shortmid"{marking}, from=3-2, to=3-3]
      \arrow[""{name=3, anchor=center, inner sep=0}, "{\id_{Fx}}", "\shortmid"{marking}, from=1-1, to=1-2]
      \arrow["Ff"', from=1-1, to=2-1]
      \arrow["{\alpha_w}"', from=2-1, to=3-1]
      \arrow[""{name=4, anchor=center, inner sep=0}, "{H\id_w}"', "\shortmid"{marking}, from=3-1, to=3-2]
      \arrow[""{name=5, anchor=center, inner sep=0}, "{K\id_z}"', "\shortmid"{marking}, from=3-3, to=3-4]
      \arrow["Gg", from=1-4, to=2-4]
      \arrow["{\beta_z}", from=2-4, to=3-4]
      \arrow[""{name=6, anchor=center, inner sep=0}, "{\id_{Gy}}", "\shortmid"{marking}, from=1-3, to=1-4]
      \arrow[""{name=7, anchor=center, inner sep=0}, "Nn"', "\shortmid"{marking}, from=4-1, to=4-4]
      \arrow[Rightarrow, no head, from=3-4, to=4-4]
      \arrow[Rightarrow, no head, from=3-1, to=4-1]
      \arrow["\res"{description, pos=0.4}, draw=none, from=1, to=0]
      \arrow["N\sigma"{description}, draw=none, from=0, to=2]
      \arrow["{\alpha_f^{-1}}"{description}, draw=none, from=3, to=4]
      \arrow["{\beta_g}"{description}, draw=none, from=6, to=5]
      \arrow["{N_{w,n,z}}"{description}, draw=none, from=2, to=7]
    \end{tikzcd}
    \quad=\quad
    \begin{tikzcd}[column sep=large]
      Fx & Gy \\
      Fw & Gz \\
      Hw & Kz
      \arrow[""{name=0, anchor=center, inner sep=0}, "{N(\alpha,\beta)(m)}", "\shortmid"{marking}, from=1-1, to=1-2]
      \arrow["Ff"', from=1-1, to=2-1]
      \arrow["Gg", from=1-2, to=2-2]
      \arrow[""{name=1, anchor=center, inner sep=0}, "{N(\alpha,\beta)(n)}", "\shortmid"{marking}, from=2-1, to=2-2]
      \arrow["{\alpha_w}"', from=2-1, to=3-1]
      \arrow["{\beta_z}", from=2-2, to=3-2]
      \arrow[""{name=2, anchor=center, inner sep=0}, "Nn"', "\shortmid"{marking}, from=3-1, to=3-2]
      \arrow["{\exists !}"{description, pos=0.4}, draw=none, from=0, to=1]
      \arrow["\res"{description}, draw=none, from=1, to=2]
    \end{tikzcd}.
  \end{equation*}
  Here $\alpha_f^{-1}$ denotes the inverse naturality comparison as in
  \cref{eq:pseudo-transformation-inverse}, hence the need for $\alpha$ to be a
  pseudo natural transformation. Finally, to define the left and right actions
  \begin{equation*}
    \begin{tikzcd}[column sep=large]
      Fx & Fy & Gz \\
      Fx && Gz
      \arrow["Fm", "\shortmid"{marking}, from=1-1, to=1-2]
      \arrow["{N(\alpha,\beta)(n)}", "\shortmid"{marking}, from=1-2, to=1-3]
      \arrow[""{name=0, anchor=center, inner sep=0}, "{N(\alpha,\beta)(m \odot n)}"', "\shortmid"{marking}, from=2-1, to=2-3]
      \arrow[Rightarrow, no head, from=1-1, to=2-1]
      \arrow[Rightarrow, no head, from=1-3, to=2-3]
      \arrow["{N(\alpha,\beta)_{m,n}^\ell}"{description}, draw=none, from=1-2, to=0]
    \end{tikzcd}
    \qquad\text{and}\qquad
    \begin{tikzcd}[column sep=large]
      Fx & Gy & Gz \\
      Fx && Gz
      \arrow[""{name=0, anchor=center, inner sep=0}, "{N(\alpha,\beta)(m \odot n)}"', "\shortmid"{marking}, from=2-1, to=2-3]
      \arrow[Rightarrow, no head, from=1-1, to=2-1]
      \arrow[Rightarrow, no head, from=1-3, to=2-3]
      \arrow["Gn", "\shortmid"{marking}, from=1-2, to=1-3]
      \arrow["{N(\alpha,\beta)(m)}", "\shortmid"{marking}, from=1-1, to=1-2]
      \arrow["{N(\alpha,\beta)_{m,n}^r}"{description}, draw=none, from=1-2, to=0]
    \end{tikzcd}
  \end{equation*}
  of the module $N(\alpha,\beta)$ for a pair of proarrows
  $x \xproto{m} y \xproto{n} z$ in $\dbl{D}$, we again use the universal
  property of restrictions in $\dbl{E}$. The left action is defined by
  \begin{equation*}
    \begin{tikzcd}
      Fx & Fy & Gz \\
      Hx & Hy & Kz \\
      Hx && Kz
      \arrow[""{name=0, anchor=center, inner sep=0}, "Fm", "\shortmid"{marking}, from=1-1, to=1-2]
      \arrow[""{name=1, anchor=center, inner sep=0}, "Hm"', "\shortmid"{marking}, from=2-1, to=2-2]
      \arrow["{\alpha_x}"', from=1-1, to=2-1]
      \arrow["{\alpha_y}", from=1-2, to=2-2]
      \arrow[""{name=2, anchor=center, inner sep=0}, "{N(\alpha,\beta)(n)}", "\shortmid"{marking}, from=1-2, to=1-3]
      \arrow[""{name=3, anchor=center, inner sep=0}, "Nn"', "\shortmid"{marking}, from=2-2, to=2-3]
      \arrow["{\beta_z}", from=1-3, to=2-3]
      \arrow[""{name=4, anchor=center, inner sep=0}, "{N(m \odot n)}"', "\shortmid"{marking}, from=3-1, to=3-3]
      \arrow[Rightarrow, no head, from=2-3, to=3-3]
      \arrow[Rightarrow, no head, from=2-1, to=3-1]
      \arrow["{\alpha_m}"{description}, draw=none, from=0, to=1]
      \arrow["\res"{description}, draw=none, from=2, to=3]
      \arrow["{N_{m,n}^\ell}"{description}, draw=none, from=2-2, to=4]
    \end{tikzcd}
    \quad=\quad
    \begin{tikzcd}
      Fx & Fy & Gz \\
      Fx && Gz \\
      Hx && Kz
      \arrow["Fm", "\shortmid"{marking}, from=1-1, to=1-2]
      \arrow["{N(\alpha,\beta)(n)}", "\shortmid"{marking}, from=1-2, to=1-3]
      \arrow[""{name=0, anchor=center, inner sep=0}, "{N(\alpha,\beta)(m \odot n)}", "\shortmid"{marking}, from=2-1, to=2-3]
      \arrow[Rightarrow, no head, from=1-1, to=2-1]
      \arrow[Rightarrow, no head, from=1-3, to=2-3]
      \arrow["{\alpha_x}"', from=2-1, to=3-1]
      \arrow["{\beta_z}", from=2-3, to=3-3]
      \arrow[""{name=1, anchor=center, inner sep=0}, "{N(m \odot n)}"', "\shortmid"{marking}, from=3-1, to=3-3]
      \arrow["{\exists !}"{description, pos=0.3}, draw=none, from=1-2, to=0]
      \arrow["\res"{description}, draw=none, from=0, to=1]
    \end{tikzcd}
  \end{equation*}
  and the right action is defined similarly. The restriction cell
  $\res: N(\alpha,\beta) \to N$ is the modulation whose component at a proarrow
  $m: x \proto y$ in $\dbl{D}$ is the corresponding restriction cell in
  $\dbl{E}$:
  \begin{equation*}
    \begin{tikzcd}
      Fx && Gy \\
      Hx && Ky
      \arrow[""{name=0, anchor=center, inner sep=0}, "{N(\alpha,\beta)(m)}", "\shortmid"{marking}, from=1-1, to=1-3]
      \arrow["{\alpha_x}"', from=1-1, to=2-1]
      \arrow["{\beta_y}", from=1-3, to=2-3]
      \arrow[""{name=1, anchor=center, inner sep=0}, "Nm"', "\shortmid"{marking}, from=2-1, to=2-3]
      \arrow["{\res_m}"{description}, draw=none, from=0, to=1]
    \end{tikzcd}
    \quad\coloneqq\quad
    \begin{tikzcd}
      Fx && Gy \\
      Hx && Ky
      \arrow[""{name=0, anchor=center, inner sep=0}, "{Nm(\alpha_x,\beta_y)}", "\shortmid"{marking}, from=1-1, to=1-3]
      \arrow["{\alpha_x}"', from=1-1, to=2-1]
      \arrow["{\beta_y}", from=1-3, to=2-3]
      \arrow[""{name=1, anchor=center, inner sep=0}, "Nm"', "\shortmid"{marking}, from=2-1, to=2-3]
      \arrow["\res"{description}, draw=none, from=0, to=1]
    \end{tikzcd}.
  \end{equation*}
  We defer for the moment the proof that these constructions are well-defined.

  To verify the universal property, suppose given pseudo transformations
  $\gamma: F_0 \to F$ and $\delta: F_k \to G$ along with a multimodulation $\mu$
  as on the left
  \begin{equation*}
    \begin{tikzcd}
      {F_0} & \cdots & {F_k} \\
      H && K
      \arrow["{\alpha \circ \gamma}"', from=1-1, to=2-1]
      \arrow["{\beta \circ \delta}", from=1-3, to=2-3]
      \arrow[""{name=0, anchor=center, inner sep=0}, "N"', "\shortmid"{marking}, from=2-1, to=2-3]
      \arrow["{M_1}", "\shortmid"{marking}, from=1-1, to=1-2]
      \arrow["{M_k}", "\shortmid"{marking}, from=1-2, to=1-3]
      \arrow["\mu"{description}, draw=none, from=1-2, to=0]
    \end{tikzcd}
    \qquad\leadsto\qquad
    \begin{tikzcd}
      {F_0} & \cdots & {F_k} \\
      F && G
      \arrow["{M_1}", "\shortmid"{marking}, from=1-1, to=1-2]
      \arrow["{M_k}", "\shortmid"{marking}, from=1-2, to=1-3]
      \arrow[""{name=0, anchor=center, inner sep=0}, "{N(\alpha,\beta)}"', "\shortmid"{marking}, from=2-1, to=2-3]
      \arrow["\gamma"', from=1-1, to=2-1]
      \arrow["\delta", from=1-3, to=2-3]
      \arrow["\nu"{description}, draw=none, from=1-2, to=0]
    \end{tikzcd},
  \end{equation*}
  and construct a new multimodulation $\nu$ as on the right, with component at a
  sequence of proarrows $x_0 \xproto{m_1} \cdots \xproto{m_k} x_k$ in $\dbl{D}$
  given by
  \begin{equation*}
    \begin{tikzcd}
      {F_0x_0} & \cdots & {F_kx_k} \\
      {Fx_0} && {Gx_k} \\
      {Hx_0} && {Kx_k}
      \arrow["{M_1m_1}", "\shortmid"{marking}, from=1-1, to=1-2]
      \arrow["{M_km_k}", "\shortmid"{marking}, from=1-2, to=1-3]
      \arrow["{\gamma_{x_0}}"', from=1-1, to=2-1]
      \arrow["{\delta_{x_k}}", from=1-3, to=2-3]
      \arrow["{\alpha_{x_0}}"', from=2-1, to=3-1]
      \arrow["{\beta_{x_k}}", from=2-3, to=3-3]
      \arrow[""{name=0, anchor=center, inner sep=0}, "{N(m_1 \odot \dots \odot m_k)}"', "\shortmid"{marking}, from=3-1, to=3-3]
      \arrow["{\mu_{m_1,\dots,m_k}}"{description}, draw=none, from=1-2, to=0]
    \end{tikzcd}
    \quad=\quad
    \begin{tikzcd}
      {F_0x_0} & \cdots & {F_kx_k} \\
      {Fx_0} && {Gx_k} \\
      {Hx_0} && {Kx_k}
      \arrow["{M_1m_1}", "\shortmid"{marking}, from=1-1, to=1-2]
      \arrow["{M_km_k}", "\shortmid"{marking}, from=1-2, to=1-3]
      \arrow["{\gamma_{x_0}}"', from=1-1, to=2-1]
      \arrow["{\delta_{x_k}}", from=1-3, to=2-3]
      \arrow["{\alpha_{x_0}}"', from=2-1, to=3-1]
      \arrow["{\beta_{x_k}}", from=2-3, to=3-3]
      \arrow[""{name=0, anchor=center, inner sep=0}, "{N(m_1 \odot \dots \odot m_k)}"', "\shortmid"{marking}, from=3-1, to=3-3]
      \arrow[""{name=1, anchor=center, inner sep=0}, "{N(\alpha,\beta)(m_1 \odot \cdots \odot m_k)}", "\shortmid"{marking}, from=2-1, to=2-3]
      \arrow["\res"{description}, draw=none, from=1, to=0]
      \arrow["{\exists !\ \nu_{m_1,\dots,m_k}}"{description, pos=0.3}, draw=none, from=1-2, to=1]
    \end{tikzcd}.
  \end{equation*}
  By construction, the multimodulation $\nu$ factors $\mu$ through the
  modulation $\res: N(\alpha,\beta) \to N$, and the uniqueness of the
  factorization is immediate from the universal property of restrictions in
  $\dbl{E}$.

  It remains to check that the restricted module $N(\alpha,\beta)$, the
  modulation $\res: N(\alpha,\beta) \to N$, and the multimodulations $\nu$ given
  by the universal property, are all well-defined. These calculations follow a
  common pattern: to prove a particular axiom for, say, the module
  $N(\alpha,\beta)$, one uses the universal property of the restrictions to
  reduce to the same axiom for $N$, possible after intermediate calculations. We
  illustrate by proving the naturality axiom for the left action of
  $N(\alpha,\beta)$. Naturality of the right action is analogous. We omit proofs
  of the other axioms.

  Given cells $\inlineCell{x}{y}{x'}{y'}{m}{m'}{f}{g}{\sigma}$ and
  $\inlineCell{y}{z}{y'}{z'}{n}{n'}{g}{h}{\tau}$ in $\dbl{D}$, we must show that
  \begin{equation*}
    \begin{dblArray}{cc}
      F\sigma & N(\alpha,\beta)\tau \\
      \Block{1-2}{N(\alpha,\beta)_{m',n'}^\ell}
    \end{dblArray}
    \quad=\quad
    \begin{dblArray}{c}
      N(\alpha,\beta)_{m,n}^\ell \\
      N(\alpha,\beta)(\sigma \odot \tau)
    \end{dblArray}.
  \end{equation*}
  By the universal property of restrictions, composing both sides of the
  equation on the bottom with the restriction cell
  $N(m' \odot n')(\alpha_{x'},\beta_{z'}) \to N(m' \odot n')$ yields an
  equivalent equation
  \begin{equation*}
    \begin{dblArray}{cc}
      F\sigma & N(\alpha,\beta)(\tau) \\
      \alpha_{m'} & \res \\
      \Block{1-2}{N_{m',n'}^\ell}
    \end{dblArray}
    \quad=\quad
    \begin{dblArray}{ccc}
      \id_{1_{Fx}} & N(\alpha,\beta)_{m,n}^\ell & \id_{1_{Gz}} \\
      \Block{2-1}{\alpha_f^{-1}} & \res & \Block{2-1}{\beta_h} \\
      & N(\sigma \odot \tau) & \\
      \Block{1-3}{N_{x',m' \odot n',z'}}
    \end{dblArray},
  \end{equation*}
  or, reducing the restrictions on both sides again,
  \begin{equation*}
    \begin{dblArray}{cccc}
      F\sigma & \Block{2-1}{\alpha_g^{-1}} & \res & \Block{2-1}{\beta_h} \\
      \alpha_{m'} & & N\tau & \\
      1_{Hm'} & \Block{1-3}{N_{y',n',z'}} \\
      \Block{1-4}{N_{m',n'}^\ell}
    \end{dblArray}
    \quad=\quad
    \begin{dblArray}{cccc}
      \Block{3-1}{\alpha_f^{-1}} & \alpha_m & \res & \Block{3-1}{\beta_h} \\
      & \Block{1-2}{N_{m,n}^\ell} & & \\
      & \Block{1-2}{N(\sigma \odot \tau)} & & \\
      \Block{1-4}{N_{x',m' \odot n',z'}} & & &
    \end{dblArray}.
  \end{equation*}
  To prove this equation, re-associate and apply the cell naturality of the
  pseudo transformation $\alpha$ to show that left-hand side is equal to
  \begin{equation*}
    \begin{dblArray}{cccc}
      F\sigma & \Block{2-1}{\alpha_g^{-1}} & \res & \Block{2-1}{\beta_h} \\
      \alpha_{m'} & & N\tau & \\
      \Block{1-4}{N_{m',y',n',z'}} & & &
    \end{dblArray}
    \quad=\quad
    \begin{dblArray}{cccc}
      \Block{2-1}{\alpha_f^{-1}} & \alpha_m & \res & \Block{2-1}{\beta_h} \\
      & H\sigma & N\tau & \\
      \Block{1-4}{N_{x',m',n',z'}} & & &
    \end{dblArray}
    \quad=\quad
    \begin{dblArray}{cccc}
      \Block{3-1}{\alpha_f^{-1}} & \alpha_m & \res & \Block{3-1}{\beta_h} \\
      & H\sigma & N\tau & \\
      & \Block{1-2}{N_{m',n'}^\ell} & & \\
      \Block{1-4}{N_{x',m' \odot n',z'}} & & &
    \end{dblArray}.
  \end{equation*}
  Finally, the right-hand side of this equation is equal to that of the previous
  equation by the naturality of the left action of the module $N$. This proves
  the naturality of the left action of restricted module $N(\alpha,\beta)$.
\end{proof}

\begin{corollary}[Virtual equipment of cartesian lax functors]
  Let $\dbl{D}$ be a cartesian double category and let $\dbl{E}$ be a cartesian
  equipment. Then:
  \begin{enumerate}[(i),noitemsep]
    \item $\vCartLaxFun_\ell(\dbl{D},\dbl{E})$ has restrictions along cartesian
      lax transformations $\alpha$ and $\beta$ whenever $\alpha$ is pseudo.
    \item $\vCartLaxFun_c(\dbl{D},\dbl{E})$ has restrictions along cartesian
      oplax transformations $\alpha$ and $\beta$ whenever $\beta$ is pseudo.
    \item $\vCartLaxFun_\pseudo(\dbl{D},\dbl{E})$ has restrictions along any
      cartesian pseudo transformations $\alpha$ and $\beta$.
  \end{enumerate}
  Consequently, the virtual double category
  $\vCartLaxFun_\pseudo(\dbl{D},\dbl{E})$ of cartesian lax functors
  $\dbl{D} \to \dbl{E}$, cartesian pseudo natural transformations, cartesian
  modules, and multimodulations is a virtual equipment.
\end{corollary}
\begin{proof}
  Since restriction cells are closed under finite products
  (\cref{lem:closure-restrictions}), the restriction of a cartesian module along
  cartesian pseudo transformations is again a cartesian module. The result thus
  follows from \cref{cor:cartesian-lax-functor-vdbl-category} and
  \cref{thm:lax-functor-veqp}.
\end{proof}

\begin{corollary}[Virtual equipment of models]
  Let $\dbl{T}$ be a cartesian double theory and let $\dbl{S}$ be a cartesian
  equipment. Then there is a virtual equipment of models of the theory $\dbl{T}$
  in $\dbl{S}$, (pseudo or strict) maps between models, bimodules between
  models, and multitransformations between those.

  Moreover, when the theory $\dbl{T}$ is generated by a 2-category or a
  bicategory (i.e., has only trivial arrows or trivial proarrows) and the
  semantics $\dbl{S}$ has local coequalizers, the virtual equipment of models is
  representable as an equipment.
\end{corollary}

\section{Conclusion}
\label{sec:conclusion}

Cartesian double categories, cartesian lax functors, and their higher morphisms
are the building blocks of a framework for categorical doctrines based on
two-dimensional functorial semantics. Every cartesian double theory can be
interpreted in a cartesian double category, such as that of spans in a finitely
complete category $\cat{S}$ or of matrices in a distributive category $\catV$,
to give a unital virtual double category of models. Under commonly satisfied
conditions, double theories in fact have genuine double categories of models.
Moreover, when interpreted in a cartesian equipment, cartesian double theories
have virtual equipments or equipments of models. Many familiar categorical
structures are models of cartesian double theories, such as categories,
presheaves, monoidal categories, braided and symmetric monoidal categories,
2-groups, multicategories, and cartesian and cocartesian categories. These
structures can be internal to $\cat{S}$ or enriched in $\catV$, depending on the
choice of semantics. Restriction sketches add further expressivity to double
theories by incorporating restrictions as in an equipment. Symmetric
multicategories and cartesian multicategories are models of restriction
sketches.

The hierarchy of categorical logics motivates the need for doctrines as a
unifying principle. But, at a higher level, doctrines themselves vary according
to the structure needed to define them. The same tradeoff presents itself for
doctrines as for theories: the richer the logic needed to define a theory or
doctrine, the fewer semantics in which it can be interpreted. Enriched category
theory illustrates this principle. The constructions possible with
$\catV$-enriched categories depend on the properties of the base of enrichment
$\catV$; the more properties that the base $\catV$ posseses, the more of
ordinary category theory that can reproduced with $\catV$-enriched categories
but the fewer examples there are of such bases to begin with. A flexible
approach to doctrines should encompass enriched category theory at different
points in this spectrum.

Through this work, we have developed two notions of doctrine in detail: simple
double theories, which are double categories without further structure, and
cartesian double theories, which are cartesian double categories. This seems a
natural place to start, just as cartesian categories were for categorical logic.
Nevertheless, many doctrines require either less or more structure than is
available in a cartesian double theory. For example, the familiar theory of
monoids is interpretable in any monoidal category or even in any multicategory;
the theory of pseudomonoids should likewise be interpretable in any monoidal
double category or even in any double multicategory. So the theory of
pseudomonoids should be a ``monoidal double theory'' or even a ``multi-double
theory,'' rather than a cartesian double theory. We could then interpret
pseudomonoids in the double category of $\catV$-matrices and recover
$\catV$-enriched monoidal categories with fewer assumptions on the base of
enrichment $\catV$. In the other direction, it seems doubtful that categories
with finite limits are models of any cartesian double theory. Determining the
additional structures on double categories needed to climb the ladder of
categorical logics is a problem for the future. In either direction, the
development of new ``double doctrines''\footnote{If double doctrines start to
  proliferate, then inevitably we will need a three-dimensional logic to
  classify them. This is an infinite regress typical of higher category theory.}
must introduce structure not only to double categories but to lax double
functors and the higher morphisms between them.

As an orthogonal direction, a variety of evidence suggests that double theories
should perhaps be \emph{virtual} double categories. We have seen that, due to
obstructions to composing modules of lax functors \cite{pare2013}, the double
category of models of a double theory is in general only virtual. In such cases,
we can achieve an increase in generality and symmetry by taking the theory and
the semantics to be virtual as well. Throughout this work we have emphasized the
fundamental role played by lax functors. Conveniently, functors between virtual
double categories are automatically lax in that they reduce to lax functors when
pseudo double categories are regarded as virtual ones \cite[\mbox{Example
  3.5}]{cruttwell2010}. But the most appealing feature of virtual double
categories for our purposes is that they do not proliferate unwanted composites
of proarrows, for the simple reason that proarrows cannot be composed at all!
The challenge of controlling the laxators for unwanted proarrow composites is
most evident in the examples of restriction sketches. A plausible alternative to
restriction sketches would make double theories be not just virtual double
categories but virtual equipments. Developing the theory of virtual double
categories far enough to meet the needs of double-categorical logic is another
task for the future.

The hierarchy of double doctrines, which evidently exists even if it is so far
mostly unexplored, should be contrasted with the prevailing situation for
2-monads. As Hyland and Power have persuasively argued \cite[\mbox{p.\
  450}]{hyland2007}, ordinary monads are isolated in the hierarchy of
categorical logic compared to Lawvere theories and finite product theories. We
expect that 2-monads will be similarly isolated in the hierarchy of
double-categorical logic. Still, since models of Lawvere theories are known to
be the same as algebras of finitary monads on $\Set$ \cite[Theorem
A.40]{adamek2010}, it is natural to ask how double theories are related to
double monads. A \emph{double monad} is a monad in the 2-category of double
categories for some choice of pseudo or (op)lax double functors
\mbox{\cite[\S{7}]{grandis2004}}. A further question would then be how such
double monads are related to the 2-monads on $\Cat$ that have been so
intensively studied.

Relatedly, we have only begun to investigate the properties possessed by the
double category of models of a cartesian double theory. Categories of models of
finite products theories, sometimes called ``algebraic categories''
\cite{adamek2010}, are very well behaved. For example, algebraic categories are
complete and cocomplete, with limits computed pointwise in $\Set$. It is natural
to wonder about analogous properties of double categories of models.

Formal category theory aims to treat the fundamental concepts and constructions
of category theory in an axiomatic or synthetic style, independent of specific
details about the 2-category of categories. A long strand of work has identified
equipments \cite{wood1982,shulman2021}, and related structures such as virtual
equipments \cite{cruttwell2010} and augmented virtual equipments
\cite{koudenburg2020}, as ideal environments for formal category theory and,
more recently, even for formal $\infty$-category theory \cite[\mbox{Chapter
  9}]{riehl2022}. Since each double theory gives a virtual equipment, and
sometimes an equipment, of models, double theories can be seen as machines that
generate new environments in which to do formal category theory. As category
theorists, we naturally wish to go further, considering not just these
environments in isolation but the passages between them. We should certainly
expect virtual equipments of models to pull back along morphisms between
theories. For this, we must go beyond the two-dimensional framework developed
here to a fully three-dimensional structure encompassing at least double
categories, lax functors, lax transformations, modules, and multimodulations. As
the proper understanding of categorical logic requires at least a 2-category of
theories, models, and model homomorphisms, so must double-categorical logic
ultimately lead to three-dimensional category theory.

\paragraph{Acknowledgments} We thank Nathanael Arkor for detailed comments on
earlier drafts of the paper, including refinements to
\cref{thm:lax-functor-veqp} and its corollaries. The last author thanks his
colleagues Kevin Arlin, Owen Lynch, Brandon Shapiro, and David Spivak for
helpful conversations during the course of this research.

\printbibliography[heading=bibintoc]

\end{document}